\documentclass[11pt]{article}
\usepackage{amsfonts}
\usepackage{amssymb}
\usepackage{amsmath}
=0pt
\def\sqr#1#2{{\vcenter{\vbox{\hrule height.#2pt
              \hbox{\vrule width.#2pt height#1pt \kern#1pt \vrule
width.#2pt}
              \hrule height.#2pt}}}}
\def\signed #1{{\unskip\nobreak\hfil\penalty50
              \hskip2em\hbox{}\nobreak\hfil#1
              \parfillskip=0pt \finalhyphendemerits=0 \par}}
\def\endpf{\signed {$\sqr69$}}
\def\dbC{{\mathbb{C}}}

\def\dbE{{\mathbb{E}}}
\def\dbF{{\mathbb{F}}}

\def\dbH{{\mathbb{H}}}

\def\dbN{{\mathbb{N}}}

\def\dbP{{\mathbb{P}}}

\def\dbR{{\mathbb{R}}}

\def\dbW{{\mathbb{W}}}

\def\a{\alpha}
\def\b{\beta}
\def\g{\gamma}
\def\d{\delta}
\def\e{\varepsilon}

\def\l{\lambda}
\def\m{\mu}

\def\si{\sigma}

\def\f{\varphi}

\def\o{\omega}

\def\3n{\negthinspace \negthinspace \negthinspace }
\def\2n{\negthinspace \negthinspace }
\def\1n{\negthinspace }
\def\ns{\noalign{\smallskip} }

\def\ds{\displaystyle}
\def\G{\Gamma}
\def\D{\Delta}

\def\L{\Lambda}
\def\Si{\Sigma}

\def\O{\Omega}
\def\cA{{\cal A}}
\def\cB{{\cal B}}

\def\cF{{\cal F}}
\def\cG{{\cal G}}
\def\cH{{\cal H}}

\def\cJ{{\cal J}}

\def\cL{{\cal L}}
\def\cM{{\cal M}}

\def\cQ{{\cal Q}}

\def\cT{{\cal T}}
\def\cU{{\cal U}}

\def\cX{{\cal X}}

\def\mE{{\mathbb{E}}}

\def\no{\noindent}

\def\ms{\medskip}
\def\bs{\bigskip}
\def\q{\quad}
\def\qq{\qquad}
\def\hb{\hbox}

\def\lan{\mathop{\langle}}
\def\ran{\mathop{\rangle}}

\def\cd{\cdot}
\def\cds{\cdots}

\def\dim{\hbox{\rm dim$\,$}}

\def\ae{\hbox{\rm a.e.{ }}}
\def\as{\hbox{\rm a.s.{ }}}

\def\span{\hbox{\rm span$\,$}}
\def\tr{\hbox{\rm tr$\,$}}

\def\deq{\mathop{\buildrel\D\over=}}
\def\Re{{\mathop{\rm Re}\,}}

\def\({\Big (}
\def\){\Big )}
\def\[{\Big[}
\def\]{\Big]}
\def\neg{\negthinspace}
\def\={\buildrel \triangle \over =}

\def\resp{{\it resp. }}

at 9pt 

\def\be{\begin{equation}}
\def\bel{\begin{equation}\label}
\def\ee{\end{equation}}
\def\bea{\begin{eqnarray}}
\def\eea{\end{eqnarray}}
\def\bt{\begin{theorem}}
\def\et{\end{theorem}}
\def\bc{\begin{corollary}}
\def\ec{\end{corollary}}
\def\bl{\begin{lemma}}
\def\el{\end{lemma}}
\def\bp{\begin{proposition}}
\def\ep{\end{proposition}}
\def\br{\begin{remark}}
\def\er{\end{remark}}
\def\ba{\begin{array}}
\def\ea{\end{array}}
\def\bd{\begin{definition}}
\def\ed{\end{definition}}

\newtheorem{lemma}{Lemma}[section]
\newtheorem{remark}{Remark}[section]

\newtheorem{theorem}{Theorem}[section]
\newtheorem{corollary}{Corollary}[section]

\newtheorem{definition}{Definition}[section]
\newtheorem{proposition}{Proposition}[section]

\makeatletter
   
   \@addtoreset{equation}{section}
\makeatother

\begin{document}

\title{\bf General Pontryagin-Type Stochastic Maximum
Principle and Backward Stochastic Evolution
Equations in Infinite Dimensions}

\author{Qi L\"{u}\thanks{School of Mathematical Sciences, University of Electronic Science and Technology of China, Chengdu 610054, China. {\small\it e-mail:} {\small\tt
luqi59@163.com}. \ms}
 ~~~
  and~~~
Xu Zhang\thanks{Yangtze Center of Mathematics,
Sichuan University, Chengdu 610064, China.
{\small\it e-mail:} {\small\tt
zhang$\_$xu@scu.edu.cn}.}}
\date{}

\maketitle

\begin{abstract}

The main purpose of this paper is to give a
solution to a long-standing unsolved problem in
stochastic control theory, i.e., to establish
the Pontryagin-type maximum principle for
optimal controls of general infinite
dimensional nonlinear stochastic evolution
equations.  Both drift and diffusion terms can
contain the control variables, and the control
domains are allowed to be nonconvex. The key to
reach it is to provide a suitable formulation
of operator-valued backward stochastic
evolution equations (BSEEs for short), as well
as a way to define their solutions. Besides,
both vector-valued and operator-valued BSEEs,
with solutions in the sense of transposition,
are studied. As a crucial preliminary, some
weakly sequential Banach-Alaoglu-type theorems
are established for uniformly bounded linear
operators between Banach spaces.

\end{abstract}

\bs

\no{\bf 2010 Mathematics Subject
Classification}.  Primary 93E20; Secondary
60G05, 60H15, 60G07.

\bs

\no{\bf Key Words}. Stochastic evolution
equation, optimal control, Pontryagin-type
maximum principle, backward stochastic
evolution equation, transposition solution.

\newpage
\tableofcontents

\newpage

\section{Introduction}\label{s1}

Let $(\O,\cF,\dbF,\dbP)$ be a complete filtered
probability space with the filtration
$\dbF=\{\cF_t\}_{t\ge0}$, on which a
one-dimensional standard Brownian motion
$\{w(t)\}_{t\ge0}$ is defined. Let $T>0$, and
let $X$ be a Banach space. For any $t\in[0,T]$
and $r\in [1,\infty)$, denote by
$L_{\cF_t}^r(\O;X)$ the Banach space of all
$\cF_t$-measurable random variables $\xi:\O\to
X$ such that $\mathbb{E}|\xi|_X^r < \infty$,
with the canonical  norm.  Also, denote by
$D_{\dbF}([0,T];L^{r}(\O;X))$ the vector space
of all $X$-valued $r$th power integrable
$\dbF$-adapted processes $\phi(\cdot)$ such
that $\phi(\cdot):[0,T] \to
L^{r}(\O,\cF_T,P;X)$ is c\`adl\`ag, i.e., right
continuous with left limits. Clearly,
$D_{\dbF}([0,T];L^{r}(\O;X))$ is a Banach space
with the following norm
 $$
|\phi(\cd)|_{D_{\dbF}([0,T];L^{r}(\O;X))} =
\sup_{t\in
[0,T)}\left[\mE|\phi(r)|_X^r\right]^{1/r}.
 $$
We denote by $C_{\dbF}([0,T];L^{r}(\O;X))$ the
Banach space of all $X$-valued $\dbF$-adapted
processes $\phi(\cdot)$ such that
$\phi(\cdot):[0,T] \to L^{r}(\O,\cF_T,P;X)$ is
continuous,  with the norm inherited from
$D_{\dbF}([0,T];$ $L^{r}(\O;X))$. Fix any
$r_1,r_2,r_3,r_4\in[1,\infty]$. Put
 $$
 \ba{ll}
\ds
L^{r_1}_\dbF(\O;L^{r_2}(0,T;X))=\Big\{\f:(0,T)\times\O\to
X\bigm|\f(\cd)\hb{
is $\dbF$-adapted and }\dbE\(\int_0^T|\f(t)|_X^{r_2}dt\)^{\frac{r_1}{r_2}}<\infty\Big\},\\
\ns\ds
L^{r_2}_\dbF(0,T;L^{r_1}(\O;X))=\Big\{\f:(0,T)\times\O\to
X\bigm|\f(\cd)\hb{ is $\dbF$-adapted and
}\int_0^T\(\dbE|\f(t)|_X^{r_1}\)^{\frac{r_2}
{r_1}}dt<\infty\Big\}.
 \ea
 $$
Clearly, both $L^{r_1}_\dbF(\O;L^{r_2}(0,T;X))$
and $L^{r_2}_\dbF(0,T;L^{r_1}(\O;X))$ are
Banach spaces with the canonical norms. If
$r_1=r_2$, we simply denote the above space by
$L^{r_1}_\dbF(0,T;X)$. Let $Y$ be another
Banach space. Denote by $\cL(X, Y)$ the
(Banach) space of all bounded linear operators
from $X$ to $Y$, with the usual operator norm
(When $Y=X$, we simply write $\cL(X)$ instead
of $\cL(X, Y)$). Further, we denote by $
\cL_{pd}\big(L^{r_1}_{\dbF}(0,T;L^{r_2}(\O;X)),\;L^{r_3}_{\dbF}(0,T;L^{r_4}(\O;Y))\big)$
(\resp $
\cL_{pd}\big(X,\;L^{r_3}_{\dbF}(0,T;L^{r_4}(\O;Y))\big)$)
the vector space of all bounded, pointwisely
defined linear operators $\cL$ from
$L^{r_1}_{\dbF}(0,T;L^{r_2}(\O;X))$ (\resp $X$)
to $L^{r_3}_{\dbF}(0,T;L^{r_4}(\O;Y))$, i.e.,
for $\ae (t,\omega)\in (0,T)\times\Omega$,
there exists an $L(t,\o)\in\cL (X,Y)$ verifying
that $\big(\cL u(\cd)\big)(t,\o)=L
(t,\o)u(t,\o), \;\, \forall\; u(\cd)\in
L^{r_1}_{\dbF}(0,T;L^{r_2}(\O;X))$ (\resp
$\big(\cL x\big)(t,\o)=L (t,\o)x, \;\,
\forall\; x\in X$). Similarly, one can define
the spaces
$\cL_{pd}\big(L^{r_2}(\O;X),\;L^{r_3}_{\dbF}(0,T;L^{r_4}(\O;Y))\big)$
and
$\cL_{pd}\big(L^{r_2}(\O;X),\;L^{r_4}(\O;Y)\big)$,
etc.

Let $H$ be a complex Hilbert space, and let $A$
be an unbounded linear operator (with domain
$D(A)$ on $H$), which is the infinitesimal
generator of a $C_0$-semigroup $\{S(t)\}_{t\geq
0}$. Denote by $A^*$ the dual operator of $A$.
Clearly, $D(A)$ is a Hilbert space with the
usual graph norm, and $A^*$ is the
infinitesimal generator of $\{S^*(t)\}_{t\geq
0}$, the dual $C_0$-semigroup of
$\{S(t)\}_{t\geq 0}$. For any $\l\in\rho(A)$,
the resolvent of $A$, denote by $A_\l$ the
Yosida approximation of $A$ and by
$\{S_\l(t)\}_{t\in\dbR}$ the $C_0$-group
generated by $A_\l$. Let $U$ be a metric space
with its metric $d(\cd,\cd)$. Put
$$\cU[0,T] \triangleq \Big\{u(\cdot):\,[0,T]\to U\;\Big|\; u(\cd) \mbox{ is $\dbF$-adapted} \Big\}.$$
Throughout this paper, we assume the following
condition.

\ms

\no{\bf (A1)} {\it Suppose that
$a(\cd,\cd,\cd):[0,T]\times H\times U\to H$ and
$b(\cd,\cd,\cd):[0,T]\times H\times U\to H$ are
two maps satisfying: i) For any $(x,u)\in
H\times U$, the maps $a(\cd,x,u):[0,T]\to H$
and $b(\cd,x,u):[0,T]\to H$ are Lebesgue
measurable; ii) For any $(t,x)\in [0,T]\times
H$, the maps $a(t,x,\cd):U\to H$ and
$b(t,x,\cd):U\to H$ are continuous; and iii)
There is a constant $C_L>0$ such that
\begin{equation}\label{ab0}
\left\{
\begin{array}{ll}\ds
|a(t,x_1,u) - a(t,x_2,u)|_H+|b(t,x_1,u) -
b(t,x_2,u)|_H \leq
C_L|x_1-x_2|_H,\\
\ns\ds |a(t,0,u)|_H +|b(t,0,u)|_H \leq C_L, \\
\ns\ds \hspace{5.5cm} \forall\;
(t,x_1,x_2,u)\in [0,T]\times H\times H\times U.
\end{array}
\right.
\end{equation}}

Consider the following controlled (forward)
stochastic evolution equation:
\begin{eqnarray}\label{fsystem1}
\left\{
\begin{array}{lll}\ds
dx = \big[Ax +a(t,x,u)\big]dt + b(t,x,u)dw(t) &\mbox{ in }(0,T],\\
\ns\ds x(0)=x_0,
\end{array}
\right.
\end{eqnarray}
where $u\in \cU[0,T]$ and $x_0\in
L^{p_0}_{\cF_0}(\O;H)$ for some given $p_0>1$.
We call $x(\cd)\equiv x(\cd\,;x_0,u)\in
C_\dbF([0,T];L^{p_0}(\O;H))$ a mild solution to
\eqref{fsystem1} if
$$
x(t)=S(t)x_0 + \int_0^t S(t-s)a(s,x(s),u(s))ds
+ \int_0^t
S(t-s)b(s,x(s),u(s))dw(s),\q\forall\;
t\in[0,T].
$$
In the sequel, we shall denote by $C$ a generic
constant, depending on $T$, $A$, $p_0$ (or $p$
to be introduced later) and $C_L$ (or $J$ and
$K$ to be introduced later), which may be
different from one place to another. Similar to
\cite[Chapter 7]{Prato}, it is easy to show the
following result:

\begin{lemma}\label{well lemma s1}
Let the assumption (A1) hold. Then, the
equation \eqref{fsystem1} is well-posed in the
sense of mild solution. Furthermore,
$$
|x(\cd)|_{C_\dbF([0,T];L^{p_0}(\O;H))} \leq
C\big(1+|x_0|_{L^{p_0}_{\cF_0}(\O;H)}\big).
$$
\end{lemma}

Also, we need the following condition:

\ms

\no{\bf (A2)} {\it Suppose that
$g(\cd,\cd,\cd):[0,T]\times H\times U\to \dbR$
and $h(\cd):H\to \dbR$ are two functions
satisfying: i) For any $(x,u)\in H\times U$,
the function $g(\cd,x,u):[0,T]\to \dbR$ is
Lebesgue measurable; ii) For any $(t,x)\in
[0,T]\times H$, the function $g(t,x,\cd):U\to
\dbR$ is continuous; and iii) There is a
constant $C_L>0$ such that
\begin{equation}\label{gh} \left\{
\begin{array}{ll}\ds
|g(t,x_1,u) - g(t,x_2,u)|_{H} +|h(x_1) -
h(x_2)|_H
 \leq C_L|x_1-x_2|_H,\\
\ns\ds |g(t,0,u)|_H +|h(0)|_H \leq C_L,\\
\ns\ds \hspace{4.5cm} \forall\;
(t,x_1,x_2,u)\in [0,T]\times H\times H\times U.
\end{array}
\right.
\end{equation}}

\ms

Define a cost functional $\cJ(\cdot)$ (for the
controlled system \eqref{fsystem1}) as follows:
 \bel{jk1}
\cJ(u(\cdot))\triangleq \dbE\Big[\int_0^T
g(t,x(t),u(t))dt + h(x(T))\Big],\qq\forall\;
u(\cdot)\in \cU[0,T],
 \ee
where $x(\cd)$ is the corresponding solution to
\eqref{fsystem1}.

Let us consider the following optimal control
problem for the system \eqref{fsystem1}:

\ms

\no {\bf Problem (P)} {\it Find a $\bar
u(\cdot)\in \cU[0,T]$ such that
 \bel{jk2}
\cJ (\bar u(\cdot)) = \inf_{u(\cdot)\in
\cU[0,T]} \cJ (u(\cdot)).
 \ee
Any $\bar u(\cdot)$ satisfying (\ref{jk2}) is
called an {\it optimal control}. The
corresponding state process $\bar x(\cdot)$ is
called an {\it optimal state (process)}, and
$(\bar x(\cdot),\bar u(\cdot))$ is called an
{\it optimal pair}.}

\ms

The main goal of this paper is to establish
some necessary conditions for optimal pairs of
Problem (P), in the spirit of the
Pontryagin-type maximum principle (\cite{PC}).
In this respect, the problem is now
well-understood in the case that $\dim
H<\infty$. We refer to \cite{Kushner} and the
references therein for early studies on the
maximum principle for controlled stochastic
differential equations in finite dimensional
spaces. Then, people established further
results on the maximum principle for stochastic
control systems under various assumptions, say,
the diffusion coefficients were non-degenerate
(e.g. \cite{Haussmann}), and/or the diffusion
coefficients were independent of the controls
(e.g. \cite{Bensoussan1,Bismut}), and/or the
control domains were convex (e.g.
\cite{Bensoussan1}). Note that, generally
speaking, many practical systems (especially in
the area of finance) do not satisfy these
assumptions. In \cite{Peng1}, a maximum
principle was obtained for general stochastic
control systems without the above mentioned
assumptions, and it was found that the
corresponding result in the general case
differs essentially from its deterministic
counterpart. As important byproducts in the
study of the above finite dimensional
stochastic control problems, one introduced
some new mathematical tools, say, backward
stochastic differential equations (BSDEs, for
short) and forward-backward stochastic
differential equations (\cite{Bismut,B1,PP} and
\cite{MY,YZ}), which are now extensively
applied to many other fields.

Let us recall here the main idea and result in
\cite{Peng1}. Suppose that $(\bar x(\cd), \bar
u(\cd))$ is a given optimal pair for the
special case that $A=0$, $H=\dbR^n$ (for some
$n\in\dbN$) and $\dbF$ is the natural
filtration $\dbW$ (generated by the Brownian
motion $\{w(\cd)\}$ and augmented by all the
$\dbP$-null sets). First, similar to the
corresponding deterministic setting, one
introduces the following first order adjoint
equation (which is however a BSDE in the
stochastic case):
\begin{equation}\label{3.8}
 \left\{
 \ba{ll}dy(t)=-\Big[a_x(t,\bar x(t),\bar
 u(t))^\top
y(t)+b_x(t,\bar x(t),\bar u(t))^\top  Y(t)\\
\ns
\qq\qq\q-g_x(t,\bar x(t),\bar
u(t))\Big]dt+Y(t)dw(t),&\mbox{ in }
[0,T),\\
y(T)=-h_x(\bar x(T)).\ea\right.
 \end{equation}
Here the unknown is a {\it pair} of
$\dbF$-adapted processes $(y(\cd),Y(\cd))\in
C_{\dbF}([0,T];L^2(\O;\dbR^n)) \times
L^2_\dbF(0,T;\dbR^n)$. Next, to establish the
desired maximum principle for stochastic
controlled systems with control-dependent
diffusion and possibly nonconvex control
domains, the author in \cite{Peng1} had the
following fundamental finding: Except for the
first order adjoint equation (\ref{3.8}), one
has to introduce an additional second order
adjoint equation as follows:
 \begin{equation}\label{3.9}\left\{\ba{ll}dP(t)=-\Big[a_x(t,\bar x(t),\bar u(t))^\top P(t)+P(t)
a_x(t,\bar x(t),\bar u(t))\\
\ns
\qq\qq+b_x(t,\bar x(t),\bar u(t))^\top
P(t)b_x(t,\bar x(t), \bar u(t))\\\ns
\qq\qq+b_x(t,\bar x(t),\bar u(t))^\top
Q(t)+Q(t) b_x(t,\bar x(t),\bar u(t))\\\ns
\qq\qq +\dbH_{xx}(t,\bar x(t),\bar
u(t),y(t),Y(t))\Big]dt+Q(t) dw(t), &\mbox{ in } [0,T)\\
\ns
P(T)=-h_{xx}(\bar x(T)).\ea\right.
\end{equation}
In (\ref{3.9}), the {\it Hamiltonian}
$\dbH(\cd,\cd,\cd,\cd,\cd)$ is defined by
 $$
 \ba{ll} \dbH(t,x,u,y_1,y_2)=\lan y_1,a(t,x,u)\ran_{\dbR^n}+\lan y_2, b(t,x,u)
\ran_{\dbR^n}-g(t,x,u),\\\ns
\hspace{4.5cm}(t,x,u,y_1,y_2)\in[0,T]\neg\times\neg\dbR^n\neg\times\neg
U\neg\times \neg\dbR^n\neg\times\neg\dbR^n.
 \ea
 $$
Clearly, the equation (\ref{3.9}) is an
$\dbR^{n\times n}$-valued BSDE in which the
unknown is a pair of processes
$(P(\cd),Q(\cd))\in
C_{\dbF}([0,T];L^2(\O;\dbR^{n\times n})) \times
L^2_\dbF(0,T;\dbR^{n\times n})$. Then,
associated with the 6-tuple $(\bar x(\cd),\bar
u(\cd), y(\cd), Y(\cd), P(\cd), Q(\cd))$,
define
$$
 \ba{ll}\ds
 \cH(t,x,u)\deq \dbH(t,x,u,y(t),Y(t))+\frac{1}{2}{\lan
P(t)b(t,x,u),b(t,x,u)\ran}_{\dbR^n}\\\ns\ds
\qq\qq\q\;-{\lan P(t)b(t,\bar x(t), \bar
u(t)),b(t,x,u)\ran}_{\dbR^n}.
 \ea
 $$
The main result in \cite{Peng1} asserts that
the optimal pair $(\bar x(\cd), \bar u(\cd))$
verifies the following stochastic maximum
principle:
 $$\cH(t,\bar x(t),\bar u(t))=\max_{u\in U}\cH(t,\bar x(t),u),\q\ae
t\in[0,T], \q\dbP\hb{-}\as
 $$

On the other hand, there exist extensive works
addressing the Pontryagin-type maximum
principle for optimal controls of deterministic
infinite dimensional controlled systems (e.g.
\cite{LY} and the rich references therein).
Naturally, one expects to extend the optimal
control theory of both stochastic finite
dimensional systems and deterministic infinite
dimensional systems to that of infinite
dimensional stochastic evolution equations. In
this respect, we refer to \cite{Bensoussan2}
for a pioneer work. Later progresses are
available in the literature
\cite{Al-H2,AGY,HP1,TL,Tu,Zhou} and so on.
Nevertheless, almost all of the existing
published works on the necessary conditions for
optimal controls of infinite dimensional
stochastic evolution equations addressed only
the case that the diffusion term does NOT
depend on the control variable (i.e., the
function $b(t,x,u)$ in \eqref{fsystem1} is
independent of $u$). As far as we know, the
stochastic maximum principle for general
infinite dimensional nonlinear stochastic
systems with control-dependent diffusion
coefficients and possibly nonconvex control
domains has been a longstanding unsolved
problem.

In this paper, we aim to give a solution to the
above mentioned unsolved problem. Inspired by
\cite{Peng1}, we will first study an $H$-valued
BSEE and an $\cL(H)$-valued BSEE, employed
accordingly as the first order adjoint equation
and the second order adjoint equation (for the
original equation (\ref{fsystem1})), and then
establish the desired necessary conditions for
optimal controls with the aid of the solutions
of these equations.

First, we need to study the following
$H$-valued BSEE:
\begin{eqnarray}\label{bsystem1}
\left\{
\begin{array}{lll}
\ds dy(t) = -  A^* y(t) dt + f(t,y(t),Y(t))dt + Y(t) dw(t) &\mbox{ in }[0,T),\\
\ns\ds y(T) = y_T.
\end{array}
\right.
\end{eqnarray}
Here $y_T \in L_{\cF_T}^{p}(\O;H)$) with $p\in
(1,2]$, $f(\cd,\cd,\cd):[0,T]\times H\times H
\to H$ satisfies
\begin{equation}\label{Lm1}
\left\{
\begin{array}{ll}\ds
f(\cd,0,0)\in
L^{1}_{\dbF}(0,T;L^{p}(\O;H)),\\\ns\ds
|f(t,x_1,y_1)-f(t,x_2,y_2)|_H\leq
C_L\big(|x_1-x_2|_H+|y_1-y_2|_H \big),\\
\ns\ds\hspace{3.5cm} \ae (t,\o)\in
[0,T]\times\O,\;\; \forall\;x_1,x_2,y_1,y_2\in
H.
\end{array}
\right.
 \end{equation}
Since neither the usual natural filtration
condition nor the quasi-left continuity is
assumed for the filtration $\dbF$ in this
paper, and because the unbounded operator $A$
is assumed to generate a general
$C_0$-semigroup, we cannot apply the existing
results on infinite dimensional BSEEs (e.g.
\cite{Al-H1,HP2,MY,MM}) to obtain the
well-posedness of the equation
\eqref{bsystem1}.

Next, it is more important that the following
$\cL(H)$-valued BSEE\footnote{Throughout this
paper, for any operator-valued process (\resp
random variable) $R$, we denote by $R^*$ its
pointwisely dual operator-valued process (\resp
random variable). For example, if $R\in
L^{r_1}_\dbF(0,T; L^{r_2}(\O; \cL(H)))$, then
$R^*\in L^{r_1}_\dbF(0,T; L^{r_2}(\O;
\cL(H)))$, and $|R|_{L^{r_1}_\dbF(0,T;
L^{r_2}(\O; \cL(H)))}=|R^*|_{L^{r_1}_\dbF(0,T;
L^{r_2}(\O; \cL(H)))}$.}:
\begin{equation}\label{op-bsystem3}
\left\{\3n
\begin{array}{ll}
\ds dP  =  - (A^*  + J^* )P dt  -  P(A + J )dt
-K^*PKdt
 - (K^* Q +  Q K)dt
  +   Fdt  +  Q dw(t) &\mbox{ in } [0,T),\\
\ns\ds P(T) = P_T
\end{array}
\right.
\end{equation}
should be studied. Here and henceforth, $F\in
L^1_\dbF(0,T;L^2(\O;\cL(H)))$, $P_T\in
L^2_{\cF_T}(\O;\cL(H))$, and $J,K\in
L^4_\dbF(0,T; L^\infty(\O; \cL(H)))$. For the
special case when $H=\dbR^n$, it is easy to see
that \eqref{op-bsystem3} is an $\dbR^{n\times
n}$ (matrix)-BSDE, and therefore, the desired
well-posedness follows from that of an
$\dbR^{n^2}$(vector)-valued BSDE. However, one
has to face a real challenge in the study of
\eqref{op-bsystem3} when $\dim H=\infty$,
without further assumption on the data $F$ and
$P_T$. Indeed, in the infinite dimensional
setting, although $\cL(H)$ is still a Banach
space, it is neither reflexive (needless to say
to be a Hilbert space) nor separable even if
$H$ itself is separable (See Problem 99 in
\cite{Halmos}). As far as we know, in the
previous literatures there exists no such a
stochastic integration/evolution equation
theory in general Banach spaces that can be
employed to treat the well-posedness of
\eqref{op-bsystem3}. For example, the existing
result on stochastic integration/evolution
equation in UMD Banach spaces (e.g.
\cite{NVW1,NVW2}) does not fit the present case
because, if a Banach space is UMD, then it is
reflexive.

The key of this work is to give accordingly
reasonable definitions of the solutions to
\eqref{bsystem1} and \eqref{op-bsystem3}, and
show the corresponding well-posedness results.
For this purpose, we employ the transposition
method developed in our previous work
\cite{LZ}, which was addressed to the BSDEs in
$\dbR^n$. Our method has several advantages.
The first one is that the usual duality
relationship is contained in our definition of
solutions, and therefore, we do NOT need to use
It\^o's formula to derive this sort of relation
as usual to obtain the desired stochastic
maximum principle. Note that, it may be very
difficult to derive the desired It\^o's formula
for the mild solutions of general stochastic
evolution equations in infinite dimensions. The
second one is that we do NOT need to use the
Martingale Representation Theorem, and
therefore we can study the problem with a
general filtration. Note that, when we deal
with BSEEs with operator unknowns, as far as we
know, there exists no Martingale Representation
Theorem (for the $\cL(H)$-valued martingale)
even if $\dbF$ is the natural filtration
$\dbW$. Thirdly, as shown in \cite{WZ} (though
it addressed only BSDEs in $\dbR^n$), similar
to the classical finite element method solving
deterministic partial differential equations,
our transposition method leads naturally
numerical schemes to solve both vector-valued
and operator-valued BSEEs (The detailed
analysis is beyond the scope of this paper and
will be presented in our forthcoming work).

In order to define the transposition solution
to \eqref{bsystem1}, we introduce the following
(forward) stochastic evolution equation:
\begin{eqnarray}\label{fsystem2}
\left\{
\begin{array}{lll}\ds
dz = (Az + v_1)ds +  v_2 dw(s) &\mbox{ in }(t,T],\\
\ns\ds z(t)=\eta,
\end{array}
\right.
\end{eqnarray}
where $t\in[0,T]$, $v_1\in
L^1_{\dbF}(t,T;L^{q}(\O;H))$, $v_2\in
L^2_{\dbF}(t,T;L^{q}(\O;H))$, $\eta\in
L^{q}_{\cF_t}(\O;H)$, and $q=\frac{p}{p-1}$
(See \cite[Chapter 6]{Prato} for the
well-posedness of \eqref{fsystem2} in the sense
of mild solution). We now introduce the
following notion.

\begin{definition}\label{definition1}
We call $(y(\cdot), Y(\cdot)) \in
D_{\dbF}([0,T];L^{p}(\O;H)) \times
L^2_{\dbF}(0,T;L^p(\O;H))$  a transposition
solution to (\ref{bsystem1}) if for any $t\in
[0,T]$, $v_1(\cdot)\in
L^1_{\dbF}(t,T;L^q(\O;H))$, $v_2(\cdot)\in
L^2_{\dbF}(t,T; L^q(\O;H))$, $\eta\in
 L^q_{\cF_t}(\O;H)$ and the corresponding solution $z\in C_{\dbF}([t,T];L^q(\O;H))$ to the equation (\ref{fsystem2}), it holds that
\begin{equation}\label{eq def solzz}
\begin{array}{ll}\ds
\q\dbE \big\langle z(T),y_T\big\rangle_{H}
- \dbE\int_t^T \big\langle z(s),f(s,y(s),Y(s) )\big\rangle_Hds\\
\ns\ds = \dbE \big\langle\eta,y(t)\big\rangle_H
+ \dbE\int_t^T \big\langle
v_1(s),y(s)\big\rangle_H ds + \dbE\int_t^T
\big\langle v_2(s),Y(s)\big\rangle_H ds.
\end{array}
\end{equation}
\end{definition}

On the other hand, to define the solution to
\eqref{op-bsystem3} in the transposition sense,
we need to introduce the following two
(forward) stochastic evolution equations:
\begin{equation}\label{op-fsystem2}
\left\{
\begin{array}{ll}
\ds dx_1 = (A+J)x_1ds + u_1ds + Kx_1 dw(s) + v_1 dw(s) &\mbox{ in } (t,T],\\
\ns\ds x_1(t)=\xi_1
\end{array}
\right.
\end{equation}
and
\begin{equation}\label{op-fsystem3}
\left\{
\begin{array}{ll}
\ds dx_2 = (A+J)x_2ds + u_2ds + Kx_2 dw(s) + v_2 dw(s) &\mbox{ in } (t,T],\\
\ns\ds x_2(t)=\xi_2.
\end{array}
\right.
\end{equation}
Here $\xi_1,\xi_2 \in L^4_{\cF_t}(\O;H)$,
$u_1,u_2\in L^2_\dbF(t,T;L^4(\O;H))$ and
$v_1,v_2\in L^2_\dbF(t,T;L^4(\O;H))$. Also, we
need to introduce the solution space for
\eqref{op-bsystem3}. For this purpose, write
 \bel{jshi1}
\begin{array}{ll}\ds
 D_{\dbF,w}([0,T];L^{2}(\O;\cL(H))\\
\ns\ds\= \Big\{P(\cd,\cd)\;\Big|\;
P(\cd,\cd)\in
\cL_{pd}\big(L^{2}_{\dbF}(0,T;L^{4}(\O;H)),\;L^2_{\dbF}(0,T;L^{\frac{4}{3}}(\O;H))\big),\\
\ns\ds\q \mbox{and for every } t\in[0,T]\hb{ and }\xi\in L^4_{\cF_t}(\O;H),\\
\ns\ds\q P(\cd,\cd)\xi\in
D_{\dbF}([t,T];L^{\frac{4}{3}}(\O;H)) \mbox{
and }
|P(\cd,\cd)\xi|_{D_{\dbF}([t,T];L^{\frac{4}{3}}(\O;H))}
\leq C|\xi|_{L^4_{\cF_t}(\O;H)} \Big\},
\end{array}
 \ee
and
 \bel{jshi1zx}
 L^2_{\dbF,w}(0,T;L^{2}(\O;\cL(H)))\=\cL_{pd}\big(L^{2}_{\dbF}(0,T;L^{4}(\O;H)),\;L^1_{\dbF}(0,T;L^{\frac{4}{3}}(\O;H))\big).
 \ee
We now define the transposition solution to
\eqref{op-bsystem3} as follows:

\begin{definition}\label{op-definition2}
We call
 $(P(\cd),Q(\cd))\in D_{\dbF,w}([0,T];
L^{2}(\O;\cL(H)))\times
L^2_{\dbF,w}(0,T;L^{2}(\O;\cL(H)))$ a
transposition solution to the equation
\eqref{op-bsystem3} if for any $t\in [0,T]$,
$\xi_1,\xi_2\in L^4_{\cF_t}(\O;H)$, $u_1(\cd),
u_2(\cd)\in L^2_{\dbF}(t,T;L^4(\O;H))$ and
$v_1(\cd),v_2(\cd)\in L^2_{\dbF}(t,T;
L^4(\O;H))$, it holds that
 \begin{equation}\label{eq def sol1}
\begin{array}{ll}\ds
\q\dbE \big\langle P_T
x_1(T),x_2(T)\big\rangle_{H}
 - \dbE\int_t^T \big\langle F(s)x_1(s),x_2(s)\big\rangle_{H}ds\\
\ns\ds = \dbE \big\langle
P(t)\xi_1,\xi_2\big\rangle_{H} + \dbE\int_t^T
\big\langle P(s)u_1(s),x_2(s)\big\rangle_{H} ds
+
\dbE\int_t^T \big\langle P(s)x_1(s),u_2(s)\big\rangle_{H} ds \\
\ns\ds \q + \dbE\int_t^T \big\langle P(s)
K(s)x_1(s), v_2(s)\big\rangle_{H} ds +
\dbE\int_t^T \big\langle P(s)v_1(s),
K(s)x_2(s)+v_2(s)\big\rangle_{H} ds\\
\ns\ds \q  + \dbE\int_t^T \big\langle
Q(s)v_1(s),x_2(s)\big\rangle_{H} ds +
\dbE\int_t^T \big\langle
Q(s)x_1(s),v_2(s)\big\rangle_{H} ds.
\end{array}
\end{equation}
Here, $x_1(\cd)$ and $x_2(\cd)$ solve
\eqref{op-fsystem2} and \eqref{op-fsystem3},
respectively.
\end{definition}

We shall derive the well-posedness of
(\ref{bsystem1}) in the sense of transposition
solution, by the method developed in \cite{LZ}.
Here, we face to another difficulty in the
study of (\ref{bsystem1}), i.e., $H$ is not
separable in our case. On the other hand, it
seems very difficult to establish the
well-posedness of transposition solutions to
the general equation \eqref{op-bsystem3}, and
therefore, in this paper we succeed in doing it
only for a particular case. Because of this,
instead, we introduce a weaker notion, i.e.,
relaxed transposition solution to
\eqref{op-bsystem3} (See Definition
\ref{op-definition2x} in Section \ref{s4}).
Nevertheless, it is still highly technical to
derive the well-posedness result for
\eqref{op-bsystem3} in the sense of relaxed
transposition solution. To do this, we need to
prove some weakly sequential compactness
results in the spirit of the classical
(sequential) Banach-Alaoglu theorem (also known
as Alaoglu's theorem, e.g. \cite{C}) but for
uniformly bounded linear operators in Banach
spaces. It seems that these sequential
compactness results have some independent
interest and may be applied in other places.
Once the well-posedness for both
(\ref{bsystem1}) and \eqref{op-bsystem3}, as
well as some properties of the relaxed
transposition solution to \eqref{op-bsystem3},
are established, we are able to derive the
desired Pontryagin-type stochastic maximum
principle for Problem (P).

In this paper, in order to present the key idea
in the simplest way, we do not pursue the full
technical generality. Firstly, we consider only
the simplest case of one dimensional standard
Brownian motion (with respect to the time $t$).
It would be interesting to extend the results
in this paper to the case of colored (infinite
dimensional) noise, or even with both time- and
space-dependent noise. Secondly, we impose
considerably strong regularity and boundedness
assumptions on the nonlinearities appeared in
the state equation (\ref{fsystem1}) and the
cost functional (\ref{jk1}) (See (A1) and (A2)
in the above, (A3) in Section \ref{s6}, and
(A4) in Section \ref{s7}). It would be quite
interesting to study the same problem but with
the minimal regularities and/or with unbounded
controls. Thirdly, we consider neither state
constraints nor partial observations in our
optimal control problem.

The rest of this paper is organized  as
follows. In Section \ref{s2}, we present some
preliminary results. Section \ref{s3} is
addressed to the well-posedness of  the
equation \eqref{bsystem1}. In Section
\ref{s4-1}, we study the well-posedness of the
equation \eqref{op-bsystem3} under some
additional assumptions. Section \ref{s34}
provides some sequential Banach-Alaoglu-type
theorems for uniformly bounded linear operators
between Banach spaces. In Section \ref{s4}, we
establish the well-posedness of the equation
\eqref{op-bsystem3} in the general case, while
Section \ref{as4} provides further properties
for solutions to this equation. Section
\ref{s6} gives the Pontryagin-type necessary
conditions for the optimal pair of Problem (P)
under the condition that $U$ is a convex subset
in some Hilbert space. Finally, in Section
\ref{s7}, we establish the Pontryagin-type
stochastic maximum principle for Problem (P)
for the general control domain $U$.

\section{Preliminaries}\label{s2}

In this section, we present some preliminary
results which will be used in the sequel.

First, we recall the following
Burkholder-Davis-Gundy inequality in infinite
dimensions (See \cite[Theorem 1.2.4]{LK}, for
example).
\begin{lemma} \label{BDG}
Let $f(\cd)\in L^2_\dbF(0,T;H)$. Then for any
$\a>0$, we have that
\begin{equation}\label{BDG1}
\mE\(\sup_{0\leq t\leq T}\Big| \int_0^t f(s)
dw(s)\Big|_H^\a\) \leq C \mE\( \int_0^T
|f(s)|_H^2 ds \)^{\frac{\a}{2}}.
\end{equation}
\end{lemma}

Next, for any given $r\in [1,\infty]$, we
denote by $C_{0,\dbF}^\infty((0,T);L^r(\O;H))$
the set of all $H$-valued $r$th power
integrable $\dbF$-adapted processes
$\phi(\cdot)$ such that $\phi(\cdot):(0,T) \to
L^r(\O,\cF_T,P;H)$ is an infinitely
differentiable (vector-valued) function and has
a compact support in $(0,T)$. We have the
following result.

\begin{lemma}\label{lemma1.1}
The space $C_{0,\dbF}^\infty((0,T);L^r(\O;H))$
is dense in $L^s_\dbF(0,T;L^r(\O;H))$ for any
$r\in [1,\infty]$ and $s\in [1,\infty)$.
\end{lemma}

{\it Proof}\,:  It suffices to show that for
any given $f\in L^s_\dbF(0,T;L^r(\O;H))$ and
each $\e>0$, there is a $g\in
C_{0,\dbF}^\infty((0,T);L^r(\O;H))$ such that
$|f-g|_{L^s_\dbF(0,T;L^r(\O;H))}<\e$.  Since
the set of simple processes is dense in
$L^s_\dbF(0,T;L^r(\O;H))$, we can find an
 $$\ds f_n =
\sum_{i=1}^n\chi_{[t_i,t_{i+1})}(t)x_i,
 $$
where $n\in\mathbb{N}$,
$0=t_1<t_2<\cds<t_n<t_{n+1}=T$ and $x_i\in
L^r_{\cF_{t_i}}(\O;H)$, such that
$$
|f-f_n|_{L^s_\dbF(0,T;L^r(\O;H))}<\frac{\e}{2}.
$$
On the other hand, for each
$\chi_{[t_i,t_{i+1})}$, we can find a $g_i\in
C_0^\infty(t_i,t_{i+1})$ such that
$$
|\chi_{[t_i,t_{i+1})} - g_i|_{L^s(0,T)}\leq
\frac{\e}{2n(1+|x_i|_{L^r(\O;H)})}.
$$
Write $\ds g=\sum_{i=1}^n g_i(t)x_i$. Then, it
is clear that $g\in
C_{0,\dbF}^\infty((0,T);L^r(\O;H))$. Moreover,
$$
\begin{array}{ll}\ds
|f-g|_{L^s_\dbF(0,T;L^r(\O;H))}\leq
|f-f_n|_{L^s_\dbF(0,T;L^r(\O;H))} +
|f_n-g|_{L^s_\dbF(0,T;L^r(\O;H))}\\
\ns\ds\hspace{3.5cm}<\frac{\e}{2} +
\sum_{i=1}^n |\chi_{[t_i,t_{i+1})}x_i -
g_ix_i|_{L^s_\dbF(0,T;L^r(\O;H))} <\e.
\end{array}
$$
This completes the proof of Lemma
\ref{lemma1.1}.
\endpf

 \vspace{0.2cm}

Fix any $t_1$ and $t_2$ satisfying $0\leq t_2 <
t_1 \leq T$, we recall the following known
Riesz-type Representation Theorem (See
\cite[Corollary 2.3 and  Remark 2.4]{LYZ}).

\begin{lemma}\label{lemma1}
Assume that $Y$ is a reflexive Banach space.
Then, for any $r,s\in [1,\infty)$, it holds
that
 $$\left(L^r_\dbF(t_2,t_1;L^s(\O;Y))\right)^*=L^{r'}_\dbF(t_2,t_1;L^{s'}(\O;Y^*)),
 $$
where $s'=s/(s-1)$ if $s\not=1$; $s'=\infty$ if
$s=1$; and $r'=r/(r-1)$ if $r\not=1$;
$r'=\infty$ if $r=1$.
\end{lemma}

Several more lemmas are in order.

\begin{lemma}\label{lemma2}
Let $q\geq 2$. For any $\big(v_1(\cdot),
v_2(\cdot),\eta\big)\in
L^1_{\dbF}(t,T;L^q(\O;H))\times
L^2_{\dbF}(t,T;L^q(\O;H))\times
L^q_{\cF_t}(\O;H)$, the mild solution
$z(\cdot)\in C_{\dbF}([t,T];L^{q}(\O;H))$ of
the equation \eqref{fsystem2}, given by
 \bel{z1x}
z(\cd) = S(\cd-t)\eta + \int_t^\cd
S(\cd-\si)v_1(\si) d\si + \int_t^\cd S(\cd-\si)
v_2(\si) dw(\si),
 \ee
satisfies
\begin{equation}\label{ine fs2}
\begin{array}{ll}
\ds \q |z(\cdot)|_{C_{\dbF}([t,T];L^{q}(\O;H))}\\
\ns\ds \leq C\left|\big(v_1(\cdot),
v_2(\cdot),\eta\big)\right|_{
L^1_\dbF(t,T;L^q(\O;H))\times
L^2_\dbF(t,T;L^q(\O;H))\times
L^q_{\cF_t}(\O;H)},\qq\forall\;t\in [0,T].
\end{array}
\end{equation}
\end{lemma}

\medskip

{\it Proof}\,:  By (\ref{z1x}), it is easy to
see that $z(\cd)\in
C_{\dbF}([t,T];L^{q}(\O;H))$. Also,  by Lemma
\ref{BDG} and Minkowski's inequality, we have
that
$$
\begin{array}{ll}\ds
\mE|z(s)|_H^q \!\!\!&=\ds \mE\Big| S(s-t)\eta + \int_t^s S(s-\si)v_1(\si) d\si + \int_t^s S(s-\si) v_2(\si) dw (\si)\Big|_H^q \\
\ns&\ds  \leq C\left\{\mE\Big| S(s-t)\eta\Big|_H^q + \mE\Big| \int_t^s S(s-\si)v_1(\si)d\si  \Big|_H^q +  \mE\[\int_t^s \Big|S(s-\si)v_2(\si)\Big|^2_H d\si\]^\frac{q}{2}\right\} \\
\ns&\ds  \leq C\left\{\mE\big|  \eta\big|_H^q + \mE\[\int_t^s \big|  v_1(\si)  \big|_Hd\si\]^q +   \mE\[\int_t^s\big| v_2(\si)\big|^2_H d\si\]^\frac{q}{2}\right\}\\
\ns&\ds   \leq C\[\mE\big|  \eta\big|_H^q + |v_1(\cd)|_{L^q_\dbF(\O;L^1(t,T;H))}+ |v_2(\cd)|^2_{L^q_\dbF(\O;L^2(t,T;H))}\] \\
\ns&\ds   \leq C\[\mE\big| \eta\big|_H^q +
|v_1(\cd)|_{L^1_\dbF(t,T;L^q(\O;H))}+
|v_2(\cd)|^2_{L^2_\dbF(t,T;L^q(\O;H))}\],
\end{array}
$$
which gives (\ref{ine fs2}). \endpf

\begin{lemma}\label{lemma2.1}
Assume that $p\in(1,\infty]$,
$q=\left\{\ba{ll}\frac{p}{p-1}&\hb{if
}\ p\in (1,\infty),\\[2mm] 1&\hb{if }\ p=\infty,\ea\right.$ $f_1\in
L^p_{\dbF}(0,T;L^2(\O;H))$ and $f_2\in
L^q_{\dbF}(0,T;L^2(\O;H))$. Then there exists a
monotonic sequence $\{h_n\}_{n=1}^\infty$ of
positive numbers such that
$\ds\lim_{n\to\infty}h_n=0$, and
\begin{equation}\label{2.21}
\lim_{n\to\infty}\frac{1}{h_n}\int_t^{t+h_n}\mathbb{E}
\langle f_1(t),f_2(\tau)\rangle_{H}
d\tau=\mathbb{E} \langle
f_1(t),f_2(t)\rangle_{H},\qq\ae\, t\in
 [0,T].
 \end{equation}
\end{lemma}

{\it Proof}\,: Write
$$\tilde{f_2} = \left\{\begin{array}{ll} \ds f_2,&  t\in [0,T],\\
\ns\ds 0, & t\in (T,2\,T].
\end{array}\right.$$ Obviously,
$\tilde{f_2}\in L^q_{\dbF}(0,2\,T;L^2(\O;H))$
and
$$|\tilde{f_2}|_{L^q_{\dbF}(0,2\,T;L^2(\O;H))}=|\tilde{f_2}|_{L^q_{\dbF}(0,T;L^2(\O;H))}=|f_2|_{L^q_{\dbF}(0,T;L^2(\O;H))}.$$ By Lemma \ref{lemma1.1}, for any
$\e>0$, one can find an $f_2^0 \in
C_{\dbF}([0,2\,T];L^2(\O;H))$ such that
 \bel{2-e1}
 |\tilde{f_2}-f_2^0|_{L^q_{\dbF}(0,2T;L^2(\O;H))}\le
 \e.
 \ee
By the uniform continuity of $f_2^0(\cd)$ in
$L^2(\O;H)$, one can find a $\d=\d(\e)>0$ such
that
 \bel{2-e2}
 |f_2^0(s_1)-f_2^0(s_2)|_{L^2_{\cF_T}(\O;H)}\le
 \e,\q \forall\;s_1,s_2\in [0,2\,T]\hb{ satisfying }|s_1-s_2|\le \d.
 \ee
By means of  \eqref{2-e2}, for each $h\le \d$,
we have
\begin{equation}\label{lm2.1eq1}
\begin{array}{ll}
\ds
\q\int_0^T\left|\frac{1}{h}\int_t^{t+h}\mathbb{E}\langle
f_1(t),f_2^0(\tau)\rangle_H
d\tau-\mathbb{E}\langle
f_1(t),f_2^0(t)\rangle_H\right|dt\\\ns
\ds=\frac{1}{h}\int_0^T\left|\int_t^{t+h}\mathbb{E}\langle
f_1(t),f_2^0(\tau)-f_2^0(t)\rangle_H d\tau
\right|dt\\\ns \ds\le
\frac{1}{h}\int_0^T\int_t^{t+h}|f_1(t)|_{L^2_{\cF_T}(\O;H)}|f_2^0(\tau)-f_2^0(t)|_{L^2_{\cF_T}(\O;H)}
d\tau dt\\\ns \ds\le
\frac{\e}{h}\int_0^T\int_t^{t+h}|f_1(t)|_{L^2_{\cF_T}(\O;H)}
d\tau dt=\e\int_0^T|f_1(t)|_{L^2_{\cF_T}(\O;H)}
dt\le C\e |f_1|_{ L^p_{\dbF}(0,T;L^2(\O;H))}.
\end{array}
\end{equation}
Owing to \eqref{2-e1}, we find that
 \begin{equation}\label{lm2.1eq2}
 \begin{array}{ll}
 \ds \q\int_0^T\left|\mathbb{E}\langle
 f_1(t),\tilde{f_2}(t)\rangle_H-\mathbb{E}\big(
 f_1(t),f_2^0(t)\rangle_H\right|dt\\\ns
  \ds\le |f_1|_{
L^p_{\dbF}(0,T;L^2(\O;H))}|\tilde{f_2}-f_2^0|_{L^q_{\dbF}(0,2T;L^2(\O;H))}\le
\e |f_1|_{ L^p_{\dbF}(0,T;L^2(\O;H))}.
 \end{array}
 \end{equation}
Further, utilizing \eqref{2-e1} again, we see
that
\begin{equation}\label{lm2.1eq3}
\begin{array}{ll}
\ds
\q\int_0^T\left|\frac{1}{h}\int_t^{t+h}\mathbb{E}
\langle  f_1(t),\tilde{f_2}(\tau)\rangle_H
d\tau-\frac{1}{h}\int_t^{t+h}\mathbb{E} \langle
f_1(t),f_2^0(\tau)\rangle_H d\tau\right|dt\\\ns
\ds=\frac{1}{h}\int_0^T\left|\int_t^{t+h}\mathbb{E}
\langle
f_1(t),\tilde{f_2}(\tau)-f_2^0(\tau)\rangle_H
d\tau \right|dt\\\ns \ds\le
\frac{1}{h}\int_0^T\int_t^{t+h}|f_1(t)|_{L^2_{\cF_T}(\O;H)}|\tilde{f_2}(\tau)-f_2^0(\tau)|_{L^2_{\cF_T}(\O;H)}
d\tau dt\\\ns \ds\le
\frac{1}{h}\left[\int_0^T\int_t^{t+h}|f_1(t)|_{L^2_{\cF_T}(\O;H)}^pd\tau
dt\right]^{1/p}\left[\int_0^T\int_t^{t+h}|\tilde{f_2}(\tau)-f_2^0(\tau)|_{L^2_{\cF_T}(\O;H)}^q
d\tau dt\right]^{1/q}\\\ns \ds=|f_1|_{
L^p_{\dbF}(0,T;L^2(\O;H))}\left[
\frac{1}{h}\int_0^T\int_0^h|\tilde{f_2}(t+\tau)-f_2^0(t+\tau)|
_{L^2_{\cF_T}(\O;H)}^q d\tau
dt\right]^{1/q}\\\ns \ds=|f_1|_{
L^p_{\dbF}(0,T;L^2(\O;H))}\left[
\frac{1}{h}\int_0^h\int_\tau^{T+\tau}|\tilde{f_2}(t)-f_2^0(t)|
_{L^2_{\cF_T}(\O;H)}^q
dtd\tau\right]^{1/q}\\\ns \ds\le |f_1|_{
L^p_{\dbF}(0,T;L^2(\O;H))}\left[
\frac{1}{h}\int_0^h\int_0^T|\tilde{f_2}(t)-f_2^0(t)|
_{L^2_{\cF_T}(\O;H)}^q  dtd\tau\right]^{1/q}\le
\e |f_1|_{ L^p_{\dbF}(0,T;L^2(\O;H))}.
\end{array}
\end{equation}

From \eqref{lm2.1eq1}, \eqref{lm2.1eq2} and
\eqref{lm2.1eq3}, we conclude that
 $$
\int_0^T\left|\frac{1}{h}\int_t^{t+h}\mathbb{E}
\langle f_1(t),\tilde{f_2}(\tau)\rangle_H
d\tau-\mathbb{E} \langle
f_1(t),\tilde{f_2}(t)\rangle_H\right|dt\leq C\e
|f_1|_{ L^p_{\dbF}(0,T;L^2(\O;H))}.
  $$
Therefore,
 $$
\lim_{h\to0}\int_0^T\left|\frac{1}{h}\int_t^{t+h}\mathbb{E}\langle
f_1(t),\tilde{f_2}(\tau)\rangle_H
d\tau-\mathbb{E}\langle
f_1(t),\tilde{f_2}(t)\rangle_H\right|dt=0.
  $$
This implies that there exists a monotonic
sequence $\{h_n\}_{n=1}^\infty$ of positive
numbers with $\ds\lim_{n\to\infty}h_n=0$, such
that
$$
\lim_{n\to
\infty}\frac{1}{h_n}\int_t^{t+h_n}\mathbb{E}\langle
f_1(t),\tilde{f_2}(\tau)\rangle_H
d\tau=\mathbb{E}\langle
f_1(t),\tilde{f_2}(t)\rangle_H,\qq \ae\, t\in
 [0,T].
$$
By this and the definition of
$\tilde{f_2}(\cdot)$, we conclude that
\begin{eqnarray*}
&\,&\lim_{n\to
\infty}\frac{1}{h_n}\int_t^{t+h_n}\mathbb{E}\langle
f_1(t),f_2(\tau)\rangle_H d\tau = \lim_{n\to
\infty}\frac{1}{h_n}\int_t^{t+h_n}\mathbb{E}\langle
f_1(t),\tilde{f_2}(\tau)\rangle_H d\tau
=\mathbb{E}\langle
f_1(t),\tilde{f_2}(t)\rangle_H
\\\ns&\,& \qq\qq\qq\qq\qq\qq\qq\;\,\;\;=\mathbb{E}\langle f_1(t),f_2(t)\rangle_H,\qq \ae \,t\in
[0,T].
\end{eqnarray*}
 This
completes the proof of Lemma
\ref{lemma2.1}.\endpf

\begin{lemma}\label{lemma5}
For each $t\in[0,T]$, the following three
conclusions hold:

i) If $u_2=v_2=0$ in the equation \eqref{op-fsystem3}, then there
exists an operator $U(\cd,t)\in \cL\big(L^4_{\cF_t}(\O;H),$ $
C_\dbF([t,T];L^4(\O;H))\big)$ such that the solution to
\eqref{op-fsystem3} can be represented as $x_2(\cd) =
U(\cd,t)\xi_2$. Further, for any $t\in [0,T)$, $\xi\in
L^4_{\cF_t}(\O;H)$ and $\e>0$, there is a $\d\in (0,T-t)$ such that
for any $s\in [t, t+\d]$, it holds that
\begin{equation}\label{10.31eq1}
|U(\cd,t)\xi-U(\cd,s)\xi|_{L^\infty_\dbF(s,T;L^4(\O;H))}
<\e.
\end{equation}

ii) If $\xi_2=0$ and $v_2=0$ in the equation
\eqref{op-fsystem3}, then there exists an
operator $V(\cd,t)\in
\cL\big(L^2_{\dbF}(t,T;L^4(\O;H)),\
C_\dbF([t,T];L^4(\O;H))\big)$ such that the
solution to \eqref{op-fsystem3} can be
represented as $x_2(\cd) = V(\cd,t)u_2$.

iii) If $\xi_2=0$ and $u_2=0$ in the equation
\eqref{op-fsystem3}, then there exists an
operator $W(\cd,t)\in
\cL\big(L^2_{\dbF}(t,T;L^4(\O;H)),\
C_\dbF([t,T];L^4(\O;H))\big)$ such that the
solution to \eqref{op-fsystem3} can be
represented as $x_2(\cd) = W(\cd,t)v_2$.
\end{lemma}

{\it Proof}\,: We prove only the first
conclusion. Define $U(\cd,t)$ as follows:
$$
\left\{
\begin{array}{ll}\ds
U(\cd,t):\ \ L^4_{\cF_t}(\O;H)\to C_\dbF([t,T];L^4(\O;H)),\\
\ns\ds U(s,t)\xi_2 = x_2(s),\qq\forall\; s\in
[t,T],
\end{array}
\right.
$$
where $x_2(\cd)$ is the mild solution to
\eqref{op-fsystem3} with $u_2 = v_2 = 0$.

By Lemma \ref{BDG} and H\"{o}lder's inequality,
and noting $J,K\in L^4_\dbF(0,T; L^\infty(\O;
\cL(H)))$, we obtain that
$$
\begin{array}{ll}\ds
\mE|x_2(s)|^4_H = \mE\Big|S(s-t) \xi_2 + \int_t^s S(s-\si)J(\si)x_2(\si)d\si + \int_t^s  S(s-\si)K(\si)x_2(\si)dw(\si) \Big|_H^4 \\
\ns\ds \hspace{1.7cm} \leq C \left\{\mE\Big|S(s-t) \xi_2 \Big|_H^4 + \mE\Big| \int_t^s S(s-\si)J(\si)x_2(\si)d\si  \Big|_H^4\right. \\
\ns\ds \hspace{1.7cm} \q+ \left.\mE \[\int_t^s  \Big|S(s-\si)K(\si)x_2(\si) \Big|_H^2d\si\]^2 \right\}\\
\ns\ds \hspace{1.7cm} \leq C \[\mE
\big|\xi_2\big|_H^4 +   \int_t^s
\(\big|J(\si)\big|_{L^\infty(\O;
\cL(H))}^4+\big|K(\si)\big|_{L^\infty(\O;
\cL(H))}^4\)\mE\big|x_2(\si)\big|_H^4d\si
\], \q\forall\;s\in [t,T].
\end{array}
$$
This, together with Gronwall's inequality,
implies that
$$
|x_2(s)|_{C_\dbF([t,T];L^4(\O;H))} \leq
C|\xi_2|_{L^4_{\cF_t}(\O;H)}.
$$
Hence, $U(\cd,t)$ is a bounded linear operator
from $L^4_{\cF_t}(\O;H)$ to $
C_\dbF([t,T];L^4(\O;H))$ and $U(\cd,t)\xi_2$
solves the equation \eqref{op-fsystem3} with
$u_2 = v_2 = 0$.

On the other hand, from the definition of $U(\cd,t)$ and $U(\cd,s)$,
for each $r\in [s,T]$, we see that
$$
U(r,t)\xi=S(r-t)\xi + \int_{t}^r
S(r-\tau)J(\tau)U(\tau,t)\xi d\tau + \int_{t}^r
S(r-\tau)K(\tau)U(\tau,t)\xi dw,
$$
and
$$
U(r,s)\xi=S(r-s)\xi + \int_{s}^r S(r-\tau)J(\tau)U(\tau,s)\xi d\tau
+ \int_{s}^r S(r-\tau)K(\tau)U(\tau,s)\xi dw.
$$
Hence,
$$
\begin{array}{ll}\ds
\q\mE|U(r,s)\xi-U(r,t)\xi|_H^4 \\
\ns\ds \leq C\mE\Big| S(r-s)\xi - S(r-t)\xi
\Big|_H^4 + C\mE\Big|\int_{s}^r
S(r-\tau)J(\tau)\big[
U(\tau,s)\xi-U(\tau,t)\xi \big]ds\Big|_H^4 +\\
\ns\ds \q + C\mE\Big|\int_{s}^r
S(r-\tau)K(\tau)\big[ U(\tau,s)\xi-U(\tau,t)\xi
\big]dw\Big|_H^4 + C\mE\Big|\int_{t}^{s}
S(r-\tau)J(\tau)
U(\tau,t)\xi d\tau\Big|_H^4 \\
\ns\ds \q + C\mE\Big|\int_{t}^{s} S(s-\tau)K(\tau)
U(\tau,t)\xi dw\Big|_H^4\\
\ns\ds \leq C\mE\Big| S(r-s)\xi - S(r-t)\xi \Big|_H^4\\
\ns\ds \q + C\int_s^r\(\big|J(\tau)\big|_{L^\infty(\O;
\cL(H))}^4+\big|K(\tau)\big|_{L^\infty(\O;
\cL(H))}^4\)\mE\Big|U(\tau,s)\xi-U(\tau,t)\xi\Big|_H^4d\tau
 \\
\ns\ds \q + C\int_{t}^s\(\big|J(\tau)\big|_{L^\infty(\O;
\cL(H))}^4+\big|K(\tau)\big|_{L^\infty(\O;
\cL(H))}^4\)\mE\Big|U(\tau,t)\xi\Big|_H^4d\tau\\
\ns\ds \leq C\int_s^r\(\big|J(\tau)\big|_{L^\infty(\O;
\cL(H))}^4+\big|K(\tau)\big|_{L^\infty(\O;
\cL(H))}^4\)\mE\Big|U(\tau,s)\xi-U(\tau,t)\xi\Big|_H^4d\tau
 \\
\ns\ds \q + C\mE\Big| S(r-s)\xi - S(r-t)\xi \Big|_H^4+
C\int_{t}^s\(\big|J(\tau)\big|_{L^\infty(\O;
\cL(H))}^4+\big|K(\tau)\big|_{L^\infty(\O;
\cL(H))}^4\)d\tau\mE\big|\xi\big|_H^4.
\end{array}
$$
Then, by Gronwall's inequality, we find that
$$
\mE\Big|U(r,s)\xi-U(r,t)\xi\Big|_H^4 \le
C\Big[h(r,s,t)+\int_s^rh(\si,s,t)d\si\Big],
$$
where
 $$
 h(r,s,t)= \mE\Big| S(r-s)\xi -
S(r-t)\xi \Big|_H^4+ \int_{t}^s\(\big|J(\tau)\big|_{L^\infty(\O;
\cL(H))}^4+\big|K(\tau)\big|_{L^\infty(\O;
\cL(H))}^4\)d\tau\mE\big|\xi\big|_H^4.
 $$
  Further, it is easy to see that
$$
\Big| \xi - S(s-t)\xi \Big|_H^4 \leq
C|\xi|_H^4.
$$
By Lebesgue's dominated convergence theorem, we have
$$
\lim_{s\to t+0}\mE\Big|\xi - S(s-t)\xi \Big|_H^4 =0.
$$
Hence, there is a $\d\in (0,T-t)$ such that  \eqref{10.31eq1} holds
for any $s\in [t, t+\d]$. This completes the proof of Lemma
\ref{lemma5}.
\endpf

\ms

For any $t\in [0,T]$ and $\l\in\rho(A)$,
consider the following two forward stochastic
differential equations:
\begin{equation}\label{op-fsystem2.1}
\left\{
\begin{array}{ll}
\ds dx_1 ^\l = (A_\l+J)x_1^\l ds + u_1ds + Kx_1^\l dw(s) + v_1 dw(s) &\mbox{ in } (t,T],\\
\ns\ds x_1^\l (t)=\xi_1
\end{array}
\right.
\end{equation}
 and
\begin{equation}\label{op-fsystem3.1}
\left\{
\begin{array}{ll}
\ds dx_2^\l = (A_\l+J)x_2^\l ds + u_2ds + Kx_2^\l dw(s) + v_2 dw(s) &\mbox{ in } (t,T],\\
\ns\ds x_2^\l (t)=\xi_2.
\end{array}
\right.
\end{equation}
 Here $(\xi_1,u_1,v_1)$ (\resp
$(\xi_2,u_2,v_2)$) is the same as that in
\eqref{op-fsystem2} (\resp
\eqref{op-fsystem3}). We have the following
result:

\bl\label{lem2.7}
The solutions of (\ref{op-fsystem2.1}) and
(\ref{op-fsystem3.1}) satisfy
\begin{equation}\label{Yo1}
\left\{
\begin{array}{ll}\ds
\lim_{\l\to\infty} x_1^\l(\cd) = x_1(\cd) \mbox{ in } C_\dbF([t,T];L^4(\O;H)),\\
\ns\ds \lim_{\l\to\infty} x_2^\l(\cd) =
x_2(\cd) \mbox{ in } C_\dbF([t,T];L^4(\O;H)).
\end{array}
\right.
\end{equation}
Here $x_1(\cd)$ and $x_2(\cd)$ are solutions of
\eqref{op-fsystem2} and \eqref{op-fsystem3},
respectively.
\el

{\it Proof}\,: Clearly, for any $s\in [t,T]$,
it holds that
$$
\begin{array}{ll}\ds
\q\mE|x_1(s)-x_1^\l(s)|_{H}^4 \\
\ns\ds  = \mE\Big| \[S(s-t) - S_\l(s-t)\]\xi_1
+ \int_t^s
\[S(s-\si)J(\si)x_1(\si)-
S_\l(s-\si)J(\si)x_1^\l(\si)\]d\si \\
\ns\ds\q + \int_t^s
\[S(s-\si)-S_\l(s-\si)\]u_1(\si)d\si + \int_t^s \[ S(s-\si)K(\si)x_1(\si) -
S_\l(s-\si)K(\si)x_1^\l(\si)\]dw(\si)\\
\ns\ds \q + \int_t^s \[
S(s-\si)-S_\l(s-\si)\]v_1(\si)dw
(\si)\Big|_H^4.
\end{array}
$$
Since $A_\l$ is the Yosida approximation of
$A$, one can find a positive constant
$C=C(A,T)$, independent of $\l$, such that
 \bel{zx-s1}
 |S_\l(\cd)|_{L^\infty(0,T;\cL(H))}\leq C.
 \ee
Hence,
$$
\begin{array}{ll}\ds
\q\mE\Big|\int_t^s \[S(s-\si)J(\si)x_1(\si)-
S_\l(s-\si)J(\si)x_1^\l(\si)\]d\si\Big|_H^4 \\
\ns\ds \leq C\mE\int_t^s \Big| \[S(s-\si)-
S_\l(s-\si)\]J(\si)x_1(\si)\Big|_H^4 d\si +
C\mE\int_t^s \Big| S_\l(s-\si)J(\si)\big[
x_1(\si) -
x_1^\l(\si)\big] \Big|_H^4d\si \\
\ns\ds \leq C\mE\Big|\int_t^s \[S(s-\si)-
S_\l(s-\si)\]J(\si)x_1(\si)d\si\Big|_H^4 +
C\mE\int_t^s\big|J(\si)\big|_{L^\infty(\O;
\cL(H))}^4 \big|x_1(\si)- x_1^\l(\si)\big|_H^4
d\si.
\end{array}
$$
It follows from Lemma \ref{BDG} that
$$
\begin{array}{ll}\ds
\q\mE \Big| \int_t^s \[ S(s-\si)K(\si)x_1(\si)
-
S_\l(s-\si)K(\si)x_1^\l(\si)\]dw(\si)  \Big|_H^4 \\
\ns\ds \leq C\mE\int_t^s \Big| \[S(s-\si)-
S_\l(s-\si)\]K(\si)x_1(\si)\Big|_H^4 d\si +
C\mE\int_t^s \Big| S_\l(s-\si)K(\si)\big[
x_1(\si) -
x_1^\l(\si)\big] \Big|_H^4d\si\\
\ns\ds \leq C\mE\int_t^s \Big| \[S(s-\si)-
S_\l(s-\si)\]K(\si)x_1(\si)\Big|_H^4 d\si +
C\mE\int_t^s\big|K(\si)\big|_{L^\infty(\O;
\cL(H))}^4 \big|x_1(\si) - x_1^\l(\si)
\big|_H^4d\si.
\end{array}
$$
Hence, for $t\le s\le T$,
$$
\mE\big|x_1(s)-x_1^\l(s)\big|_{H}^4 \leq
\L(\l,s) +
C\mE\int_t^s\(\big|J(\si)\big|_{L^\infty(\O;
\cL(H))}^4+\big|K(\si)\big|_{L^\infty(\O;
\cL(H))}^4\) \big|x_1(\si) - x_1^\l(\si)
\big|_H^4d\si.
$$
Here
$$
\begin{array}{ll}\ds
\L(\l,s) = C\mE\Big| \[S(s-t) -
S_\l(s-t)\]\xi_1\Big|_H^4+C\mE\Big|\int_t^s
\[S(s-\si)-S_\l(s-\si)\]u_1(\si)d\si\Big|_H^4 \\
\ns\ds \hspace{1.5cm} + C\mE\int_t^s \Big|\[
S(s-\si)-S_\l(s-\si)\]v_1(\si)\Big|_H^4d\si\\
\ns\ds \hspace{1.5cm} + C\mE\Big|\int_t^s
\[S(s-\si)- S_\l(s-\si)\]J(\si)x_1(\si)d\si\Big|_H^4\\
\ns\ds \hspace{1.5cm} + C\mE\int_t^s \Big|
\[S(s-\si)-
S_\l(s-\si)\]K(\si)x_1(\si)\Big|_H^4 d\si.
\end{array}
$$
By Gronwall's inequality, it follows that
$$
\mE\big|x_1(s)-x_1^\l(s)\big|_{H}^4 \leq
\L(\l,s)+C\int_t^se^{C(s-\tau)}\L(\l,\tau)d\tau,\qq
t\le s\le T.
$$
Since $A_\l$ is the Yosida approximation of
$A$, we see that
$\ds\lim_{\l\to\infty}\L(\l,s)=0$, which
implies that
$$
\lim_{\l\to\infty}\big|x_1^\l(\cd)-x_1(\cd)\big|_{C_\dbF([t,T]:L^4(\O;H))}=0.
$$
This leads to the first equality in
\eqref{Yo1}. The second equality in \eqref{Yo1}
can be proved similarly. This completes the
proof of Lemma \ref{lem2.7}.\endpf

\begin{lemma}\label{lemma8}
Let $H$ be a separable Hilbert space. Then, for
any $r\ge 1$, $\xi\in L^r_{\cF_T}(\O;H)$ and
$t\in [0,T)$, it holds that
\begin{equation}\label{lemma8eq1}
\lim_{s\to t^+} \big|\mE(\xi\,|\,\cF_s) -
\mE(\xi\,|\,\cF_t)\big|_{L^r_{\cF_T}(\O;H)} =
0.
\end{equation}
\end{lemma}

{\it Proof}\,:   Assume that
$\ds\xi=\sum_{i=1}^\infty\xi_i e_i$, where
$\{e_i\}_{i=1}^\infty$ is an orthonormal basis
of $H$. It is easy to see that
 $$
 \mE\Big(\sum_{i=1}^\infty |\xi_i|^2\Big)^{r/2} < \infty.
 $$
Hence, for any $\e>0$, there exists a $N>0$
such that $\ds \mE\Big(\sum_{i=N+1}^\infty
|\xi_i|^2\Big)^{r/2} < \frac{\e^r}{3^r}$.
Obviously,
$\ds\mE(\xi\;|\;\cF_t)=\sum_{i=1}^\infty
\mE(\xi_i\;|\;\cF_t)e_i$ for any $t\in [0,T)$.
By
 $$
 \ba{ll}\ds
 \Big(\sum_{i=N+1}^\infty
|\mE(\xi_i\;|\;\cF_{t})|^2\Big)^{r/2}=\Big|\mE
\Big(\sum_{i=N+1}^\infty
\xi_ie_i\;\Big|\;\cF_{t}\Big)\Big|_H^r\\\ns\ds\le
\mE \Big(\Big|\sum_{i=N+1}^\infty
\xi_ie_i\Big|_H^r\;\Big|\;\cF_{t}\Big)=\mE
\Big(\Big(\sum_{i=N+1}^\infty
|\xi_i|^2\Big)^{r/2}\;\Big|\;\cF_{t}\Big),
\q\hbox{a.s.},
 \ea
 $$
we see that
 $$
  \mE\Big(\sum_{i=N+1}^\infty
|\mE(\xi_i\;|\;\cF_{t})|^2\Big)^{r/2}\le
\mE\Big(\sum_{i=N+1}^\infty
|\xi_i|^2\Big)^{r/2} <
\frac{\e^r}{3^r},\qq\forall\;t\in [0,T].
 $$
On the other hand, since
$\big\{\mE(\xi_i\,|\,\cF_t)\big\}_{t\in[0,T]}$
is an $H$-valued
$\{\cF_t\}_{t\in[0,T]}$-martingale for each
$i\in\dbN$, we conclude that there is an
$H$-valued c\'{a}dl\'{a}g process
$\big\{x_i(t)\big\}_{t\in[0,T]}$ such that
$x_i(t)=\mE(\xi_i\,|\,\cF_t)$, $P$-a.s. Now,
for each $i\in \{1,2,\cdots,N\}$, by the fact
that the family
$\{\mE(\xi_i\,|\,\cF_t)\}_{t\in[0,T]}$ is
uniformly $r$th power integrable, we can find a
$\d>0$ such that for any $t\leq s\leq t+\d$, it
holds that $\ds \mE|x_i(t)-x_i(s)|^r <
\frac{\e^r}{3^rN^r}$. Therefore, for any $t\leq
s\leq t+\d$, it holds that
$$
\ba{ll}\ds \big[\mE\big|\mE(\xi\;|\;\cF_s) -
\mE(\xi\;|\;\cF_t)\big|_H^r\big]^{1/r} \\\ns\ds
\leq \Big[\mE\Big(\sum_{i=N+1}^\infty
|\mE(\xi_i\;|\;\cF_{s})|^2\Big)^{r/2}\Big]^{1/r}
+ \Big[\mE\Big(\sum_{i=N+1}^\infty
|\mE(\xi_i\;|\;\cF_{t})|^2\Big)^{r/2}\Big]^{1/r}
+ \sum_{i=1}^N\Big(\mE
|x_i(t)-x_i(s)|^r\Big)^{1/r}\\\ns\ds <\e,
 \ea
$$
which completes the proof.
\endpf

\begin{lemma}\label{lemma4}
Assume that  $H_1$ is a Hilbert space, and $U$
is a nonempty subset of $H_1$. If  $F(\cdot)\in
L_{\dbF}^2(0,T;H_1)$ and $\bar u(\cdot)\in
\cU[0,T]$ such that
\begin{equation}\label{lemma4 ine1}
\Re\dbE \int_0^T \big\langle F(t,\cd), u(t,\cd)
- \bar u(t,\cd) \big\rangle_{H_1} dt \leq 0,
\end{equation}
holds for any $u(\cdot)\in \cU[0,T]$ satisfying
$u(\cd)-\bar u(\cd)\in
L^2_\dbF(0,T;L^2(\O;H_1))$, then, for any point
$u\in U$, the following pointwise inequality
holds:
\begin{equation}\label{lemma4 ine2}
\Re\big\langle F(t,\o), u- \bar
u(t,\o)\big\rangle_{H_1} \leq 0, \,\ \ae
(t,\omega)\in [0,T]\times\Omega.
\end{equation}
\end{lemma}

\medskip

{\it Proof}\,: We use the contradiction
argument. Suppose that the inequality
\eqref{lemma4 ine2} did not hold. Then, there
would exist a $u_0\in U$ and an $\e>0$ such
that
$$
\a_\e\=\int_\O\int_0^T \chi_{\L_\e}(t,\o)dtdP
>0,
$$
where $ \L_\e \triangleq \Big\{ (t,\o)\in
[0,T]\times \O\;\Big|\; \Re\big\langle F(t,\o),
u_0 - \bar u(t,\o)\big\rangle_{H_1} \geq \e
\Big\}$, and $\chi_{\L_\e}$ is the
characteristic function of $\L_\e$. For any
$m\in \dbN$, define
 $
 \L_{\e,m}\=\L_\e\cap \big\{(t,\omega)\in [0,T]\times\O\;\big|\;|\bar u(t,\o)|_{H_1}\le m\big\}$.
It is clear that
$\ds\lim_{m\to\infty}\L_{\e,m}=\L_\e$. Hence,
there is an $m_\e\in\dbN$ such that
 $$
\int_\O\int_0^T \chi_{\L_{\e,m}}(t,\o)dtdP
>\frac{\a_\e}{2}>0, \qq \forall\;m\ge m_\e.
$$

Since $\big\langle F(\cdot), u_0 - \bar
u(\cdot)\big\rangle_{H_1}$ is $\{ {\cal F}_t
\}$-adapted, so is the process
$\chi_{\L_{\e,m}}(\cdot)$. Define
$$
\hat u_{\e,m}(t,\o) = u_0
\chi_{\L_{\e,m}}(t,\o)+ \bar
u(t,\o)\chi_{\L_{\e,m}^c}(t,\o),\q (t,\o)\in
[0,T]\times \O.
$$
Noting that $|\bar u(\cd)|_{H_1}\le m$ on
$\L_{\e,m}$, we see that $\hat
u_{\e,m}(\cdot)\in \cU[0,T]$ and satisfies
$\hat u_{\e,m}(\cd)-\bar u(\cd)\in
L^2_\dbF(0,T;H)$. Hence, for any $m\ge m_\e$,
we obtain that
\begin{equation}
\begin{array}{ll}\ds
\Re\dbE\int_0^T \big\langle F(t), \hat
u_{\e,m}(t) -\bar u(t) \big\rangle_{H_1} dt
&=\ds \int_\O\int_0^T
\chi_{\L_{\e,m}}(t,\o)\Re\big\langle F(t,\o),
u_0 -
\bar u(t,\o) \big\rangle_{H_1}\, dtdP \nonumber\\
\ns&\ds \geq \e\int_\O\int_0^T \chi_{\L_{\e,m}}(t,\o)dtdP\nonumber\\
\ns&\ds \geq \frac{\e\a_\e}{2}
>0,\nonumber
\end{array}
\end{equation}
which contradicts \eqref{lemma4 ine1}. This
completes the proof of Lemma
\ref{lemma4}.\endpf

\section{Well-posedness of the vector-valued BSEEs}\label{s3}

This section is devoted to proving the
following result.
\begin{theorem}\label{the1}
For any $p\in (1,2]$, $y_T \in
L^p_{\cF_T}(\O;H)$, $f(\cd,\cd,\cd):[0,T]\times
H\times H \to H$ satisfying (\ref{Lm1}), the
equation \eqref{bsystem1} admits one and only
one transposition solution $(y(\cdot),
Y(\cdot)) \in D_{\dbF}([0,T];L^{p}(\O; $ $H))
\times L^2_{\dbF}(0,T;L^{p}(\O;H))$.
Furthermore,
 \begin{equation}\label{ine the1zz}
 \begin{array}{ll}\ds
 \q|(y(\cdot), Y(\cdot))|_{
D_{\dbF}([t,T];L^{p}(\O;H)) \times L^2_{\dbF}(t,T;L^{p}(\O;H))}\\
\ns\ds\leq C\left[
 |f(\cd,0,0)|_{ L^1_{\dbF}(t,T;L^p(\O;H))} +|y_T|_{
L^p_{\cF_T}(\O;H)}\right], \q\forall\;t\in
[0,T].
\end{array}
\end{equation}
\end{theorem}

{\it Proof:} We borrow some ideas from the
proof of \cite[Theorem 3.1]{LZ}. The proof is
divided into five steps.  In the first four
steps, we study \eqref{bsystem1} for a special
case, in which $f(\cd,\cd,\cd)$ is independent
of $y$ and $Y$. More precisely, for any $y_T
\in L^p_{\cF_T}(\O;H)$ and $f(\cdot)\in
L^1_{\dbF}(0,T; L^p(\O;H))$, we consider first
the following equation:
\begin{eqnarray}\label{bsystem1-zx}
\left\{
\begin{array}{lll}
\ds dy (t)= -  A^* y(t) dt + f(t)dt + Y (t)dw (t)&\mbox{ in }[0,T),\\
\ns\ds y(T) = y_T.
\end{array}
\right.
\end{eqnarray}
In the last step, we deal with \eqref{bsystem1}
for the general case by the fixed point
technique.

\ms

{\bf Step 1.} For any $t\in [0,T]$, we define a
linear functional $\ell$ (depending on $t$) on
the Banach space
$L^1_{\dbF}(t,T;L^q(\O;H))\times
L^2_{\dbF}(t,T;L^q(\O;H))\times
L^q_{\cF_t}(\O;H)$ as follows (Recall that
$q=\frac{p}{p-1}$): \bel{xx1}\ba{ll}
\ds\ell\big(v_1(\cdot), v_2(\cdot),\eta\big) =
\mathbb{E}\big\langle z(T),y_T\big\rangle_H -
\mathbb{E}\int_t^T \big\langle
z(s),f(s)\big\rangle_H ds,\\\ns\ds
\qq\qq\qq\forall\; \big(v_1(\cdot),
v_2(\cdot),\eta\big)\in
L^1_{\dbF}(t,T;L^q(\O;H))\times
L^2_{\dbF}(t,T;L^q(\O;H))\times
L^q_{\cF_t}(\O;H),
\ea\ee
where $z(\cdot)\in C_{\dbF}([t,T];L^{q}(\O;H))$
solves the equation (\ref{fsystem2}).

By means of the H\"older inequality and Lemma
\ref{lemma2},  it is easy to show that
\begin{equation}\label{bound 1}
\begin{array}{ll}\ds
 \q\left|\ell\big(v_1(\cdot),
v_2(\cdot),\eta\big)\right| \\
\ns\ds \leq |z(T)|_{ L^q_{\cF_T}(\O;H)}|y_T|_{
L^p_{\cF_T}(\O;H)} +
|z(\cdot)|_{C_{\dbF}([t,T];L^{q}(\O;H))}
|f|_{L^1_{\dbF}(t,T;L^p(\O;H))} \\
\ns\ds\leq C\left[
 |f(\cdot)|_{L^1_{\dbF}(t,T;L^p(\O;H))}+|y_T|_{
L^p_{\cF_T}(\O;H)}\right] \\
 \ns\ds \qq\times \left|\big(v_1(\cdot),
v_2(\cdot),\eta\big)\right|_{
L^1_{\dbF}(t,T;L^q(\O;H))\times
L^2_{\dbF}(t,T;L^q(\O;H))\times
L^q_{\cF_t}(\O;H)}, \q\forall\;t\in [0,T],
\end{array}
\end{equation}
where the positive constant $C=C(T,A)$ is
independent of $t$. From (\ref{bound 1}), it
follows that $\ell$ is a bounded linear
functional on $L^1_{\dbF}(t,T;L^q(\O;H))\times
L^2_{\dbF}(t,T;L^q(\O;H))\times
L^q_{\cF_t}(\O;H)$. By Lemma \ref{lemma1},
there exist $y^t(\cdot)\in
L^\infty_{\dbF}(t,T;L^p(\O;H))$, $Y^t(\cdot)
\in L^2_{\dbF}(t,T;L^p(\O;H))$ and $\xi^t \in
L^p_{\cF_t}(\O;H)$ such that
\begin{eqnarray}\label{th1eq1}
&\,&\q\mathbb{E}\big\langle
z(T),y_T\big\rangle_H -
\mathbb{E}\int_t^T \big\langle z(\tau),f(\tau)\big\rangle_H \,d\tau\nonumber\\
&\,&
 =  \mathbb{E}\int_t^T
\big\langle v_1(\tau),y^t(\tau)\big\rangle_H
\,d\tau + \mathbb{E} \int_t^T\big\langle
v_2(\tau),Y^t(\tau)\big\rangle_H \,d\tau
+\mathbb{E} \big\langle\eta,\xi^t\big\rangle_H.
\end{eqnarray}
It is clear that $\xi^T=y_T$. Furthermore,
there is a positive constant $C=C(T,A)$,
independent of $t$, such that
\begin{eqnarray}\label{th1eq2}
&\,&\q |(y^t(\cdot), Y^t(\cdot),\xi^t)|_{
 L^\infty_{\dbF}(t,T;L^p(\O;H)) \times L^2_{\dbF}(t,T;L^p(\O;H))\times L^p_{\cF_t}(\O;H)}\nonumber\\
&\,& \leq C\left[
 |f(\cdot)|_{ L^1_{\dbF}(t,T;L^p(\O;H))}+|y_T|_{
L^p_{\cF_T}(\O;H)}\right], \qq\forall\;t\in
[0,T].
\end{eqnarray}

\ms

{\bf Step 2.} Note that the function
$(y^t(\cdot), Y^t(\cdot))$ obtained in Step 1
may depend on $t$. In this step, we show the
time consistency of $(y^t(\cdot), Y^t(\cdot))$,
that is, for any $t_1$ and $t_2$ satisfying
$0\leq t_2 \leq t_1 \leq T$, it holds that
\begin{equation}\label{th1eq3}
\big(y^{t_2} (\tau,\o),Y^{t_2}
(\tau,\o)\big)=\big( y^{t_1}(\tau,\o),
Y^{t_1}(\tau,\o)\big),\qq \ae (\tau,\o) \in
[t_1,T]\times\O.
\end{equation}

Since the solution $z(\cdot)$ of
(\ref{fsystem2}) depends on $t$, we also denote
it by $z^t(\cdot)$ whenever there exists a
possible confusion. To show (\ref{th1eq3}), we
fix arbitrarily $\varrho (\cdot)\in
L^1_{\dbF}(t_1,T;L^q(\O;H))$ and $\varsigma
(\cdot)\in L^2_{\dbF}(t_1,T;L^q(\O;H))$, and
choose first $t=t_1$, $\eta = 0$,
$v_1(\cdot)=\varrho (\cdot)$ and $v_2(\cdot) =
\varsigma (\cdot)$ in (\ref{fsystem2}). From
(\ref{th1eq1}), we obtain that
\begin{equation}\label{th1eq4}
 \ba{ll}\ds
 \q\mathbb{E}\big\langle z^{t_1}(T),y_T\big\rangle_H -
\mathbb{E}\int_{t_1}^T \big\langle
z^{t_1}(\tau),f(\tau)\big\rangle_H
d\tau\\\ns\ds
 =  \mathbb{E}\int_{t_1}^T
\big\langle\varrho
(\tau),y^{t_1}(\tau)\big\rangle_H
d\tau+\mathbb{E}\int_{t_1}^T
\big\langle\varsigma
(\tau),Y^{t_1}(\tau)\big\rangle_H d\tau.
 \ea
\end{equation}
Then, choose $t=t_2$, $\eta = 0$,
$v_1(t,\omega) = \chi_{[t_1,T]}(t) \varrho
(t,\omega)$ and $v_2(t,\omega) =
\chi_{[t_1,T]}(t) \varsigma (t,\omega)$ in
(\ref{fsystem2}). It is clear that
$$z^{t_2}(\cdot) =\left\{
\begin{array}{ll}\ds  z^{t_1}(\cdot),& t \in [t_1,T],\\
\ns\ds 0, & t \in [t_2,t_1).
\end{array}\right.
 $$
From the equality \eqref{th1eq1}, it follows
that
\begin{equation}\label{th1eq5}
\ba{ll}\ds \q\mathbb{E}\big\langle
z^{t_1}(T),y_T\big\rangle_H -
\mathbb{E}\int_{t_1}^T \big\langle
z^{t_1}(\tau),f(\tau)\big\rangle_H
d\tau\\\ns\ds
 =  \mathbb{E}\int_{t_1}^T
\big\langle\varrho
(\tau),y^{t_2}(\tau)\big\rangle_H
d\tau+\mathbb{E}\int_{t_1}^T
\big\langle\varsigma
(\tau),Y^{t_2}(\tau)\big\rangle_H d\tau. \ea
\end{equation}
Combining  (\ref{th1eq4}) and (\ref{th1eq5}),
we obtain that
 $$
 \ba{ll}\ds
 \mathbb{E}\int_{t_1}^T
\big\langle\varrho
(\tau),y^{t_1}(\tau)-y^{t_2}(\tau)\big\rangle_H
d\tau +\mathbb{E}\int_{t_1}^T
\big\langle\varsigma
(\tau),Y^{t_1}(\tau)-Y^{t_2}(\tau)\big\rangle_H
d\tau=0,\\\ns\ds\qq \q\forall\; \varrho
(\cdot)\in L^1_{\dbF}(t_1,T;L^q(\O;H)),\q
\varsigma (\cdot)\in
L^2_{\dbF}(t_1,T;L^q(\O;H)). \ea
 $$
This yields the desired equality
(\ref{th1eq3}).

Put
\begin{equation}\label{th1eq6}
 y(t,\o)=y^0(t,\o),\qq Y (t,\o)= Y^0(t,\o),\qq \forall\;(t,\o) \in
[0,T]\times\O.
\end{equation}
From  \eqref{th1eq3}, it follows that
\begin{equation}\label{th1eq7}
 \big(y^t (\tau,\o),Y^t
(\tau,\o)\big)=\big( y(\tau,\o),Y
(\tau,\o)\big), \qq \ae(\tau,\o) \in
[t,T]\times\O.
\end{equation}
Combining (\ref{th1eq1}) and (\ref{th1eq7}), we
end up with
 \be\label{th1eq8}
 \ba{ll}
 \ds \q\mathbb{E}\big\langle
z(T),y_T\big\rangle_H - \mathbb{E}
\big\langle\eta,\xi^t \big\rangle_H\\\ns\ds
 =\mathbb{E}\int_t^T \big\langle
z(\tau),f(\tau)\big\rangle_H d\tau+
\mathbb{E}\int_t^T \big\langle
v_1(\tau),y(\tau)\big\rangle_H d\tau
+\mathbb{E} \int_t^T \big\langle
v_2(\tau),Y(\tau)\big\rangle_H d\tau,\\\ns\ds
\qq\q\forall\; \big(v_1(\cdot),
v_2(\cdot),\eta\big)\in
L^1_{\dbF}(t,T;L^q(\O;H))\times
L^2_{\dbF}(t,T;L^q(\O;H))\times
L^q_{\cF_t}(\O;H). \ea
 \ee

\ms

{\bf Step 3.} We show in this step that $\xi^t$
has a c\`adl\`ag modification.

First of all, we claim that, for each $t\in
[0,T]$,
 \begin{equation}\label{6e1}
  \dbE\Big(S^*(T-t) y_T -
\int_t^T S^*(s-t) f(s)ds \;\Big|\;\cF_t\Big) =
\xi^t,\ \ \dbP\mbox{-}\as
 \end{equation}
To prove this, we note that for any $\eta\in
L^q_{\cF_t}(\O;H)$, $v_1=0$ and $ v_2=0$, the
corresponding solution to \eqref{fsystem2} is
given by $z(s)=S(s-t)\eta$ for $s\in [t,T]$.
Hence, by (\ref{th1eq8}), we obtain that
 \begin{equation}\label{M2}
\dbE\big\langle S(T-t)\eta,y_T\big\rangle_H -
\dbE \langle \eta,\xi^t\rangle_H = \dbE\int_t^T
\big\langle S(s-t)\eta,f(s)\big\rangle_H ds.
 \ee
Noting that
$$
\dbE\big\langle
S(T-t)\eta,y_T\big\rangle_H=\dbE\big\langle
\eta,S^*(T-t)y_T\big\rangle_H=\dbE
\big\langle\eta,\dbE(S^*(T-t)y_T\;|\;\cF_t)\big\rangle_H
$$
and
$$
\dbE\int_t^T \big\langle
S(s-t)\eta,f(s)\big\rangle_H ds=\dbE\Big\langle
\eta,\int_t^TS^*(s-t) f(s)ds\Big\rangle_H =\dbE
\Big\langle \eta,\dbE\Big(\int_t^T S^*(s-t)
f(s)ds\;\Big|\;\cF_t\Big) \Big\rangle_H,
$$
by (\ref{M2}), we conclude that
\begin{equation}\label{M1}
\dbE
\Big\langle\eta,\dbE\Big(S^*(T-t)y_T-\int_t^T
S^*(s-t) f(s)ds\;\Big|\;\cF_t\Big)-\xi^t
\Big\rangle_H=0,\qq\forall\;\eta\in
L^q_{\cF_t}(\O;H).
\end{equation}
Clearly, (\ref{6e1}) follows from (\ref{M1})
immediately.

In the rest of this step, we show that the
process
$$\left\{\dbE\Big(S^*(T-t) y_T - \int_t^T S^*(s-t) f(s)ds
\;\Big|\;\cF_t\Big)\right\}_{t\in[0,T]}$$ has a
c\`adl\`ag modification. Unlike the case that
$H$ is a finite dimensional space, the proof of
this fact (in the infinite dimensional space)
is quite technical.

Noting that $H$ is not assumed to be separable
(in this section), we are going to construct a
separable subspace of $H$ as our working space.
For this purpose, noting that the set of simple
functions is dense in $L_{\cF_T}^p(\O;H)$, we
conclude that there exists a sequence
$\{y^m\}_{m=1}^\infty \subset
L_{\cF_T}^2(\O;H)$ satisfying the following two
conditions:

1)  \ $\ds
y^m=\sum_{k=1}^{N_m}\a_k^m\chi_{\O_k^m}(\o)$,
where $N_m\in\dbN$, $\a_k^m\in H$ and
$\O_k^m\in \cF_T$ with $\{\O_k^m\}_{k=1}^{N_m}$
to be a partition of $\O$; and

2) \ $\ds\lim_{m\to \infty}|y^m -
y_T|_{L_{\cF_T}^p(\O;H)}=0$.

\no Likewise, since the set of simple adapted
processes is dense in
$L^p_\dbF(\O;L^1(0,T;H))$, there exists a
sequence $\{f^m\}_{m=1}^\infty \!\subset
L^1_{\dbF}(0,T;L^p(\O;H))$ satisfying the
following two conditions:

i) \ $\ds
f^m=\sum_{j=1}^{L_m}\sum_{k=1}^{M_j^m}\a_{j,k}^m\chi_{\O_{j,k}^m}(\o)\chi_{[t_j^m,t_{j+1}^m)}(t)$,
where $L_m\in\dbN$, $M_j^m\in\dbN$,
$\a_{j,k}^m\in H$, $\O_{j,k}^m\in \cF_{t_j^m}$
with $\{\O_{j,k}^m\}_{k=1}^{M_j^m}$ being a
partition of $\O$, and $0=t_1^m< t_2^m \cds <
t_{J_m}^m < t_{J_m+1}^m=T $; and

\vspace{0.3cm}

ii) \ $\ds\lim_{m\to \infty}|f^m -
f|_{L^1_{\dbF}(0,T;L^p(\O;H))}=0$.

\vspace{0.3cm}

\no Denote by $\Xi$ the set of all the above
elements $\a_k^m$ ($k=1,2,\cdots, N_m;\
m=1,2,\cdots$) and $\a_{j,k}^m$ ($k=1,2,\dots,
M_j^m; \ j=1,2,\cdots, L_m; \ m=1,2,\cdots$) in
$H$, and by $\widetilde H$ the closure of
$\span \Xi$ under the topology of $H$. Clearly,
$\widetilde H$ is a separable closed subspace
of $H$, and hence, $\widetilde H$ itself is
also a Hilbert space.

Recall that for any $\l\in\rho(A)$, the bounded
operator $A_\l$ (resp. $A^*_\l$) generates a
$C_0$-group $\{S_\l(t)\}_{t\in\dbR}$ (resp.
$\{S^*_\l(t)\}_{t\in\dbR}$) on $H$.

For each $m\in\dbN$ and $t\in [0,T]$, put
 \begin{equation}\label{6e1.1z}
 \xi_{\l,m}^t\=\dbE\Big(S^*_\l(T-t) y^m -
\int_t^T S^*_\l(s-t) f^m(s)ds
\;\Big|\;\cF_t\Big)
 \end{equation}
and
 \begin{equation}\label{X}
 X^m_\l(t)\=S^*_\l(t)\xi_{\l,m}^t - \int_0^t S^*_\l(s) f^m(s)ds.
 \end{equation}
We claim that $\{ X^m_\l(t)\}$ is an
$\widetilde H$-valued $\{\cF_t\}$-martingale.
In fact, for any $\tau_1, \tau_2 \in [0,T]$
with $\tau_1 \leq \tau_2$, by  \eqref{6e1.1z}
and (\ref{X}), it follows that
 \begin{equation}
\begin{array}{ll}\ds
\q\dbE(X^m_\l(\tau_2)\;|\;\cF_{\tau_1})\\
\ns\ds =
\dbE\Big(S^*_\l(\tau_2)\xi_{\l,m}^{\tau_2} -
\int_0^{\tau_2}S^*_\l(s)
f^m(s)ds\;\Big|\;\cF_{\tau_1}\Big) \\
\ns\ds = \dbE\left.\left[\dbE\Big(S^*_\l(T)y^m
- \int_{\tau_2}^T
S^*_\l(s)f^m(s)ds\;\Big|\;\cF_{\tau_2}\Big) -
\int_0^{\tau_2}
S^*_\l(s)f^m(s)ds\;\right|\;\cF_{\tau_1}\right]\\
\ns\ds  = \dbE\Big(S^*_\l(T)y^m-\int_0^{T}
S^*_\l(s)f^m(s)ds\;\Big|\;\cF_{\tau_1}\Big)\\
\ns\ds =
S^*_\l(\tau_1)\dbE\Big(S^*_\l(T-\tau_1)y^m-\int_{\tau_1}^{T}
S^*_\l(s-\tau_1)f^m(s)ds\;\Big|\;\cF_{\tau_1}\Big)-\int_0^{\tau_1}S^*_\l(s)f^m(s)ds\\
\ns\ds = S^*_\l(\tau_1)\xi_{\l,m}^{\tau_1}- \int_0^{\tau_1}S^*_\l(s)f^m(s)ds \\
\ns\ds = X^m_\l(\tau_1), \q \dbP\mbox{-a.s.},
\end{array}
\end{equation}
as desired.

Now, since $\{X^m_\l(t)\}_{0\leq t\leq T}$ is
an $\widetilde H$-valued $\dbF$-martingale, it
enjoys a c\`adl\`ag modification, and hence so
does the following process
 $$\{\xi_{\l,m}^t\}_{0\leq t \leq
T}=\left\{S^*_\l(-t)\left[X^m_\l(t)+ \int_0^t
S^*_\l(s) f^m(s)ds\right]\right\}_{0\leq t \leq
T}.
 $$
Here we have used the fact that
$\{S_\l^*(t)\}_{t\in\dbR}$ is a $C_0$-group on
$H$. We still use $\{\xi_{\l,m}^t\}_{0\leq t
\leq T}$ to stand for its c\`adl\`ag
modification.

From \eqref{6e1} and \eqref{6e1.1z}, it follows
that
\begin{equation}\label{th1eq3.1}
\begin{array}{ll}
\ds \q\lim_{m\to\infty}\lim_{\l\to\infty} |\xi^\cd - \xi_{\l,m}^{\cd}|_{L^\infty_{\dbF}(0,T;L^p(\O;H))} \\
\ns\ds = \lim_{m\to\infty}\lim_{\l\to\infty}
\Big| \dbE\Big(S^*(T-\cd)) y_T -
\int_\cd^T S^*(s-\cd)) f(s)ds \;\Big|\;\cF_\cd\Big) \\
\ns\ds \q - \dbE\Big(S^*_\l(T-\cd)) y^m-
\int_\cd^T S^*_\l(s-\cd) f^m(s)ds \;\Big|\;\cF_\cd\Big) \Big|_{L^\infty_{\dbF}(0,T;L^p(\O;H))}\\
\ns\ds \leq \lim_{m\to\infty}\lim_{\l\to\infty}   \Big| S^*(T-\cd) y_T - S^*_\l(T-\cd) y^m \Big|_{L^\infty_{\dbF}(0,T;L^p(\O;H))} \\
\ns\ds \q + \lim_{m\to\infty}\lim_{\l\to\infty}
\Big| \int_\cd^T S^*(s-\cd) f(s)ds - \int_\cd^T
S^*_\l(s-\cd) f^m(s)ds
\Big|_{L^\infty_{\dbF}(0,T;L^p(\O;H))}.
\end{array}
\end{equation}

Let us prove the right hand side of
\eqref{th1eq3.1} equals zero. First, we prove
 \bel{z-x1}
\lim_{m\to\infty}\lim_{\l\to\infty}   \Big|
S^*(T-\cd) y_T - S^*_\l(T-\cd) y^m
\Big|_{L^\infty_{\dbF}(0,T;L^p(\O;H))} = 0.
 \ee
Since $\{S(t)\}_{t\geq 0}$ is a
$C_0$-semigroup, for any $\e>0$, there is an
$M>0$ such that for any $m>M$, it holds that
$$
|S^*(T-\cd)y_T -
S^*(T-\cd)y^m|_{L^\infty_{\dbF}(0,T;L^p(\O;H))}
< \frac{\e}{2}.
$$
On the other hand, by the property of Yosida
approximations, we deduce that for any $\a\in
H$, it holds that
$\ds\lim_{\l\to\infty}|S^*(T-\cd)\a -
S^*_\l(T-\cd)\a|_{L^\infty(0,T;H)}=0$. Thus,
there is a $\L=\L(m)>0$ such that for any
$\l>\L$, it holds that
$$|S^*(T-\cd)\a_k^m -
S^*_\l(T-\cd)\a_k^m|_{L^\infty(0,T;H)} <
\frac{\e}{2N_m},\q k=1,2,\cdots, N_m,$$ which
implies that
$$
|S^*(T-\cdot)y^m -
S^*_\l(T-\cdot)y^m|_{L^\infty_{\dbF}(0,T;L^p(\O;H))}
\leq \sum_{k=1}^{N_m} |S^*(T-\cd)\a_k^m -
S^*_\l(T-\cd)\a_k^m|_{L^\infty(0,T;H)} <
\frac{\e}{2}.
$$
Therefore, for each $m>M$, there is a
$\L=\L(m)$ such that when $\l>\L(m)$, it holds
that
$$
\begin{array}{ll}
\ds \q \Big| S^*(T-\cd) y_T - S^*_\l(T-\cd) y^m \Big|_{L^\infty_{\dbF}(0,T;L^p(\O;H))}\\
\ns\ds \leq |S^*(T-\cd)y_T - S^*(T-\cd)y^m|_{L^\infty_{\dbF}(0,T;L^p(\O;H))} + |S^*(T-\cd)y^m - S^*_\l(T-\cd)y^m|_{L^\infty_{\dbF}(0,T;L^p(\O;H))}\\
\ns\ds < \frac{\e}{2} + \frac{\e}{2}  =\e.
\end{array}
$$
This gives (\ref{z-x1}).

Further, we show that
 \bel{z-x2}
\lim_{m\to\infty}\lim_{\l\to\infty}\Big|
\int_\cd^T S^*(s-\cd) f(s)ds - \int_\cd^T
S^*_\l(s-\cd) f^m(s)ds
\Big|_{L^\infty_{\dbF}(0,T;L^p(\O;H))}=0.
 \ee
For any $\e>0$, there is a $M^*>0$ such that
for any $m>M^*$,
$$
\Big|\int_\cd^T S^*(s-\cd)f(s)ds - \int_\cd^T
S^*(s-\cd)f^m(s)ds
\Big|_{L^\infty_{\dbF}(0,T;L^p(\O;H))} <
\frac{\e}{2}.
$$
By the property of Yosida approximations again,
we  know that for any $\a\in H$, it holds that
 $$\ds\lim_{\l\to\infty}\Big|\int_\cd^T
S^*(s-\cd)\a ds - \int_\cd^T S^*_\l(s-\cd)\a
ds\Big|_{L^\infty(0,T;H)}=0.
 $$
Thus, there is a $\L^*=\L^*(m)>0$ such that for
any $\l>\L^*$,
 $$
 \ba{ll}\ds
 \Big|\int_\cd^T S^*(s-\cd) \a_{j,k}^mds -
\int_\cd^T S^*_\l(s-\cd) \a_{j,k}^mds
\Big|_{L^\infty(0,T;H)} &\ds <
\frac{\e}{2 J_m\max(M_1^m,M_2^m,\cdots,M_{J_m}^m)},\\
 \ns \ds & \ds j=1,2,\cdots,L_m; \ k=1,2,\cdots, M_j^m.\ea$$
This implies that
$$
\begin{array}{ll}\ds
\q\Big|\int_\cd^T S^*(s-\cd)f^m(s) ds - \int_\cd^T S^*_\l(s-\cd)f^m(s)ds \Big|_{L^\infty_{\dbF}(0,T;L^p(\O;H))} \\
\ns\ds  \leq \sum_{j=1}^{L_m}\sum_{k=1}^{M_j^m}
\Big|\int_\cd^T S^*(s-\cd)\a_{j,k}^mds -
\int_\cd^T S^*_\l(s-\cd)\a_{j,k}^mds
\Big|_{L^\infty(0,T;H)} < \frac{\e}{2}.
\end{array}
$$
Therefore, for any $m>M^*$ and
$\l>\L^*=\L^*(m)$, we have
$$
\begin{array}{ll}
\ds \q\Big| \int_\cd^T S^*(s-\cd) f(s)ds - \int_\cd^T S^*_\l(s-\cd) f^m(s)ds \Big|_{L^\infty_{\dbF}(0,T;L^p(\O;H))}\\
\ns\ds \leq \Big|\int_\cd^T S^*(s-\cd)f(s)ds - \int_\cd^T S^*(s-\cd)f^m(s)ds \Big|_{L^\infty_{\dbF}(0,T;L^p(\O;H))} \\
\ns\ds \q + \Big|\int_\cd^T S^*(s-\cd)f^m(s) ds - \int_\cd^T S^*_\l(s-\cd)f^m(s)ds \Big|_{L^\infty_{\dbF}(0,T;L^p(\O;H))}\\
\ns\ds < \frac{\e}{2} + \frac{\e}{2}   =\e.
\end{array}
$$
This gives (\ref{z-x2}).

By \eqref{th1eq3.1}, (\ref{z-x1}) and
(\ref{z-x2}), we obtain that
$\ds\lim_{m\to\infty}\lim_{\l\to\infty}
|\xi^\cd -
\xi_{\l,m}^{\cd}|_{L^\infty_{\dbF}(0,T;L^p(\O;H))}=0$.
Recalling that $\xi_{\l,m}^{\cd}\in
D_{\dbF}([0,T];L^p(\O;H))$, we deduce that
$\xi^\cd$ enjoys a c\'{a}dl\'{a}g modification.

\ms

{\bf Step 4.} In this step, we show that, for
a.e. $t\in [0,T]$,
 \begin{equation}\label{6e3}
 \xi^t= y(t),\ \ \dbP\mbox{-}\as
 \end{equation}
We consider first the case that $p=2$ and fix
any $\gamma \in L^2_{\cF_{t_2}}(\O;H)$.
Choosing $t=t_2$, $v_1(\cd) = 0$, $v_2(\cd) =
0$ and $\eta = (t_1-t_2)\gamma$ in
(\ref{fsystem2}), utilizing \eqref{th1eq8}, we
obtain that
\begin{equation}\label{eq6xz}
\mathbb{E}\big\langle S(T-t_2)(t_1-t_2)\gamma,
y_T \big\rangle_H - \mathbb{E}\big\langle
(t_1-t_2)\gamma, \xi^{t_2} \big\rangle_H=
\mathbb{E}\int_{t_2}^T \big\langle S(\tau-t_2)
(t_1-t_2)\gamma, f(\tau)\big\rangle_H d\tau .
\end{equation}
Choosing $t=t_2$, $v_1(\tau,\o) =
\chi_{[t_2,t_1]}(\tau)\gamma(\o)$, $v_2(\cd) =
0$ and $\eta = 0$ in (\ref{fsystem2}),
utilizing \eqref{th1eq8} again, we find that
\begin{equation}\label{eq7zx}
 \begin{array}{ll}\ds
\q\mathbb{E}\Big\langle \int_{t_2}^T S(T-s) \chi_{[t_2,t_1]}(s)\gamma ds, y_T \Big\rangle_H \\
\ns \ds= \mathbb{E}\int_{t_2}^{t_1}\Big\langle
\int_{t_2}^\tau S(\tau-s)\gamma ds,
f(\tau)\Big\rangle_H d\tau +
\mathbb{E}\int_{t_1}^T \Big\langle S(\tau-t_1)
\int_{t_2}^{t_1}S(t_1-s)\gamma ds,
f(\tau)\Big\rangle_H d\tau \\
\ns\ds \q + \mathbb{E}\int_{t_2}^{t_1}\langle
\gamma,y(\tau)\rangle_H d\tau.
 \end{array}
\end{equation}
From (\ref{eq6xz}) and (\ref{eq7zx}), we find
\begin{equation}\label{eq7.1}
\begin{array}{ll}\ds
 \q\mathbb{E}\langle \gamma, \xi^{t_2}\rangle_H \\
 \ns\ds =
\frac{1}{t_1-t_2}\int_{t_2}^{t_1}\mathbb{E}\lan
\gamma,y(\tau) \rangle_H
d\tau \! +\! \mathbb{E}\big\langle S(T-t_2) \gamma, y_T \big\rangle_H \!-\! \frac{1}{t_1-t_2}\mathbb{E}\Big\langle \int_{t_2}^T S(T-\tau) \chi_{[t_2,t_1]}(\tau)\gamma d\tau, y_T \Big\rangle_H \\
\ns\ds \q  - \mathbb{E}\int_{t_2}^T \langle
S(\tau-t_2)
 \gamma, f(\tau)\rangle_H d\tau + \frac{1}{t_1-t_2}\mathbb{E}\int_{t_2}^{t_1}\Big\langle \int_{t_2}^\tau S(\tau-s)\gamma,
f(\tau)\Big\rangle_H d\tau\\
\ns\ds \q +
\frac{1}{t_1-t_2}\mathbb{E}\int_{t_1}^T
\Big\langle S(\tau-t_1)
\int_{t_2}^{t_1}S(t_1-s)\gamma ds,
f(\tau)\Big\rangle_H d\tau.
\end{array}
\end{equation}
Now we analyze the terms in the right hand side
of \eqref{eq7.1} one by one. First, it is easy
to show that
\begin{equation}\label{eq7.11}
 \lim_{t_1\to t_2+0} \frac{1}{t_1-t_2}\mathbb{E}\int_{t_2}^{t_1}\Big\langle \int_{t_2}^\tau S(s-t_2)\gamma,
f(\tau)\Big\rangle_H d\tau
=0,\qq\forall\;\gamma \in
L^2_{\cF_{t_2}}(\O;H).
\end{equation}
 Further,
\begin{equation}\label{eq7.13}
 \begin{array}{ll}\ds
\q\lim_{t_1\to t_2+0}   \frac{1}{t_1-t_2}\mathbb{E}\Big\langle \int_{t_2}^T S(T-\tau) \chi_{[t_2,t_1]}(\tau)\gamma d\tau, y_T \Big\rangle_H \\
\ns\ds = \lim_{t_1\to t_2+0}  \frac{1}{t_1-t_2} \mathbb{E} \Big\langle \int_{t_2}^{t_1} S(T-\tau) \gamma d\tau, y_T \Big\rangle_H \\
\ns\ds = \mathbb{E}\big\langle S(T-t_2) \gamma,
y_T \big\rangle_H.
\end{array}
\end{equation}
Utilizing the semigroup property of
$\{S(t)\}_{t\geq 0}$, we have
\begin{equation}\label{eq7.14}
\lim_{t_1\to t_2+0}
\frac{1}{t_1-t_2}\mathbb{E}\int_{t_1}^T
\Big\langle S(\tau-t_1)
\int_{t_2}^{t_1}S(t_1-s)\gamma ds,
f(\tau)\Big\rangle_H d\tau  =
\mathbb{E}\int_{t_2}^T \big\langle S(\tau-t_2)
 \gamma, f(\tau)\big\rangle_H d\tau.
\end{equation}
From \eqref{eq7.1}, \eqref{eq7.11},
\eqref{eq7.13} and \eqref{eq7.14}, we arrive at
\begin{equation}\label{eq8xz}
\lim_{t_1\to t_2+0}
\frac{1}{t_1-t_2}\int_{t_2}^{t_1}\mathbb{E}\langle
\gamma,y(\tau) \rangle_H
d\tau=\mathbb{E}\big\langle \gamma,
\xi^{t_2}\big\rangle_H,\qq\forall\;\gamma \in
L^2_{\cF_{t_2}}(\O;H), \ t_2\in [0,T).
\end{equation}

Now, by \eqref{eq8xz}, we conclude that, for
$\ae t_2\in (0,T)$
\begin{equation}\label{ezq10}
\lim_{t_1\to t_2+0}
\frac{1}{t_1-t_2}\int_{t_2}^{t_1}\mathbb{E}\big\langle
\xi^{t_2}-y(t_2),y(\tau) \big\rangle_H
d\tau=\mathbb{E}\big\langle \xi^{t_2}-y(t_2),
\xi^{t_2}\big\rangle_H.
\end{equation}
By Lemma \ref{lemma2.1}, we can find a
monotonic sequence $\{h_n\}_{n=1}^\infty$ of
positive numbers with
$\ds\lim_{n\to\infty}h_n=0$, such that
\begin{equation}\label{eq12}
\lim_{n\to\infty}
\frac{1}{h_n}\int_{t_2}^{t_2+h_n}\!\mathbb{E}\langle
\xi^{t_2}-y(t_2),y(\tau) \rangle_H
d\tau=\mathbb{E}\langle
\xi^{t_2}\!-\!y(t_2),y(t_2)\rangle_H,\q\ \ae
t_2\in [0,T).
\end{equation}
By (\ref{ezq10})--(\ref{eq12}), we arrive at
\begin{equation}\label{eq14}
\mathbb{E}\langle \xi^{t_2}-y(t_2),
\xi^{t_2}\rangle_H=\mathbb{E}\langle
\xi^{t_2}-y(t_2),y(t_2)\rangle_H,\q\ \ae t_2\in
[0,T].
\end{equation}
By (\ref{eq14}), we find that $\mathbb{E}\left|
\xi^{t_2}-y(t_2)\right|_H^2=0$ for $ t_2\in
[0,T]$ a.e., which implies (\ref{6e3}) for
$p=2$ immediately.

When $p\in (1,2]$, we choose
$\{y_T^n\}_{n=1}^\infty\subset
L^2_{\cF_T}(\O;H)$ and
$\{f_n\}_{n=1}^\infty\subset
L^1_\dbF(0,T;L^2(\O;H))$ such that
\begin{equation}\label{6e1.01}
\left\{
\begin{array}{ll}
\ds \lim_{n\to\infty} y_T^n = y_T \mbox{ in }
L^p_{\cF_T}(\O;H),\\
\ns\ds \lim_{n\to\infty} f_n = f \mbox{ in }
L^1_\dbF(0,T;L^p(\O;H)).
\end{array}
\right.
\end{equation}
We replace $y_T$ (\resp $f$) by $y_T^n$ (\resp
$f_n$) in the definition of the functional
$\ell$ (See (\ref{xx1})) and denote by
$(y_n(\cdot), Y_n(\cdot),\xi_n^t)$ the
corresponding triple satisfying (\ref{th1eq8}).
By the definition of $(y(\cd),Y(\cd),\xi^t)$
and $(y_n(\cdot), Y_n(\cdot),\xi_n^t)$, it is
easy to see that $\big(y(\cd)-y_n(\cdot),
Y(\cd) - Y_n(\cd)\big)$, $n=1,2,\cdots$,
satisfy the following:
\begin{eqnarray}\label{bsystem11}
\ba{ll}
 \ds \q\mathbb{E}\big\langle
z(T),y_T-y_T^n\big\rangle_H - \mathbb{E}
\big\langle\eta,\xi^t-\xi_n^t
\big\rangle_H\\\ns\ds
 =\mathbb{E}\int_t^T \big\langle
z(\tau),f(\tau)-f_n(\tau)\big\rangle_H d\tau+
\mathbb{E}\int_t^T \big\langle
v_1(\tau),y(\tau)-y_n(\tau)\big\rangle_H
d\tau\\\ns\ds\q +\mathbb{E} \int_t^T
\big\langle
v_2(\tau),Y(\tau)-Y_n(\tau)\big\rangle_H
d\tau,\\\ns\ds \qq\q\forall\; \big(v_1(\cdot),
v_2(\cdot),\eta\big)\in
L^1_{\dbF}(t,T;L^q(\O;H))\times
L^2_{\dbF}(t,T;L^q(\O;H))\times
L^q_{\cF_t}(\O;H). \ea
\end{eqnarray}
Hence,
\begin{equation}\label{s4ine1zz}
\begin{array}{ll}\ds
\q |y(\cd) -
y_n(\cd)|_{L^\infty_\dbF(0,T;L^p(\O;H))} +
|\xi^t -
\xi^t_n|_{L^p_{\cF_t}(\O;H)} \\
\ns\ds \leq
C\big(|f-f_n|_{L^1_\dbF(0,T;L^p(\O;H))} + |y_T
- y_T^n|_{L^p_{\cF_T}(\O;H)}\big).
\end{array}
\end{equation}
Here the constant $C$ is independent of $n$.
From the above inequality, we conclude
$$
 \lim_{n\to\infty} y_n(\cd) = y(\cd) \mbox{ in }L^\infty_\dbF(0,T;L^p(\O;H)) \mbox{ and } \lim_{n\to\infty} \xi^t_n = \xi^t \mbox{ in
} L^p_{\cF_t}(\O;H).
$$
Therefore,
$$
|y(t)-\xi^t|_{L^p_{\cF_t}(\O;H)} \leq
\lim_{n\to\infty}|y_n(t)-\xi_n^t|_{L^p_{\cF_t}(\O;H)}\leq
\lim_{n\to\infty}|y_n(t)-\xi_n^t|_{L^2_{\cF_t}(\O;H)}
= 0,\q \ae t\in [0,T],
$$
which implies \eqref{6e3} immediately.

Finally, by (\ref{6e3}) and recalling that
$\xi^t$ has a c\`adl\`ag modification, we see
that there is a c\`adl\`ag $H$-valued process
$\{\tilde y(t)\}_{t\in [0,T]}$ such that
$y(\cd)=\tilde y(\cd)$ in $[0,T]\times \Omega$
\ae It is easy to check that $(\tilde
y(\cdot),Y(\cdot))$ is a transposition solution
to the equation (\ref{bsystem1-zx}). To
simplify the notation, we still use $y$
(instead of $\tilde y$) to denote the first
component of the solution. Clearly, $(y(\cdot),
Y(\cdot)) \in D_{\dbF}([0,T];L^{p}(\O;H))
\times L^2_{\dbF}(0,T;L^{p}(\O;H))
 $ satisfies that
\begin{eqnarray}\label{s4eq1z}
&\,&\q |(y(\cdot), Y(\cdot))|_{
L^\infty_{\dbF}(t,T;L^p(\O;H)) \times L^2_{\dbF}(t,T;L^p(\O;H))}\nonumber\\
&\,& \leq C\left[
 |f(\cdot)|_{ L^1_{\dbF}(t,T;L^p(\O;H))}+|y_T|_{
L^p_{\cF_T}(\O;H)}\right], \qq\forall\;t\in
[0,T].
\end{eqnarray}
Also, the uniqueness of the transposition
solution to (\ref{bsystem1-zx}) is obvious.

\ms

{\bf Step 5.} In this step, we consider the
equation (\ref{bsystem1}) for the general case.

Fix any $T_1\in [0,T]$. For any $(\si(\cdot),
\Si(\cdot)) \in D_{\dbF}([T_1,T];L^{p}(\O;H))
\times L^2_{\dbF}(T_1,T;L^{p}(\O;H))$, we
consider the following equation:
\begin{equation}\label{s5eq1}
\left\{
\begin{array}{lll}
\ds dy_1 = - A^*y dt +  f(t,\si(t),\Si(t)) dt + Y_1dw(t) & \mbox{ in } [T_1,T),\\
\ns\ds y_1(T)=y_T.
\end{array}
\right.
\end{equation}
By the condition (\ref{Lm1}) and the result
obtained in the above, the equation
\eqref{s5eq1} admits a unique transposition
solution $(y_1(\cdot), Y_1(\cdot)) \in
D_{\dbF}([T_1,T];L^{p}(\O; H)) \times
L^2_{\dbF}(T_1,T;L^{p}(\O;H))$. This defines a
map $\cJ$ as
$$
\left\{ \!\!\begin{array}{ll}\ds \cJ: \;
D_{\dbF}([T_1,T];L^{p}(\O; H)) \times
L^2_{\dbF}(T_1,T;L^{p}(\O;H))\to
D_{\dbF}([T_1,T];L^{p}(\O; H)) \times
L^2_{\dbF}(T_1,T;L^{p}(\O;H)),\\
\ns\ds  \cJ(p(\cdot), P(\cdot)) = (y_1(\cdot),
Y_1(\cdot)).
\end{array}
\right.
$$

Now we show that the map $\cJ$ is contractive
provided that $T-T_1$ is small enough. Indeed,
for another $( \theta(\cdot), \Theta(\cdot))
\in D_{\dbF}([T_1,T];L^{p}(\O; H)) \times
L^2_{\dbF}(T_1,T;L^{p}(\O;H))$, we define
$(y_2(\cdot), Y_2(\cdot))=\cJ( \theta(\cdot),
\Theta(\cdot))$. Put
 $$
 y_3(\cd)=y_1(\cd)- y_2(\cd),\q Y_3(\cd)=Y_1(\cd)-
 Y_2(\cd),\q  f_3(\cd)=f(\cd,p(\cd),P(\cd))-f(\cd,q(\cd),Q(\cd)).
 $$
Clearly, $( y_3(\cdot),  Y_3(\cdot))$ solves
the following equation
 \begin{equation}\label{s5eq2}
\left\{
\begin{array}{lll}
\ds d y_3 = - A^*y_3 dt +  f_3(t) dt + Y_3dw(t) & \mbox{ in } [T_1,T),\\
\ns\ds y_3(T)=0.
\end{array}
\right.
\end{equation}
By the condition (\ref{Lm1}), it is easy to see
that $ f_3(\cd)\in L^1_{\dbF}(T_1,T;L^p(\O;H))$
and
 \begin{equation}\label{s5eq3}
 \begin{array}{ll}\ds
 | f_3(\cd)|_{L^1_{\dbF}(T_1,T;L^p(\O;H))}\\
  \ns\ds\le C_L\left[|\si(\cd)-\theta(\cd)|_{L^1_{\dbF}(T_1,T;L^p(\O;H))}+|\Si(\cd)-\Theta(\cd)|_{L^1_{\dbF}(T_1,T;L^p(\O;H))}\right]\\
 \ns
 \ds\le C_L\big(T-T_1+\sqrt{T-T_1}\big)\left[|\si(\cd)-\theta(\cd)|_{D_{\dbF}([T_1,T];L^{p}(\O;H)) } + |\Si(\cd)-\Theta(\cd)|_{L^2_{\dbF}(T_1,T;L^p(\O;H))}\right].
 \end{array}
 \end{equation}
By \eqref{s4eq1z}, it follows that
  \begin{equation}\label{s5eq4}
 \begin{array}{ll}\ds
\q |( y_3(\cdot), Y_3(\cdot))|_{
D_{\dbF}([T_1,T];L^{p}(\O;H)) \times
L^2_{\dbF}(T_1,T;L^p(\O;H))}\leq
C | f_3(\cd)|_{L^1_{\dbF}(T_1,T;L^p(\O;H))}\\
\ns \ds\le
C\big(T-T_1+\sqrt{T-T_1}\big)\left[|\si(\cd)-\theta(\cd)|_{D_{\dbF}([T_1,T];L^{p}(\O;H))}+|\Si(\cd)-\Theta(\cd)|_{L^2_{\dbF}(T_1,T;L^p(\O;H))}\right].
 \end{array}
 \end{equation}
Choose $T_1$ so that
$C\big(T-T_1+\sqrt{T-T_1}\big)<1$. Then, $\cJ$
is a contractive map.

By means of the Banach fixed point theorem,
$\cJ$ enjoys a unique fixed point $(y(\cdot),
Y(\cdot))\in D_{\dbF}([T_1,T];L^{p}(\O;H))
\times L^2_{\dbF}(T_1,T;L^{p}(\O;H))$. It is
clear that $(y(\cdot), Y(\cdot))$ is a
transposition solution to the following
equation:
 \begin{equation}\label{s5eq5}
\left\{
\begin{array}{lll}
\ds dy(t) = - A^*y (t)dt + f(t,y(t),Y(t)) dt + Y(t)dw(t) & \mbox{ in } [T_1,T),\\
\ns\ds y(T)=y_T.
\end{array}
\right.
\end{equation}
Using again (\ref{Lm1}) and similar to
\eqref{s5eq3}, we see that $
f(\cd,y(\cd),Y(\cd) )\in
L^1_{\dbF}(T_1,T;L^p(\O;H))$ and
  \begin{equation}\label{s5eq6}
 \begin{array}{ll}\ds
 \q |f(\cd,y(\cd),Y(\cd) )|_{ L^1_{\dbF}(T_1,T;L^p(\O;H))}\\
  \ns\ds\le |f(\cd,0,0 )|_{ L^1_{\dbF}(T_1,T;L^p(\O;H))}+C_L\left[|y(\cd)|_{L^1_{\dbF}(T_1,T;L^p(\O;H))}
  +|Y(\cd)|_{L^1_{\dbF}(T_1,T;L^p(\O;H))}\right]\\
\ns \ds\le |f(\cd,0,0 )|_{
L^1_{\dbF}(T_1,T;L^p(\O;H))}\!+\!
C_L\big(T\!\!-\!\!T_1\!+\!\sqrt{T\!-\!T_1}\big)\left[|y(\cd)|_{D_{\dbF}([T_1,T];L^{p}(\O;H))
}\!+\!|Y(\cd)|_{L^2_{\dbF}(T_1,T;L^p(\O;H))}\right].
 \end{array}
 \end{equation}
Therefore, we find that
  \begin{equation}\label{s5eq7}
 \begin{array}{ll}\ds
 \q |(y(\cdot),   Y(\cdot))|_{
D_{\dbF}([T_1,T];L^{p}(\O;H)) \times L^2_{\dbF}(T_1,T;L^p(\O;H))}\\
\ns \ds\le C\left[
 |f(\cd,y(\cd),Y(\cd) )|_{ L^1_{\dbF}(T_1,T;L^p(\O;H))}+|y_T|_{
L^p_{\cF_T}(\O;H)}\right]\\
\ns \ds\le
C\[\big(T-T_1+\sqrt{T-T_1}\big)|(y(\cdot),
Y(\cdot))|_{ D_{\dbF}([T_1,T];L^{p}(\O;H))
\times
L^2_{\dbF}(T_1,T;L^p(\O;H))}\\
\ns\ds\qq+|f(\cd,0,0 )|_{
 L^1_{\dbF}(T_1,T;L^p(\O;H))}+|y_T|_{
L^p_{\cF_T}(\O;H)}].
 \end{array}
 \end{equation}
Since $C\big(T-T_1+\sqrt{T-T_1}\big)<1$, it
follows from \eqref{s5eq7} that
\begin{equation}\label{s5eq8}
\begin{array}{ll}\ds
\3n\3n\1n\big|(y(\cdot), Y(\cdot))\big|_{
D_{\dbF}([T_1,T];L^{p}(\O;H)) \times
L^2_{\dbF}(T_1,T;L^p(\O;H))}\leq
C\big[\,|f(\cd,0,0 )|_{
L^1_{\dbF}(T_1,T;L^p(\O;H))}+|y_T|_{
L^p_{\cF_T}(\O;H)}\big].
\end{array}
\end{equation}

Repeating the above argument, we obtain the
transposition solution of the equation
\eqref{bsystem1}. The uniqueness of such
solution to \eqref{bsystem1} is obvious. The
desired estimate \eqref{ine the1zz} follows
from \eqref{s5eq8}. This completes the proof of
Theorem \ref{the1}.\endpf

\section{Well-posedness result for the operator-valued BSEEs with special data}\label{s4-1}

This section is addressed to proving a
well-posedness result for the transposition
solutions of the operator-valued BSEEs with
special data $P_T$ and $F$.

We begin with the following uniqueness result
for the transposition solution to
\eqref{op-bsystem3}.

\begin{theorem}\label{op well th4-1}
If $P_T\in L^2_{\cF_T}(\O;\cL(H))$, $F\in
L^1_\dbF(0,T;L^2(\O; \cL(H)))$ and $J,K\in
L^4_\dbF(0,T; L^\infty(\O; $ $\cL(H)))$, then
\eqref{op-bsystem3} admits at most one
transposition solution $(P(\cd),Q(\cd))\in
D_{\dbF,w}([0,T]; L^{2}(\O; $ $\cL(H)))\times
L^2_{\dbF,w}(0,T;L^{2}(\O;\cL(H)))$.
\end{theorem}

{\it Proof}\,: Assume that $(\overline
P(\cd),\overline Q(\cd))$ is another
transposition solution to the equation
\eqref{op-bsystem3}. Then, by Definition
\ref{op-definition2}, it follows that
\begin{equation}\label{ueq1as}
\begin{array}{ll}
\ds 0=\mE\Big\langle \(\overline P(t) - P(t)\)
\xi_1,\xi_2 \Big\rangle_{H} + \mE \int_t^T
\Big\langle \(\overline P(s) - P(s)\)u_1(s),
x_2(s)\Big\rangle_{H}ds \\
\ns\ds \qq + \mE \int_t^T \Big\langle
\(\overline P(s) - P(s)\)x_1(s),
u_2(s)\Big\rangle_{H}ds  + \mE \int_t^T
\Big\langle \(\overline P(s) - P(s)\) K
(s)x_1 (s), v_2 (s)\Big\rangle_{H}ds \\
\ns\ds \qq  +  \mE \int_t^T \Big\langle
\(\overline P(s) - P(s)\)v_1
(s), K (s)x_2 (s)+v_2(s)  \Big\rangle_{H}ds \\
\ns\ds \qq + \mE \int_t^T
\Big\langle\(\overline Q(s) - Q(s)\) v_1(s),
x_2(s)\Big\rangle_{H}ds  + \mE \int_t^T
\Big\langle \(\overline Q(s) - Q(s)\)
x_1(s), v_2(s) \Big\rangle_{H}ds,\\
\ns\ds\hspace{13cm}\forall\; t\in [0,T].
\end{array}
\end{equation}

Choosing $u_1=v_1=0$ and $u_2=v_2=0$ in
equations \eqref{op-fsystem2} and
\eqref{op-fsystem3}, respectively, by
\eqref{ueq1as}, we obtain that, for any $t\in
[0,T]$,
$$
0=\mE\Big\langle \(\overline P(t) - P(t)\)
\xi_1,\xi_2 \Big\rangle_{H},\qq\forall\;
\xi_1,\; \xi_2\in L^4_{\cF_t}(\O;H).
$$
Hence, we find that $\overline P(\cd)=P(\cd)$.
By this, it is easy to see that (\ref{ueq1as})
becomes the following
 \begin{equation}\label{u111ezq1as}
\begin{array}{ll}
\ds 0=\mE \int_t^T \Big\langle\(\overline Q(s)
- Q(s)\) v_1(s), x_2(s)\Big\rangle_{H}ds  + \mE
\int_t^T \Big\langle \(\overline Q(s) - Q(s)\)
x_1(s), v_2(s) \Big\rangle_{H}ds,\\
\ns\ds\hspace{13cm}\forall\; t\in [0,T].
\end{array}
\end{equation}
Choosing $t=0$, $\xi_2=0$ and $v_2=0$ in the
equation  \eqref{op-fsystem3}, we see that
(\ref{u111ezq1as}) becomes
 \begin{equation}\label{uezq1as}
0=\mE \int_0^T \Big\langle\(\overline Q(s) -
Q(s)\) v_1(s), x_2(s)\Big\rangle_{H}ds.
\end{equation}

We claim that the set
$$
\Xi\=\Big\{x_2(\cd)\;\Big|\; x_2(\cd)\hb{
solves } \eqref{op-fsystem3} \hb{ with
}t=0,\;\xi_2=0,\; v_2=0\hb{ and some }u_2\in
L^4_{\dbF}(0,T;H) \Big\}
$$
is dense in $L^4_{\dbF}(0,T;H)$. Indeed,
arguing by contradiction, if this was not true,
then there would be a nonzero $r\in
L^{\frac{4}{3}}_{\dbF}(0,T;H)$ such that
\begin{equation}\label{r}
\mE\int_0^T \big\langle r,x_2 \big\rangle_H ds
= 0,\q\mbox{ for any } x_2\in \Xi.
\end{equation}
Let us consider the following $H$-valued BSEE:
\begin{equation}\label{xz2}
\left\{
\begin{array}{ll} \ds
dy=-A^*ydt+\big(r-J(t)^*y-K(t)^*Y\big)dt+Ydw(t),\qq\hb{in
} [0,T),\\
\ns\ds y(T)=0.
\end{array}
\right.
\end{equation}
The solution to the equation (\ref{xz2}) is
understood in the transposition sense. By
Theorem \ref{the1}, the BSEE (\ref{xz2}) admits
one and only one transposition solution
$(y(\cdot), Y(\cdot)) \in
D_{\dbF}([0,T];L^{\frac{4}{3}}(\O;H)) \times
L^2_{\dbF}(0,T;L^{\frac{4}{3}}(\O;H))$. Hence,
for any $\phi_1(\cdot)\in
L^1_{\dbF}(0,T;L^4(\O;H))$ and
$\phi_2(\cdot)\in L^2_{\dbF}(0,T;L^4(\O;H))$,
it holds that
\begin{equation}\label{zx3}
\begin{array}{ll}\ds
\q - \dbE\int_0^T \big\langle z(s),r(s)-J(s)^*y(s)-K(s)^*Y(s)\big\rangle_Hds\\
\ns\ds = \dbE\int_0^T \big\langle
\phi_1(s),y(s)\big\rangle_H ds + \dbE\int_0^T
\big\langle \phi_2(s),Y(s)\big\rangle_H ds,
\end{array}
\end{equation}
where $z(\cd)$ solves
 \bel{xz4}
 \left\{
 \ba{ll}
 \ds dz=(Az+\phi_1)dt+\phi_2dw(t),\qq\hb{in } (0,T],\\
 \ns\ds
 z(0)=0.
 \ea
 \right.
 \ee
In particular, for any $x_2(\cd)$ solving
\eqref{op-fsystem3} with $t=0$, $\xi_2=0$,
$v_2=0$ and an arbitrarily given $u_2\in
L^4_{\dbF}(0,T;H)$, we choose $z=x_2$,
$\phi_1=Jx_2+u_2$ and $\phi_2=Kx_2$. By
(\ref{zx3}), it follows that
\begin{equation}\label{zx-5}
- \dbE\int_0^T \big\langle
x_2(s),r(s)\big\rangle_Hds= \dbE\int_0^T
\big\langle u_2(s),y(s)\big\rangle_H
ds,\q\forall\; u_2\in L^4_{\dbF}(0,T;H).
\end{equation}
By (\ref{zx-5}) and recalling \eqref{r}, we
conclude that $y(\cd)=0$. Hence, (\ref{zx3}) is
reduced to
 \begin{equation}\label{zx-6}
- \dbE\int_0^T \big\langle
z(s),r(s)-K(s)^*Y(s)\big\rangle_Hds =
\dbE\int_0^T \big\langle
\phi_2(s),Y(s)\big\rangle_H ds.
\end{equation}
Choosing $\phi_2(\cd)=0$ in (\ref{xz4}) and
(\ref{zx-6}), we obtain that
  \begin{equation}\label{zx-7}
 \dbE\int_0^T \Big\langle \int_0^sS(s-\si)\phi_1(\si)d\si,r(s)-K(s)^*Y(s)\Big\rangle_Hds =
0,\q\forall\; \phi_1(\cdot)\in
L^1_{\dbF}(0,T;L^4(\O;H)).
\end{equation}
Hence,
  \begin{equation}\label{zx-8}
 \int_\si^TS(s-\si)\big[r(s)-K(s)^*Y(s)\big]ds =0,\qq\forall\;\si\in [0,T].
\end{equation}
Then, for any given $\l_0\in\rho(A)$ and
$\si\in [0,T]$, we have
\begin{equation}\label{zx-9}
\begin{array}{ll}\ds
\q\int_\si^TS(s-\si)(\l_0-A)^{-1}\big[r(s)-K(s)^*Y(s)\big]ds\\
\ns\ds
=(\l_0-A)^{-1}\int_\si^TS(s-\si)\big[r(s)-K(s)^*Y(s)\big]ds
=0.
\end{array}
\end{equation}
Differentiating the equality (\ref{zx-9}) with
respect to $\si$, and noting (\ref{zx-8}), we
see that
 $$
\ba{ll}\ds
(\l_0-A)^{-1}\big[r(\si)-K(\si)^*Y(\si)\big]&\ds=-\int_\si^TS(s-\si)A(\l_0-A)^{-1}\big[r(s)-K(s)^*Y(s)\big]ds\\\ns\ds
 &\ds=\int_\si^T S(s-\si)\big[r(s)-K(s)^*Y(s)\big]ds\\\ns\ds
 &\ds\q-\l_0\int_\si^TS(s-\si)(\l_0-A)^{-1}\big[r(s)-K(s)^*Y(s)\big]ds\\\ns\ds
 &\ds=0,\qq\forall\;\si\in [0,T].
 \ea
 $$
Therefore,
 \bel{zx-10}
r(\cd)=K(\cd)^*Y(\cd).
 \ee
By (\ref{zx-10}), the equation (\ref{xz2}) is
reduced to
 \bel{xoz2}
 \left\{
 \ba{ll}
 \ds dy=-A^*ydt-J(s)^*ydt+Ydw(t),\qq\hb{in } [0,T),\\\ns\ds
 y(T)=0.
 \ea
 \right.
 \ee
It is clear that the unique transposition of
(\ref{xoz2}) is $(y(\cd),Y(\cd))=(0,0)$. Hence,
by (\ref{zx-10}), we conclude that $r(\cd)=0$,
which is a contradiction. Therefore, $\Xi$ is
dense in $L^4_{\dbF}(0,T;H)$. This, combined
with (\ref{uezq1as}), yields that
 $$
\(\overline Q(\cd) - Q(\cd)\)
v_1(\cd)=0,\qq\forall\; v_1(\cdot)\in
L^4_{\dbF}(0,T;H).
 $$
Hence $\overline Q(\cd) = Q(\cd)$. This
completes the proof of Theorem \ref{op well
th4-1}.
\endpf

\ms \ms

In the rest of this section, we assume that $H$
is a separable Hilbert space. Denote by
$\cL_2(H)$ the Hilbert space of all
Hilbert-Schmidt operators on $H$. We have the
following well-posedness result.

\begin{theorem}\label{op well th0}
If $P_T\in L^2_{\cF_T}(\O;\cL_2(H))$, $F\in
L^1_\dbF(0,T;L^2(\O;\cL_2(H)))$ and $J,K\in
L^4_\dbF(0,T; L^\infty(\O; $ $\cL(H)))$, then
the equation \eqref{op-bsystem3} admits one and
only one transposition solution $\big(P(\cd),
Q(\cd)\big)$ with the following regularity:
$$
\big(P(\cd), Q(\cd)\big) \in
D_\dbF([0,T];L^{2}(\O;\cL_2(H)))\times
L^2_{\dbF}(0,T;\cL_2(H)).
$$
Furthermore,
 \begin{equation}\label{ine the1z0}
 \begin{array}{ll}\ds
 |(P, Q)|_{
D_\dbF([0,T];L^{2}(\O;\cL_2(H)))\times
L^2_{\dbF}(0,T;\cL_2(H))}\leq C\left[
 |F|_{L^1_\dbF(0,T;L^2(\O;\cL_2(H)))} + |P_T|_{L^2_{\cF_T}(\O;\cL_2(H))}\right].
\end{array}
\end{equation}
\end{theorem}

\medskip

{\it Proof}\,: We divide the proof into several
steps.

\ms

{\bf Step 1.} Define a family of operators
$\{\cT (t)\}_{t\geq 0}$ on $\cL_2(H)$ as
follows:
$$
\cT (t)O = S(t)OS^*(t), \q\forall\; O\in
\cL_2(H).
$$
We claim that $\{\cT (t)\}_{t\geq 0}$ is a
$C_0$-semigroup on $\cL_2(H)$. Indeed, for any
nonnegative $s$ and $t$, we have $$ \cT (t+s)O
= S(t+s)OS^*(t+s) = S(t)S(s)OS^*(s)S^*(t)=\cT
(t)\cT (s)O, \q\forall\; O\in \cL_2(H).
$$
Hence, $\{\cT (t)\}_{t\geq 0}$ is a semigroup
on $\cL_2(H)$. Next, we choose an orthonormal
basis $\{e_i\}_{i=1}^\infty$ of $H$. For any
$O\in \cL_2(H)$ and $t\in [0,\infty)$,
\begin{equation}\label{s5eq.1}
\begin{array}{ll}\ds
\q\lim_{s\to t^+} \big|\cT (s)O-
\cT (t)O\big|^2_{\cL_2(H)}\\
\ns\ds \leq |\!|S(t)|\!|^2_{\cL(H)}\lim_{s\to
t^+}\big|S(s-t)OS^*(s-t)-
O\big|^2_{\cL_2(H)}|\!|S^*(t)|\!|^2_{\cL(H)}\\
\ns\ds \leq |\!|S(t)|\!|_{\cL(H)}^4 \lim_{s\to
t^+}\sum_{i=1}^\infty\big|S(s-t)OS^*(s-t)e_i-
Oe_i\big|^2_{H} \\
\ns\ds \leq 2|\!|S(t)|\!|_{\cL(H)}^4
\lim_{s\to
t^+}\sum_{i=1}^\infty\[\big|S(s-t)OS^*(s-t)e_i-
S(s-t)Oe_i\big|_{H}^2+ \big|S(s-t)Oe_i-
Oe_i\big|^2_{H} \].
\end{array}
\end{equation}
For the first series in the right hand side of
(\ref{s5eq.1}), we have
$$
\begin{array}{ll}\ds
\q \sum_{i=1}^\infty\big|S(s-t)OS^*(s-t)e_i-
S(s-t)Oe_i\big|^2_{H}\\
\ns\ds \leq C
\sum_{i=1}^\infty\big|OS^*(s-t)e_i-
Oe_i\big|^2_{H} =C
\big|OS^*(s-t)-O\big|^2_{\cL_2(H)} =C
\big|\big[OS^*(s-t)-O\big]^*\big|^2_{\cL_2(H)}\\
\ns\ds
=C\sum_{i=1}^\infty\big|\big[S(s-t)O^*-O^*\big]
e_i\big|^2_H.
\end{array}
$$
For each $i\in\dbN$,
$$
\big|\big[S(s-t)O^*-O^*\big] e_i\big|^2_H \leq
2\big| S(s-t)O^* e_i\big|^2_H + \big|O^*
e_i\big|^2_H \leq C\big|O^* e_i\big|^2_H.
$$
It is clear that
$$
\sum_{i=1}^\infty\big| O^*e_i\big|^2_{H} =
|O^*|_{\cL_2(H)}^2= |O|_{\cL_2(H)}^2.
$$
Hence, by Lebesgue's dominated convergence
theorem, it follows that
\begin{equation}\label{s5eq.2}
\begin{array}{ll}\ds
\q\lim_{s\to t^+}
\sum_{i=1}^\infty\big|S(s-t)OS^*(s-t)e_i-
S(s-t)Oe_i\big|_{H}^2\\
\ns\ds \leq C \lim_{s\to
t^+}\sum_{i=1}^\infty\big|OS^*(s-t)e_i-
Oe_i\big|^2_{H}=C \sum_{i=1}^\infty\lim_{s\to
t^+}\big|OS^*(s-t)e_i- Oe_i\big|^2_{H}=0.
\end{array}
\end{equation}
By a similar argument, it follows that
\begin{equation}\label{s5eq.3}
\lim_{s\to t^+}
\sum_{i=1}^\infty\big|S(s-t)Oe_i-
Oe_i\big|_{H}^2=0.
\end{equation}
From \eqref{s5eq.1}--\eqref{s5eq.3}, we find
that
$$
\lim_{s\to t^+} \big|\cT (s)O- \cT
(t)O\big|^2_{\cL_2(H)}=0,\q\forall\; t\in
[0,\infty)\mbox{ and } O\in\cL_2(H).
$$
Similarly,
$$
\lim_{s\to t^-} \big|\cT (s)O- \cT
(t)O\big|^2_{\cL_2(H)}=0,\q\forall\; t\in
(0,\infty)\mbox{ and } O\in\cL_2(H).
$$
Hence, $\{\cT (t)\}_{t\geq 0}$ is a
$C_0$-semigroup on $\cL_2(H)$.

\ms

{\bf Step 2.} Denote by $\cA$ the infinitesimal
generater of $\{\cT (t)\}_{t\geq 0}$. We
consider the following $\cL_2(H)$-valued BSEE:
\begin{equation}\label{qqqs5f}
\left\{
\begin{array}{ll}\ds
dP = -\cA^* Pdt + f(t,P,Q)dt + Qdw &\mbox{ in
}[0,T),\\
\ns\ds P(T)=P_T,
\end{array}
\right.
\end{equation}
where
\begin{equation}\label{s5f}
f(t,P,Q) = -J^*P-PJ-K^*PK - K^*Q - QK + F.
\end{equation}
Noting that $J,K\in L^4_\dbF(0,T;
L^\infty(\O;\cL(H)))$ and $F\in
L^1_\dbF(0,T;L^2(\O;\cL_2(H)))$, we see that
$f(\cd,\cd,\cd)$ satisfies \eqref{Lm1}. By
Theorem \ref{the1}, we conclude that there
exists a $ (P, Q) \in
D_\dbF([0,T];L^{2}(\O;\cL_2(H)))\times
L^2_{\dbF}(0,T;\cL_2(H)) $ solves
\eqref{qqqs5f} in the sense of Definition
\ref{definition1}, where the Hilbert space $H$
is replaces by $\cL_2(H)$. Further, $(P,Q)$
satisfies \eqref{ine the1z0}.

Denote by $O(\cd)$ the tensor product of
$x_1(\cd)$ and $x_2(\cd)$, where $x_1$ and
$x_2$ solve respectively \eqref{op-fsystem2}
and \eqref{op-fsystem3}. As usual, $O(t,\o)x =
\langle x,x_1 \rangle_H x_2$ for a.e.
$(t,\o)\in [0,T]\times \O$ and $x\in H$. Hence,
$O(t,\o)\in \cL_2(H)$. For any $\l\in\rho(A)$,
define a family of operators $\{\cT_\l
(t)\}_{t\geq 0}$ on $\cL_2(H)$ as follows:
$$
\cT_\l (t)O = S_\l(t)OS_\l^*(t), \q\forall\;
O\in \cL_2(H).
$$
By the result proved in Step 1, it follows that
$\{\cT_\l (t)\}_{t\geq 0}$ is a $C_0$-semigroup
on $\cL_2(H)$. Further, for any $O\in\cL_2(H)$,
we have
$$
\begin{array}{ll}\ds
\lim_{t\to 0^+}\frac{\cT_\l (t)O-O}{t}\!\!\!&
\ds = \lim_{t\to
0^+}\frac{S_\l(t)OS_\l^*(t)-O}{t} = \lim_{t\to
0^+}\frac{S_\l(t)OS_\l^*(t)-OS_\l^*(t)+OS_\l^*(t)-O}{t}\\
\ns&\ds = A_\l O +  O A_\l^*.
\end{array}
$$
Hence, the infinitesimal generater $\cA_\l$ of
$\{\cT_\l (t)\}_{t\geq 0}$ is as follows:
$$
\cA_\l O = A_\l O +  O A_\l^*,\q\mbox{ for
every }O\in \cL_2(H).
$$

Now, for any $O\in\cL_2(H)$, it holds that
$$
\begin{array}{ll}\ds
\q\lim_{\l\to\infty}\Big|\cT(t)O- \cT_\l (t)O \Big|_{\cL_2(H)}\\
\ns\ds = \lim_{\l\to\infty} \Big| S(t)OS^*(t) -
S_\l(t)OS_\l^*(t) \Big|_{\cL_2(H)}\\
\ns\ds \leq\lim_{\l\to\infty} \Big| S(t)OS^*(t)
- S(t)OS_\l^*(t)\Big|_{\cL_2(H)} +
\lim_{\l\to\infty}\Big| S(t)OS_\l^*(t) -
S_\l(t)OS_\l^*(t) \Big|_{\cL_2(H)}.
\end{array}
$$
Let us compute the value of each term in the
right hand side of the above inequality. First,
\begin{equation}\label{HS2}
\begin{array}{ll}\ds
\Big| S(t)OS^*(t) -
S(t)OS_\l^*(t)\Big|^2_{\cL_2(H)}\!\!\!&\ds\leq
C\Big| OS^*(t) -  OS_\l^*(t)\Big|^2_{\cL_2(H)}
= C\Big|
S(t)O^* -  S_\l (t)O^*\Big|^2_{\cL_2(H)}\\
\ns&\ds = C\sum_{i=1}^\infty \Big| \(S(t) -
S_\l (t)\)O^* e_i\Big|^2_{H}.
\end{array}
\end{equation}
Since
$$
\Big| \(S(t) - S_\l (t)\)O^* e_i\Big|^2_{H}\leq
C|O^* e_i|^2_{H}\q\mbox{ and }\sum_{i=1}^\infty
\Big| O^* e_i\Big|^2_{H} =
|O|^2_{\cL_2(H)}<\infty,
$$
by means of Lebesgue's dominated theorem and
\eqref{HS2}, we find that
$$
\lim_{\l\to\infty}\Big| S(t)OS^*(t) -
S(t)OS_\l^*(t)\Big|^2_{\cL_2(H)}=0.
$$
By a similar argument, we find that
$$
\lim_{\l\to\infty}\Big| S(t)OS_\l^*(t) -
S_\l(t)OS_\l^*(t) \Big|^2_{\cL_2(H)}=0.
$$
Hence,
\begin{equation}\label{HS2.1}
\lim_{\l\to\infty}\Big|\cT(t)O- \cT_\l (t)O
\Big|_{\cL_2(H)}=0, \q\mbox{ for any }t\geq 0.
\end{equation}

Write $O^\l = x_1^\l\otimes x_2^\l$, where
$x_1^\l$ and $x_2^\l$ solve accordingly
\eqref{op-fsystem2.1} and
\eqref{op-fsystem3.1}. Then,
\begin{equation}\label{HS3}
dO^\l = (A_\l x_1^\l) \otimes x_2^\l ds +
x_1^\l\otimes (A_\l x_2^\l)ds + u^{\l}ds + v^\l
dw,
\end{equation}
where
\begin{equation}\label{qz11}
\left\{
\begin{array}{ll}\ds
u^\l = \big(J x_1^\l\big)\otimes x_2^\l +
x_1^\l \otimes \big(J x_2^\l\big) + u_1\otimes
x_2^\l + x_1^\l\otimes u_2 + \big(K
x_1^\l\big)\otimes \big(K x_2^\l\big)  +
\big(K x_1^\l\big)
\otimes v_2 \\
\ns\ds \qq + v_1 \otimes \big( K x_2^\l\big) +
v_1 \otimes
v_2 , \\
\ns\ds v^\l = \big(K x_1^\l\big)\otimes x_2^\l
+ x_1^\l \otimes \big(K x_2^\l\big) + v_1
\otimes x_2^\l + x_1^\l \otimes v_2.
\end{array}
\right.
\end{equation}
Further, for any $h\in H$, we find
$$
\big((A_\l x_1^\l) \otimes
x_2^\l\big)(h)=\big\langle h,A_\l x_1^\l
\big\rangle_H x_2^\l = \big\langle A_\l^*h,
x_1^\l \big\rangle_H x_2^\l =  \big(x_1^\l
\otimes x_2^\l\big) A_\l^*h.
$$
Thus,
\begin{equation}\label{HS4}
\big(A_\l x_1^\l\big) \otimes x_2^\l = O^\l
A_\l^*.
\end{equation}
Similarly, we have the following equalities:
\begin{equation}\label{HS5}
\left\{
\begin{array}{ll}\ds
x_1^\l\otimes \big(A_\l x_2^\l\big) = A_\l O^\l,\\
\ns\ds \big(J x_1^\l\big)\otimes x_2^\l +
x_1^\l \otimes \big(J x_2^\l\big) = O^\l J^* +
J
O^\l,\\
\ns\ds \big(K x_1^\l\big)\otimes \big(K
x_2^\l\big)  = K
O^\l K^*,\\
\ns\ds \big(K x_1^\l\big) \otimes v_2   + v_1
\otimes \big( K x_2^\l\big)=\big(x_1^\l \otimes
v_2 \big)K^*   +  K^* \big(v_1 \otimes
x_2^\l\big),\\
\ns\ds \big(K x_1^\l\big)\otimes x_2^\l +
x_1^\l \otimes \big(K x_2^\l\big)= O^\l K^* + K
O^\l.
\end{array}
\right.
\end{equation}
By \eqref{qz11}--(\ref{HS5}), we obtain that
\begin{equation}\label{qz111}
\left\{
\begin{array}{ll}\ds
u^\l = J O^\l + O^\l J^* + u_1\otimes x_2^\l +
x_1^\l\otimes u_2 + K O^\l K^*  +  K x_1^\l
\otimes v_2 + v_1 \otimes K x_2^\l + v_1
\otimes
v_2 , \\
\ns\ds v^\l = K O^\l + O^\l K^* + v_1 \otimes
x_2^\l + x_1^\l \otimes v_2.
\end{array}
\right.
\end{equation}
From \eqref{HS3}, \eqref{HS4}, the first
equality in \eqref{HS5} and \eqref{qz111}, we
see that $O^\l$ solves
\begin{equation}\label{HS6}
\left\{\3n
\begin{array}{ll}
\ds dO^\l  = \cA_\l O^\l ds + u^\l ds + v^\l dw(s) &\mbox{ in } (t,T],\\
\ns\ds O^\l(t) = \xi_1\otimes \xi_2.
\end{array}
\right.
\end{equation}
Hence,
\begin{equation}\label{HS7}
O^\l(s)=\cT_{\l}(s-t)(\xi_1\otimes \xi_2) +
\int_t^s \cT_\l(\tau-t) u^{\l}(\tau)d\tau +
\int_t^s \cT_\l(\tau-t) v^{\l}(\tau)dw(\tau),
\q\forall\;s\in [t,T].
\end{equation}

We claim that
\begin{equation}\label{HS8}
\lim_{\l\to\infty}|O^\l(\cd)-O(\cd)|_{C_\dbF([t,T];L^4(\O;\cL_2(H)))}=0,\qq\forall\;t\in[0,T].
\end{equation}
Indeed, for any $s\in [t,T]$, we have
$$
\begin{array}{ll}\ds
\q|O^\l(s)-O(s)|^2_{\cL_2(H)} =
\sum_{i=1}^\infty \big|O^\l(s) e_i -
O(s) e_i\big|_H^2 \\
\ns\ds = \sum_{i=1}^\infty \big|\langle
e_i,x_1^\l(s)\rangle_H x_2^\l(s) - \langle e_i,
x_1(s)\rangle_H x_2(s)\big|_H^2 \\
\ns \ds \leq 2\sum_{i=1}^\infty \big|\langle
e_i,x_1^\l(s)\rangle_H x_2^\l(s) - \langle e_i,
x_1^\l(s)\rangle_H x_2(s)\big|_H^2 +
2\sum_{i=1}^\infty \big|\langle
e_i,x_1^\l(s)\rangle_H x_2 (s) - \langle e_i,
x_1(s)\rangle_H x_2(s)\big|_H^2\\
\ns\ds \leq 2 \big|x_2^\l(s)-x_2(s)
\big|_H^2\sum_{i=1}^\infty\big|\langle
e_i,x_1^\l(s)\rangle_H\big|^2 +
2\big|x_2(s)\big|_H^2\sum_{i=1}^\infty
\big|\langle e_i,x_1^\l(s)-x_1(s)\rangle_H
\big|^2\\
\ns\ds =2 \big|x_1(s)\big|_H^2
\big|x_2^\l(s)-x_2 (s) \big|_H^2 + 2
\big|x_2(s)\big|_H^2 \big|x_1^\l(s)-x_1 (s)
\big|_H^2.
\end{array}
$$
This, together Lemma \ref{lem2.7}, implies that
\eqref{HS8} holds.

By a similar argument and noting \eqref{HS2.1},
using Lemma \ref{BDG} and \eqref{HS8}, we can
show that, for any $t\in [0,T]$, it holds that
\begin{equation}\label{HS9}
\left\{
\begin{array}{ll}\ds
\lim_{\l\to\infty}\Big|\int_t^\cd
\cT_\l(\tau-t) u^{\l}(\tau)d\tau - \int_t^\cd
\cT(\tau-t)
u(\tau)d\tau\Big|_{C_\dbF([t,T];L^4(\O;\cL_2(H)))}=0,\\
\ns\ds \lim_{\l\to\infty}\Big|\int_t^\cd
\cT_\l(\tau-t) v^{\l}(\tau)dw(\tau) -
\int_t^\cd \cT(\tau-t)
v(\tau)dw(\tau)\Big|_{C_\dbF([t,T];L^4(\O;\cL_2(H)))}=0,
\end{array}
\right.
\end{equation}
where
 \bel{qz1}
\left\{
\begin{array}{ll}\ds
u = J O(\cd) + O(\cd)J^* + u_1\otimes x_2 +
x_1\otimes u_2 + K O(\cd) K^* + (K x_1) \otimes
v_2 + v_1 \otimes (K x_2) + v_1 \otimes
v_2 , \\
\ns\ds v = K O(\cd) + O(\cd)K^* + v_1 \otimes
x_2 + x_1 \otimes v_2.
\end{array}
\right.
 \ee
From \eqref{HS7}--\eqref{HS9}, we obtain that
\begin{equation}\label{HS10}
O(s)=\cT(s-t)(\xi_1\otimes \xi_2) + \int_t^s
\cT(\tau-t) u(\tau)d\tau + \int_t^s \cT(\tau-t)
v(\tau)dw(\tau), \q\forall\;s\in [t,T].
\end{equation}
Hence, $O(\cd)$ verifies that
\begin{equation}\label{op-fsystem4x}
\left\{\3n
\begin{array}{ll}
\ds dO(s)  = \cA O(s) ds + uds + vdw(s) &\mbox{ in } (t,T],\\
\ns\ds O(t) = \xi_1\otimes \xi_2.
\end{array}
\right.
\end{equation}

\ms

{\bf Step 3.} Since $(P,Q)$ solves
\eqref{qqqs5f} in the transposition sense and
by (\ref{op-fsystem4x}), it follows that
\begin{equation}\label{eq def solzzz}
\begin{array}{ll}\ds
\q\dbE \big\langle
O(T),P_T\big\rangle_{\cL_2(H)}
- \dbE\int_t^T \big\langle O(s),f(s,P(s),Q(s) )\big\rangle_{\cL_2(H)}ds\\
\ns\ds = \dbE
\big\langle\xi_1\otimes\xi_2,P(t)\big\rangle_{\cL_2(H)}
+ \dbE\int_t^T \big\langle
u(s),P(s)\big\rangle_{\cL_2(H)} ds +
\dbE\int_t^T \big\langle
v(s),Q(s)\big\rangle_{\cL_2(H)} ds.
\end{array}
\end{equation}
By (\ref{s5f}) and recalling that
$O(\cd)=x_1(\cd)\otimes x_2(\cd)$, we find that
\begin{equation}\label{eq def solzzz1}
\begin{array}{ll}\ds
\q\dbE\int_t^T \overline{\big\langle
O(s),f(s,P(s),Q(s)
)\big\rangle_{\cL_2(H)}}ds\\
\ns\ds =\dbE\int_t^T \big\langle
\big(-J(s)^*P(s)-P(s)J(s)-K^*(s)P(s)K(s)-K^*(s)Q(s)-Q(s)K(s)
\\
\ns\ds \qq\qq + F(s) \big)x_1(s),x_2(s)
\big\rangle_{H}ds.
\end{array}
\end{equation}
Further, by (\ref{qz1}), we have
\begin{equation}\label{eq def solzzz2}
\begin{array}{ll}\ds
\q\dbE\int_t^T \overline{\big\langle
u(s),P(s)\big\rangle_{\cL_2(H)}} ds\\
\ns\ds = \dbE\int_t^T \big\langle
J^*(s)P(s)x_1(s),x_2(s)\big\rangle_{H} ds +
\dbE\int_t^T \big\langle
P(s)J(s)x_1(s),x_2(s)\big\rangle_{H} ds
\\
\ns\ds \q +\dbE\int_t^T \big\langle
P(s)u_1(s),x_2(s)\big\rangle_{H}
ds+\dbE\int_t^T \big\langle
P(s)x_1(s),u_2(s)\big\rangle_{H} ds
\\
\ns\ds \q + \dbE\int_t^T \big\langle
K^*(s)P(s)K(s)x_1(s),x_2(s)\big\rangle_{H} ds
+\dbE\int_t^T \big\langle
P(s)K(s)x_1(s),v_2(s)\big\rangle_{H} ds + \\
\ns\ds\q+ \dbE\int_t^T \big\langle
K(s)^*P(s)v_1(s),x_2(s)\big\rangle_{H} ds +
\dbE\int_t^T \big\langle P(s)
v_1(s),v_2(s)\big\rangle_{H} ds,
\end{array}
\end{equation}
and
\begin{equation}\label{eq def solzzz3}
\begin{array}{ll}\ds
\q\dbE\int_t^T \overline{\big\langle
v(s),Q(s)\big\rangle_{\cL_2(H)}} ds\\
\ns\ds  =\dbE\int_t^T \big\langle
K^*(s)Q(s)x_1(s),x_2(s)\big\rangle_{H} ds +
\dbE\int_t^T \big\langle
Q(s)K(s)x_1(s),x_2(s)\big\rangle_{H} ds \\
\ns\ds \q + \dbE\int_t^T \big\langle
Q(s)v_1(s),x_2(s)\big\rangle_{H} ds +
\dbE\int_t^T \big\langle
x_1(s),Q(s)v_2(s)\big\rangle_{H} ds.
\end{array}
\end{equation}
From \eqref{eq def solzzz}--\eqref{eq def
solzzz3}, we see that $(P(\cd),Q(\cd))$
satisfies \eqref{eq def sol1}. Hence,
$(P(\cd),Q(\cd))$ is a transposition solution
of \eqref{op-bsystem3} (in the sense of
Definition \ref{op-definition2}). The
uniqueness of $(P(\cd),Q(\cd))$ follows from
Theorem \ref{op well th4-1}. This concludes the
proof of Theorem \ref{op well th0}.
\endpf

\br
Theorems \ref{op well th4-1}--\ref{op well th0}
indicate that, in some sense, the transposition
solution introduced in Definition
\ref{op-definition2} is a reasonable notion for
the solution to \eqref{op-bsystem3}.
Unfortunately, we are unable to prove the
existence of transposition solution to
\eqref{op-bsystem3} in the general case though
a weak version, i.e., the relaxed transposition
solution to this equation, introduced/studied
in the next three sections, suffices to
establish the desired Pontryagin-type
stochastic maximum principle for Problem (P) in
the general setting.
\er

\section{Sequential Banach-Alaoglu-type theorems in the operator
version}\label{s34}

The classical Banach-Alaoglu Theorem (e.g.
\cite[p. 130]{C}) states that the closed unit
ball of the dual space of a normed vector space
is compact in the weak* topology. This theorem
has an important special (sequential) version,
asserting that the closed unit ball of the dual
space of a separable normed vector space ({\it
resp.}, the closed unit ball of a reflexive
Banach space) is sequentially compact in the
weak* topology ({\it resp.}, the weak
topology). In this section, we shall present
several sequential Banach-Alaoglu-type theorems
for uniformly bounded linear operators (between
suitable Banach spaces). These results will
play crucial roles in the study of the
well-posedness of (\ref{op-bsystem3}) in the
general case.

Let $\{y_n\}_{n=1}^\infty\subset Y$ and $y\in
Y$. Let $\{z_n\}_{n=1}^\infty\subset Y^*$ and
$z\in Y^*$. In the sequel, we denote by
$$
\mbox{(w)-}\lim_{n\to\infty}y_n = y \ \mbox{ in
} Y
$$
when $\{y_n\}_{n=1}^\infty $ weakly converges
to $y$ in $Y$; and by
$$
\mbox{(w*)}\mbox{-}\lim_{n\to\infty}z_n = z\
\mbox{ in } Y^*
$$
when $\{z_n\}_{n=1}^\infty $ weakly$^*$
converges to $z$ in $Y^*$. Let us show first
the following result (It seems for us that this
is a known result. However we have not found it
in any reference):

\begin{lemma}\label{lemma6}
Let $X$ be a separable Banach space and let $Y$
be a reflexive Banach space. Assume that
$\{G_n\}_{n=1}^\infty\subset \cL(X,Y)$ is a
sequence of bounded linear operators such that
$\{G_nx\}_{n=1}^\infty$ is bounded for any
given $x\in X$. Then, there exist a subsequence
$\{G_{n_k}\}_{k=1}^\infty$ and a bounded linear
operator $G$ from $X$ to $Y$ such that
$$
\mbox{\rm(w)}\mbox{-}\lim_{k\to\infty} G_{n_k}x
= Gx \ \mbox{ in }\, Y,\qq \forall\; x\in X,
$$
$$
\mbox{\rm(w*)}\mbox{-}\lim_{k\to\infty}
G^*_{n_k}y^* = G^*y^* \ \mbox{ in }\,
X^*,\qq\forall\: y^*\in Y^*,
$$
and
 \bel{z03}
 |\!|G|\!|_{\cL(X,Y)}\le \sup_{n\in\dbN}|\!|G_n|\!|_{\cL(X,Y)}(<\infty).
 \ee
\end{lemma}

\br
Lemma \ref{lemma6} is not a direct consequence
of the classical sequential Banach-Alaoglu
Theorem. Indeed, as we mentioned before, the
Banach space $\cL(X,Y)$ is neither reflexive
nor separable even if both $X$ and $Y$ are
(infinite dimensional) separable Hilbert
spaces.
\er

{\it Proof of Lemma \ref{lemma6}}\,:  Noting
that $X$ is separable, we can find a countable
subset $\{x_i\}_{i=1}^\infty$ of $X$ such that
$\{x_1,x_2,\cdots\}$ is dense in $X$. Since
$\{G_n x_1\}_{n=1}^\infty$ is bounded in $Y$
and $Y$ is reflexive, there exists a
subsequence $\{n_k^{(1)}\}_{k=1}^\infty \subset
\{n\}_{n=1}^\infty$ such that $\ds
\mbox{\rm(w)}\mbox{-}\lim_{k\to\infty}G_{n_k^{(1)}}
x_1 = y_1$. Now, the sequence
$\big\{G_{n_k^{(1)}} x_2\big\}_{k=1}^\infty$ is
still bounded in $Y$, one can find a
subsequence $\{n_k^{(2)}\}_{k=1}^\infty \subset
\{n_k^{(1)}\}_{k=1}^\infty$ such that $\ds
\mbox{\rm(w)}\mbox{-}\lim_{k\to\infty}G_{n_k^{(2)}}
x_2 = y_2$. By the induction, for any $m\in
\dbN$, we can find a subsequence
$\{n_k^{(m+1)}\}_{k=1}^\infty \subset
\{n^{(m)}_k\}_{k=1}^\infty\subset\cds\subset
\{n_k^{(1)}\}_{k=1}^\infty \subset
\{n\}_{n=1}^\infty$ such that
$\ds\mbox{\rm(w)}\mbox{-}\lim_{k\to\infty}G_{n_k^{(m+1)}}
x_{m+1}=y_{m+1}$. We now use the classical
diagonalisation argument. Write $n_m =
n^{(m)}_m$, $m=1,2,\cds$. Then, it is clear
that $\{G_{n_m }x_i \}_{m=1}^\infty$ converges
weakly to $y_i$ in $Y$.

Let us define an operator $G$ (from $X$ to $Y$)
as follows: For any $x\in X$,
$$
G x = \lim_{k\to\infty} y_{i_k}=
\lim_{k\to\infty}\(\mbox{(w)-}\1n\lim_{m\to\infty}G_{n_m}
x_{i_k}\),
$$
where $\{x_{i_k}\}_{k=1}^\infty$ is any
subsequence of $ \{x_i\}_{i=1}^\infty$ such
that $\ds\lim_{k\to\infty}x_{i_k}=x$ in $X$. We
shall show below that $G\in \cL(X,Y)$.

First, we show that $G$ is well-defined. By the
Principle of Uniform Boundedness,  it is clear
that $\{G_n\}_{n=1}^\infty$ is uniformly
bounded in $\cL(X,Y)$. We choose $M>0$ such
that $|G_n|_{\cL(X,Y)}\leq M$ for all
$n\in\dbN$. Since $\{x_{i_{k}}\}_{k=1}^\infty$
is a Cauchy sequence in $X$, for any $\e>0$,
there is a $N>0$ such that
$|x_{i_{k_1}}-x_{i_{k_2}}|<\frac{\e}{M}$ when
$k_1,k_2>N$. Hence,
$|G_n(x_{i_{k_1}}-x_{i_{k_2}})|_Y<\e$ for any
$n\in\dbN$. Then, by the weakly sequentially
lower semicontinuity (of Banach spaces), we
deduce that
$$|y_{i_{k_1}}-y_{i_{k_2}}|_Y\leq \mathop{\underline{\hspace*{2pt}\lim}}\limits_{m\to\infty}
|G_{n_m}(x_{i_{k_1}}-x_{i_{k_2}})|_Y<\e,$$
which implies that $\{y_{i_k}\}_{k=1}^\infty$
is a Cauchy sequence in $Y$. Therefore, we see
that $\ds\lim_{k\to\infty}y_{i_k}$ exists in
$Y$. On the other hand, assume that there is
another subsequence
$\{x'_{i_{k}}\}_{k=1}^\infty\subset
\{x_i\}_{i=1}^\infty$ such that
$\ds\lim_{k\to\infty}x'_{i_{k}} = x$. Let
$y'_{i_k}$ be the corresponding weak limit of
$G_{n_m}x'_{i_k}$ in $Y$ for $m\to\infty$. Then
we find that
$$
\begin{array}{ll}\ds
|\lim_{k\to\infty}y_{i_k} -
\lim_{k\to\infty}y'_{i_k}|_Y \leq
\lim_{k\to\infty}\mathop{\underline{\hspace*{2pt}\lim}}\limits_{m\to\infty}|G_{n_m}(x_{i_k}-x'_{i_k})|_Y
\leq M\lim_{k\to\infty}| x_{i_k}-x'_{i_k}
|_X \\
\ns\ds \hspace{3.7cm} \leq M\lim_{k\to\infty}|
x_{i_k}-x |_X + M\lim_{k\to\infty}| x-x'_{i_k}
|_X = 0.
\end{array}
$$
Hence,  $G$ is well-defined.

Next, we prove that $G$ is a bounded linear
operator. For any $x\in X$ and the above
sequence $\{x_{i_{k}}\}_{k=1}^\infty$, it
follows that
$$
|Gx|_Y = \lim_{k\to\infty}|y_{i_k}|_Y\leq
\lim_{k\to\infty}\mathop{\underline{\hspace*{2pt}\lim}}\limits_{m\to\infty}|G_{n_m}x_{i_k}|_Y
\leq M\lim_{k\to\infty}|x_{i_k}|_X\leq M|x|_X.
$$
Hence, $G$ is a bounded operator. Further, for
any $x^{(1)}, x^{(2)} \in X$, $\a\in \dbC$ and
$\b \in \dbC$, we choose
$\{x^{(j)}_{i_k}\}_{k=1}^\infty\subset
\{x_i\}_{i=1}^\infty$, $j=1,2$, such that
$\ds\lim_{k\to\infty}x^{(j)}_{i_k} = x^{(j)}$,
and denote by $y^{(j)}_{i_k}$ the weak limit of
$G_{n_m}x^{(j)}_{i_k}$ in $Y$ for $m\to\infty$.
Hence $\ds
Gx^{(j)}=\lim_{k\to\infty}y^{(j)}_{i_k}$. Then,
$$
\a x^{(1)} + \b x^{(2)} = \lim_{k\to\infty}(\a
x^{(1)}_{i_k} + \b x^{(2)}_{i_k}) =
\a\lim_{k\to\infty} x^{(1)}_{i_k} +
 \b \lim_{k\to\infty}x^{(2)}_{i_k}
$$
and
$$
\mbox{(w)-}\1n\lim_{m\to\infty} G_{n_m}(\a
x^{(1)}_{i_k} + \b x^{(2)}_{i_k})
=\a\(\mbox{(w)-}\1n\lim_{m\to\infty} G_{n_m}
x^{(1)}_{i_k} \)+
\b\(\mbox{(w)-}\1n\lim_{m\to\infty} G_{n_m}
x^{(2)}_{i_k}\).
$$
Hence,
$$
\begin{array}{ll}\ds
G(\a x^{(1)} + \b x^{(2)}) \3n&\ds=
\lim_{k\to\infty}
\(\mbox{(w)-}\1n\lim_{m\to\infty} G_{n_m}(\a
x^{(1)}_{i_k} +
\b x^{(2)}_{i_k})\)\\
\ns&\ds  =
\a\lim_{k\to\infty}\(\mbox{(w)-}\1n\lim_{m\to\infty}
G_{n_m} x^{(1)}_{i_k}\)
+ \b\lim_{k\to\infty}\(\mbox{(w)-}\1n\lim_{m\to\infty} G_{n_m} x^{(2)}_{i_k}\) \\
\ns&\ds   = \a Gx^{(1)} + \b G x^{(2)}.
\end{array}
$$
Therefore, $G\in \cL(X,Y)$.

Also, for any $x\in X$ and $y^*\in Y^*$, it
holds that
$$
(x,G^*y^*)_{X,X^*} =(Gx,y^*)_{Y,Y^*} =
\lim_{k\to\infty}(G_{n_k}x,y^*)_{Y,Y^*} =
\lim_{k\to\infty}(x,G_{n_k}^*y^*)_{X,X^*}.
$$
Hence,
$$
\mbox{(w*)-}\lim_{k\to\infty} G^*_{n_k} y^* =
G^* y^* \mbox{ in } X^*.
$$

Finally, from the above proof, (\ref{z03}) is
obvious. This completes the proof of Lemma
\ref{lemma6}.
\endpf

\medskip

Let us introduce the following set class:
 \bel{x1z}
\cM=\Big\{O\in (0,T)\times
\O\;\Big|\;\{\chi_O(\cd)\} \mbox{ is an }
\dbF\hb{-adapted process}\Big\}.
 \ee
This set class will be used several times in
the sequel. We now show the following
``stochastic process" version of Lemma
\ref{lemma6}.

\begin{theorem}\label{lemma6.1}
Let $X$ and $Y$ be respectively a separable and
a reflexive Banach space, and let
$L^p(\O,\cF_T,\dbP)$, with $1\leq p <\infty$,
be separable. Let $1\leq p_1, p_2 < \infty$ and
$1<q_1, q_2 < \infty$. Assume that
$\{\cG_n\}_{n=1}^\infty$ is a sequence of
uniformly bounded, pointwisely defined linear
operators from
$L^{p_1}_{\dbF}(0,T;L^{p_2}(\O;X))$ to
$L^{q_1}_{\dbF}(0,T;L^{q_2}(\O;Y))$. Then,
there exist a subsequence
$\{\cG_{n_k}\}_{k=1}^\infty\subset
\{\cG_n\}_{n=1}^\infty$ and a
 $\cG\in\cL_{pd}\big(L^{p_1}_{\dbF}(0,T;L^{p_2}(\O;X)),\;
L^{q_1}_{\dbF}(0,T;L^{q_2}(\O;Y))\big)$ such
that
$$
\cG u(\cd) = {\mbox{\rm(w)-}}\lim_{k\to\infty}
\cG_{n_k}u(\cd) \ \mbox{ in }
L^{q_1}_{\dbF}(0,T;L^{q_2}(\O;Y)), \qq
\forall\; u(\cd)\in
L^{p_1}_{\dbF}(0,T;L^{p_2}(\O;X)).
$$
Moreover,
 $\ds|\!|\cG|\!|_{\cL(L^{p_1}_{\dbF}(0,T;L^{p_2}(\O;X)), \
L^{q_1}_{\dbF}(0,T;L^{q_2}(\O;Y)))}\le
\sup_{n\in\dbN}|\!|\cG_n|\!|_{\cL(L^{p_1}_{\dbF}(0,T;L^{p_2}(\O;X)),
\ L^{q_1}_{\dbF}(0,T;L^{q_2}(\O;Y))}$.
\end{theorem}

\br\label{r5.2}
i) As we shall see later, the most difficult
part in the proof of Theorem \ref{lemma6.1} is
to show that the weak limit operator $\cG$ is a
bounded, pointwisely defined linear operators
from $L^{p_1}_{\dbF}(0,T;L^{p_2}(\O;X))$ to
$L^{q_1}_{\dbF}(0,T;L^{q_2}(\O;Y))$. Note that,
a simple application of Lemma \ref{lemma6} to
the operators $\{\cG_n\}_{n=1}^\infty$ does not
guarantee this point but only that $\cG\in
\cL\big(L^{p_1}_{\dbF}(0,T;L^{p_2}(\O;X)),\,L^{q_1}_{\dbF}(0,T;
$ $L^{q_2}(\O;Y))\big)$.

ii) Theorem \ref{lemma6.1} indicates that
$\cL_{pd}\big(L^{p_1}_{\dbF}(0,T;L^{p_2}(\O;X)),\,L^{q_1}_{\dbF}(0,T;L^{q_2}(\O;Y))\big)$
is a closed linear subspace of the Banach space
$\cL\big(L^{p_1}_{\dbF}(0,T;L^{p_2}(\O;X)),\,L^{q_1}_{\dbF}(0,T;L^{q_2}(\O;Y))\big)$.
\er

{\it Proof of Theorem \ref{lemma6.1}}\,:  We
divide the proof into several steps.

\ms

{\bf Step 1.} Since $\{\cG_n\}_{n=1}^\infty$ is
a sequence of uniformly bounded, pointwisely
defined linear operators from
$L^{p_1}_{\dbF}(0,T;L^{p_2}(\O;X))$ to
$L^{q_1}_{\dbF}(0,T;L^{q_2}(\O;Y))$, for each
$n\in\dbN$ and $\ae (t,\omega)\in
(0,T)\times\Omega$, there exists an
$G_n(t,\o)\in\cL (X,Y)$ verifying that
 \bel{z20}
\big(\cG_nu(\cd)\big)(t,\o)=G_n(t,\o)u(t,\o),
\;\, \forall\; u(\cd)\in
L^{p_1}_{\dbF}(0,T;L^{p_2}(\O;X)).
 \ee
Write
 $$
 M=\sup_{n\in\dbN}|\!|\cG_n|\!|_{\cL(L^{p_1}_{\dbF}(0,T;L^{p_2}(\O;X)), \ L^{q_1}_{\dbF}(0,T;L^{q_2}(\O;Y))}.
 $$
By Lemma \ref{lemma6}, we conclude that there
exist a bounded linear operator $\cG$ from
$L^{p_1}_{\dbF}(0,T;L^{p_2}(\O;X))$ to
$L^{q_1}_{\dbF}(0,T;L^{q_2}(\O;Y))$ and a
subsequence $\{\cG_{n_k}\}_{k=1}^\infty\subset
\{\cG_n\}_{n=1}^\infty$ such that
 \bel{z21}
\cG u(\cd) = \mbox{(w)-}\lim_{k\to\infty}
\cG_{n_k}u(\cd) \ \mbox{ in }
L^{q_1}_{\dbF}(0,T;L^{q_2}(\O;Y)),
 \ee
and
 \bel{z01}
\left|\cG
u(\cd)\right|_{L^{q_1}_{\dbF}(0,T;L^{q_2}(\O;Y))}\le
M|u(\cd)|_{ L^{p_1}_{\dbF}(0,T;L^{p_2}(\O;X))},
\q\forall\; u(\cd)\in
L^{p_1}_{\dbF}(0,T;L^{p_2}(\O;X)).
 \ee

We claim that
\begin{equation}\label{lm6.1eq1}
 \sum_{i=1}^mf_i\cG u_i =\cG \Big( \sum_{i=1}^mf_iu_i\Big)=
  \mbox{(w)}\mbox{-}\lim_{k\to\infty}
 \sum_{i=1}^mf_i  \cG_{n_k}u_i\
\mbox{ in } L^{q_1}_{\dbF}(0,T;L^{q_2}(\O;Y))
\end{equation}
and
 \bel{z04}
\Big|\sum_{i=1}^mf_i\cG
u_i\Big|_{L^{q_1}_{\dbF}(0,T;L^{q_2}(\O;Y))}\le
M\Big|\sum_{i=1}^mf_i u_i\Big|_{
L^{p_1}_{\dbF}(0,T;L^{p_2}(\O;X))},
 \ee
where $m\in\dbN$, $f_i\in L^\infty_{\dbF}(0,T)
$ and $u_i\in
L^{p_1}_{\dbF}(0,T;L^{p_2}(\O;X))$,
$i=1,2,\cdots,m$. To show this, write $q_1' =
\frac{q_1}{q_1-1}$ and
$q_2'=\frac{q_2}{q_2-1}$. It follows from
(\ref{z21}) and (\ref{z20}) that for any
$v(\cd)\in
L^{q_1'}_{\dbF}(0,T;L^{q_1'}(\O;Y^*))$,
 \bel{z07}
\begin{array}{ll}\ds
\int_0^T \mE\Big\langle \sum_{i=1}^mf_i(s)(\cG u_i)(s), v(s) \Big\rangle_{Y,Y^*}ds = \int_0^T \mE\sum_{i=1}^m\big\langle  \cG u_i)(s),  f_i(s) v(s) \big\rangle_{Y,Y^*} ds \\
\ns\ds   =  \lim_{k\to\infty} \int_0^T\mE
\sum_{i=1}^m\big\langle
\big(\cG_{n_k}u_i\big)(s),  f_i(s) v(s)
\big\rangle_{Y,Y^*} ds= \lim_{k\to\infty}
\int_0^T\mE \sum_{i=1}^m\big\langle
(f_i\cG_{n_k}u_i)(s),   v(s)
\big\rangle_{Y,Y^*} ds,
\end{array}
 \ee
 and
 \bel{oq1-1}
 \ba{ll}
 \ds\lim_{k\to\infty} \int_0^T\mE \sum_{i=1}^m\big\langle
\cG_{n_k}(s)u_i(s),  f_i(s) v(s)
\big\rangle_{Y,Y^*} ds=\lim_{k\to\infty}
\int_0^T\mE \sum_{i=1}^m\big\langle
G_{n_k}(s)u_i(s),  f_i(s) v(s)
\big\rangle_{Y,Y^*} ds
\\\ns\ds=\lim_{k\to\infty} \int_0^T\mE \big\langle
G_{n_k}(s)\Big(\sum_{i=1}^mf_i(s)u_i(s)\Big),
v(s) \big\rangle_{Y,Y^*}
ds\\\ns\ds=\lim_{k\to\infty} \int_0^T\mE
\big\langle
\Big(\cG_{n_k}\big(\sum_{i=1}^mf_iu_i\big)\Big)(s),
v(s) \big\rangle_{Y,Y^*} ds=\int_0^T\mE
\big\langle
\Big(\cG\big(\sum_{i=1}^mf_iu_i\big)\Big)(s),
v(s) \big\rangle_{Y,Y^*} ds.
 \ea
 \ee
By \eqref{z07}--\eqref{oq1-1}, we obtain
(\ref{lm6.1eq1})--(\ref{z04}).

\ms

{\bf Step 2.} Each $x\in X$ can be regarded as
an element (i.e., $\chi_{(0,T)\times\O}(\cd)
x$) in $L^{p_1}_{\dbF}(0,T;L^{p_2}(\O;X))$.
Hence, $\cG x$ makes sense and belongs to
$L^{q_1}_{\dbF}(0,T;L^{q_2}(\O;Y))$. It is easy
to see that $\cL$ is a bounded linear operator
from $X$ to
$L^{q_1}_{\dbF}(0,T;L^{q_1}(\O;Y))$. By
(\ref{z01}), we find that
 \bel{z02} \left|\big(\cG
x\big)(\cd)\right|_{L^{q_1}_{\dbF}(0,T;L^{q_2}(\O;Y))}\le
MT^{1/p_2}|x|_X, \qq \forall\; x\in X.
 \ee

Write $B_X=\big\{x\in X\;\big|\;|x|_X\le
1\big\}$. By the separability of $X$, we see
that $\ds\Big\{\sup_{x\in B_X}\big|\big(\cG
x\big)(\cd)\big|_Y\Big\}$ is an $\dbF$-adapted
process. We claim that
 \bel{z11}
 \sup_{x\in B_X}\big|\big(\cG x\big)(t,\o)\big|_Y <\infty,\q\ae \,(t,\o)\in (0,T)\times \O.
 \ee
In the rest of this step, we shall prove
(\ref{z11}) by the contradiction argument.

Assume that (\ref{z11}) was not true. Then,
thanks to the adaptedness of $\ds\Big\{\sup_{x
\in B_X}\big|\big(\cG x\big)(\cd)\big|_Y\Big\}$
with respect to $\dbF$, there would be a set
$A\in \cM$, defined by (\ref{x1z}), such that
$\mu(A)>0$ (Here $\mu$ stands for the product
measure of the Lebesgue measure (on $[0,T]$)
and the probability measure $\dbP$) and that
$$
\sup_{x\in B_X}\big|\big(\cG
x\big)(t,\o)\big|_Y=\infty, \mbox{ for
}(t,\o)\in A.
$$
Let $\{x_i\}_{i=1}^\infty$ be a sequence  in
$B_X$ such that it is dense in $B_X$.  Then
$$
\sup_{i\in \dbN}\big|\big(\cG
x_i\big)(t,\o)\big|_Y = \sup_{x\in
B_X}\big|\big(\cG x\big)(t,\o)\big|_Y=\infty,\q
\mbox{ for } (t,\o)\in A.
$$
For any $n\in \dbN$, we define a sequence of
subsets of $(0,T)\times \O$ in the following
way.
\begin{equation}\label{Min}
\left\{
\begin{array}{ll}\ds
A_1^{(n)} = \Big\{(t,\o)\in(0,T)\times \O\;\Big|\; \big|\big(\cG x_1\big)(t,\o)\big|_Y \geq n \Big\},\\
\ns\ds  A_i^{(n)} = \Big\{(t,\o)\in
\big((0,T)\times \O\big)\setminus
\big(\bigcup_{k=1}^{i-1}A_k^{(n)}\big)\;\Big|\;\big|\big(\cG
x_i\big)(t,\o)\big|_Y \geq n \Big\},\; \mbox{
if } i>1.
\end{array}
\right.
\end{equation}
It follows from the adaptedness of
$\big|\big(\cG x\big)(\cd)\big|_Y$ that
$A_i^{(n)}\in\cM$ for every $i\in\mathbb{N}$
and $n\in\mathbb{N}$. It is clear that $\ds
A\subset \bigcup_{i=1}^\infty A_i^{(n)}$  for
any $n\in\mathbb{N}$ and $A_i^{(n)}\cap
A_j^{(n)} =\emptyset$ for $i\neq j$. Hence, we
see that
 $$
\sum_{i=1}^\infty \mu(A_i^{(n)}) = \mu\big(
\bigcup_{i=1}^\infty A_i^{(n)} \big) \geq
\mu(A)>0, \mbox{ for all } n\in \mathbb{N}.
 $$
Thus, for each $n\in\dbN$, there is a
$N_n\in\dbN$ such that
\begin{equation}\label{muMin}
\sum_{i=1}^{N_n}\mu(A_i^{(n)}) = \mu\big(
\bigcup_{i=1}^{N_n} A_i^{(n)} \big) \geq
\frac{\mu(A)}{2}>0.
\end{equation}

Write
 \bel{z06} x^{(n)}(t,\o) = \sum_{i=1}^{N_n} \chi_{A_i^{(n)}}(t,\o)x_i.
 \ee
Clearly, $\{x^{(n)}(t)\}_{t\in [0,T]}$ is an
adapted process. By $\cG\in
\cL\big(L^{p_1}_{\dbF}(0,T;L^{p_2}(\O;X)),\;L^{q_1}_{\dbF}(0,T;L^{q_2}(\O;Y))\big)$
and $|x^{(n)}(t,\o)|_X\le 1$ for $\ae (t,\o)\in
(0,T)\times\O$, we find that
\begin{equation}\label{lm6.q1}
 \ba{ll}\ds
|\cG
x^{(n)}|_{L^{q_1}_{\dbF}(0,T;L^{q_2}(\O;Y))}\!\!\!&\ds\le
M\left\{\int_0^T\left[\int_\O\left|x^{(n)}(t,\o)\right|_X^{p_1}
\dbP(d\o)\right]^{\frac{p_2}{p_1}}dt\right\}^{\frac{1}{p_2}}
 \\
\ns&\ds \le M T^{1/p_2},\qq \mbox{ for all }
n\in\mathbb{N}.
 \ea
\end{equation}

On the other hand, let us choose a $n >
\frac{2M}{\mu(A)}T^{\frac{1}{p_1}+\frac{1}{q'_1}}$.
From (\ref{Min})--(\ref{z06}), and noting
\eqref{lm6.1eq1}, we obtain that
$$
\begin{array}{ll}\ds
|\cG
x^{(n)}|_{L^{q_1}_{\dbF}(0,T;L^{q_2}(\O;Y))}
\ge T^{-\frac{1}{q'_1}}|\cG
x^{(n)}|_{L^1_{\dbF}(0,T;L^1
(\O;Y))}\\
\ns\ds\hspace{3.8cm} =
T^{-\frac{1}{q'_1}}\sum_{i=1}^{N_n}
\int_{A_i^{(n)}}\big|  \cG x_i \big|_Y dtd\dbP
 =T^{-\frac{1}{q'_1}}\sum_{i=1}^{N_n} \int_{A_i^{(n)}}\big|  \cG x_i \big|_Y d\mu   \\
\ns\ds\hspace{3.8cm}\geq
T^{-\frac{1}{q'_1}}n\sum_{i=1}^{N_n}\mu(A_i^{(n)})
\geq \frac{\mu(A)}{2}T^{-\frac{1}{q'_1}}n >
MT^{\frac{1}{p_1}},
\end{array}
$$
which contradicts the inequality
\eqref{lm6.q1}. Therefore, we conclude that
\eqref{z11} holds.

\ms

{\bf Step 3.} By (\ref{z11}), for $\ae
(t,\o)\in (0,T)\times\O$, we may define an
operator $ G(t,\o)\in\cL (X,Y)$ by
 \bel{z15}
X\ni x\mapsto G(t,\o)x=\big(\cG x\big)(t,\o).
 \ee
Further, we introduce the following subspace of
$L^{p_1}_{\dbF}(0,T;L^{p_2}(\O;X))$:
$$
\cX=\Big\{u(\cd)=\sum_{i=1}^m
\chi_{A_i}(\cd)h_i\;\Big|\; m\in\dbN,\,
A_i\in\cM, h_i\in X\Big\}.
$$
It is clear that $\cX$ is dense in
$L^{p_1}_{\dbF}(0,T;L^{p_2}(\O;X))$. We now
define a linear operator $\widetilde \cG$ from
$\cX$ to $L^{q_1}_{\dbF}(0,T;L^{q_2}(\O;Y))$ by
\bel{z16} \cX\ni u(\cd)=\sum_{i=1}^m
\chi_{A_i}(\cd)h_i\mapsto (\widetilde \cG
u)(t,\o)=\sum_{i=1}^m
\chi_{A_i}(t,\o)G(t,\o)h_i. \ee

We claim that
 \bel{z17}
(\widetilde \cG u)(\cd)=(\cG
u)(\cd),\qq\forall\; u(\cd)\in \cX.
 \ee
Indeed, it follows from (\ref{z20}) that for
any $v(\cd)\in
L^{q_1'}_{\dbF}(0,T;L^{q_2'}(\O;Y^*))$, and
$u(\cd)$ to be of the form in (\ref{z16}),
 $$
\begin{array}{ll}\ds
\mE\int_0^T \big\langle  (\widetilde \cG u)(s), v(s) \big\rangle_{Y,Y^*}ds  = \mE\int_0^T \big\langle  \sum_{i=1}^m \chi_{A_i}(s) G (s)h_i,  v(s) \big\rangle_{Y,Y^*} ds\\
\ns\ds = \mE\int_0^T \big\langle  \sum_{i=1}^m \chi_{A_i}(s)\big(\cG h_i\big)(s),  v(s) \big\rangle_{Y,Y^*} ds  = \sum_{i=1}^m\mE\int_0^T \big\langle  \big(\cG h_i\big)(s), \chi_{A_i}(s)v(s) \big\rangle_{Y,Y^*} ds \\
\ns\ds   = \sum_{i=1}^m \lim_{k\to\infty} \mE\int_0^T \big\langle  G_{n_k}(s)h_i,  \chi_{A_i}(s) v(s) \big\rangle_{Y,Y^*} ds \\
\ns\ds   =  \lim_{k\to\infty} \mE\int_0^T
\big\langle
G_{n_k}(s)\Big(\sum_{i=1}^m\chi_{A_i}(s)h_i\Big),
v(s) \big\rangle_{Y,Y^*} ds  = \mE\int_0^T
\big\langle (\cG u)(s), v(s)
\big\rangle_{Y,Y^*} ds.
\end{array}
 $$
This gives (\ref{z17}).

Recall that $\cG$ is a bounded linear operator
from $L^{p_1}_{\dbF}(0,T;L^{p_2}(\O;X))$ to
$L^{q_1}_{\dbF}(0,T;L^{q_2}(\O;Y))$. Hence, it
is also a bounded linear operator from $\cX$ to
$L^{q_1}_{\dbF}(0,T;L^{q_2}(\O;Y))$. By
(\ref{z17}), we conclude that $\widetilde \cG$
is a bounded linear operator from $\cX$ to
$L^{q_1}_{\dbF}(0,T;L^{q_2}(\O;Y))$. Since
$\cX$ is dense in
$L^{p_1}_{\dbF}(0,T;L^{p_2}(\O;X))$, it is
clear that $\widetilde \cG$ can be uniquely
extended as a bounded linear operator from
$L^{p_1}_{\dbF}(0,T;L^{p_2}(\O;X))$ to
$L^{q_1}_{\dbF}(0,T;L^{q_2}(\O;Y))$ (We still
denote by $\widetilde \cG$ its extension). By
(\ref{z17}) again, we conclude that
  \bel{z18}
\widetilde \cG =\cG .
 \ee

It remains to show that
 \bel{zz15}\big(\widetilde \cG
u(\cd)\big)(t,\omega) = G(t,\o)u(t,\omega), \qq
\ae (t,\omega)\in (0,T)\times \Omega,
 \ee
for all $u\in L^{p_1}_\dbF(0,T;L^{p_2}(\O;X))$.
For this purpose, by the fact that $\cX$ is
dense in $L^{p_1}_{\dbF}(0,T;L^{p_2}(\O;X))$,
we may assume that
 \bel{z22}
 u(\cd)=\sum_{i=1}^\infty \chi_{A_i}(\cd)h_i,
 \ee
for some $A_i\in\cM$ and $h_i\in X$,
$i=1,2,\cdots$ (Note that here we assume
neither  $A_i\bigcap A_j=\emptyset$ nor
$h_i\not=h_j$ for $i,j=1,2,\cdots$). For each
$n\in \dbN$, write $\ds u^n(\cd)=\sum_{i=1}^n
\chi_{A_i}(\cd)h_i$. From (\ref{z22}), it is
clear that
 \bel{z23}
 u(\cd)=\lim_{n\to\infty}u^n(\cd),\qq \hb{in
 }L^{p_1}_\dbF(0,T;L^{p_2}(\O;X)).
 \ee
By (\ref{z04}), (\ref{z15}), (\ref{z16}) and
(\ref{z23}), it is easy to see that
 \bel{z24}\big(\widetilde \cG u^n(\cd)\big)(t,\o)=\sum_{i=1}^n
\chi_{A_i}(t,\o)G(t,\o)h_i
 \ee
is a Cauchy sequence in
$L^{q_1}_{\dbF}(0,T;L^{q_2}(\O;Y))$. Hence, by
(\ref{z24}) and recalling that $\widetilde \cG$
is a bounded linear operator from
$L^{p_1}_{\dbF}(0,T;L^{p_2}(\O;X))$ to
$L^{q_1}_{\dbF}(0,T;L^{q_2}(\O;Y))$, we
conclude that
 \bel{z25}\big(\widetilde \cG u(\cd)\big)(t,\o)=\sum_{i=1}^\infty
\chi_{A_i}(t,\o)G(t,\o)h_i.
 \ee
Combining (\ref{z22}) and (\ref{z25}), we
obtain (\ref{zz15}).

Finally, by (\ref{z18}) and (\ref{zz15}), the
desired result follows. This completes the
proof of Theorem \ref{lemma6.1}.
\endpf

\ms

From the proof of Theorem \ref{lemma6.1}, it is
easy to deduce the following result.

\begin{corollary}\label{cor1}
Let $X$ and $Y$ be respectively a separable and
a reflexive Banach space, and let
$L^p(\O,\cF_T,\dbP)$, with $1\leq p <\infty$,
be separable. Let $1<q_1,q_2 < \infty$. Assume
that $\{\cG_n\}_{n=1}^\infty$ is a sequence of
uniformly bounded, pointwisely defined linear
operators from $X$ to
$L^{q_1}_{\dbF}(0,T;L^{q_2}(\O;Y))$. Then,
there exist a subsequence
$\{\cG_{n_k}\}_{k=1}^\infty\subset
\{\cG_n\}_{n=1}^\infty$ and an
$\cG\in\cL_{pd}\big(X,\;
L^{q_1}_{\dbF}(0,T;L^{q_2}(\O;Y))\big)$ such
that
$$
\cG x = {\mbox{\rm(w)-}}\lim_{k\to\infty}
\cG_{n_k}x \ \mbox{ in }
L^{q_1}_{\dbF}(0,T;L^{q_2}(\O;Y)), \qq
\forall\; x\in X.
$$
Moreover,
 $\ds|\!|\cG|\!|_{\cL(X,
L^{q_1}_{\dbF}(0,T;L^{q_2}(\O;Y)))}\le
\sup_{n\in\dbN}|\!|\cG_n|\!|_{\cL(X, \
L^{q_1}_{\dbF}(0,T;L^{q_2}(\O;Y))}$.
\end{corollary}

\ms

Proceeding exactly as in the proof of Theorem
\ref{lemma6.1}, we can show the following
``random variable" and ``random
variable-stochastic process" versions of Lemma
\ref{lemma6} (Hence, the detailed proof will be
omitted).

\begin{theorem}\label{lemma6--1}
Let $X$ and $Y$ be accordingly a separable and
a reflexive Banach space, and let
$L^p(\O,\cF_T,\dbP)$, with $1\leq p <\infty$,
be separable. Let $1\leq p_1< \infty$ and
$1<q_1 < \infty$. Assume that
$\{\cG_n\}_{n=1}^\infty$ is a sequence of
uniformly bounded, pointwisely defined linear
operators from $L^{p_1}_{\cF_T}(\O;X)$ to
$L^{q_1}_{\cF_T}(\O;Y)$. Then, there exist a
subsequence $\{\cG_{n_k}\}_{k=1}^\infty\subset
\{\cG_n\}_{n=1}^\infty$ and an
$\cG\in\cL_{pd}\big(L^{p_1}_{\cF_T}(\O;X),\;$ $
L^{q_1}_{\cF_T}(\O;Y)\big)$ such that
$$
\cG u(\cd) = {\mbox{\rm(w)-}}\lim_{k\to\infty}
\cG_{n_k}u(\cd) \ \mbox{ in }
L^{q_1}_{\cF_T}(\O;Y), \qq \forall\; u(\cd)\in
L^{p_1}_{\cF_T}(\O;X).
$$
Moreover,
 $\ds|\!|\cL|\!|_{\cL(L^{p_1}_{\cF_T}(\O;X), \
L^{q_1}_{\cF_T}(\O;Y))}\le
\sup_{n\in\dbN}|\!|\cG_n|\!|_{\cL(L^{p_1}_{\cF_T}(\O;X),\
L^{q_1}_{\cF_T}(\O;Y))}$.
\end{theorem}

\begin{theorem}\label{lemma6---1}
Let $X$ and $Y$ be respectively a separable and
a reflexive Banach space, and let
$L^p(\O,\cF_T,\dbP)$, with $1\leq p <\infty$,
be separable. Let $1\leq p_1 < \infty$, $1<q_1,
q_2 < \infty$ and $0\le t_0\le T$. Assume that
$\{\cG_n\}_{n=1}^\infty$ is a sequence of
uniformly bounded, pointwisely defined linear
operators from $L^{p_1}_{\cF_{t_0}}(\O;X)$ to
$L^{q_1}_{\dbF}(t_0,T;L^{q_2}(\O;Y))$. Then,
there exist a subsequence
$\{\cG_{n_k}\}_{k=1}^\infty\subset
\{\cG_n\}_{n=1}^\infty$  and an
$\cG\in\cL_{pd}\big(L^{p_1}_{\cF_{t_0}}(\O;X),\;L^{q_1}_{\dbF}(t_0,T;L^{q_2}(\O;Y))\big)$
such that
$$
\cG u(\cd) = {\mbox{\rm(w)-}}\lim_{k\to\infty}
\cG_{n_k}u(\cd) \ \mbox{ in }
L^{q_1}_{\dbF}(t_0,T;L^{q_2}(\O;Y)), \q
\forall\; u(\cd)\in L^{p_1}_{\cF_{t_0}}(\O;X).
$$
Moreover,
 $\ds|\!|\cG|\!|_{\cL(L^{p_1}_{\cF_{t_0}}(\O;X), \
L^{q_1}_{\dbF}(t_0,T;L^{q_2}(\O;Y)))}\le
\sup_{n\in\dbN}|\!|\cG_n|\!|_{\cL(L^{p_1}_{\cF_{t_0}}(\O;X),
\ L^{q_1}_{\dbF}(t_0,T;L^{q_2}(\O;Y))}$.
\end{theorem}

\br
Similar to Remark \ref{r5.2} ii), we see that
$\cL_{pd}\big(X,\;
L^{q_1}_{\dbF}(0,T;L^{q_2}(\O;Y))\big)$,
$\cL_{pd}\big(L^{p_1}_{\cF_T}(\O;X),\;$
$L^{q_1}_{\cF_T}(\O;Y)\big)$ and
$\cL_{pd}\big(L^{p_1}_{\cF_{t_0}}(\O;X),\;L^{q_1}_{\dbF}(t_0,T;L^{q_2}(\O;Y))\big)$
(for any given $t_0\in [0,T]$) are accordingly
closed linear subspaces of the Banach spaces
$\cL\big(X,\;
L^{q_1}_{\dbF}(0,T;L^{q_2}(\O;Y))\big)$,
$\cL\big(L^{p_1}_{\cF_T}(\O;X),\;L^{q_1}_{\cF_T}(\O;Y)\big)$
and
$\cL\big(L^{p_1}_{\cF_{t_0}}(\O;X),\;L^{q_1}_{\dbF}(t_0,T;L^{q_2}(\O;Y))\big)$.
\er

It is clear that the probability space
$(\O,\cF,\dbP)$ plays no special role in the
above Theorems
\ref{lemma6.1}--\ref{lemma6---1}. For possible
applications in other places, we give below a
``deterministic" modification of Theorem
\ref{lemma6.1}.

Let $(\O_1,\cM_1,\m_1)$ and $(\O_2,\cM_2,\m_2)$
be two finite measure spaces. Let $\cM$ be a
sub-$\si$-field of the $\si$-field generated by
$\cM_1\times\cM_2$, and for any $1\le
p,q<\infty$, let
$$\ba{ll}
\ns\ds
L^p_\cM(\O_1;L^q(\O_2;X))=\Big\{\f:\O_1\times
\O_2\to X\;\Big|\;\f(\cd)\hb{ is
$\cM$-measurable and}\\
\ns\ds\qq\qq\qq\qq\qq\qq\qq\qq\int_{\O_1}\left(\int_{\O_2}|\f(\o_1,\o_2)|_H^qd\m_2(\o_2)\right)^{p\over
q}d\m_1(\o_1)<\infty\Big\}.\ea$$
It is easy to show that $
L^p_\cM(\O_1;L^q(\O_2;X))$ is a Banach space
with the canonical norm. Similar to the proof
of Theorem \ref{lemma6.1}, one can prove the
following result:

\begin{theorem}\label{lma6.1}
Let $X$ and $Y$ be respectively a separable and
a reflexive Banach space. Let $1\leq p_1, p_2 <
\infty$ and $1<q_1, q_2 < \infty$, and let
$L^{p_1}_\cM(\O_1;L^{p_2}(\O_2;\dbC))$ be
separable. Assume that $\{\cG_n\}_{n=1}^\infty$
is a sequence of uniformly bounded, pointwisely
defined linear operators from
$L^{p_1}_\cM(\O_1;L^{p_2}(\O_2;X))$ to
$L^{q_1}_\cM(\O_1;L^{q_2}(\O_2;Y))$. Then,
there exist a subsequence
$\{\cG_{n_k}\}_{k=1}^\infty\subset
\{\cG_n\}_{n=1}^\infty$ and a
 $\cG\in\cL_{pd}\big(L^{p_1}_\cM(\O_1;$ $L^{p_2}(\O_2;X)),\,
L^{q_1}_\cM(\O_1;L^{q_2}(\O_2;Y))\big)$
(defined similarly as
$\cL_{pd}\big(L^{p_1}_{\dbF}(0,T;L^{p_2}(\O;X)),\;
L^{q_1}_{\dbF}(0,T;L^{q_2}(\O;$ $Y))\big)$)
such that
$$
\cG u(\cd) = {\mbox{\rm(w)-}}\lim_{k\to\infty}
\cG_{n_k}u(\cd) \ \mbox{ in }
L^{q_1}_\cM(\O_1;L^{q_2}(\O_2;Y)), \qq
\forall\; u(\cd)\in
L^{p_1}_\cM(\O_1;L^{p_2}(\O_2;X)).
$$
Moreover,
 $\ds|\!|\cG|\!|_{\cL(L^{p_1}_\cM(\O_1;L^{p_2}(\O_2;X)), \
L^{q_1}_\cM(\O_1;L^{q_2}(\O_2;Y)))}\le
\sup_{n\in\dbN}|\!|\cG_n|\!|_{\cL(L^{p_1}_\cM(\O_1;L^{p_2}(\O_2;X)),
\ L^{q_1}_\cM(\O_1;L^{q_2}(\O_2;Y))}$.
\end{theorem}

\section{Well-posedness of the operator-valued BSEEs in the general case}\label{s4}

This section is addressed to proving the
well-posedness result for the equation
(\ref{op-bsystem3}) with general data in the
sense of relaxed transposition solution, to be
defined later.

Write
 \bel{jshi2}
\3n\begin{array}{ll}\ds
\cQ[0,T]\!\=\!\Big\{\big(Q^{(\cd)},\widehat
Q^{(\cd)}\big)\;\Big|\;\mbox{For any } t\in
[0,T], \mbox{ both }Q^{(t)}\mbox{ and }\widehat
Q^{(t)}\mbox{ are bounded
linear operators}\\
\ns\ds\hspace{1.7cm}\mbox{ from
}L^4_{\cF_t}(\O;H)\times
L^2_\dbF(t,T;L^4(\O;H))\times
L^2_\dbF(t,T;L^4(\O;H)) \mbox{ to }
L^{2}_\dbF(t,T;L^{\frac{4}{3}}(\O;H))\\
\ns\ds \hspace{1.7cm} \mbox{ and
}Q^{(t)}(0,0,\cd)^*=\widehat
Q^{(t)}(0,0,\cd)\Big\}.
\end{array}
 \ee
We now define the relaxed transposition
solution to \eqref{op-bsystem3} as follows:

\begin{definition}\label{op-definition2x}
We call $\big(P(\cd),Q^{(\cd)},\widehat
Q^{(\cd)}\big)\in D_{\dbF,w}([0,T];
L^{\frac{4}{3}}(\O;\cL(H)))\times \cQ[0,T]$ a
relaxed transposition solution to
\eqref{op-bsystem3} if for any $t\in [0,T]$,
$\xi_1,\xi_2\in L^4_{\cF_t}(\O;H)$, $u_1(\cd),
u_2(\cd)\in L^2_{\dbF}(t,T;L^4(\O;H))$ and
$v_1(\cd),v_2(\cd)\in L^2_{\dbF}(t,T;
L^4(\O;H))$, it holds that
 \begin{equation}\label{6.18eq1}
\begin{array}{ll}
\ds \q\mE\big\langle P_T x_1(T), x_2(T)
\big\rangle_{H} - \mE \int_t^T \big\langle
F(s) x_1(s), x_2(s) \big\rangle_{H}ds\\
\ns\ds =\mE\big\langle P(t) \xi_1,\xi_2
\big\rangle_{H} + \mE \int_t^T \big\langle
P(s)u_1(s), x_2(s)\big\rangle_{H}ds + \mE
\int_t^T \big\langle P(s)x_1(s),
u_2(s)\big\rangle_{H}ds \\
\ns\ds \q  + \mE \int_t^T \big\langle P(s)K
(s)x_1 (s), v_2 (s)\big\rangle_{H}ds +  \mE
\int_t^T \big\langle  P(s)v_1 (s), K (s)x_2 (s)+ v_2(s)\big\rangle_{H}ds\\
\ns\ds \q + \mE \int_t^T \big\langle v_1(s),
\widehat
Q^{(t)}(\xi_2,u_2,v_2)(s)\big\rangle_{H}ds+ \mE
\int_t^T \big\langle Q^{(t)}(\xi_1,u_1,v_1)(s),
v_2(s) \big\rangle_{H}ds.
\end{array}
\end{equation}
Here, $x_1(\cd)$ and $x_2(\cd)$ solve
\eqref{op-fsystem2} and \eqref{op-fsystem3},
respectively.
\end{definition}

\br
It is easy to see that, if
$\big(P(\cd),Q(\cd)\big)$ is a transposition
solution to \eqref{op-bsystem3}, then
$\big(P(\cd), Q^{(\cd)},\widehat
Q^{(\cd)}\big)$ is a relaxed transposition
solution to the same equation, where (Recall
Lemma \ref{lemma5} for $U(\cd,t)$, $V(\cd,t)$
and $W(\cd,t)$)
 $$
 \left\{
 \ba{ll}
 \ds Q^{(t)}(\xi,u,v)= Q(\cd)U(\cd,t)\xi+Q(\cd)V(\cd,t)u+\frac{Q(\cd)W(\cd,t)+\big(Q(\cd)^*W(\cd,t)\big)^*}{2}v,\\\ns\ds
 \widehat Q^{(t)}(\xi,u,v)=Q(\cd)^*U(\cd,t)\xi+Q(\cd)^*V(\cd,t)u+\frac{Q(\cd)^*W(\cd,t)+\big(Q(\cd)W(\cd,t)\big)^*}{2}v, \ea
 \right.
$$
for any $(\xi,u,v)\in L^4_{\cF_t}(\O;H)\times
L^2_\dbF(t,T;L^4(\O;H))\times
L^2_\dbF(t,T;L^4(\O;H))$. However, it is
unclear how to obtain a transposition solution
$\big(P(\cd),Q(\cd)\big)$  to
\eqref{op-bsystem3} by means of its relaxed
transposition solution $\big(P(\cd),
Q^{(\cd)},\widehat Q^{(\cd)}\big)$. It seems
that this is possible but we cannot do it at
this moment.
\er

We have the following well-posedness result for
the equation (\ref{op-bsystem3}).

\begin{theorem}\label{op well th}
Assume that $H$ is a separable Hilbert space,
and $L^p_{\cF_T}(\O;\dbC)$ ($1\leq p < \infty$)
is a separable Banach space. Then, for any
$P_T\in L^2_{\cF_T}(\O;\cL(H))$, $F\in
L^1_\dbF(0,T;L^2(\O;\cL(H)))$ and $J,K\in
L^4_\dbF(0,T; L^\infty(\O; \cL(H)))$, the
equation \eqref{op-bsystem3} admits one and
only one relaxed transposition solution
$\big(P(\cd),Q^{(\cd)},\widehat
Q^{(\cd)}\big)\in D_{\dbF,w}([0,T];
L^{\frac{4}{3}}(\O;\cL(H)))\times \cQ[0,T]$.
Furthermore,
\begin{equation}\label{ine the1z}
\begin{array}{ll}\ds
|\!|P|\!|_{ \cL(L^2_\dbF(0,T;L^4(\O;H)),\;
L^{2}_{\dbF}(0,T;L^{\frac{4}{3}}(\O;H)))}\\\ns\ds\q+\sup_{t\in
[0,T]}\big|\!\big|\big(Q^{(t)},\widehat
Q^{(t)}\big)\big|\!\big|_{\big(\cL(L^4_{\cF_t}(\O;H)\times
L^2_\dbF(t,T;L^4(\O;H))\times
L^2_\dbF(t,T;L^4(\O;H)),\;
L^{2}_\dbF(t,T;L^{\frac{4}{3}}(\O;H))\big)^2}\\
\ns\ds\leq C\left[
|F|_{L^1_\dbF(0,T;\;L^2(\O;\cL(H)))} +
|P_T|_{L^2_{\cF_T}(\O;\;\cL(H))}\right].
\end{array}
\end{equation}
\end{theorem}

{\it Proof}\,: We consider only the case that
$H$ is a real Hilbert space (The case of
complex Hilbert spaces can be treated
similarly). The proof is divided into several
steps.

\ms

{\bf Step 1}. In this step, we introduce a
suitable approximation to the equation
\eqref{op-bsystem3}.

Let $\{e_n\}_{n=1}^\infty$ be an orthonormal
basis of $H$ and $\{\G_n\}_{n=1}^\infty$ be the
standard projection operator from $H$ onto its
subspace $\span \{e_1,e_2,\cdots,e_n\}$, that
is, $\ds \G_n x = \sum_{i=1}^n x_i e_i$ for any
$\ds x=\sum_{i=1}^\infty x_i e_i\in H$. Write
$H_n = \G_n H$. It is clear that, for each
$n\in\dbN$, $H_n$ is isomorphic to the
$n$-dimensional Euclidean space $\dbR^n$. In
the sequel, we identify $H_n$ by $\dbR^n$, and
hence $\cL(H_n)=\cL(\dbR^n)$ is the set of all
$n\times n$ (real) matrices. For any $M_1,
M_2\in \cL(\dbR^n)$, put $\big\langle M_1,
M_2\big\rangle_{\cL(\dbR^n)} = \tr(M_1
M^\top_2)$. It is easy to check that
$\big\langle \cd,\cd\big\rangle_{\cL(\dbR^n)}$
is an inner product on $\cL(\dbR^n)$, and
$\cL(\dbR^n)$ is a Hilbert space with this
inner product.

Consider the following matrix-valued BSDE:
\begin{equation}\label{op-bsystem4}
\left\{
\begin{array}{ll}
\ds dP^{n,\l} = - (A_{\l,n}^* + J_n^*)P^{n,\l}
dt - P^{n,\l}(A_{\l,n} + J_n)dt -
K^*_nP^{n,\l}K_n dt
 \\
 \ns\ds \hspace{1.46cm} - (K_n^* Q^{n,\l} + Q^{n,\l} K_n)dt +  F_ndt + Q^{n,\l} dw(t) &\mbox{ in } [0,T),\\
\ns\ds P^{n,\l}(T) =  P^n_{T},
\end{array}
\right.
\end{equation}
where $\l\in\rho(A)$, $A_{\l,n} = \G_n A_{\l}
\G_n$, $A_{\l}$ (as before) stands for the
Yosida approximation of $A$, $J_n = \G_n J
\G_n$, $K_n = \G_n K \G_n$, $F_n = \G_n F \G_n$
and $P^n_{T} = \G_n P_T\G_n $.

The solution to \eqref{op-bsystem4} is
understood in the transposition sense.
According to Theorem \ref{the1} (or
\cite[Theorem 4.1]{LZ}), the equation
\eqref{op-bsystem4} admits a unique
transposition solution
$(P^{n,\l}(\cd),Q^{n,\l}(\cd))\in
D_{\dbF}([0,T];L^2(\O;\cL(\dbR^n)))\times
L^2_\dbF(0,T;L^2(\O;\cL(\dbR^n))) $ such that,
for any $t\in [0,T]$, $U^n_1(\cdot)\in
L^1_{\dbF}(t,T;$ $L^2(\O;\cL(\dbR^n)))$,
$V^n_1(\cdot)\in
 L^2_{\dbF}(t,T;L^2(\O;\cL(\dbR^n)))$ and $\eta\in
 L^2_{\cF_t}(\O;\cL(\dbR^n))$, and the corresponding solution $X^n(\cd)\in  C_{\dbF}([t,T];L^2(\O;\cL(\dbR^n)))$ of the following equation:
\begin{equation}\label{op-fsystem1}
\left\{
\begin{array}{ll}
\ds dX^n = U_1^n ds + V_1^n dw(s) &\mbox{ in } (t,T],\\
\ns\ds X^n(t) = \eta,
\end{array}
\right.
\end{equation}
it holds that
\begin{equation}\label{eq def solz}
\begin{array}{ll}\ds
\q\dbE \big\langle
X^n(T),P^n_T\big\rangle_{\cL(\dbR^n)}
 -\dbE\int_t^T \big\langle X^n(s),\Phi^{n,\l}(s)\big\rangle_{\cL(\dbR^n)}ds\\
\ns\ds = \dbE
\big\langle\eta,P^{n,\l}(t)\big\rangle_{\cL(\dbR^n)}
+ \dbE\int_t^T \big\langle
U^n_1(s),P^{n,\l}(s)\big\rangle_{\cL(\dbR^n)}
ds + \dbE\int_t^T \big\langle V^n_1(s),
Q^{n,\l}(s)\big\rangle_{\cL(\dbR^n)} ds,
\end{array}
\end{equation}
where
 \begin{equation}\label{x1x} \Phi^{n,\l}= - \big(A_{\l,n}^* +
J_n^*\big)P^{n,\l} - P^{n,\l}\big(A_{\l,n} +
J_n\big)  - K^*_nP^{n,\l}K_n -  K_n^* Q^{n,\l}
- Q^{n,\l} K_n + F_n.
 \end{equation}

Clearly, \eqref{op-bsystem4} can be regarded as
finite dimensional approximations of the
equation \eqref{op-bsystem3}.  In the rest of
the proof, we shall construct the desired
solution to the equation \eqref{op-bsystem3} by
means of the solutions to \eqref{op-bsystem4}.

\ms

{\bf Step 2}.  This step is devoted to
introducing   suitable finite approximations of
the equations \eqref{op-fsystem2} and
\eqref{op-fsystem3}.

We approximate accordingly \eqref{op-fsystem2}
and \eqref{op-fsystem3} by the following finite
dimensional systems:
\begin{equation}\label{op-fsystem5}
\left\{
\begin{array}{ll}
\ds dx^{n,\l}_1 = (A_{\l,n}+J_n)x^{n,\l}_1ds + u^n_1ds + K_nx^{n,\l}_1 dw(s) + v^n_1 dw(s) &\mbox{ in } (t,T],\\
\ns\ds x^{n,\l}_1(t)=\xi^n_1
\end{array}
\right.
\end{equation}
 and
\begin{equation}\label{op-fsystem6}
\left\{
\begin{array}{ll}
\ds dx^{n,\l}_2 = (A_{\l,n} + J_n)x^{n,\l}_2ds + u^n_2ds + K_n x^{n,\l}_2 dw(s) + v^n_2 dw(s) &\mbox{ in } (t,T],\\
\ns\ds x^{n,\l}_2(t)=\xi^n_2.
\end{array}
\right.
\end{equation}
Here $\xi_1^n=\G_n\xi_1$, $\xi_2^n =
\G_n\xi_2$, $u_1^n(\cd) = \G_n u_1(\cd)$,
$u_2^n(\cd)=\G_n u_2(\cd)$, $v^n_1(\cd) = \G_n
v_1(\cd)$ and $v^n_2(\cd) = \G_n v_2(\cd)$. It
is easy to see that both (\ref{op-fsystem5})
and (\ref{op-fsystem6}) are stochastic
differential equations. Obviously,
$\xi^{n}_1,\xi^{n}_2\in
L^4_{\cF_t}(\O;\dbR^n)$, $u_1^n, u_2^n \in
L^2_\dbF(t,T;L^4(\O;\dbR^n))$ and $v_1^n, v_2^n
\in L^2_\dbF(t,T;L^4(\O;\dbR^n))$. One can
easily check that, for $ k=1,2$,
\begin{equation}\label{eq3.1}
\left\{
\begin{array}{ll}\ds
\lim_{n\to \infty}\xi_k^{n} = \xi_k \mbox{ in
}L^4_{\cF_t}(\O;H),\\
\ns\ds \lim_{n\to \infty}u_k^n = u_k \mbox{ in
}L^2_\dbF(t,T;L^4(\O;H)),\\
\ns\ds \lim_{n\to \infty}v_k^n = v_k \mbox{ in
}L^2_\dbF(t,T;L^4(\O;H)).
\end{array}
\right.
\end{equation}
Then, similar to Lemma \ref{lem2.7}, one can
show that
\begin{equation}\label{eq4xz}
\begin{array}{ll}\ds
\lim_{n\to\infty} x_k^{n,\l} = x_k^{\l} \mbox{
in }L^4_\dbF(\O;C([t,T];H)),\qq k=1,2.
\end{array}
\end{equation}
Hence, by Lemma \ref{lem2.7}, we obtain that
\begin{equation}\label{eq4.1}
\begin{array}{ll}\ds
\lim_{\l\to\infty}\lim_{n\to\infty}  x_k^{n,\l}
= x_k \mbox{ in }C_\dbF([t,T];L^4(\O;H)),\qq
k=1,2.
\end{array}
\end{equation}

Denote by $U_{n,\l}(\cd,\cd)$ the bounded
linear operator such that $
U_{n,\l}(\cd,t)\xi_1^n $ solves the equation
\eqref{op-fsystem5} with $u_1^n=v_1^n=0$.
Clearly, $ U_{n,\l}(\cd,t)\xi_2^n $ solves the
equation \eqref{op-fsystem6} with
$u_2^n=v_2^n=0$. We claim that that for any
$\l\in\rho(A)$, there is a constant $C(\l)>0$
such that for all $n\in\dbN$, it holds that
\begin{equation}\label{op-fsystem5
enest}
|U_{n,\l}(\cd,t)\xi_1^n|_{L^\infty_\dbF(t,T;L^{4}(\O;H))}\leq
C(\l)|\xi_1|_{L^4_{\cF_t}(\O;H)}.
\end{equation}
Indeed, by
$$
x^{n,\l}_1(s) = S_{n,\l}(s-t)\xi_1^n + \int_t^s
S_{n,\l}(s-\tau)J_n(\tau)x_1^{n,\l}(\tau)d\tau
+ \int_t^s
S_{n,\l}(s-\tau)K_n(\tau)x_1^{n,\l}(\tau)dw(\tau),
$$
noting that for all $n\in\dbN$,
$|S_{n,\l}(\cd)|_{L^\infty(0,T;\cL(H))}\leq
e^{|\!|A_\l|\!|_{\cL(H)}T}$ and utilizing Lemma
\ref{BDG}, we obtain that
$$
\begin{array}{ll}\ds
\q\mE|U_{n,\l}(s,t)\xi_1^n|^4_H \\
\ns\ds = \mE\Big| S_{n,\l}(s-t)\xi_1^n +
\int_t^s
S_{n,\l}(s-\tau)J(\tau)U_{n,\l}(\tau,t)\xi_1^n
d\tau + \int_t^s
S_{n,\l}(s-\tau)K(\tau)U_{n,\l}(\tau,t)\xi_1^n dw\Big|^4_{H}\\
\ns\ds  \leq C(\l)\mE|\xi_1|^4_H +
C(\l)\int_t^s\[|J(\tau)|^4_{ L^\infty(\O;
\cL(H))}+|K(\tau)|^4_{ L^\infty(\O; \cL(H))}\]
|U_{n,\l}(\tau,t)\xi_1^n|_H^4d\tau.
\end{array}
$$
This, together with Gronwall's inequality,
implies \eqref{op-fsystem5 enest}.

Also, denote by $U_{\l}(\cd,\cd)$ the bounded
linear operator such that $ U_{\l}(\cd,t)\xi_1
$ solves the equation \eqref{op-fsystem2.1}
with $u_1=v_1=0$. Clearly, $ U_{\l}(\cd,t)\xi_2
$ solves the equation \eqref{op-fsystem3.1}
with $u_2=v_2=0$. We claim that that there is a
constant $C>0$ such that for any $\l\in\rho(A)$
it holds that
\begin{equation}\label{op-fsystem5
enest x}
|U_{\l}(\cd,t)\xi_1|_{L^\infty_\dbF(t,T;L^{4}(\O;H))}\leq
C|\xi_1|_{L^4_{\cF_t}(\O;H)}.
\end{equation}
Indeed, similar to the above proof of
\eqref{op-fsystem5 enest}, by (\ref{zx-s1}) and
utilizing Lemma \ref{BDG}, we obtain that
$$
\begin{array}{ll}\ds
\q\mE|U_{\l}(s,t)\xi_1|^4_H \\
\ns\ds = \mE\Big| S_{\l}(s-t)\xi_1 + \int_t^s
S_{\l}(s-\tau)J(\tau)U_{\l}(\tau,t)\xi_1 d\tau
+ \int_t^s
S_{\l}(s-\tau)K(\tau)U_{\l}(\tau,t)\xi_1 dw\Big|^4_{H}\\
\ns\ds  \leq C\mE|\xi_1|^4_H +
C\int_t^s\[|J(\tau)|^4_{ L^\infty(\O;
\cL(H))}+|K(\tau)|^4_{ L^\infty(\O; \cL(H))}\]
|U_{\l}(\tau,t)\xi_1|_H^4d\tau.
\end{array}
$$
Hence \eqref{op-fsystem5 enest x} follows from
Gronwall's inequality.

\ms

{\bf Step 3}. In this step, we show that
$(P^{n,\l}(\cd),Q^{n,\l}(\cd))$ satisfies a
variational equality, which can be viewed as an
approximation of (\ref{eq def sol1}).

Denote by $X^{n,\l}$ the tensor product of
$x_1^{n,\l}$ and $x_2^{n,\l}$, i.e., $X^{n,\l}
= x_1^{n,\l}\otimes x_2^{n,\l}$.
 Since
$$
\begin{array}{ll}\ds
\q d(x_1^{n,\l}\otimes x_2^{n,\l})\\
\ns\ds = \big(dx_1^{n,\l}\big)\otimes x_2^{n,\l} + x_1^{n,\l}\otimes d\big( x_2^{n,\l}\big) +  \big(dx_1^{n,\l}\big)\otimes d\big( x_2^{n,\l}\big)\\
\ns\ds  = \[\big(A_{\l,n}+J_n\big) x_{1}^{n,\l}\]\otimes x_2^{n,\l}ds + x_1^{n,\l}\otimes \[\big(A_{\l,n}+J_n\big)x_2^{n,\l}\]ds \\
\ns\ds \q + \!\[  u_1^n\otimes x_2^{n,\l}\! +\! x_1^{n,\l}\otimes u_2^n \!+\! \big(K_n x_1^{n,\l}\big)\otimes \big(K_n x_2^{n,\l}\big)\! + \!\big(K_n x_1^{n,\l}\big)\otimes v_2^n \!+\! v_1^n\otimes \big(K_nx_2^{n,\l}\big) \!+\! v_1^n\otimes v_2^n \]ds \\
\ns\ds \q + \[ K_nx_1^{n,\l}\otimes x_2^{n,\l}
+ x_1^{n,\l}\otimes \big(K_n x_2^{n,\l}\big) +
v_1^n\otimes x_2^{n,\l} + x_1^{n,\l}\otimes
v_2^n \]dw(s),
\end{array}
$$
we see that $X^{n,\l}$ solves the following
equation:
\begin{equation}\label{op-fsystem4}
\left\{\3n
\begin{array}{ll}
\ds dX^{n,\l} = \a^{n,\l}ds + \b^{n,\l}dw(s) &\mbox{ in } (t,T],\\
\ns\ds X^{n,\l}(t) = \xi_1^{n}\otimes
\xi_2^{n},
\end{array}
\right.
\end{equation}
where
$$
\left\{
\begin{array}{ll}\ds
\a^{n,\l} = \[\big(A_{\l,n}+J_n\big)
x_{1}^{n,\l}\]\otimes x_2^{n,\l}  +
x_1^{n,\l}\otimes
\[\big(A_{\l,n}+J_n\big)x_2^{n,\l}\]
+   u_1^n\otimes x_2^{n,\l} + x_1^{n,\l}\otimes u_2^n \\
\ns\ds \hspace{1.5cm} + \big(K_n
x_1^{n,\l}\big)\otimes \big(K_n x_2^{n,\l}\big)
+ \big(K_n x_1^{n,\l}\big)\otimes v_2^n +
v_1^n\otimes \big(K_nx_2^{n,\l}\big) +
v_1^n\otimes
v_2^n, \\
\ns\ds \b^{n,\l} = K_nx_1^{n,\l}\otimes
x_2^{n,\l} + x_1^{n,\l}\otimes \big(K_n
x_2^{n,\l}\big) + v_1^n\otimes x_2^{n,\l} +
x_1^{n,\l}\otimes v_2^n.
\end{array}
\right.
$$

Recalling that $(P^{n,\l}(\cd),Q^{n,\l}(\cd))$
is the transposition solution to
(\ref{op-bsystem4}), by (\ref{eq def solz}) and
(\ref{op-fsystem4}), we obtain that
\begin{equation}\label{s5eq01}
\begin{array}{ll}\ds
\ds \q\mE\big\langle x_1^{n,\l}(T)\otimes
x_2^{n,\l}(T),P^{n}_T \big\rangle_{\cL(\dbR^n)}
- \mE \int_t^T \big\langle x_1^{n,\l}(s)\otimes
x_2^{n,\l}(s),\Phi^{n,\l}(s)
\big\rangle_{\cL(\dbR^n)}ds
 \\
 \ns\ds = \big\langle \xi_1\otimes \xi_2,P^{n,\l}(t)
\big\rangle_{\cL(\dbR^n)} + \mE \int_t^T
\big\langle
\a^{n,\l}(s),P^{n,\l}(s) \big\rangle_{\cL(\dbR^n)}ds\\
 \ns\ds\q + \mE \int_t^T
\big\langle \b^{n,\l}(s),Q^{n,\l}(s)
\big\rangle_{\cL(\dbR^n)}ds,
\end{array}
\end{equation}
where $\Phi^{n,\l}(\cd)$ is given by
(\ref{x1x}).

A direct computation shows that
\begin{equation}\label{s5eq02}
\begin{array}{ll}\ds
\q\mE \int_t^T \big\langle x_1^{n,\l}(s)\otimes
x_2^{n,\l}(s),\Phi^{n,\l}(s)
\big\rangle_{\cL(\dbR^n)}ds  = \mE \int_t^T
\big\langle \Phi^{n,\l}(s)x_1^{n,\l}(s),
x_2^{n,\l}(s) \big\rangle_{\dbR^n}ds \\
\ns\ds = -\mE \int_t^T \big\langle P^{n,\l}
x_1^{n,\l}(s), \big(A_{\l,n}\! +
\!J_n\big)x_2^{n,\l}(s) \big\rangle_{\dbR^n}ds
\! -\!\mE \int_t^T \big\langle
P^{n,\l}\big(A_{\l,n}\!+ \! J_n\big)
x_1^{n,\l}(s), x_2^{n,\l}(s)
\big\rangle_{\dbR^n}ds \\
\ns\ds \q - \mE \int_t^T \big\langle
P^{n,\l}K_n x_1^{n,\l}(s), K_nx_2^{n,\l}(s)
\big\rangle_{\dbR^n}ds - \mE \int_t^T
\big\langle Q^{n,\l} x_1^{n,\l}(s),
K_nx_2^{n,\l}(s) \big\rangle_{\dbR^n}ds \\
\ns\ds \q - \mE \int_t^T \big\langle Q^{n,\l}
K_n x_1^{n,\l}(s), x_2^{n,\l}(s)
\big\rangle_{\dbR^n}ds + \mE \int_t^T
\big\langle F_n x_1^{n,\l}(s), x_2^{n,\l}(s)
\big\rangle_{\dbR^n}ds.
\end{array}
\end{equation}
 Next,
\begin{equation}\label{s5eq03}
\begin{array}{ll}\ds
\q\mE \int_t^T \big\langle
\a^{n,\l}(s),P^{n,\l}(s)
\big\rangle_{\cL(\dbR^n)}ds \\
\ns\ds = \mE \int_t^T \big\langle P^{n,\l}(s)
\big(A_{\l,n}+J_n\big) x_{1}^{n,\l}(s),
x_2^{n,\l}(s) \big\rangle_{\dbR^n}ds
\\\ns\ds\q+ \mE \int_t^T \big\langle
P^{n,\l}(s)
 x_{1}^{n,\l}(s),
\big(A_{\l,n} +J_n\big)x_2^{n,\l}(s)
\big\rangle_{\dbR^n}ds \\
\ns\ds \q + \mE \int_t^T \big\langle
P^{n,\l}(s)
 u_1^n(s), x_2^{n,\l}(s)
\big\rangle_{\dbR^n}ds + \mE \int_t^T
\big\langle P^{n,\l}(s) x_1^{n,\l}(s), u_2^n(s)
\big\rangle_{\dbR^n}ds \\
\ns\ds \q + \mE \int_t^T \big\langle
P^{n,\l}(s)K_n x_1^{n,\l}(s), K_n x_2^{n,\l}(s)
\big\rangle_{\dbR^n}ds   + \mE \int_t^T
\big\langle P^{n,\l}(s)K_n x_1^{n,\l}(s),
v_2^n(s)
\big\rangle_{\dbR^n}ds \\
\ns\ds \q + \mE \int_t^T \big\langle
P^{n,\l}(s)v_1^n(s), K_nx_2^{n,\l}
(s)\big\rangle_{\dbR^n}ds   + \mE \int_t^T
\big\langle P^{n,\l}(s)v_1^n(s), v_2^n(s)
\big\rangle_{\dbR^n}ds.
\end{array}
\end{equation}
  Further,
\begin{equation}\label{s5eq04}
\begin{array}{ll}\ds
\q \mE\int_t^T \big\langle
\b^{n,\l}(s),Q^{n,\l}(s)
\big\rangle_{\cL(\dbR^n)}ds \\
\ns\ds = \mE\int_t^T \big\langle
Q^{n,\l}(s)K_nx_1^{n,\l}(s),x_2^{n,\l}
(s)\big\rangle_{\dbR^n}ds + \mE\int_t^T
\big\langle Q^{n,\l}(s) x_1^{n,\l}(s),K_n
x_2^{n,\l}(s) \big\rangle_{\dbR^n}ds \\
\ns\ds \q + \mE\int_t^T \big\langle Q^{n,\l}(s)
v_1^n(s),x_2^{n,\l}(s) \big\rangle_{\dbR^n}ds
+ \mE\int_t^T \big\langle Q^{n,\l}(s)
x_1^{n,\l}(s),v_2^n(s) \big\rangle_{\dbR^n}ds.
\end{array}
\end{equation}

From \eqref{s5eq01}--\eqref{s5eq04}, we arrive
at
\begin{equation}\label{eq2}
\begin{array}{ll}
\ds \q\mE\big\langle P^{n}_T x_1^{n,\l}(T),
x_2^{n,\l}(T) \big\rangle_{\dbR^n} - \mE
\int_t^T \big\langle
F_n(s) x_1^{n,\l}(s), x_2^{n,\l}(s) \big\rangle_{\dbR^n}ds\\
\ns\ds =\mE\big\langle P^{n,\l}(t)
\xi_1^n,\xi_2^n \big\rangle_{\dbR^n} + \mE
\int_t^T \big\langle P^{n,\l}(s)u_1^n(s),
x_2^{n,\l}(s)\big\rangle_{\dbR^n}ds \\
\ns\ds \q   + \mE \int_t^T \big\langle
P^{n,\l}(s)x_1^{n,\l}(s),
u_2^n(s)\big\rangle_{\dbR^n}ds + \mE \int_t^T
\big\langle
P^{n,\l}(s)K_n(s)x_1^{n,\l}(s), v_2^n(s)\big\rangle_{\dbR^n}ds  \\
\ns\ds \q + \mE \int_t^T \big\langle
P^{n,\l}(s)v_1^{n,\l}(s),
K_n(s)x_2^{n,\l}(s)+v_2^n(s) \big\rangle_{\dbR^n}ds  \\
\ns\ds \q + \mE \int_t^T \big\langle
Q^{n,\l}(s) v_1^n(s),
x_2^{n,\l}(s)\big\rangle_{\dbR^n}ds+ \mE
\int_t^T \big\langle Q^{n,\l}(s) x_1^{n,\l}(s),
v_2^n(s) \big\rangle_{\dbR^n}ds.
\end{array}
\end{equation}

From the above $\dbR^{n\times n}$-valued
processes $P^{n,\l}(\cd)$ and $Q^{n,\l}(\cd)$,
one obtains two $\cL(H)$-valued processes
$P^{n,\l}(\cd)\G_n$ and $Q^{n,\l}(\cd)\G_n$. To
simply the notations, we simply identify
$P^{n,\l}(\cd)$ (\resp $Q^{n,\l}(\cd)$) and
$P^{n,\l}(\cd)\G_n$ (resp.
$Q^{n,\l}(\cd)\G_n$).

\ms

{\bf Step 4}. In this step, we take
$n\to\infty$ in (\ref{eq2}) with $t\in
\{r_j\}_{j=1}^\infty$. Here
$\{r_j\}_{j=1}^\infty$ stands for the subset of
all rational numbers in $[0,T]$.

In the sequel, we fix a sequence
$\{\l_m\}_{m=1}^\infty\subset \rho(A)$ such
that $\l_m\to\infty$ as $m\to\infty$.

Choose $u^n_1=v^n_1=0$ and $u^n_2=v^n_2=0$ in
\eqref{op-fsystem5} and \eqref{op-fsystem6},
respectively. From the equality \eqref{eq2}, it
follows that
\begin{equation}\label{eq201}
\mE\big\langle P^{n}_T x_1^{n,\l_m}(T),
x_2^{n,\l_m}(T) \big\rangle_{\dbR^n} - \mE
\int_t^T \big\langle F_n(s) x_1^{n,\l_m}(s),
x_2^{n,\l_m}(s) \big\rangle_{\dbR^n}ds
=\mE\big\langle P^{n,\l_m}(t) \xi_1^n,\xi_2^n
\big\rangle_{\dbR^n}.
\end{equation}
Combing \eqref{op-fsystem5 enest} and
\eqref{eq201}, we find that
\begin{equation}\label{eq202}
\ba{ll} \ds\q\Big|\mE\big\langle P^{n,\l_m}(t)
\xi_1^n,\xi_2^n
\big\rangle_H\Big|=\Big|\mE\big\langle
P^{n,\l_m}(t) \xi_1^n,\xi_2^n
\big\rangle_{\dbR^n}\Big|\\\ns\ds\leq
C(\l_m)\big(|P_T|_{L^2_{\cF_T}(\O;\cL(H))} +
|F|_{L^1_\dbF(0,T;\;L^2(\O;\cL(H)))}
\big)|\xi_1|_{L^4_{\cF_t}(\O;H)}|\xi_2|_{L^4_{\cF_t}(\O;H)}.
 \ea
\end{equation}
Here and henceforth $C(\l_m)$ denotes a generic
constant depending only on $\l_m$, independent
of $n$.

For $P^{n,\l_m}(t)$, we can find a
$\xi_{1,n,m}\in L^4_{\cF_t}(\O;H)$ with
$|\xi_{1,n,m}|_{L^4_{\cF_t}(\O;H)}=1$ such that
\begin{equation}\label{eq203}
\big|P^{n,\l_m}(t)\xi_{1,n,m}\big|_{L^{\frac{4}{3}}_{\cF_t}(\O;H)}\geq
\frac{1}{2}\big|P^{n,\l_m}(t)
\big|_{L^{2}_{\cF_t}(\O;\cL(H))}.
\end{equation}
 Moreover, we
can find a $\xi_{2,n,m}\in L^4_{\cF_t}(\O;H)$
with $|\xi_{2,n,m}|_{L^4_{\cF_t}(\O;H)}=1$ such
that
\begin{equation}\label{eq204}
\mE\big\langle P^{n,\l_m}(t)
\xi_{1,n,m},\xi_{2,n,m}
\big\rangle_{\dbR^n}\geq
\frac{1}{2}\big|P^{n,\l_m}(t)\xi_1^{n,m}\big|_{L^{\frac{4}{3}}_{\cF_t}(\O;H)}.
\end{equation}
From \eqref{eq202}--\eqref{eq204}, we obtain
that
\begin{equation}\label{eq205}
|P^{n,\l_m}|_{L^\infty_\dbF(0,T;L^2(\O;\cL(H)))}
\leq C(\l_m)\big(|P_T|_{L^2_{\cF_T}(\O;\cL(H))}
+ |F|_{L^1_\dbF(0,T;L^2(\O;\cL(H)))}
\big),\q\;\forall\;n\in\dbN.
\end{equation}
Thanks to Theorem \ref{lemma6.1}, one can find
a $
P^{\l_m}\in\cL_{pd}\big(L^2_\dbF(0,T;L^4(\O;H)),\;L^2_\dbF(0,T;L^{\frac{4}{3}}(\O;H))\big)$
such that
\begin{equation}\label{eq204*}
|\!|
P^{\l_m}|\!|_{\cL(L^2_\dbF(0,T;L^4(\O;H)),\;
L^2_\dbF(0,T;L^{\frac{4}{3}}(\O;H)))} \leq
C(\l_m)\big(|P_T|_{L^2_{\cF_T}(\O;\cL(H))} +
|F|_{L^1_\dbF(0,T;L^2(\O;\cL(H)))} \big),
 \ee
and a subsequence
$\{n^{(1)}_{k}\}_{k=1}^\infty\subset
\{n\}_{n=1}^\infty$  such that
\begin{equation}\label{s5oeq3}
\mbox{(w)}\mbox{-}\lim_{k\to\infty}P^{n^{(1)}_k,\l_{m}}u
= P^{\l_m} u \ \mbox{ in }\
L^2_\dbF(0,T;L^{\frac{4}{3}}(\O;H)),\qq
\forall\; u\in L^2_\dbF(0,T;L^4(\O;H)).
\end{equation}
Note that, by means of the standard
diagonalisation argument, one can choose the
subsequence $\{n^{(1)}_{k}\}_{k=1}^\infty$ to
be independent of $\l_m$.

Next, by Theorem \ref{lemma6--1}, for each
$r_j$ and $\l_m$, there exist an $
R^{(r_j,\l_m)}\in\cL_{pd}\big(\!L^4_{\cF_{r_j}}\!(\O;H),\,L^{\frac{4}{3}}_{\cF_{r_j}}\!(\O;H)\!\big)$
and a subsequence
$\{n^{(2)}_{k}\}_{k=1}^\infty\subset
\{n^{(1)}_{k}\}_{k=1}^\infty$ such that
\begin{equation}\label{s5oeq3.1}
\mbox{(w)}\mbox{-}\lim_{k\to\infty}P^{n^{(2)}_k,\l_m}(r_j)
\xi = R^{(r_j,\l_m)}\xi \ \ \ \mbox{ in }
 L^{\frac{4}{3}}_{\cF_{r_j}}(\O;H),\qq \forall\;
\xi\in L^4_{\cF_{r_j}}(\O;H).
\end{equation}
Here, again, by the diagonalisation argument,
one can choose the subsequence
$\{n^{(2)}_{k}\}_{k=1}^\infty$ to be
independent of $r_j$ and $\l_m$.

Let $u^n_1=v^n_1=0$ and $\xi^n_2=0$, $u^n_2=0$
in \eqref{op-fsystem5} and \eqref{op-fsystem6},
respectively. From \eqref{eq2}, we find that
$$
\begin{array}{ll}\ds
\q\mE \int_t^T \big\langle
Q^{n,\l_m}(s)U_{n,\l_m}(s,t) \xi_1^n, v_2^n(s)
\big\rangle_Hds \\
\ns\ds =\mE \int_t^T \big\langle
Q^{n,\l_m}(s)U_{n,\l_m}(s,t) \xi_1^n, v_2^n(s)
\big\rangle_{\dbR^n}ds
\\
\ns\ds =\mE\big\langle P^{n}_T x_1^{n,\l_m}(T),
x_2^{n,\l_m}(T) \big\rangle_{\dbR^n} - \mE
\int_t^T \big\langle F_n(s)
x_1^{n,\l_m}(s), x_2^{n,\l_m}(s) \big\rangle_{\dbR^n}ds\\
\ns\ds\q- \mE \int_t^T \big\langle
P^{n,\l_m}(s)K_n(s)x_1^{n,\l_m}(s),
v_2^n(s)\big\rangle_{\dbR^n}ds.
\end{array}
$$
This implies that
\begin{equation}\label{eq301}
\begin{array}{ll}\ds
\q\mE \int_t^T \big\langle
Q^{n,\l_m}(s)U_{n,\l_m}(s,t) \xi_1^n, v_2^n(s)
\big\rangle_Hds
\\
\ns\ds\leq
C(\l_m)\big(|P_T|_{L^2_{\cF_T}(\O;\cL(H))}+|F|_{L^1_\dbF(0,T;L^2(\O;\cL(H)))}\big)
|\xi_1|_{L^4_{\cF_t}(\O;H)}|v_2|_{L^2_\dbF(t,T;L^4(\O;H))}.
\end{array}
\end{equation}

We define two operators $Q_1^{n,\l_m,t}$ and
$\widehat Q_1^{n,\l_m,t}$ from
$L^4_{\cF_t}(\O;H)$ to
$L^{2}_\dbF(t,T;L^{\frac{4}{3}}(\O;H))$ as
follows:
$$
\left\{
\begin{array}{ll}\ds
Q_1^{n,\l_m,t}\xi=Q^{n,\l_m}(\cd)U_{n,\l_m}(\cd,t)
\xi^n, \q \forall\,\xi\in L^4_{\cF_t}(\O;H);\\
\ns\ds \widehat  Q_1^{n,\l_m,t}\xi =
Q^{n,\l_m}(\cd)^*U_{n,\l_m}(\cd,t) \xi^n, \q
\forall\,\xi\in L^4_{\cF_t}(\O;H).
\end{array}
\right.
$$
Here $\xi^n = \G_n \xi$. It is clear that
$Q_1^{n,\l_m,t},\widehat Q_1^{n,\l_m,t}\in
\cL(L^4_{\cF_t}(\O;H),\;
L^{2}_\dbF(t,T;L^{\frac{4}{3}}(\O;H)))$. For
any given $n$, $\l_m$ and $t$, we can find a
$\xi_1^{n,m,t}\in L^4_{\cF_t}(\O;H)$ with
$|\xi_1^{n,m,t}|_{L^4_{\cF_t}(\O;H)}=1$, such
that
\begin{equation}\label{eq302}
\big|Q_1^{n,\l_m,t}\xi_1^{n,m,t}\big|_{L^{2}_\dbF(t,T;L^{\frac{4}{3}}(\O;H))}\geq
\frac{1}{2}\big|\!\big| Q_1^{n,\l_m,t}
\big|\!\big|_{\cL(L^4_{\cF_t}(\O;H),\;
L^{2}_\dbF(t,T;L^{\frac{4}{3}}(\O;H)))}.
\end{equation}
Furthermore, we can find a $v_2^{n,m,t}(\cd)\in
L^2_\dbF(t,T;L^4(\O;H))$ with
$|v_2^{n,m,t}(\cd)|_{
L^2_\dbF(t,T;L^4(\O;H))}=1$ such that
\begin{equation}\label{eq303}
\mE \int_t^T \big\langle
Q^{n,\l_m,t}(s)U_{n,\l_m}(s,t) \xi_1^{n,m,t},
v_2^{n,m,t}(s) \big\rangle_Hds \geq
\frac{1}{2}\big| Q_1^{n,\l_m,t}
\xi_1^{n,m,t}\big|_{L^{2}_\dbF(t,T;L^{\frac{4}{3}}(\O;H))}.
\end{equation}
Hence, combining (\ref{eq301}), (\ref{eq302})
and (\ref{eq303}), it follows that
\begin{equation}\label{eq304}
\big|\!\big| Q_1^{n,\l_m,t}
\big|\!\big|_{\cL(L^4_{\cF_t}(\O;H),\;
L^{2}_\dbF(t,T;L^{\frac{4}{3}}(\O;H)))} \leq
C(\l_m)
\big(|P_T|_{L^2_{\cF_T}(\O;\cL(H))}+|F|_{L^1_\dbF(0,T;L^2(\O;\cL(H)))}\big).
\end{equation}
Similarly,
\begin{equation}\label{eq304x}
\big|\!\big|\widehat
Q_1^{n,\l_m,t}\big|\!\big|_{\cL(L^4_{\cF_t}(\O;H),\;
L^{2}_\dbF(t,T;L^{\frac{4}{3}}(\O;H)))} \leq
C(\l_m)
\big(|P_T|_{L^2_{\cF_T}(\O;\cL(H))}+|F|_{L^1_\dbF(0,T;L^2(\O;\cL(H)))}\big).
\end{equation}
By Lemma \ref{lemma6}, for each $r_j$ and
$\l_m$, there  exist two bounded linear
operators $Q^{\l_m,r_j}_1$ and $\widehat
Q^{\l_m,r_j}_1$, from $L^4_{\cF_{r_j}}(\O;H)$
to
$L^{\frac{4}{3}}_\dbF(r_j,T;L^{\frac{4}{3}}(\O;H))$,
and a subsequence
$\{n^{(3)}_{k}\}_{k=1}^\infty\subset
\{n^{(2)}_{k}\}_{n=1}^\infty$, independent of
$r_j$ and $\l_m$, such that
\begin{equation}\label{eq305}
\left\{
\begin{array}{ll}\ds
\mbox{(w)}\mbox{-}\lim_{k\to\infty}Q_1^{n^{(3)}_k,\l_{m},r_j}\xi
= Q^{\l_m,r_j}_1 \xi \ \ \ \mbox{ in }
L^{2}_\dbF({r_j},T;L^{\frac{4}{3}}(\O;H)),\qq
\forall\; \xi\in L^4_{\cF_{r_j}}(\O;H),\\
\ns\ds
\mbox{(w)}\mbox{-}\lim_{k\to\infty}\widehat
Q_1^{n^{(3)}_k,\l_{m},r_j}\xi = \widehat
Q^{\l_m,r_j}_1 \xi \ \ \ \mbox{ in }
L^{2}_\dbF({r_j},T;L^{\frac{4}{3}}(\O;H)),\qq
\forall\; \xi\in L^4_{\cF_{r_j}}(\O;H).
\end{array}
\right.
\end{equation}

Next, we choose $\xi_1^n=0$, $v_1^n=0$ in
\eqref{op-fsystem5} and $\xi_2^n=0$, $u_2^n=0$
in \eqref{op-fsystem6}. From \eqref{eq2}, we
obtain that
\begin{equation}\label{eq401}
\begin{array}{ll}
\ds \q\mE\big\langle P^{n}_T x_1^{n,\l_m}(T),
x_2^{n,\l_m}(T) \big\rangle_{\dbR^n} - \mE
\int_t^T \big\langle
F_n(s) x_1^{n,\l_m}(s), x_2^{n,\l_m}(s) \big\rangle_{\dbR^n}ds\\
\ns\ds =  \mE \int_t^T \big\langle
P^{n,\l_m}(s)u_1^n(s),
x_2^{n,\l_m}(s)\big\rangle_{\dbR^n}ds + \mE
\int_t^T \big\langle
P^{n,\l_m}(s)K_n(s)x_1^{n,\l_m}(s), v_2^n(s)\big\rangle_{\dbR^n}ds  \\
\ns\ds \q   + \mE \int_t^T \big\langle
Q^{n,\l_m}(s) x_1^{n,\l_m}(s), v_2^n(s)
\big\rangle_{\dbR^n}ds.
\end{array}
\end{equation}
Define an operator $Q_2^{n,\l_m,t}$ from
$L^2_\dbF(t,T;L^4(\O;H))$ to
$L^{2}_{\dbF}(t,T;L^{\frac{4}{3}}(\O;H))$ as
follows:
$$
\big(Q_2^{n,\l_m,t} u\big)(\cd) =
Q^{n,\l_m}(\cd)\int_t^{\cd}U_{n,\l_m}(\cd,\tau)u^n(\tau)d\tau,
\q \forall\,u\in L^2_\dbF(t,T;L^4(\O;H)),
$$
where $u^n = \G_n u$. From \eqref{eq401}, we
get that
\begin{equation}\label{eq402}
\begin{array}{ll}\ds
\q\mE \int_t^T \big\langle \big(Q_2^{n,\l_m,t}
u_1^n\big)(s),v_2^n(s) \big\rangle_Hds=\int_t^T
\big\langle \big(Q_2^{n,\l_m,t}
u_1^n\big)(s),v_2^n(s) \big\rangle_{\dbR^n}ds
\\\ns \ds= \mE\big\langle P^{n}_T
x_1^{n,\l_m}(T), x_2^{n,\l_m}(T)
\big\rangle_{\dbR^n} - \mE \int_t^T \big\langle
F_n(s) x_1^{n,\l_m}(s),
x_2^{n,\l_m}(s) \big\rangle_{\dbR^n}ds\\
\ns\ds \q - \mE \int_t^T \big\langle
P^{n,\l_m}(s)u_1^n(s),
x_2^{n,\l_m}(s)\big\rangle_{\dbR^n}ds -\mE
\int_t^T \big\langle
P^{n,\l_m}(s)K_n(s)x_1^{n,\l_m}(s),
v_2^n(s)\big\rangle_{\dbR^n}ds
\\
\ns\ds \leq C(\l_m)\big(
|P_T|_{L^2_{\cF_T}(\O;\cL(H))} +
|F|_{L^1_\dbF(0,T;L^2(\O;\cL(H)))}\big)|u_1|_{L^2_\dbF(t,T;L^4(\O;H))}
|v_2|_{L^2_\dbF(t,T;L^4(\O;H))}.
\end{array}
\end{equation}

Let us choose a $u_1^{n,m,t}\in
L^2_\dbF(t,T;L^4(\O;H))$ satisfying
$|u_1^{n,m,t}|_{L^2_\dbF(t,T;L^4(\O;H))}=1$,
such that
\begin{equation}\label{eq403}
\big|Q_2^{n,\l_m,t}
u_1^{n,m,t}\big|_{L^{2}_{\dbF}(t,T;L^{\frac{4}{3}}(\O;H))}
\geq
\frac{1}{2}\big|\!\big|Q_2^{n,\l_m,t}\big|\!\big|_{\cL(L^2_\dbF(t,T;L^4(\O;H)),\
L^{2}_{\dbF}(t,T;L^{\frac{4}{3}}(\O;H)))}.
\end{equation}
Then we choose a $v_2^{n,m,t}\in
L^2_\dbF(t,T;L^4(\O;H))$ satisfying
$|v_2^{n,m,t}|_{L^2_\dbF(t,T;L^4(\O;H))}=1$,
such that
\begin{equation}\label{eq404}
\mE \int_t^T \big\langle \big(Q_2^{n,\l_m,t}
u_1^{n,m,t}\big)(s),v_2^{n,m,t}(s)
\big\rangle_Hds\geq
\frac{1}{2}\big|Q_2^{n,\l_m,t}
u_1^{n,m,t}\big|_{L^{2}_{\dbF}(t,T;L^{\frac{4}{3}}(\O;H))}.
\end{equation}
From \eqref{eq402}--\eqref{eq404}, we see that
\begin{equation}\label{eq405}
\big|\!\big|Q_2^{n,\l_m,t}\big|\!\big|_{\cL(L^2_\dbF(t,T;L^4(\O;H)),\
L^{\frac{4}{3}}_{\dbF}(t,T;L^{\frac{4}{3}}(\O;H)))}\leq
C(\l_m)\big( |P_T|_{L^2_{\cF_T}(\O;\cL(H))} +
|F|_{L^1_\dbF(0,T;L^2(\O;\cL(H)))}\big).
\end{equation}
Also, we define an operator $\widehat
Q_2^{n,\l_m,t}$ from $L^2_\dbF(t,T;L^4(\O;H))$
to $L^{2}_{\dbF}(t,T;L^{\frac{4}{3}}(\O;H))$ as
follows:
$$
\big(\widehat Q_2^{n,\l_m,t} u\big)(\cd) =
Q^{n,\l_m}(\cd)^*\int_t^{\cd}U_{n,\l_m}(\cd,\tau)u^n(\tau)d\tau,
\q \forall\,u\in L^2_\dbF(t,T;L^4(\O;H)),
$$
where $u^n = \G_n u$. By a similar argument to
derive the inequality \eqref{eq405}, we find
that
\begin{equation}\label{eq405x}
\big|\!\big|\widehat
Q_2^{n,\l_m,t}\big|\!\big|_{\cL(L^2_\dbF(t,T;L^4(\O;H)),\
L^{\frac{4}{3}}_{\dbF}(t,T;L^{\frac{4}{3}}(\O;H)))}\leq
C(\l_m)\big( |P_T|_{L^2_{\cF_T}(\O;\cL(H))} +
|F|_{L^1_\dbF(0,T;L^2(\O;\cL(H)))}\big).
\end{equation}

By Lemma \ref{lemma6}, we conclude that, for
each $r_j$ and $\l_m$, there exist two bounded
linear operators $Q^{\l_m,r_j}_2$ and $\widehat
Q^{\l_m,r_j}_2$ from
$L^2_\dbF(r_j,T;L^4(\O;H))$ to
$L^{2}_\dbF(r_j,T;L^{\frac{4}{3}}(\O;H))$ and a
subsequence
$\{n^{(4)}_{k}\}_{k=1}^\infty\subset
\{n^{(3)}_{k}\}_{n=1}^\infty$, independent of
$r_j$ and $\l_m$, such that
\begin{equation}\label{eq406}
\left\{
\begin{array}{ll}\ds
\mbox{(w)}\mbox{-}\lim_{k\to\infty}Q_2^{n^{(4)}_k,\l_{m},r_j}u
= Q^{\l_m,r_j}_2 u  \ \mbox{ in }
L^{2}_\dbF(r_j,T;L^{\frac{4}{3}}(\O;H)),\;\
\forall\; u\in L^2_\dbF(r_j,T;L^4(\O;H)),\\
\ns\ds
\mbox{(w)}\mbox{-}\lim_{k\to\infty}\widehat
Q_2^{n^{(4)}_k,\l_{m},r_j}u = \widehat
Q^{\l_m,r_j}_2 u \ \mbox{ in }
L^{2}_\dbF(r_j,T;L^{\frac{4}{3}}(\O;H)),\;\
\forall\; u\in L^2_\dbF(r_j,T;L^4(\O;H)).
\end{array}
\right.
\end{equation}

Now, we choose $\xi_1^n=0$ and $u_1^n=0$ in
\eqref{op-fsystem5}, and $\xi_2^n=0$ and
$u_2^n=0$ in \eqref{op-fsystem6}. From
\eqref{eq2}, we obtain that
\begin{equation}\label{eq501x}
\begin{array}{ll}
\ds \q\mE\big\langle P^{n}_T x_1^{n,\l_m}(T),
x_2^{n,\l_m}(T) \big\rangle_{\dbR^n} - \mE
\int_t^T \big\langle
F_n(s) x_1^{n,\l_m}(s), x_2^{n,\l_m}(s) \big\rangle_{\dbR^n}ds\\
\ns\ds =   \mE \int_t^T \!\!\big\langle
P^{n,\l_m}(s)K_n(s)x_1^{n,\l_m}(s),
v_2^n(s)\big\rangle_{\dbR^n}ds\! + \!\mE
\int_t^T \!\!\big\langle
P^{n,\l_m}(s)v_1^{n}(s),
K_n(s)x_2^{n,\l_m}(s)\!+\!v_2^n(s) \!\big\rangle_{\dbR^n}ds  \\
\ns\ds \q +  \mE \int_t^T \big\langle
Q^{n,\l_m}(s)v_1^n(s),
x_2^{n,\l_m}(s)\big\rangle_{\dbR^n}ds+ \mE
\int_t^T \big\langle Q^{n,\l_m}(s)
x_1^{n,\l_m}(s), v_2^n(s)
\big\rangle_{\dbR^n}ds.
\end{array}
\end{equation}
Let us define a bilinear functional
$\cB_{n,\l_{m},t}(\cd,\cd)$ on
$L^2_\dbF(t,T;L^4(\O;H))\times
L^2_\dbF(t,T;L^4(\O;H))$ as follows:
 \bel{jiushi3}
\begin{array}{ll}\ds
\q \cB_{n,\l_{m},t}(v_1,v_2)\\
\ns\ds =\mE \int_t^T \big\langle
Q^{n,\l_{m}}(s)v_1^n(s),
x_2^{n,\l_{m}}(s)\big\rangle_{\dbR^n}ds + \mE
\int_t^T \big\langle Q^{n,\l_{m}}(s)
x_1^{n,\l_{m}}(s), v_2^n(s)
\big\rangle_{\dbR^n}ds,\\
\ns\ds \hspace{9cm} \forall\, v_1,v_2\in
L^2_\dbF(t,T;L^4(\O;H)).
\end{array}
 \ee
It is easy to check that
$\cB_{n,\l_{m},t}(\cd,\cd)$ is a bounded
bilinear functional. From \eqref{eq501x}, it
follows that
\begin{equation}\label{eq502}
\begin{array}{ll}
\ds \q
\cB_{{n_{k}^{(4)}},\l_{m},t}(v_1,v_2)\\
\ns\ds = \mE\Big\langle P^{{n_{k}^{(4)}}}_T
x_1^{{n_{k}^{(4)}},\l_{m}}(T),
x_2^{{n_{k}^{(4)}},\l_{m}}(T)
\Big\rangle_{\dbR^{n_{k}^{(4)}}} - \mE \int_t^T
\Big\langle
F_{n_{k}^{(4)}}(s) x_1^{{n_{k}^{(4)}},\l_{m}}(s), x_2^{{n_{k}^{(4)}},\l_{m}}(s) \Big\rangle_{\dbR^{n_{k}^{(4)}}}ds\\
\ns\ds \q -  \mE \int_t^T  \Big\langle
P^{{n_{k}^{(4)}},\l_{m}}(s)K_{n_{k}^{(4)}}(s)x_1^{{n_{k}^{(4)}},\l_{m}}(s),
v_2^{n_{k}^{(4)}}(s)\Big\rangle_{\dbR^{n_{k}^{(4)}}}ds\\
\ns\ds \q - \mE \int_t^T \Big\langle
P^{{n_{k}^{(4)}},\l_{m}}(s)v_1^{{n_{k}^{(4)}}}(s),
K_{n_{k}^{(4)}}(s)x_2^{{n_{k}^{(4)}},\l_{m}}(s)
+ v_2^{n_{k}^{(4)}}(s)
\Big\rangle_{\dbR^{n_{k}^{(4)}}}ds.
\end{array}
\end{equation}

From the definition of $P^{\l_m}$,
$x_1^{{n_{k}^{(4)}}}$ and
$x_2^{{n_{k}^{(4)}}}$, we find that
$$
\left\{
\begin{array}{ll}\ds
\lim_{k\to\infty} \mE\big\langle
P^{{n_{k}^{(4)}}}_T
x_1^{n_{k}^{(4)},\l_{m}}(T),
x_2^{n_{k}^{(4)},\l_{m}}(T)
\big\rangle_{\dbR^{n_{k}^{(3)}}} =
\mE\big\langle P_T x_1^{\l_m} (T), x_2^{\l_m}
(T)
\big\rangle_{H},\\
\ns\ds \lim_{k\to\infty}  \mE \int_t^T
\big\langle F_{n_{k}^{(3)}}(s)
x_1^{n_{k}^{(4)},\l_{m}}(s),
x_2^{n_{k}^{(4)},\l_{m}}(s)
\big\rangle_{\dbR^{n_{k}^{(4)}}}ds = \mE
\int_t^T \big\langle F(s) x_1^{\l_m} (s),
x_2^{\l_m} (s)
\big\rangle_{H}ds,\\
\ns\ds \!\!\lim_{k\to\infty}  \mE
\!\int_t^T\!\! \big\langle
P^{{n_{k}^{(4)}},\l_{m}}(s)K_{n_{k}^{(4)}}(s)x_1^{{n_{k}^{(4)}},\l_{m}}(s),\!
v_2^{n_{k}^{(4)}}(s)\big\rangle_{\dbR^{n_{k}^{(4)}}}ds
\!= \! \mE \!\int_t^T \!\!\big\langle
P^{\l_m}(s)K(s)x_1^{\l_m} (s),\!
v_2(s)\big\rangle_{H}ds,\\
\ns\ds \lim_{k\to\infty}  \mE \int_t^T
 \big\langle P^{{n_{k}^{(4)}},\l_{m}}(s)v_1^{{n_{k}^{(4)}}}(s),
K_{n_{k}^{(4)}}(s)x_2^{{n_{k}^{(4)}},\l_{m}}(s)
+ v_2^{n_{k}^{(4)}}(s)
 \big\rangle_{\dbR^{n_{k}^{(4)}}}ds \\
 \ns\ds  \q=  \mE \int_t^T
 \big\langle P^{\l_m}(s)v_1 (s),
K(s)x_2^{\l_m}(s) + v_2(s)
 \big\rangle_{H}ds,
\end{array}
\right.
$$
where $x_1 $ (\resp $x_2 $) solves the equation
\eqref{op-fsystem2} (\resp \eqref{op-fsystem3})
with $\xi_1=0$ and $u_1=0$ (\resp $\xi_2=0$ and
$u_2=0$). This, together with \eqref{eq502},
implies that
\begin{equation}\label{eq503}
\3n\!  \begin{array}{ll} \ds \q \cB^t_{\l_{m}}
(v_1,v_2)\=
\lim_{k\to\infty}\cB_{{n_{k}^{(4)}},\l_{m},t}(v_1,v_2)\\
\ns\ds = \mE\big\langle P_T x^{\l_{m}}_1 (T),
x^{\l_{m}}_2 (T) \big\rangle_{H} - \mE \int_t^T
\big\langle
F(s) x^{\l_{m}}_1 (s), x^{\l_{m}}_2 (s) \big\rangle_{H}ds\\
\ns\ds \q -  \mE \int_t^T \big\langle
P^{\l_{m}}(s)K(s)x^{\l_{m}}_1(s),
v_2(s)\big\rangle_{H}ds  - \mE \int_t^T
\big\langle P^{\l_{m}}(s)v_1(s),
K(s)x^{\l_{m}}_2(s) + v_2(s) \Big\rangle_{H}ds.
\end{array}
\end{equation}

Noting that the solution of
(\ref{op-fsystem2.1}) (with $\xi_1=0$, $u_1=0$
and $\l$ replaced by $\l_m$) satisfies
$$
\begin{array}{ll}\ds
x_1^{\l_{m}}(s)\!=\!\int_t^s
S_{\l_{m}}(s-\tau)J(\tau)x_1^{\l_{m}}(\tau)d\tau
\!+\! \int_t^s
S_{\l_{m}}(s-\tau)K(\tau)x_1^{\l_{m}}(\tau)d\tau
\!+\! \int_t^s
S_{\l_{m}}(s-\tau)v(\tau)dw(\tau),
\end{array}
$$
by means of Lemma \ref{BDG} and Gronwall's
inequality, we conclude that
\begin{equation}\label{eq504}
|x_1^{\l_{m}} |_{L^\infty_\dbF(t,T;L^4(\O;H))}
\leq C|v_1|_{L^2_\dbF(t,T;L^4(\O;H))}.
\end{equation}
Similarly,
\begin{equation}\label{eq505}
|x_2^{\l_{m}}|_{L^\infty_\dbF(t,T;L^4(\O;H))}
\leq C|v_2|_{L^2_\dbF(t,T;L^4(\O;H))}.
\end{equation}
Combining \eqref{eq503}, \eqref{eq504},
\eqref{eq505} and \eqref{eq204*}, we obtain
that
$$
\begin{array}{ll}\ds
|\cB^t_{\l_{m}}(v_1,v_2)|\leq C(\l_m)\big(
|P_T|_{L^2_{\cF_T}(\O;\cL(H))} +
|F|_{L^1_\dbF(0,T;L^2(\O;\cL(H)))}\big)|v_1|_{L^2_\dbF(t,T;L^4(\O;H))}
|v_2|_{L^2_\dbF(t,T;L^4(\O;H))}.
\end{array}
$$
Hence, $\cB^t_{\l_{m}}(\cd,\cd)$ is a bounded
bilinear functional on
$L^2_\dbF(t,T;L^4(\O;H))\times
L^2_\dbF(t,T;L^4(\O;H))$. Now, for any fixed
$v_2\in L^2_\dbF(t,T;L^4(\O;H))$,
$\cB^t_{\l_{m}}(\cd,v_2)$ is a bounded  linear
functional on $L^2_\dbF(t,T;L^4(\O;H))$.
Therefore, by Lemma \ref{lemma1}, we can find a
unique $\tilde v_1\in
L^{2}_\dbF(t,T;L^{\frac{4}{3}}(\O;H))$ such
that
$$
\cB^t_{\l_{m}}(v_1,v_2)=\big\langle \tilde v_1,
v_2
\big\rangle_{L^{2}_\dbF(t,T;L^{\frac{4}{3}}(\O;H)),\,L^2_\dbF(t,T;L^4(\O;H))},\qq
\forall\; v_2\in L^2_\dbF(t,T;L^4(\O;H)).
$$
Define an operator $\widehat Q_3^{\l_m,t} $
from $L^2_\dbF(t,T;L^4(\O;H))$ to
$L^{2}_\dbF(t,T;L^{\frac{4}{3}}(\O;H))$ as
follows:
$$
\widehat Q_3^{\l_m,t}\, v_1 = \tilde v_1.
$$
From the uniqueness of $\tilde v_1$, it is
clear that $\widehat Q_3^{\l_m,t}$ is
well-defined. Further,
$$
\begin{array}{ll}\ds
\q|\widehat Q_3^{\l_m,t}
v_1|_{L^{2}_\dbF(t,T;L^{\frac{4}{3}}(\O;H))}
=|\tilde
v_1|_{L^{2}_\dbF(t,T;L^{\frac{4}{3}}(\O;H))}
\\
\ns\ds\leq C(\l_m)\big(
|P_T|_{L^2_{\cF_T}(\O;\cL(H))} +
|F|_{L^1_\dbF(0,T;L^2(\O;\cL(H)))}\big)
|v_1|_{L^2_\dbF(t,T;L^4(\O;H))}.
\end{array}
$$
This shows that $\widehat Q_3^{\l_m,t}$ is a
bounded operator. For any $\a,\b\in\dbR$ and
$v_2,v_3,v_4\in L^2_\dbF(t,T;L^4(\O;H))$,
$$
\begin{array}{ll}\ds
\q\big\langle\widehat Q_3^{\l_m,t}(\a v_3 + \b
v_4),v_2\big\rangle_{L^{2}_\dbF(t,T;L^{\frac{4}{3}}(\O;H)),L^2_\dbF(t,T;L^4(\O;H))}\\\ns\ds
= \cB^t_{\l_{m}}(\a v_3 + \b v_4,v_2) =
\a\cB^t_{\l_{m}}( v_3 ,v_2) + \b
\cB^t_{\l_{m}}(v_4,v_2),
\end{array}
$$
which indicates that $ \widehat Q_3^{\l_m,t}(\a
v_3 + \b v_4) = \a \widehat Q_3^{\l_m,t}v_3 +
\b \widehat Q_3^{\l_m,t} v_4$. Hence, $\widehat
Q_3^{\l_m,t}$ is a bounded linear operator from
$L^2_\dbF(t,T;L^4(\O;H))$ to
$L^{2}_\dbF(t,T;L^{\frac{4}{3}}(\O;H))$. Put
$Q_3^{\l_m,t} = \frac{1}{2}\widehat
Q_3^{\l_m,t}$. Then, for any $v_1, v_2\in
L^2_\dbF(t,T;L^4(\O;H))$, it holds that
\begin{equation}\label{eqbilinear1}
\begin{array}{ll}\ds
\q\cB^t(v_1,v_2)\\
\ns\ds = \big\langle Q_3^{\l_m,t} v_1, v_2
\big\rangle_{L^{2}_\dbF(t,T;L^{\frac{4}{3}}(\O;H)),L^2_\dbF(t,T;L^4(\O;H))}
+ \big\langle  v_1, \big(Q_3^{\l_m,t}\big)^*
v_2
\big\rangle_{L^2_\dbF(t,T;L^4(\O;H)),L^{2}_\dbF(t,T;L^{\frac{4}{3}}(\O;H))}.
\end{array}
\end{equation}

From \eqref{eq2}, \eqref{s5oeq3},
\eqref{s5oeq3.1}, \eqref{eq305}, \eqref{eq406},
\eqref{jiushi3}--\eqref{eq503} and
\eqref{eqbilinear1}, we obtain that
\begin{equation}\label{eq602}
\begin{array}{ll}
\ds \q\mE\big\langle P_T x_1^{\l_m}(T),
x_2^{\l_m}(T) \big\rangle_{H} - \mE
\int_{r_j}^T \big\langle
F(s) x_1^{\l_m}(s), x_2^{\l_m}(s) \big\rangle_{H}ds\\
\ns\ds =\mE\big\langle R^{(r_j,\l_m)}
\xi_1,\xi_2 \big\rangle_{H} + \mE \int_{r_j}^T
\big\langle P^{\l_m}(s)u_1(s),
x_2^{\l_m}(s)\big\rangle_{H}ds + \mE
\int_{r_j}^T \big\langle
P^{\l_m}(s)x_1^{\l_m}(s),
u_2(s)\big\rangle_{H}ds\\
\ns\ds \q   + \mE \int_{r_j}^T \big\langle
P^{\l_m }(s)K(s)x_1^{\l_m}(s),
v_2(s)\big\rangle_{H}ds + \mE \int_{r_j}^T
\big\langle P^{\l_m}(s)v_1(s),
K(s)x_2^{\l_m}(s)+v_2(s)
\big\rangle_{H}ds  \\
\ns\ds \q  + \mE \int_{r_j}^T \big\langle
v_1(s), \widehat Q^{\l_m,r_j}_{1}(\xi_2)(s) +
\widehat Q^{\l_m,r_j}_{2}(u_2)(s) +
\big(Q_3^{\l_m,r_j}\big)^*(v_2)(s)
\big\rangle_{H}ds \\
\ns\ds \q  + \mE \int_{r_j}^T \big\langle
Q^{\l_m,r_j}_{1}(\xi_1)(s) +
Q^{\l_m,r_j}_{2}(u_1)(s) +
Q_3^{\l_m,r_j}(v_1)(s), v_2(s)
\big\rangle_{H}ds,\\
\ns\ds\q
\forall\;(\xi_1,u_1,v_1),(\xi_2,u_2,v_2) \in
L^4_{\cF_{r_j}}(\O;H)\times
L^2_\dbF(r_j,T;L^4(\O;H)) \times
L^2_\dbF(r_j,T;L^4(\O;H)),\;j\in\dbN.
\end{array}
\end{equation}

{\bf Step 5}. In this step, we take
$n\to\infty$ in (\ref{eq2}) for all $t\in
[0,T]$.

Let $u_1=v_1=0$ in \eqref{op-fsystem2} and
$u_2=v_2=0$ in \eqref{op-fsystem3}. By
\eqref{eq602}, we obtain that
$$
 \ba{ll}\ds
 \mE\big\langle P_T U_{\l_m}(T,r_j)\xi_1, U_{\l_m}(T,{r_j})\xi_2 \big\rangle_{H} - \mE
\int_{r_j}^T \!\big\langle F(s)
U_{\l_m}(s,{r_j})\xi_1, U_{\l_m}(s,{r_j})\xi_2
\big\rangle_{H}ds \\\ns\ds
 =\mE\big\langle
R^{(r_j,\l_m)}\xi_1,\xi_2 \big\rangle_{H}.
 \ea
 $$
Hence, for any $\xi_1,\xi_2\in
L^4_{\cF_{r_j}}(\O;H)$, it holds that
$$
\mE\Big\langle U^*_{\l_m}(T,{r_j})P_T
U_{\l_m}(T,{r_j})\xi_1 - \int_{r_j}^T
U^*_{\l_m}(s,{r_j})F(s)
U_{\l_m}(s,{r_j})\xi_1ds, \xi_2 \Big\rangle_{H}
=\mE\big\langle R^{(r_j,\l_m)}\xi_1,\xi_2
\big\rangle_{H}.
$$
This leads to
\begin{equation}\label{10.29eq1}
\mE \(U^*_{\l_m}(T,{r_j})P_T
U_{\l_m}(T,{r_j})\xi_1 - \int_{r_j}^T
U^*_{\l_m}(s,t)F(s)
U_{\l_m}(s,{r_j})\xi_1ds\;\Big|\; \cF_{r_j}\) =
R^{(r_j,\l_m)}\xi_1.
\end{equation}

For any $t\in [0,T]$, $h\in [t,T]$ and $\xi\in
L^4_{\cF_t}(\O;H)$, let us define
 \bel{jiushi4}
 R^{(h,\l_m)}\xi\=\mE \(U^*_{\l_m}(T,h)P_T
U_{\l_m}(T,h)\xi - \int_h^T U^*_{\l_m}(s,h)F(s)
U_{\l_m}(s,h)\xi ds\;\Big|\; \cF_h\).
 \ee
For any $t \leq h_1 \leq h_2 \leq T$ and
$\xi\in L^4_{\cF_t}(\O;H)$, by (\ref{jiushi4}),
it follows that
\begin{equation}\label{eq6.1-c}
\begin{array}{ll}\ds
\q\mE\big| R^{(h_2,\l_m)}\xi - R^{(h_1,\l_m)}\xi\big|_H^{\frac{4}{3}} \\
\ns\ds \le C\Big[\mE\Big| \mE \(U^*_{\l_m}(T,h_2)P_T U_{\l_m}(T,h_2)\xi - \int_{h_2}^T  U^*_{\l_m}(s,h_2)F(s) U_{\l_m}(s,h_2)\xi ds\;\Big|\; \cF_{h_2}\)  \\
\ns\ds \q - \mE \(U^*_{\l_m}(T,h_1)P_T U_{\l_m}(T,h_1)\xi - \int_{h_1}^T  U^*_{\l_m}(s,h_1)F(s) U_{\l_m}(s,h_1)\xi ds\;\Big|\; \cF_{h_2}\) \Big|_H^{\frac{4}{3}}\\
\ns\ds \q + \mE\Big| \mE \(U^*_{\l_m}(T,h_1)P_T
U_{\l_m}(T,h_1)\xi - \int_{h_1}^T
U^*_{\l_m}(s,h_1)F(s) U_{\l_m}(s,h_1)\xi
ds\;\Big|\; \cF_{h_2}\)
\\
\ns\ds \q - \mE \(U^*_{\l_m}(T,h_1)P_T
U_{\l_m}(T,h_1)\xi - \int_{h_1}^T
U^*_{\l_m}(s,h_1)F(s) U_{\l_m}(s,h_1)\xi
ds\;\Big|\;
\cF_{h_1}\)\Big|_H^{\frac{4}{3}}\Big].
\end{array}
\end{equation}
By Lemma \ref{lemma8}, it is easy to show that
\begin{equation}\label{eq6.2-c}
\begin{array}{ll}\ds
\lim_{h_2\to h_1^+}\mE\Big|\mE
\(U^*_{\l_m}(T,h_1) P_T U_{\l_m}(T,h_1)\xi -
\int_{h_1}^T U^*_{\l_m}(s,h_1)F(s)
U_{\l_m}(s,h_1)\xi ds\;\Big|
\; \cF_{h_2}\)  \\
\ns\ds \qq\;\, - \mE \(U^*_{\l_m}(T,h_1)P_T
U_{\l_m}(T,h_1)\xi - \int_{h_1}^T
U^*_{\l_m}(s,h_1)F(s) U_{\l_m}(s,h_1)\xi
ds\;\Big|\;
\cF_{h_1}\)\Big|_H^{\frac{4}{3}}\\
\ns\ds=0.
\end{array}
\end{equation}
On the other hand,
\begin{equation}\label{eq6.3-c}
\begin{array}{ll}\ds
\q \mE\Big|   \mE \(U^*_{\l_m}(T,h_2)P_T
U_{\l_m}(T,h_2)\xi - \int_{h_2}^T
U^*_{\l_m}(s,h_2)F(s) U_{\l_m}(s,h_2)\xi
ds\;\Big|\; \cF_{h_2}\)
\\
\ns\ds \q - \mE \(U^*_{\l_m}(T,h_1)P_T
U_{\l_m}(T,h_1)\xi
- \int_{h_1}^T  U^*_{\l_m}(s,h_1)F(s) U_{\l_m}(s,h_1)\xi ds\;\Big|\; \cF_{h_2}\)\Big|_H^{\frac{4}{3}}\\
\leq C\mE\Big|U^*_{\l_m}(T,h_2)P_T U_{\l_m}(T,h_2)\xi - U^*_{\l_m}(T,h_1)P_T U_{\l_m}(T,h_1)\xi \Big|_H^{\frac{4}{3}} \\
\ns\ds \q + C\mE\Big|\int_{h_2}^T  \[U^*_{\l_m}(s,h_2)F(s) U_{\l_m}(s,h_2)\xi - U^*_{\l_m}(s,h_1)F(s) U_{\l_m}(s,h_1)\xi\]ds \Big|_H^{\frac{4}{3}}\\
\ns\ds \q  + C\mE\Big|\int_{h_1}^{h_2}
U^*_{\l_m}(s,h_1)F(s) U_{\l_m}(s,h_1)\xi
ds\Big|_H^{\frac{4}{3}}.
\end{array}
\end{equation}
Hence, noting that, for each $m\in\dbN$,  $A_{\l_m}$ is a bounded
linear operator on $H$, we obtain that
\begin{equation}\label{eq6.4-c}
\begin{array}{ll}\ds
\lim_{h_2\to h_1^+}\mE\Big|   \mE
\(U^*_{\l_m}(T,h_2)P_T U_{\l_m}(T,h_2)\xi -
\int_{h_2}^T U^*_{\l_m}(s,h_2)F(s)
U_{\l_m}(s,h_2)\xi ds\Big| \cF_{h_2}\)
\\
\ns\ds \qq\;\, - \mE \(U^*_{\l_m}(T,h_1)P_T
U_{\l_m}(T,h_1)\xi - \int_{h_1}^T
U^*_{\l_m}(s,h_1)F(s) U_{\l_m}(s,h_1)\xi
ds\;\Big|\;
\cF_{h_2}\)\Big|_H^{\frac{4}{3}}\\
\ns\ds=0.
\end{array}
\end{equation}
Put
 \bel{zq1}\widehat P^{\l_m}(\cd)\xi \= R^{(\cd,\l_m)}\xi.
 \ee
By (\ref{eq6.1-c})--(\ref{eq6.2-c}) and
(\ref{eq6.4-c})--(\ref{zq1}), it is easy to see
that $\widehat P^{\l_m}(\cd)\xi $ is right
continuous in $L^{\frac{4}{3}}_{\cF_T}(\O;H)$
on $[t,T]$.

For any $t\in [0,T)\setminus
\{r_j\}_{j=1}^\infty$, we can find a
subsequence $\{r_{j_k}\}_{k=1}^\infty\subset
\{r_j\}_{j=1}^\infty$ such that $r_{j_k}>t$ and
$\ds\lim_{k\to\infty}r_{j_k}=t$. Letting
$u_1=v_1=0$ in the equation
\eqref{op-fsystem2}, and letting $\xi_2=0$ and
$u_2=0$ in the equation \eqref{op-fsystem3}, by
\eqref{eq602}, we find that
\begin{equation}\label{6.6eq3}
\begin{array}{ll}
\ds \q\mE\big\langle P_T x_1^{\l_m}(T),
x_2^{\l_m}(T) \big\rangle_{H} - \mE
\int_{r_{j_k}}^T \big\langle
F(s) x_1^{\l_m}(s), x_2^{\l_m}(s) \big\rangle_{H}ds\\
\ns\ds = \mE \int_{r_{j_k}}^T \big\langle
P^{\l_m}(s) K(s)x_1^{\l_m}(s),
 v_2(s)\big\rangle_{H}ds + \mE
\int_{t}^T \big\langle
\chi_{[r_{j_k},T]}Q^{{\l_m},r_{j_k}}_{1}(\xi_1)
(s), v_2(s) \big\rangle_{H}ds.
\end{array}
\end{equation}
Let us choose $\xi_{1,j_k}\in
L^4_{\cF_{r_{j_k}}}(\O;H)$ such that
$|\xi_{1,j_k}|_{L^4_{\cF_{r_{j_k}}}(\O;H)}=1$
and
\begin{equation}\label{6.6eq4}
\begin{array}{ll}\ds
\big|\chi_{[r_{j_k},T]}Q^{{\l_m},r_{j_k}}_{1}(\xi_{1,j_k})
\big|_{L^{2}_\dbF(r_{j_k},T;L^{\frac{4}{3}}(\O;H)))}=\big|\chi_{[r_{j_k},T]}Q^{{\l_m},r_{j_k}}_{1}(\xi_{1,j_k})
\big|_{L^{2}_\dbF(t,T;L^{\frac{4}{3}}(\O;H)))}\\\ns\ds
 \geq
\frac{1}{2}\big|\!\big|\chi_{[r_{j_k},T]}
Q^{{\l_m},r_{j_k}}_{1}\big|\!\big|_{\cL(L^4_{\cF_{r_{j_k}}}\!\!(\O;H),\;
L^{2}_\dbF(t,T;L^{\frac{4}{3}}(\O;H)))}\geq
\frac{1}{2}\big|\!\big|\chi_{[r_{j_k},T]}
Q^{{\l_m},r_{j_k}}_{1}\big|\!\big|_{\cL(L^4_{\cF_{t}}(\O;H),\;
L^{2}_\dbF(t,T;L^{\frac{4}{3}}(\O;H)))}.
\end{array}
\end{equation}
Then, we choose $v_{2,j_k}\in
L^4_\dbF(r_{j_k},T;L^4(\O;H))$ with
$|v_{2,j_k}|_{
L^4_\dbF(r_{j_k},T;L^4(\O;H))}=1$ such that
\begin{equation}\label{6.6eq5}
\mE \int_{t}^T \big\langle
\chi_{[r_{j_k},T]}Q^{{\l_m},r_{j_k}}_{1}(\xi_{1,j_k})
(s), v_{2,j_k}(s) \big\rangle_{H}ds \geq
\frac{1}{2}\big|\chi_{[r_{j_k},T]}Q^{{\l_m},r_{j_k}}_{1}(\xi_{1,j_k})
\big|_{L^{2}_\dbF(r_{j_k},T;L^{\frac{4}{3}}(\O;H)))}.
\end{equation}
From \eqref{6.6eq3}--\eqref{6.6eq5}, we get
that
$$
\big|\!\big|\chi_{[r_{j_k},T]}
Q^{{\l_m},r_{j_k}}_{1}\big|\!\big|_{\cL(L^4_{\cF_{t}}(\O;H),\;
L^{2}_\dbF(t,T;L^{\frac{4}{3}}(\O;H)))}\leq
C(\l_m)\big( |P_T|_{L^2_{\cF_T}(\O;\cL(H))} +
|F|_{L^1_\dbF(0,T;\cL(H))} \big),
$$
where the constant $C(\l_m)$ is independent of
$r_{j_k}$. Similarly,
$$
\big|\!\big|\chi_{[r_{j_k},T]} \widehat
Q^{{\l_m},r_{j_k}}_{1}\big|\!\big|_{\cL(L^4_{\cF_{t}}(\O;H),\;
L^{2}_\dbF(t,T;L^{\frac{4}{3}}(\O;H)))}\leq
C(\l_m)\big( |P_T|_{L^2_{\cF_T}(\O;\cL(H))} +
|F|_{L^1_\dbF(0,T;\cL(H))} \big).
$$
From Lemma \ref{lemma6}, we conclude that
there exist two bounded linear operators $
Q^{{\l_m},t}_{1}$ and $\widehat
Q^{{\l_m},t}_{1}$, which are from
$L^4_{\cF_t}(\O;H)$ to
$L^{2}_\dbF(t,T;L^{\frac{4}{3}}(\O;H))$, and a
subsequence $\{j^{(1)}_k\}_{k=1}^\infty\subset
\{j_k\}_{k=1}^\infty$  such that
\begin{equation}\label{6.6eq6}
\left\{
\begin{array}{ll}\ds
\mbox{(w)}\mbox{-}
\lim_{k\to\infty}\chi_{[t_{j_k^{(1)}},T]}Q^{\l_m,r_{j_k}}_{1}\xi
= Q^{{\l_m},t}_{1} \xi \mbox{ in }
L^{2}_\dbF(t,T;L^{\frac{4}{3}}(\O;H)),\qq
\forall\; \xi\in L^4_{\cF_t}(\O;H),\\
\ns\ds \mbox{(w)}\mbox{-}
\lim_{k\to\infty}\chi_{[t_{j_k^{(1)}},T]}\widehat
Q^{\l_m,r_{j_k}}_{1}\xi =\widehat
Q^{{\l_m},t}_{1} \xi \mbox{ in }
L^{2}_\dbF(t,T;L^{\frac{4}{3}}(\O;H)),\qq
\forall\; \xi\in L^4_{\cF_t}(\O;H).
\end{array}
\right.
\end{equation}

Letting $\xi_1 =0$, $v_1 =0$ in
\eqref{op-fsystem2} and $\xi_2 =0$, $u_2 =0$ in
\eqref{op-fsystem3}, by \eqref{eq602}, we
obtain that
\begin{equation}\label{6.6eq7}
\begin{array}{ll}
\ds \q\mE\big\langle P_T x_1^{\l_m}(T),
x_2^{\l_m}(T) \big\rangle_{H} - \mE
\int_{r_{j_k}}^T \big\langle
F(s) x_1^{\l_m}(s), x_2^{\l_m}(s) \big\rangle_{H}ds\\
\ns\ds = \mE \int_{r_{j_k}}^T \big\langle
P^{\l_m}(s)u_1(s),
x_2^{\l_m}(s)\big\rangle_{H}ds   + \mE
\int_{r_{j_k}}^T \big\langle
P^{\l_m}(s)K(s)x_1^{\l_m}(s),v_2(s)\big\rangle_{H}ds    \\
\ns\ds \q   + \mE \int_t^T
\big\langle\chi_{[r_{j_k},T]}
Q^{\l_m,r_{j_k}}_{2}(u_1)(s), v_2(s)
\big\rangle_{H}ds,\qq \forall\;u_1,v_2 \in
L^2_\dbF(r_{j_k},T;L^4(\O;H)),\; k\in\dbN.
\end{array}
\end{equation}
We choose $u_1^{(r_{j_k})}\in
L^2_\dbF(r_{j_k},T;L^4(\O;H))$ satisfying
$|u_1^{(r_{j_k})}|_{L^2_\dbF(r_{j_k},T;L^4(\O;H))}=1$,
and
\begin{equation}\label{6.6eq8}
\begin{array}{ll}\ds
\big|\chi_{[r_{j_k},T]}Q^{\l_m,r_{j_k}}_{2}\big(u_1^{(r_{j_k})}\big)
\big|_{L^{2}_{\dbF}(r_{j_k},T;L^{\frac{4}{3}}(\O;H))}\\\ns\ds
\geq
\frac{1}{2}\big|\!\big|\chi_{[r_{j_k},T]}Q^{\l_m,r_{j_k}}_{2}
\big|\!\big|_{\cL(L^2_\dbF(r_{j_k},T;L^4(\O;H)),\; L^{2}_{\dbF}(r_{j_k},T;L^{\frac{4}{3}}(\O;H)))}\\
\ns    \geq \ds\frac{1}{2}\big|\!\big|
\chi_{[r_{j_k},T]}Q^{\l_m,r_{j_k}}_{2}
\circ\chi_{[r_{j_k},T]}
\big|\!\big|_{\cL(L^2_\dbF(t,T;L^4(\O;H)),\;
L^{2}_{\dbF}(t,T;L^{\frac{4}{3}}(\O;H)))}.
\end{array}
\end{equation}
Then we choose $v_2^{(r_{j_k})}\in
L^2_\dbF(r_{j_k},T;L^4(\O;H))$ satisfying
$\big|v_2^{(r_{j_k})}\big|_{L^2_\dbF(r_{j_k},T;L^4(\O;H))}=1$,
and
\begin{equation}\label{6.6eq9}
\mE \int_t^T \Big\langle
\chi_{[r_{j_k},T]}Q^{\l_m,r_{j_k}}_{2}
\big(u_1^{(r_{j_k})}\big)(s),v_2^{(r_{j_k})}(s)
\Big\rangle_{H}ds \geq
\frac{1}{2}\big|\chi_{[r_{j_k},T]}
Q^{\l_m,r_{j_k}}_{2}
\big(u_1^{(r_{j_k})}\big)\big|_{L^{2}_{\dbF}(r_{j_k},T;L^{\frac{4}{3}}(\O;H))}.
\end{equation}
From \eqref{6.6eq7}--\eqref{6.6eq9}, we obtain
that
$$
 \ba{ll}\ds
\big|\!\big|
\chi_{[r_{j_k},T]}Q^{\l_m,r_{j_k}}_{2}\circ\chi_{[r_{j_k},T]}\big|\!\big|_{\cL(L^2_\dbF(t,T;L^4(\O;H)),\;
L^{2}_{\dbF}(t,T;L^{\frac{4}{3}}(\O;H)))}\\\ns\ds\leq
C(\l_m)\big( |P_T|_{L^2_{\cF_T}(\O;\cL(H))} +
|F|_{L^1_\dbF(t,T;L^2(\O;\cL(H)))}\big),
 \ea
$$
where the constant $C(\l_m)$ is independent of
$r_{j_k}$. By a similar argument, we obtain
that
$$
 \ba{ll}\ds
\big|\!\big|\chi_{[r_{j_k},T]}\widehat
Q^{\l_m,r_{j_k}}_{2}\circ\chi_{[r_{j_k},T]}\big|\!\big|_{\cL(L^2_\dbF(t,T;L^4(\O;H)),\;
L^{2}_{\dbF}(t,T;L^{\frac{4}{3}}(\O;H)))}\\\ns\ds\leq
C(\l_m)\big( |P_T|_{L^2_{\cF_T}(\O;\cL(H))} +
|F|_{L^1_\dbF(t,T;L^2(\O;\cL(H)))}\big).
 \ea
$$
By means of Lemma \ref{lemma6}, we see that
there exist two bounded linear operators
$Q^{\l_m,t}_{2}$ and $\widehat
Q^{\l_m,t}_{2}$, from $L^2_\dbF(t,T;L^4(\O;H))$
to $L^{2}_\dbF(t,T;L^{\frac{4}{3}}(\O;H))$, and
a subsequence
$\{r_{j_k}^{(2)}\}_{k=1}^\infty\subset
\{r_{j_k}^{(1)}\}_{k=1}^\infty$ such that, for
any $ u\in L^2_\dbF(t,T;L^4(\O;H))$,
\begin{equation}\label{6.6eq11}
\left\{
\begin{array}{ll}\ds
\mbox{(w)}\mbox{-}\lim_{k\to\infty}\(\chi_{[r_{j_k},T]}Q^{\l_m,r_{j_k}}_{2}\circ\chi_{[r_{j_k},T]}\)u
= Q^{\l_m,t}_{2} u \mbox{ in }
L^{2}_\dbF(t,T;L^{\frac{4}{3}}(\O;H)),\\
\ns\ds
\mbox{(w)}\mbox{-}\lim_{k\to\infty}\(\chi_{[r_{j_k},T]}\widehat
Q^{\l_m,r_{j_k}}_{2}\circ\chi_{[r_{j_k},T]}\)u
= \widehat Q^{\l_m,t}_{2} u \mbox{ in }
L^{2}_\dbF(t,T;L^{\frac{4}{3}}(\O;H)).
\end{array}
\right.
\end{equation}
For any $t\in [0,T]$, we define two operators
$Q^{(\l_m,t)}$ and $\widehat Q^{(\l_m,t)}$ on
$L^4_{\cF_{t}}(\O;H)\times
L^2_\dbF({t},T;L^4(\O;H)) \times
L^2_\dbF({t},T;L^4(\O;H))$ as follows:
\begin{equation}\label{3Q}
\left\{
\begin{array}{ll}\ds
Q^{(\l_m,t)}(\xi,u,v) = Q^{\l_m,t}_{1}\xi +
Q^{\l_m,t}_{2} u + Q_3^{\l_m,t} v,\\
\ns\ds \widehat Q^{(\l_m,t)}(\xi,u,v) =
\widehat Q^{\l_m,t}_{1}\xi + \widehat
Q^{\l_m,t}_{2} u +\big(Q_3^{\l_m,t}\big)^*
v,\\
\ns\ds \qq\qq\qq\q\;\;\forall\; (\xi,u,v)\in
L^4_{\cF_{t}}(\O;H)\times
L^2_\dbF(t,T;L^4(\O;H)) \times
L^2_\dbF(t,T;L^4(\O;H)).
\end{array}
\right.
\end{equation}
From the definition of $Q^{\l_m,t}_{1}$,
$Q^{\l_m,t}_{2}$ and $Q_3^{\l_m,t}$(\resp
$\widehat Q^{\l_m,t}_{1}$, $\widehat
Q^{\l_m,t}_{2}$ and
$\big(Q_3^{\l_m,t}\big)^*$), we find that
$Q^{(\l_m,t)}(\cd,\cd,\cd)$ (\resp $\widehat
Q^{(\l_m,t)}(\cd,\cd,\cd)$) is a bounded linear
operator from $L^4_{\cF_{t}}(\O;H)\times
L^2_\dbF(t,T;L^4(\O;H)) $ $\times
L^2_\dbF(t,T;L^4(\O;H))$ to
$L^{2}_\dbF(t,T;L^{\frac{4}{3}}(\O;H))$ and
$Q^{(\l_m,t)}(0,0,\cd)^*=\widehat
Q^{(\l_m,t)}(0,0,\cd)$.

For any $t\in[0,T]$, from \eqref{eq602},
(\ref{zq1}), \eqref{6.6eq6}, \eqref{6.6eq11}
and \eqref{3Q}, we obtain that
\begin{equation}\label{6.6eq15}
\begin{array}{ll}
\ds \q\mE\big\langle P_T x_1^{\l_m}(T),
x_2^{\l_m}(T) \big\rangle_{H} - \mE \int_t^T
\big\langle
F(s) x_1^{\l_m}(s), x_2^{\l_m}(s) \big\rangle_{H}ds\\
\ns\ds =\mE\big\langle \widehat P^{\l_m}(t)
\xi_1,\xi_2 \big\rangle_{H} + \mE \int_t^T
\big\langle  P^{\l_m}(s)u_1(s),
x_2^{\l_m}(s)\big\rangle_{H}ds  + \mE \int_t^T
\big\langle P^{\l_m}(s)x_1^{\l_m}(s),
u_2(s)\big\rangle_{H}ds\\
\ns\ds \q   + \mE \int_t^T \big\langle
P^{\l_m}(s)K(s)x_1^{\l_m}(s),
v_2(s)\big\rangle_{H}ds + \mE \int_t^T
\big\langle P^{\l_m}(s)v_1(s),
K(s)x_2^{\l_m}(s)+v_2(s)
\big\rangle_{H}ds  \\
\ns\ds \q  + \mE \int_t^T \big\langle v_1(s),
\widehat Q^{(\l_m,t)}(\xi_2,u_2,v_2)(s)
\big\rangle_{H}ds   + \mE \int_t^T \big\langle
 Q^{(\l_m,t)}(\xi_1,u_1,v_1)(s),
v_2(s) \big\rangle_{H}ds,\\
\ns\ds\qq\forall\;(\xi_1,u_1,v_1),(\xi_2,u_2,v_2)
\in L^4_{\cF_{t}}(\O;H)\times
L^2_\dbF(t,T;L^4(\O;H)) \times
L^2_\dbF(t,T;L^4(\O;H)).
\end{array}
\end{equation}

We claim that
\begin{equation}\label{P1}
P^{\l_m}(t)\= \widehat P^{\l_m}(t), \ \ \ae
t\in [0,T].
\end{equation}
To show this, for any $0\leq t_1<t_2<T$, we
choose $x_1(t_1)=\eta_1\in
L^4_{\cF_{t_1}}(\O;H)$ and $u_1=v_1=0$ in the
equation \eqref{op-fsystem2.1}, and
$x_2(t_1)=0$,
$u_2(\cd)=\frac{\chi_{[t_1,t_2]}}{t_2-t_1}\eta_2$
with $\eta_2\in L^4_{\cF_{t_1}}(\O;H)$ and
$v_2=0$ in the equation \eqref{op-fsystem3.1},
by \eqref{6.6eq15} and recalling the definition
of the evolution operator $U_{\l_m}(\cd,\cd)$,
we see that
\begin{equation}\label{s5eq11}
\begin{array}{ll}\ds
\q\frac{1}{t_2-t_1}\mE\int_{t_1}^{t_2}\Big\langle
P^{\l_m}(s)
U_{\l_m}(s,t_1)\eta_1, \eta_2 \Big\rangle_H ds\\
\ns\ds =\mE\Big\langle P_T
U_{\l_m}(T,t_1)\eta_1, x_{2,t_2}^{\l_m}(T)
\Big\rangle_H   - \mE\int_{t_1}^T \Big\langle
F(s)U_{\l_m}(s,t_1)\eta_1, x_{2,t_2}^{\l_m}(s)
\Big\rangle_H ds,
\end{array}
\end{equation}
where $x_{2,t_2}^{\l_m}(\cd)$ stands for the
solution to the equation \eqref{op-fsystem3.1}
with the above choice of $\xi_2$, $u_2$ and
$v_2$. It is clear that
 \bel{vz1}
 \ds
x_{2,t_2}^{\l_m}(s)=\left\{
\ba{ll}\ds\int_{t_1}^{s}S_{\l_m}(s-\tau)J(\tau)x_{2,t_2}^{\l_m}(\tau)d\tau
+\int_{t_1}^{s}S_{\l_m}(s-\tau)K(\tau)x_{2,t_2}^{\l_m}(\tau)dw(\tau)\\\ns\ds\q
+\frac{1}{t_2-t_1}\int_{t_1}^{s}S_{\l_m}(s-\tau)\eta_2
d\tau,\qq s\in [t_1,t_2],\\\ns\ds
U_{\l_m}(s,t_2)x_{2,t_2}^{\l_m}(t_2),
\qq\qq\qq\q\;\; s\in [t_2,T].
 \ea\right.
 \ee
Then, by Lemma \ref{BDG}, we see that
$$
\begin{array}{ll}\ds
\q\mE\big|x_{2,t_2}^{\l_m}(s)\big|_H^4
\\\ns\ds\leq C(\l_m)\Big\{\int_{t_1}^s\[|J(\tau)|^4_{ L^\infty(\O;
\cL(H))}+|K(\tau)|^4_{ L^\infty(\O; \cL(H))}\]
\mE\big|x_{2,t_2}^{\l_m}(\tau)\big|_H^4d\tau
+\mE|\eta_2|_H^4\Big\},\q \forall\;s\in
[t_1,t_2].
\end{array}
$$
By Gronwall's inequality, it follows that
\begin{equation}\label{10.20eq1}
\big|x_{2,t_2}^{\l_m}\big|_{L^\infty_\dbF(t_1,t_2;L^4(\O;H))}\leq
C(\l_m)|\eta_2|_{L^4_{\cF_{t_1}}(\O;H)},
\end{equation}
where the constant $C(\l_m)$ is independent of
$t_2$. On the other hand, by \eqref{vz1}, we
have
$$
\begin{array}{ll}\ds
\q\mE\big|x_{2,t_2}^{\l_m}(t_2)-\eta_2\big|_H^4\\
\ns\ds \leq
C(\l_m)\Big\{\int_{t_1}^s\[|J(\tau)|^4_{
L^\infty(\O; \cL(H))}+|K(\tau)|^4_{
L^\infty(\O; \cL(H))}\]
\mE\big|x_{2,t_2}^{\l_m}(\tau)\big|_H^4d\tau\\\ns\ds\qq\qq\q
+
\mE\Big|\frac{1}{t_2-t_1}\int_{t_1}^{t_2}S(t_2-\tau)\eta_2
d\tau -\eta_2\Big|_H^4\Big\}.
\end{array}
$$
This, together with \eqref{10.20eq1}, implies
that
$$
\begin{array}{ll}\ds
\lim_{t_2\to
t_1+0}\mE\big|x_{2,t_2}^{\l_m}(t_2)-\eta_2\big|_H^4\3n&\ds\leq
C(\l_m)\lim_{t_2\to
t_1+0}\mE\Big|\frac{1}{t_2-t_1}\int_{t_1}^{t_2}S(t_2-\tau)\eta_2
d\tau -\eta_2\Big|_H^4=0.
\end{array}
$$
Therefore, for any $s\in [t_2,T]$,
$$
\begin{array}{ll}\ds
\q\lim_{t_2\to
t_1+0}\mE\big|U_{\l_m}(s,t_2)x_{2,t_2}^{\l_m}(t_2)-U_{\l_m}(s,t_1)\eta_2\big|_H^4\\
\ns\ds \leq 8\lim_{t_2\to
t_1+0}\[\mE\big|U_{\l_m}(s,t_2)x_{2,t_2}^{\l_m}(t_2)-U_{\l_m}(s,t_2)\eta_2\big|_H^4
+\mE\big|U_{\l_m}(s,t_2)\eta_2-U_{\l_m}(s,t_1)\eta_2\big|_H^4\]\\
\ns\ds \leq C(\l_m)\lim_{t_2\to
t_1+0}\[\mE\big| x_{2,t_2}^{\l_m}(t_2)-
\eta_2\big|_H^4
+\mE\big|U_{\l_m}(s,t_2)\eta_2-U_{\l_m}(s,t_1)\eta_2\big|_H^4\]=0.
\end{array}
$$
Hence, we obtain that
 \bel{vz1q}
\lim_{t_2\to t_1+0} x_{2,t_2}^{\l_m}(s) =
U_{\l_m}(s,t_1)\eta_2 \ \mbox{ in }
L^4_{\cF_s}(\O;H),\qq\forall\;s\in [t_2,T].
 \ee
By \eqref{10.20eq1} and \eqref{vz1q}, we
conclude that
\begin{equation}\label{10.20eq2}
\begin{array}{ll}\ds
\q\lim_{t_2\to t_1+0}\[\mE\Big\langle P_T
U_{\l_m}(T,t_1)\eta_1, x_{2,t_2}^{\l_m}(T)
\Big\rangle_H   - \mE\int_{t_1}^T \Big\langle
F(s)U_{\l_m}(s,t_1)\eta_1, x_{2,t_2}^{\l_m}(s)
\Big\rangle_H ds\] \\
\ns\ds = \mE\Big\langle P_T
U_{\l_m}(T,t_1)\eta_1, U_{\l_m}(T,t_1)\eta_2
\Big\rangle_H   - \mE\int_{t_1}^T \Big\langle
F(s)U_{\l_m}(s,t_1)\eta_1,
U_{\l_m}(s,t_1)\eta_2 \Big\rangle_H ds.
\end{array}
\end{equation}

 By choosing $x_1(t_1)=\eta_1$ and
$u_1=v_1=0$ in \eqref{op-fsystem2}, and
$x_2(t_1)=\eta_2$ and $u_2=v_2=0$ in
\eqref{op-fsystem3}, by \eqref{6.6eq15}, we
find that
\begin{equation}\label{s5eq12}
\begin{array}{ll}\ds
\q\mE\Big\langle \widehat P^{\l_m}(t_1) \eta_1,
 \eta_2
\Big\rangle_H\\
\ns\ds =\mE\Big\langle P_T
U_{\l_m}(T,t_1)\eta_1, U_{\l_m}(T,t_1) \eta_2
\Big\rangle_H  -  \mE \int_{t_1}^T  \Big\langle
F(s)U_{\l_m}(s,t_1)\eta_1, U_{\l_m}(s,t_1)
\eta_2 \Big\rangle_H ds.
\end{array}
\end{equation}

Combining \eqref{s5eq11}, \eqref{10.20eq2} and
\eqref{s5eq12}, we obtain that
\begin{equation}\label{s5eq15}
\lim_{t_2\to t_1+0} \frac{1}{t_2 - t_1}\mE
\int_{t_1}^{t_2} \Big\langle P^{\l_m}(s)
U_{\l_m}(s,t_1)\eta_1,  \eta_2 \Big\rangle_H ds
= \mE\Big\langle \widehat P^{\l_m}(t_1) \eta_1,
\eta_2 \Big\rangle_H.
\end{equation}
By Lemma \ref{lemma2.1}, we see that there is a
monotonically decreasing sequence
$\{t_2^{(n)}\}_{n=1}^\infty$ with
$t_2^{(n)}>t_1$ for every $n$, such that
$$
\lim_{t_2^{(n)}\to
t_1+0}\frac{1}{t_2^{(n)}-t_1}\mE\int_{t_1}^{t_2^{(n)}}\Big\langle
P^{\l_m}(s) U_{\l_m}(s,t_1)\eta_1, \eta_2
\Big\rangle_H ds = \mE\Big\langle
P^{\l_m}(t_1)\eta_1, \eta_2 \Big\rangle_H,\q\ae
t_1\in [0,T).
$$
This, together with \eqref{s5eq15}, implies
that
$$
\mE\Big\langle \widehat P^{\l_m}(t_1) \eta_1,
\eta_2 \Big\rangle_H = \mE\Big\langle
P^{\l_m}(t_1)\eta_1, \eta_2 \Big\rangle_H,
\;\mbox{ for }\ae t_1\in [0,T).
$$
Since $\eta_1$ and $\eta_2$ are arbitrary
elements in $L^4_{\cF_{t_1}}(\O;H)$, we
conclude \eqref{P1}.

By \eqref{6.6eq15} and \eqref{P1}, we end up
with
\begin{equation}\label{eq5}
\begin{array}{ll}
\ds \q\mE\big\langle P_T x_1^{\l_m}(T),
x_2^{\l_m}(T) \big\rangle_{H} - \mE \int_t^T
\big\langle
F(s) x_1^{\l_m}(s), x_2^{\l_m}(s) \big\rangle_{H}ds\\
\ns\ds =\mE\big\langle P^{\l_m}(t) \xi_1,\xi_2
\big\rangle_{H} + \mE \int_t^T \big\langle
P^{\l_m}(s)u_1(s),
x_2^{\l_m}(s)\big\rangle_{H}ds + \mE \int_t^T
\big\langle P^{\l_m}(s)x_1^{\l_m}(s),
u_2(s)\big\rangle_{H}ds \\
\ns\ds \q  + \mE \int_t^T \big\langle
P^{\l_m}(s)K (s)x_1^{\l_m} (s), v_2
(s)\big\rangle_{H}ds +  \mE \int_t^T
\big\langle
P^{\l_m}(s)v_1 (s), K (s)x_2^{\l_m} (s)+v_2(s)\big\rangle_{H}ds\\
\ns\ds \q  + \mE \int_t^T \big\langle v_1(s),
\widehat
Q^{(\l_m,t)}(\xi_2,u_2,v_2)(s)\big\rangle_{H}ds
+ \mE \int_t^T \big\langle
Q^{(\l_m,t)}(\xi_1,u_1,v_1)(s),
v_2(s) \big\rangle_{H}ds,\\
\ns\ds\qq\forall\;(\xi_1,u_1,v_1),(\xi_2,u_2,v_2)
\in L^4_{\cF_{t}}(\O;H)\times
L^2_\dbF(t,T;L^4(\O;H)) \times
L^2_\dbF(t,T;L^4(\O;H)).
\end{array}
\end{equation}

{\bf Step 6}. In this step, we show the
well-posedness of the relaxed transposition
solution to \eqref{op-bsystem3}.

Similar to the argument in Steps 4--5, thanks
to the uniform estimate \eqref{op-fsystem5
enest x} (with respect to $\l_m$), we conclude
that there exist a subsequence
$\{\l^{(1)}_{m_j}\}_{j=1}^\infty\subset
\{\l_m\}_{m=1}^\infty$, a
$P(\cd)\in\cL_{pd}\big(L^2_\dbF(0,T;L^4(\O;H)),\;$
$ L^2_\dbF(0,T;L^{\frac{4}{3}}(\O;H))\big)$, an
$
R^{(t)}\in\cL_{pd}\big(L^4_{\cF_{t}}(\O;H),\,L^{\frac{4}{3}}_{\cF_{t}}(\O;H)\big)$,
and two bounded linear operators $Q^{(t)}$ and
$\widehat Q^{(t)}$ from
$L^4_{\cF_t}(\O;H)\times
L^2_\dbF(t,T;L^4(\O;H))\times
L^2_\dbF(t,T;L^4(\O;H))$ to
$L^2_\dbF(t,T;L^{\frac{4}{3}}(\O;H))$
satisfying $Q^{(t)}(0,0,\cd)^*=\widehat
Q^{(t)}(0,0,\cd)$, such that
\begin{equation}\label{6.18eq2}
\begin{array}{ll}\ds
\mbox{(w)-}\lim_{j\to\infty}P^{\l^{(1)}_{m_j}}(\cd)u(\cd)
= P(\cd)u(\cd) \ \ \mbox{ in }
L^2_\dbF(0,T;L^{\frac{4}{3}}(\O;H)),\q
\forall\; u(\cd)\in L^2_\dbF(0,T;L^{4}(\O;H)),
\end{array}
\end{equation}
\begin{equation}\label{s5oeq3w1}
\mbox{(w)}\mbox{-}\lim_{j\to\infty}P^{\l^{(1)}_{m_j}}(t)
\xi = R^{(t)}\xi \ \ \ \mbox{ in }
 L^{\frac{4}{3}}_{\cF_{t}}(\O;H),\qq \forall\;
\xi\in L^4_{\cF_{t}}(\O;H),
\end{equation}
and
\begin{equation}\label{6.18eq3}
\left\{
\begin{array}{ll}\ds
\mbox{(w)-}\lim_{j\to\infty}Q^{(\l^{(1)}_{m_j},t)}
(\xi,u(\cd),v(\cd)) =
Q^{(t)}(\xi,u(\cd),v(\cd)) \ \ \mbox{ in }
L^2_\dbF(t,T;L^{\frac{4}{3}}(\O;H)),
\\
\ns\ds\mbox{(w)-}\lim_{j\to\infty}\widehat
Q^{(\l^{(1)}_{m_j},t)} (\xi,u(\cd),v(\cd)) =
\widehat Q^{(t)}(\xi,u(\cd),v(\cd)) \ \ \mbox{
in }
L^2_\dbF(t,T;L^{\frac{4}{3}}(\O;H)),\\
\ns\ds \qq \forall\; (\xi,u(\cd),v(\cd))\in
L^4_{\cF_t}(\O;H)\times
L^2_\dbF(t,T;L^{4}(\O;H))\times
L^2_\dbF(t,T;L^{4}(\O;H)).
\end{array}
\right.
\end{equation}

By (\ref{eq5}), and noting
(\ref{6.18eq2})--(\ref{6.18eq3}), we find that
$\big(P(\cd), R^{(\cd)},Q^{(\cd)},\widehat
Q^{(\cd)}\big)$ satisfies the following
variational equality:
 \begin{equation}\label{6.1-8eq1}
\begin{array}{ll}
\ds \q\mE\big\langle P_T x_1(T), x_2(T)
\big\rangle_{H} - \mE \int_t^T \big\langle
F(s) x_1(s), x_2(s) \big\rangle_{H}ds\\
\ns\ds =\mE\big\langle R^{(t)} \xi_1,\xi_2
\big\rangle_{H} + \mE \int_t^T \big\langle
P(s)u_1(s), x_2(s)\big\rangle_{H}ds + \mE
\int_t^T \big\langle P(s)x_1(s),
u_2(s)\big\rangle_{H}ds \\
\ns\ds \q  + \mE \int_t^T \big\langle P(s)K
(s)x_1 (s), v_2 (s)\big\rangle_{H}ds +  \mE
\int_t^T \big\langle  P(s)v_1 (s), K (s)x_2 (s)+ v_2(s)\big\rangle_{H}ds\\
\ns\ds \q + \mE \int_t^T \big\langle v_1(s),
\widehat
Q^{(t)}(\xi_2,u_2,v_2)(s)\big\rangle_{H}ds+ \mE
\int_t^T \big\langle
Q^{(t)}(\xi_1,u_1,v_1)(s), v_2(s) \big\rangle_{H}ds,\\
\ns\ds\qq\forall\;(\xi_1,u_1,v_1),(\xi_2,u_2,v_2)
\in L^4_{\cF_{t}}(\O;H)\times
L^2_\dbF(t,T;L^4(\O;H)) \times
L^2_\dbF(t,T;L^4(\O;H)).
\end{array}
\end{equation}

Now we show that $R^{(\cd)}\xi$ is right continuous in
$L^{\frac{4}{3}}_{\cF_T}(\O;H)$ on $[t,T]$ for any $\xi\in
L^4_{\cF_t}(\O;H)$. Since $A$ is usually an unbounded linear
operator on $H$, here we cannot employ the same method for treating
$R^{(\cd,\l_m)}$. Let $(P^n(\cd),Q^n(\cd))$ be the transposition
solution to \eqref{op-bsystem3} with the final datum
$P_T^n(=\G_nP_T\G_n)$ and the nonhomogeneous term $F_n(=\G_nF\G_n)$.
By Theorem \ref{op well th0}, it follows that
 \bel{jkj11}
P^n(\cd)\in
D_\dbF(0,T;L^2(\O;\cL_2(H)))\subset
D_\dbF(0,T;L^2(\O;\cL(H))).
 \ee
Further, similar to the derivation of the equality \eqref{10.29eq1},
we conclude that for any $t\in [0,T]$, $\tau\in [t,T]$  and $\xi\in
L^4_{\cF_t}(\O;H)$, it holds that
$$
R^{(\tau)}\xi =\mE \(U^*(T,\tau)P_T
U(T,\tau)\xi - \int_\tau^T U^*(s,\tau)F(s)
U(s,\tau)\xi ds\;\Big|\; \cF_\tau\),
$$
and
$$
P^n(\tau)\xi =\mE \(U^*(T,\tau)P_T^n
U(T,\tau)\xi - \int_\tau^T U^*(s,\tau)F_n(s)
U(s,\tau)\xi ds\;\Big|\; \cF_\tau\).
$$
By \eqref{jkj11}, in order to prove the right continuity of
$R^{(\cd)}\xi$ in $L^{\frac{4}{3}}_{\cF_T}(\O;H)$, it remains to
show that
\begin{equation}\label{10.30eq1}
\lim_{n\to\infty}\big|R^{(\cd)}\xi-P^n(\cd)\xi\big|_{L^\infty_\dbF(t,T;L^{\frac{4}{3}}(\O;H))}=0.
\end{equation}
For this purpose, for any $\tau\in [t,T]$, we see that
$$
\begin{array}{ll}\ds
\mE\big|R^{(\tau)}\xi -
P^n(\tau)\xi\big|^{\frac{4}{3}}_H\3n&\ds \leq
C\mE\Big| \int_\tau^T
U^*(s,\tau)\big[F(s)-F_n(s)\big] U(s,\tau)\xi
ds \Big|_H^{\frac{4}{3}}\\
\ns&\ds \q + C\mE\Big| U^*(T,\tau)\big(P_T-P_T^n\big) U(T,\tau)\xi
\Big|_H^{\frac{4}{3}}.
\end{array}
$$
By the first conclusion in Lemma \ref{lemma5}, we deduce that for
any $\e_1>0$, there is a $\d_1>0$ so that for all $\tau\in [t,T]$
and $\tau\leq\si\leq\tau+\d_1$,
\begin{equation}\label{10.31eq2}
\mE\big|U(r,\tau)\xi -
U(r,\si)\xi\big|_H^{\frac{4}{3}}<\e_1,\qq\forall\; r\in [\si,T].
\end{equation}
Now, we choose a monotonicity increasing sequence
$\{\tau_i\}_{i=1}^{N_1}\subset [0,T]$ for $N_1$ sufficiently large
such that $\tau_{i+1}-\tau_i \leq \d_1$ with $\tau_1=t$ and
$\tau_{N_1}=T$, and that
\begin{equation}\label{10.31eq2ggg}
\(\int_{\tau_i}^{\tau_{i+1}}\mE|F(s)|^2_{\cL(H)}ds\)^{\frac{2}{3}}
<\e_1,\;\mbox{ for all } i=1,\cds,N_1-1.
\end{equation}
For any $\tau_i<\tau\leq\tau_{i+1}$, recalling $F_n=\G_nF\G_n$, we
conclude that
\begin{equation}\label{10.31eq3}
\begin{array}{ll}\ds
\q\mE\Big| \int_\tau^T
U^*(s,\tau)\big[F(s)-F_n(s)\big] U(s,\tau)\xi
ds \Big|_H^{\frac{4}{3}}\\
\ns\ds \leq C\mE\Big| \int_{\tau_i}^T
U^*(s,\tau)\big[F(s)-F_n(s)\big] U(s,\tau_i)\xi ds
\Big|_H^{\frac{4}{3}} + C\mE\Big| \int_\tau^{\tau_i}
U^*(s,\tau)\big[F(s)-F_n(s)\big] U(s,\tau_i)\xi
ds \Big|_H^{\frac{4}{3}}\\
\ns\ds \q + C\mE\Big| \int_\tau^T
U^*(s,\tau)\big[F(s)-F_n(s)\big]
\big[U(s,\tau_i)-U(s,\tau)\big]\xi ds
\Big|_H^{\frac{4}{3}}\\
\ns\ds \leq C \int_\tau^T
\mE\Big|\big[F(s)-F_n(s)\big] U(s,\tau_i)\xi
\Big|_H^{\frac{4}{3}}ds +
C\(\int_{\tau_i}^{\tau_{i+1}}\mE|F|^2_{\cL(H)}ds\)^{\frac{2}{3}}
\\
\ns\ds \q + C\max_{s\in [\tau,T]}
\mE\Big|\big[U(s,\tau_i)-U(s,\tau)\big]\xi\Big|_H^{\frac{4}{3}}
ds.
\end{array}
\end{equation}
By the choice of $F_n$, there is an integer $N_2(\e_1)>0$ so that
for all $n>N_2$ and $i=1,\cds,N_1-1$,
\begin{equation}\label{10.31eq4}
\int_{\tau_i}^T \mE\Big|\big[F(s)-F_n(s)\big]
U(s,\tau_i)\xi \Big|_H^{\frac{4}{3}}ds \leq
\e_1.
\end{equation}
Combing \eqref{10.31eq2}--\eqref{10.31eq4}, we conclude that for all
$n>N_2$ and $\tau \in [t,T]$,
\begin{equation}\label{10.31eq3ggg}
\mE\Big| \int_\tau^T
U^*(s,\tau)\big[F(s)-F_n(s)\big] U(s,\tau)\xi
ds \Big|_H^{\frac{4}{3}} \leq C_1\e_1.
\end{equation}
Here the constant $C_1$ is independent of $\e_1$, $n$ and $\tau$.
Similarly, there is an integer $N_3(\e_1)>0$ such that for every
$n>N_3$,
\begin{equation}\label{10.31eq5}
\mE\Big| U^*(T,\tau)\big[P_T-P_T^n\big]
U(T,\tau)\xi \Big|_H^{\frac{4}{3}}\leq C_2\e_1,
\end{equation}
for the constant $C_2$ which is independent of $\e_1$, $n$ and
$\tau$. Now for any $\e>0$, let us choose $\e_1=\frac{\e}{C_1 +
C_2}$. Then, for all $n>\max\{N_2(\e_1),\,N_3(\e_1)\}$ and $\tau\in
[t,T]$,
$$
\mE\big|R^{(\cd)}\xi-P^n(\cd)\xi\big|_{H}^{\frac{4}{3}}<\e.
$$
Therefore, we obtain the desired result \eqref{10.30eq1}.

Further, similar to \eqref{P1}, one can show
that
\begin{equation}\label{PPP1}
P(t)\= R^{(t)} \ \mbox{ in }\
L^{\frac{4}{3}}_{\cF_t}(\O;H), \ \ \ae t\in
[0,T].
\end{equation}
Combining \eqref{6.1-8eq1} and \eqref{PPP1}, we
see that $\big(P(\cd), Q^{(\cd)},\widehat
Q^{(\cd)}\big)$ satisfies (\ref{6.18eq1}).
Hence, $\big(P(\cd), Q^{(\cd)},\widehat
Q^{(\cd)}\big)$ is a relaxed transposition
solution to \eqref{op-bsystem3}, and satisfies
the estimate (\ref{ine the1z}).

Finally, we show the uniqueness of the relaxed
transposition solution to \eqref{op-bsystem3}.
Assume that $\big(\overline P(\cd),\overline
Q^{(\cd)},\overline{\widehat Q}^{(\cd)}\big)\in
D_{\dbF,w}([0,T];
L^{\frac{4}{3}}(\O;\cL(H)))\times \cQ[0,T]$ is
another relaxed transposition solution to the
equation \eqref{op-bsystem3}. Then, by
Definition \ref{op-definition2x}, it follows
that
\begin{equation}\label{ueq1}
\begin{array}{ll}
\ds 0=\mE\Big\langle \(\overline P(t) - P(t)\)
\xi_1,\xi_2 \Big\rangle_{H} + \mE \int_t^T
\Big\langle \(\overline P(s) - P(s)\)u_1(s),
x_2(s)\Big\rangle_{H}ds \\
\ns\ds \qq + \mE \int_t^T \Big\langle
\(\overline P(s) - P(s)\)x_1(s),
u_2(s)\Big\rangle_{H}ds  + \mE \int_t^T
\Big\langle \(\overline P(s) - P(s)\) K
(s)x_1 (s), v_2 (s)\Big\rangle_{H}ds \\
\ns\ds \qq  +  \mE \int_t^T \Big\langle
\(\overline P(s) - P(s)\)v_1
(s), K (s)x_2 (s) +v_2(s) \Big\rangle_{H}ds \\
\ns\ds \qq + \mE \!\int_t^T \!\Big\langle
v_1(s), \!\(\overline{\widehat Q}^{(t)} \!-
\!\widehat Q^{(t)}
\)(\xi_2,u_2,v_2)(s)\Big\rangle_{H}ds \!+\! \mE
\int_t^T\! \Big\langle \(\overline Q^{(t)} \!
-\! Q^{(t)} \)
(\xi_1,u_1,v_1)(s), v_2(s) \!\Big\rangle_{H}ds,\\
\ns\ds\hspace{13cm}\forall\; t\in [0,T].
\end{array}
\end{equation}

Choosing $u_1=u_2=0$ and $v_1=v_2=0$ in the
test equations \eqref{op-fsystem2} and
\eqref{op-fsystem3}, by (\ref{ueq1}), we obtain
 that, for any $t\in [0,T]$,
$$
0=\mE\Big\langle \(\overline P(t) - P(t)\)
\xi_1,\xi_2 \Big\rangle_{H},\qq\forall\;
\xi_1,\; \xi_2\in L^4_{\cF_t}(\O;H).
$$
Hence, we find that $\overline P(\cd)=P(\cd)$.
By this, it is easy to see that (\ref{ueq1})
becomes that
\begin{equation}\label{u111ezq1}
\begin{array}{ll}
\ds 0=\mE \!\int_t^T \!\Big\langle v_1(s),
\!\(\overline{\widehat Q^{(t)}} \!- \!\widehat
Q^{(t)} \)(\xi_2,u_2,v_2)(s)\Big\rangle_{H}ds
\!+\! \mE \int_t^T\! \Big\langle \(\overline
Q^{(t)} \! -\! Q^{(t)} \)
(\xi_1,u_1,v_1)(s), v_2(s) \!\Big\rangle_{H}ds,\\
\ns\ds\hspace{13cm}\forall\; t\in [0,T].
\end{array}
\end{equation}
Choosing $v_2=0$ in the test equation
\eqref{op-fsystem3}, we see that
\eqref{u111ezq1} becomes
 \begin{equation}\label{uezq1}
0=\mE \int_t^T  \Big\langle v_1(s),
 \(\overline{\widehat Q^{(t)}} \!- \!\widehat
Q^{(t)} \)(\xi_2,u_2,0)(s)\Big\rangle_{H}ds.
\end{equation}
Noting that $v_1$ is arbitrarily in
$L^2_\dbF(0,T;L^4(\O;H))$, we conclude from
\eqref{uezq1} that $ \overline{\widehat
Q^{(t)}}(\cd,\cd,0) = \widehat
Q^{(t)}(\cd,\cd,0)$. Similarly, $
Q^{(t)}(\cd,\cd,0) = Q^{(t)}(\cd,\cd,0)$.
Hence,
\begin{equation}\label{uezq2}
\begin{array}{ll}
\ds 0=\mE \int_t^T \!\Big\langle v_1(s),
\(\overline{\widehat Q^{(t)}} \!- \!\widehat
Q^{(t)} \)(0,0,v_2)(s)\Big\rangle_{H}ds \! +\!
\mE \int_t^T \!\Big\langle \(\overline Q^{(t)}
\! - \! Q^{(t)} \!\) (0,0,v_1)(s), \!v_2(s)
\!\Big\rangle_{H}ds.
\end{array}
\end{equation}
Since $\overline
Q^{(t)}(0,0,\cd)^*=\overline{\widehat
Q^{(t)}}(0,0,\cd)$ and
$Q^{(t)}(0,0,\cd)^*=\widehat Q^{(t)}(0,0,\cd)$,
from \eqref{uezq2}, we find that
\begin{equation}\label{uezq3}
\begin{array}{ll}
\ds 0=2\mE \int_t^T \!\Big\langle v_1(s),
\(\overline{\widehat Q^{(t)}} \!- \!\widehat
Q^{(t)} \)(0,0,v_2)(s)\Big\rangle_{H}ds,
\end{array}
\end{equation}
which implies that $ \overline
Q^{(t)}(0,0,\cd)= Q^{(t)} (0,0,\cd)$ and $
\overline{\widehat Q^{(t)}}(0,0,\cd)= \widehat
Q^{(t)} (0,0,\cd)$. Hence $\overline
Q^{(t)}(\cd,\cd,\cd) = Q^{(t)}(\cd,\cd,\cd)$
and $\overline{\widehat Q^{(t)}}(\cd,\cd,\cd) =
\widehat Q^{(t)}(\cd,\cd,\cd)$. This completes
the proof of Theorem \ref{op well th}.
\endpf

\br
1) From the variational identity \eqref{eq2},
it is quite easy to obtain an {\it a priori}
estimate for $P^{n,\l_m}$ with respect to $n$
(See \eqref{eq205}). However, from the same
identity, it is clear that $Q^{n,\l_m}$ is not
coercive, and therefore, it is very hard to
derive any {\it a priori} estimate for
$Q^{n,\l_m}$. This is the main obstacle to
prove the existence of transposition solution
to the equation \eqref{op-bsystem3} in the
general case. As a remedy, we introduce four
operators $Q_1^{n,\l_m,t}$, $\widehat
Q_1^{n,\l_m,t}$, $Q_2^{n,\l_m,t}$ and $\widehat
Q_2^{n,\l_m,t}$ and the bilinear functional
$\cB_{n,\l_{m},t}(\cd,\cd)$ so that one can
obtain suitable {\it a priori} estimates and
take limit in some sense, and via which we are
able to establish the existence of relaxed
transposition solution to \eqref{op-bsystem3}
with general data.

\ms

2) Alternatively, one may use Theorem \ref{op
well th0} (instead of Theorem \ref{the1} (or
\cite[Theorem 4.1]{LZ})) to prove Theorem
\ref{op well th} (by approximating the data
$P_T$ and $F$ respectively by a sequence of
$\{P_T^k\}_{k=1}^\infty\subset
L^2_{\cF_T}(\O;\cL_2(H))$ and
$\{F^k\}_{k=1}^\infty\subset
L^1_\dbF(0,T;L^2(\O;\cL_2(H)))$ in the strong
operator topology). Nevertheless, in some
sense, our present proof seems to be more close
to the numerical approach to solve the equation
\eqref{op-bsystem3}.
\er

\section{Some properties of the relaxed transposition solutions to the operator-valued BSEEs}\label{as4}

In this section, we shall derive some
properties for the relaxed transposition
solutions to the equation \eqref{op-bsystem3}.
These properties will play key roles in the
proof of our general Pontryagin-type stochastic
maximum principle, presented in Section
\ref{s7}.

The following result shows the local Lipschitz
continuity of the relaxed transposition
solution to \eqref{op-bsystem3} with respect to
its coefficient $K$.

\begin{theorem}\label{10.2th1}
Let the assumptions in Theorem \ref{op well th}
hold and let $(P(\cd),Q^{(\cd)},\widehat
Q^{(\cd)})$ be the relaxed transposition
solution to \eqref{op-bsystem3}. Let
$K^\vartriangle\in
L^4_\dbF(0,T;L^\infty(\O;\cL(H)))$ and let
$(P^\vartriangle(\cd),Q^{(\cd,\vartriangle)},\widehat
Q^{(\cd,\vartriangle)})$ be the relaxed
transposition solution to the equation
\eqref{op-bsystem3} with $K$ replaced by
$K^\vartriangle$. Then,
\begin{equation}\label{10.2th1 eq1}
\begin{array}{ll}\ds
\big|\!\big|Q^{(0)}(0,0,\cd)-Q^{(0,\vartriangle)}(0,0,\cd)\big|\!\big|_{\cL(L^2_\dbF(0,T;L^4(\O;H)),\;L^2_\dbF(0,T;L^{\frac{4}{3}}(\O;H)))}\\
\ns\ds \q+\big|\!\big|\widehat
Q^{(0)}(0,0,\cd)-\widehat
Q^{(0,\vartriangle)}(0,0,\cd)\big|\!\big|_{\cL(L^2_\dbF(0,T;L^4(\O;H)),\;L^2_\dbF(0,T;L^{\frac{4}{3}}(\O;H)))}\\\ns\ds
\leq C(K^\vartriangle)|K
-K^\vartriangle|_{L^4_\dbF(0,T;L^\infty(\O;\cL(H)))}.
\end{array}
\end{equation}
Here the positive constant $C(K^\vartriangle)$
depends on $A$, $T$,
$|J|_{L^4_\dbF(0,T;L^\infty(\O;\cL(H)))}$,
$|K|_{L^4_\dbF(0,T;L^\infty(\O;\cL(H)))}$,
$|K^\vartriangle|_{L^4_\dbF(0,T;L^\infty(\O;\cL(H)))}$,
$|P_T|_{L^2_{\cF_T}(\O;\cL(H))}$ and
$|F|_{L^1_\dbF(0,T;L^2(\O;\cL(H)))}$.
\end{theorem}

{\it Proof}\,: The proof is divided into
several steps.

\ms

{\bf Step 1}. For any $t\in[0,T]$, consider the
following two equations:
\begin{equation}\label{op-fsystem2l}
\left\{
\begin{array}{ll}
\ds dx_1^\vartriangle  = (A+J)x_1^\vartriangle ds + u_1ds + K^\vartriangle x_1^\vartriangle  dw(s) + v_1 dw(s) &\mbox{ in } (t,T],\\
\ns\ds x_1^\vartriangle (t)=\xi_1
\end{array}
\right.
\end{equation}
and
\begin{equation}\label{op-fsystem3l}
\left\{
\begin{array}{ll}
\ds dx_2^\vartriangle  = (A+J)x_2^\vartriangle ds + u_2ds + K^\vartriangle x_2^\vartriangle  dw(s) + v_2 dw(s) &\mbox{ in } (t,T],\\
\ns\ds x_2^\vartriangle (t)=\xi_2.
\end{array}
\right.
\end{equation}
Here $\xi_1,\xi_2 \in L^4_{\cF_t}(\O;H)$,
$u_1,u_2\in L^2_\dbF(t,T;L^4(\O;H))$ and
$v_1,v_2\in L^2_\dbF(t,T;L^4(\O;H))$ are the
same as that in
\eqref{op-fsystem2}--\eqref{op-fsystem3}.
Clearly, for any $s\in [t,T]$, it holds that
 \begin{equation}\label{10.3eq2}
\3n\begin{array}{ll}\ds \mE|x_1^\vartriangle
(s)|_H^4 \3n&\ds= \mE\Big| S(s-t)\xi_1 +
\int_t^s S(s-\tau)J(\tau)x_1^\vartriangle
(\tau)d\tau + \int_t^s
S(s-\tau)u_1(\tau)d\tau \\
\ns&\ds\q + \int_t^s
S(s-\tau)K^\vartriangle(\tau)x_1^\vartriangle
(\tau)dw + \int_t^s
S(s-\tau)v_1(\tau)dw \Big|_H^4\\
\ns&\ds \leq\! C\!\[\mE|S(s-t)\xi_1|_H^4 \!+
\!\mE\Big|\int_t^s
S(s-\tau)J(\tau)x_1^\vartriangle (\tau)d\tau
\Big|_H^4 + \mE\Big|\int_t^s
S(s-\tau)u_1(\tau)d\tau
\Big|_H^4 \\
\ns&\ds \q + \mE\Big|\int_t^s
S(s-\tau)K^\vartriangle (\tau)x_1^\vartriangle
(\tau)dw\Big|_H^4 + \mE\Big|\int_t^s
S(s-\tau)v_1(\tau)dw \Big|_H^4\].
\end{array}
\end{equation}
By Lemma \ref{BDG}, it is easy to see that
$$
\begin{array}{ll}\ds
\mE\Big|\int_t^s S(s-\tau)K^\vartriangle
(\tau)x_1^\vartriangle (\tau)dw\Big|_H^4 \leq
C\mE\[\int_t^s
|S(s-\tau)K^\vartriangle(\tau)x_1^\vartriangle (\tau)|_H^2d\tau \]^2\\
\ns\ds\leq C\int_t^s\mE\[ | K^\vartriangle
(\tau)x_1^\vartriangle (\tau)|_H\]^4d\tau\ds
\leq C  \int_t^s|K^\vartriangle(\tau)
|^4_{L^\infty(\O;\cL(H))}\mE|x_1^\vartriangle
(\tau)|^4_{H}d\tau.
\end{array}
$$
This, together with \eqref{10.3eq2}, implies
that
 \bel{gogi1}
\begin{array}{ll}\ds
\mE|x_1^\vartriangle (s)|_H^4  &\ds\leq
C\big[|\xi_1|_{L^4_{\cF_0}(\O;H)}^4
+|u_1|^4_{L^2_\dbF(0,T;L^4(\O;H))}+|v_1|^4_{L^2_\dbF(0,T;L^4(\O;H))}\big]\\
\ns&\ds \q + C\int_t^s[|J(\tau)
|^4_{L^\infty(\O;\cL(H))}+|K^\vartriangle(\tau)
|^4_{L^\infty(\O;\cL(H))}]\mE|x_1^\vartriangle
(\tau)|^4_{H}d\tau.
\end{array}
 \ee
By Gronwall's  inequality, we obtain that
\begin{equation}\label{10.3eq1}
|x_1^\vartriangle
|_{L^\infty_\dbF(t,T;L^4(\O;H))} \leq
C(K^\vartriangle)\big(|\xi_1|_{L^4_{\cF_t}(\O;H)}
+ |u_1|_{L^2_\dbF(t,T;L^4(\O;H))} +
|v_1|_{L^2_\dbF(t,T;L^4(\O;H))} \big).
\end{equation}
Similarly,
\begin{equation}\label{10.3eq3}
|x_2^\vartriangle
|_{L^\infty_\dbF(t,T;L^4(\O;H))} \leq
C(K^\vartriangle)\big(|\xi_2|_{L^4_{\cF_t}(\O;H)}
+ |u_2|_{L^2_\dbF(t,T;L^4(\O;H))} +
|v_2|_{L^2_\dbF(t,T;L^4(\O;H))} \big).
\end{equation}

Let $y_1^\vartriangle = x_1-x_1^\vartriangle $
and $y_2^\vartriangle = x_2-x_2^\vartriangle $.
From \eqref{op-fsystem2} and
\eqref{op-fsystem2l}, we see that
$y_1^\vartriangle $ solves
\begin{equation}\label{op-fsystem2ll}
\left\{
\begin{array}{ll}
\ds dy_1^\vartriangle  = (A+J)y_1^\vartriangle ds  + K y_1^\vartriangle  dw(s) + (K-K^\vartriangle )x_1^\vartriangle dw(s)   &\mbox{ in } (t,T],\\
\ns\ds y_1^\vartriangle (t)=0.
\end{array}
\right.
\end{equation}
Then, similar to \eqref{gogi1} and by
(\ref{10.3eq1}), we have
$$
\begin{array}{ll}\ds
\mE|y_1^\vartriangle (s)|^4_{H} \\\ns\ds\ds\leq
C|(K-K^\vartriangle
)x_1^\vartriangle|^4_{L^2_\dbF(0,T;L^4(\O;H))}
+ C\int_t^s[|J(\tau)
|^4_{L^\infty(\O;\cL(H))}+|K(\tau)
|^4_{L^\infty(\O;\cL(H))}]\mE|y_1^\vartriangle
(\tau)|^4_{H}d\tau
\\
\ns\ds \leq C(K^\vartriangle)|K-K^\vartriangle
|^4_{L^4_\dbF(0,T;L^\infty(\O;\cL(H)))}\big(|\xi_1|_{L^4_{\cF_t}(\O;H)}
+ |u_1|_{L^2_\dbF(t,T;L^4(\O;H))} +
|v_1|_{L^2_\dbF(t,T;L^4(\O;H))}
\big)^4\\
\ns\ds \q + C\int_t^s[|J(\tau)
|^4_{L^\infty(\O;\cL(H))}+|K(\tau)
|^4_{L^\infty(\O;\cL(H))}]\mE|y_1^\vartriangle
(\tau)|^4_{H}d\tau.
\end{array}
$$
This, together with the Gronwall's inequality,
implies that
 \bel{go1}
 \ba{ll}\ds
\sup_{s\in [t,T]}\mE|x_1(s)-x_1^\vartriangle
(s)|_{H}=\sup_{s\in [t,T]}\mE|y_1^\vartriangle
(s)|_{H}\\\ns\ds\leq
C(K^\vartriangle)|K-K^\vartriangle
|_{L^4_\dbF(0,T;L^\infty(\O;\cL(H)))}\big(|\xi_1|_{L^4_{\cF_t}(\O;H)}
+ |u_1|_{L^2_\dbF(t,T;L^4(\O;H))} +
|v_1|_{L^2_\dbF(t,T;L^4(\O;H))} \big).
 \ea
 \ee
Similarly,
 \bel{go2}
 \!\!\ba{ll}
 \ds
\sup_{s\in [t,T]}\mE|x_2(s)-x_2^\vartriangle
(s)|_{H}\\\ns\ds \leq
\!C(K^\vartriangle)|K-K^\vartriangle
|_{L^4_\dbF(0,T;L^\infty(\O;\cL(H)))}\!\big(|\xi_2|_{L^4_{\cF_t}(\O;H)}
+ |u_2|_{L^2_\dbF(t,T;L^4(\O;H))} +
|v_2|_{L^2_\dbF(t,T;L^4(\O;H))} \big).
 \ea
 \ee

\ms

{\bf Step 2}. By Definition
\ref{op-definition2x}, it follows that
\begin{equation}\label{10.2eq2}
\begin{array}{ll}
\ds \q\mE\big\langle P_T x_1(T), x_2(T)
\big\rangle_{H} - \mE\big\langle P_T
x_1^\vartriangle (T), x_2^\vartriangle (T)
\big\rangle_{H} - \mE \int_t^T \big\langle
F(s) x_1(s), x_2(s) \big\rangle_{H}ds\\
\ns\ds \q + \mE \int_t^T \big\langle F(s)
x_1^\vartriangle (s), x_2^\vartriangle (s) \big\rangle_{H}ds\\
\ns\ds = \mE\big\langle
\big(P(t)-P^\vartriangle (t)\big)\xi_1,
\xi_2\big\rangle_{H}   + \mE \int_t^T
\big\langle P(s)u_1(s), x_2(s)\big\rangle_{H}ds
- \mE \int_t^T \big\langle P^\vartriangle
(s)u_1(s),
x_2^\vartriangle (s)\big\rangle_{H}ds \\
\ns\ds \q + \mE \int_t^T \big\langle P(s)
x_1(s), u_2(s) \big\rangle_{H}ds - \mE \int_t^T
\big\langle P^\vartriangle (s)
x_1^\vartriangle(s), u_2(s)
\big\rangle_{H}ds \\
\ns\ds\q + \mE \int_t^T \big\langle P(s) K
(s)x_1 (s), v_2 (s)\big\rangle_{H}ds - \mE
\int_t^T \big\langle P^\vartriangle (s)
K^\vartriangle (s)x_1^\vartriangle (s), v_2
(s)\big\rangle_{H}ds \\
\ns\ds \q +  \mE \int_t^T \big\langle  P(s)v_1
(s), K (s)x_2 (s)+v_2(s)\big\rangle_{H}ds-\mE
\int_t^T \big\langle  P^\vartriangle (s)v_1 (s), K^\vartriangle (s)x_2^\vartriangle (s)+v_2(s)\big\rangle_{H}ds\\
\ns\ds \q + \mE \int_t^T \big\langle v_1(s),
\widehat
Q^{(t)}(\xi_2,u_2,v_2)(s)\big\rangle_{H}ds -\mE
\int_t^T \big\langle v_1(s), \widehat Q^{(t,\vartriangle )}(\xi_2,u_2,v_2)(s)\big\rangle_{H}ds\\
\ns\ds \q + \mE \int_t^T \big\langle
Q^{(t)}(\xi_1,u_1,v_1)(s), v_2(s)
\big\rangle_{H}ds - \mE \int_t^T \big\langle
Q^{(t,\vartriangle )}(\xi_1,u_1,v_1)(s),
v_2(s) \big\rangle_{H}ds,\\
\ns\ds\qq\q\forall\;(\xi_1,u_1,v_1),(\xi_2,u_2,v_2)
\in L^4_{\cF_{t}}(\O;H)\times
L^2_\dbF(t,T;L^4(\O;H)) \times
L^2_\dbF(t,T;L^4(\O;H)).
\end{array}
\end{equation}

Letting  $u_1=u_2=v_1=v_2=0$ in the test
equations \eqref{op-fsystem2} and
\eqref{op-fsystem3}, respectively, from
\eqref{10.2eq2}, we find that
\begin{equation}\label{10.2eq3}
\begin{array}{ll}
\ds\q \mE\big\langle \big(P(t)-P^\vartriangle
(t)\big)\xi_1,
\xi_2\big\rangle_{H}\\
\ns\ds = \mE\big\langle P_T x_1(T), x_2(T)
\big\rangle_{H} - \mE\big\langle
P_T x_1^\vartriangle (T), x_2^\vartriangle (T) \big\rangle_{H}\\
\ns\ds\q - \mE \int_0^T \big\langle F(s)
x_1(s), x_2(s) \big\rangle_{H}ds + \mE \int_t^T
\big\langle F(s) x_1^\vartriangle (s),
x_2^\vartriangle (s)
\big\rangle_{H}ds\\
\ns\ds = \mE\big\langle P_T
[x_1(T)-x_1^\vartriangle (T)], x_2(T)
\big\rangle_{H}  + \mE\big\langle P_T
x_1^\vartriangle (T), x_2(T)-x_2^\vartriangle
(T)
\big\rangle_{H} \\
\ns\ds\q - \mE \int_0^T \big\langle F(s)
\big[x_1(s)-x_1^\vartriangle (s)\big], x_2(s)
\big\rangle_{H}ds + \mE \int_t^T \big\langle
F(s) x_1^\vartriangle (s), x_2(s) -
x_2^\vartriangle (s) \big\rangle_{H}ds.
\end{array}
\end{equation}
In (\ref{10.2eq3}), we choose $\xi_1,\xi_2\in
L^4_{\cF_t}(\O;H)$ with
$|\xi_1|_{L^4_{\cF_t}(\O;H)}=|\xi_2|_{L^4_{\cF_t}(\O;H)}=1$,
such that
$$
\mE\big\langle \big(P(t)-P^\vartriangle
(t)\big)\xi_1, \xi_2\big\rangle_{H} \geq
\frac{1}{2}\big|\!\big|P(t)-P^\vartriangle
(t)\big|\!\big|_{\cL(L^4_{\cF_t}(\O;H),\;L^{\frac{4}{3}}_{\cF_t}(\O;H))}.
$$
On the other hand, by
(\ref{10.3eq1})--(\ref{10.3eq3}) and
(\ref{go1})--(\ref{go2}), we have
\begin{equation}\label{10.2eq5}
\begin{array}{ll}\ds
\q\Big|\mE\big\langle P_T
[x_1(T)-x_1^\vartriangle (T)], x_2(T)
\big\rangle_{H} + \mE\big\langle P_T
x_1^\vartriangle (T), x_2(T)-x_2^\vartriangle
(T)
\big\rangle_{H} \\
\ns\ds \q - \mE \int_0^T \big\langle F(s)
\big[x_1(s)-x_1^\vartriangle (s)\big], x_2(s)
\big\rangle_{H}ds + \mE \int_t^T \big\langle
F(s) x_1^\vartriangle (s), x_2(s) -
x_2^\vartriangle (s)
\big\rangle_{H}ds\Big|\\
\ns\ds \leq
C(K^\vartriangle)|P_T|_{L^2_{\cF_T}(\O;\cL(H))}\big(|x_1(T)-x_1^\vartriangle
(T)|_{L^4_{\cF_T}(\O;H)}
+|x_2(T)-x_2^\vartriangle (T)|_{L^4_{\cF_T}(\O;H)}\big)\\
\ns\ds\q +
C(K^\vartriangle)|F|_{L^1_\dbF(0,T;L^2(\O;\cL(H)))}\big(|x_1-x_1^\vartriangle |_{L^\infty_\dbF(t,T;L^4(\O;H))}+|x_2-x_2^\vartriangle |_{L^\infty_\dbF(t,T;L^4(\O;H))}\big)\\
\ns\ds \leq C(K^\vartriangle)|K-K^\vartriangle
|_{L^4_\dbF(0,T;L^\infty(\O;\cL(H)))}.
\end{array}
\end{equation}
Hence,
\begin{equation}\label{10.2eq6}
\big|\!\big|P(t)-P^\vartriangle
(t)\big|\!\big|_{\cL(L^4_{\cF_t}(\O;H),\;L^{\frac{4}{3}}_{\cF_t}(\O;H))}
\leq C(K^\vartriangle)|K-K^\vartriangle
|_{L^4_\dbF(0,T;L^\infty(\O;\cL(H)))},\q\forall\;t\in
[0,T].
\end{equation}

\ms

{\bf Step 3}. Letting  $\xi_1=\xi_2=0$ and
$u_1=u_2=0$ in the test equations
\eqref{op-fsystem2} and \eqref{op-fsystem3}
respectively, from \eqref{10.2eq2} and noting
that
$$
\mE \int_0^T \big\langle v_1(s), \widehat
Q^{(t)}(0,0,v_2)(s)\big\rangle_{H}ds = \mE
\int_0^T \big\langle Q^{(t)}(0,0,v_1)(s),
v_2(s) \big\rangle_{H}ds
$$
and
$$
\mE \int_0^T \big\langle v_1(s), \widehat
Q^{(0,\vartriangle
)}(0,0,v_2)(s)\big\rangle_{H}ds = \mE \int_0^T
\big\langle Q^{(0,\vartriangle )}(0,0,v_1)(s),
v_2(s) \big\rangle_{H}ds,
$$
we find that, for any $v_1,v_2\in L^2_\dbF(0,T;L^4(\O;H))$, it holds
that
\begin{equation}\label{10.2eq4}
\begin{array}{ll}
\ds \q\mE\big\langle P_T x_1(T), x_2(T)
\big\rangle_{H} - \mE\big\langle P_T
x_1^\vartriangle (T), x_2^\vartriangle (T)
\big\rangle_{H} - \mE \int_0^T \big\langle
F(s) x_1(s), x_2(s) \big\rangle_{H}ds\\
\ns\ds \q + \mE \int_0^T \big\langle F(s)
x_1^\vartriangle (s), x_2^\vartriangle (s) \big\rangle_{H}ds\\
\ns\ds =   \mE \int_0^T \big\langle P(s) K
(s)x_1 (s), v_2 (s)\big\rangle_{H}ds - \mE
\int_0^T \big\langle P^\vartriangle (s)
K^\vartriangle (s)x_1^\vartriangle (s), v_2
(s)\big\rangle_{H}ds \\
\ns\ds \q +  \mE \int_0^T \big\langle  P(s)v_1
(s), K (s)x_2 (s)+v_2(s)\big\rangle_{H}ds-\mE
\int_0^T \big\langle  P^\vartriangle (s)v_1 (s), K^\vartriangle (s)x_2^\vartriangle (s)+v_2(s)\big\rangle_{H}ds\\
\ns\ds \q + 2\mE \int_0^T \big\langle Q^{(0)}(0,0,v_1)(s), v_2(s)
\big\rangle_{H}ds - 2\mE \int_0^T \big\langle Q^{(0,\vartriangle
)}(0,0,v_1)(s), v_2(s) \big\rangle_{H}ds.
\end{array}
\end{equation}
We choose $v_1,v_2\in L^2_\dbF(0,T;L^4(\O;H))$
with
$|v_1|_{L^2_\dbF(0,T;L^4(\O;H))}=|v_2|_{L^2_\dbF(0,T;L^4(\O;H))}=1$,
such that
\begin{equation}\label{10.2eq7}
\begin{array}{ll}\ds
\q 2\mE \int_0^T \big\langle
Q^{(0)}(0,0,v_1)(s), v_2(s) \big\rangle_{H}ds -
2\mE \int_0^T \big\langle Q^{(0,\vartriangle
)}(0,0,v_1)(s), v_2(s)
\big\rangle_{H}ds\\
\ns\ds \geq
\big|\!\big|Q^{(0)}(0,0,\cd)-Q^{(0,\vartriangle
)}(0,0,\cd)\big|\!\big|_{\cL(L^2_\dbF(0,T;L^4(\O;H)),\;L^2_\dbF(0,T;L^{\frac{4}{3}}(\O;H)))}.
\end{array}
\end{equation}
By the above choice of $v_1$ and $v_2$, similar
to \eqref{10.2eq5}, we have
\begin{equation}\label{10.2eq8}
\begin{array}{ll}\ds
\ds\Big|\mE\big\langle P_T x_1(T), x_2(T)
\big\rangle_{H} - \mE\big\langle P_T
x_1^\vartriangle (T), x_2^\vartriangle (T)
\big\rangle_{H} - \mE \int_0^T \big\langle
F(s) x_1(s), x_2(s) \big\rangle_{H}ds\\
\ns\ds  \q + \mE \int_0^T \big\langle F(s)
x_1^\vartriangle (s),
x_2^\vartriangle (s) \big\rangle_{H}ds\Big|\\
\ns\ds \leq C(K^\vartriangle)|K-K^\vartriangle
|_{L^4_\dbF(0,T;L^\infty(\O;\cL(H)))}.
\end{array}
\end{equation}
By \eqref{10.2eq6}, it follows that
\begin{equation}\label{10.2eq9}
\begin{array}{ll}\ds
\q\Big|\mE \int_0^T \big\langle P(s) K (s)x_1
(s), v_2 (s)\big\rangle_{H}ds - \mE \int_0^T
\big\langle P^\vartriangle (s) K^\vartriangle
(s)x_1^\vartriangle (s), v_2
(s)\big\rangle_{H}ds\Big|\\
\ns\ds \leq \Big|\mE \int_0^T \big\langle P(s)
K (s)\big[x_1(s)-x_1^\vartriangle (s)\big], v_2
(s)\big\rangle_{H}ds \Big|  + \Big|\mE \int_0^T
\big\langle P(s) \big[K(s)-K^\vartriangle
(s)\big]x_1^\vartriangle (s),
v_2 (s)\big\rangle_{H}ds\Big|\\
\ns\ds \q+ \Big|\mE \int_0^T \big\langle
\big[P(s) -P^\vartriangle (s)\big]
K^\vartriangle (s)x_1^\vartriangle (s), v_2
(s)\big\rangle_{H}ds\Big|\\
\ns\ds \leq C(K^\vartriangle)|K-K^\vartriangle
|_{L^4_\dbF(0,T;L^\infty(\O;\cL(H)))}.
\end{array}
\end{equation}
Similarly,
\begin{equation}\label{10.2eq11}
 \ba{ll}\ds
\Big|\mE \int_0^T \big\langle  P(s)v_1 (s), K
(s)x_2 (s)+v_2(s)\big\rangle_{H}ds-\mE \int_0^T
\big\langle P^\vartriangle (s)v_1 (s),
K^\vartriangle (s)x_2^\vartriangle
(s)+v_2(s)\big\rangle_{H}ds \\
\ns\ds \leq C(K^\vartriangle)|K-K^\vartriangle
|_{L^4_\dbF(0,T;L^\infty(\O;\cL(H)))}.
\end{array}
\end{equation}
From \eqref{10.2eq4}--\eqref{10.2eq11}, we
obtain that
 $$
 \begin{array}{ll}\ds
\big|\!\big|Q^{(0)}(0,0,\cd)-Q^{(0,\vartriangle)}(0,0,\cd)\big|\!\big|_{\cL(L^2_\dbF(0,T;L^4(\O;H)),\;L^2_\dbF(0,T;L^{\frac{4}{3}}(\O;H)))}\\\ns\ds
\leq C(K^\vartriangle)|K
-K^\vartriangle|_{L^4_\dbF(0,T;L^\infty(\O;\cL(H)))}.
 \ea
 $$
Similarly,
 $$
 \ba{ll}
 \ds\big|\!\big|\widehat Q^{(0)}(0,0,\cd)-\widehat
Q^{(0,\vartriangle)}(0,0,\cd)\big|\!\big|_{\cL(L^2_\dbF(0,T;L^4(\O;H)),\;L^2_\dbF(0,T;L^{\frac{4}{3}}(\O;H)))}\\\ns\ds
\leq C(K^\vartriangle)|K
-K^\vartriangle|_{L^4_\dbF(0,T;L^\infty(\O;\cL(H)))}.
 \ea
 $$
Hence, we obtain the desired estimate
(\ref{10.2th1 eq1}). This completes the proof
of Theorem \ref{10.2th1}.
\endpf

\ms

Next, we shall show a property of the relaxed
transposition solution to \eqref{op-bsystem3}
when the coefficient $K$ is piecewisely
constant with respect to the time variable. For
this purpose, we introduce the following
subspace of $L^{2}_\dbF(0,T;L^{4}(\O;H))$
(Recall (\ref{x1z}) for the definition of
$\cM$):
\begin{equation}\label{cH}
\cH=\Big\{\sum_{i=1}^\ell
\chi_{O_i}(\cd)h_i\;\Big|\; \ell\in\dbN,\,
O_i\in\cM, \, h_i\in D(A)\Big\}.
\end{equation}
It is clear that $\cH$ is dense in
$L^{2}_\dbF(0,T;L^{4}(\O;H))$. We have the
following result.

\begin{theorem}\label{10.1th}
Suppose that the assumptions in Theorem \ref{op
well th} hold and $\ds K=\sum_{i=1}^{n_0}
\chi_{[t_i,t_{i+1})}(t)K_i$ for some
$n_0\in\dbN$,
$0=t_1<t_2<\cdots<t_{n_0}<t_{n_0+1}=T$, and
$K_i\in L^\infty_{\cF_{t_i}}(\O;\cL(D(A)))$,
$i=1,\cds,n_0$. Let $(P(\cd),Q^{(\cd)},\widehat
Q^{(\cd)})$ be the relaxed transposition
solution to \eqref{op-bsystem3}. Then, there
exist two pointwisely defined linear operators
$Q$ and $\widehat Q$, both of which are from
$\cH$ to $L^2_\dbF(0,T;L^{\frac{4}{3}}(\O;H))$,
such that
\begin{equation}\label{10.9eq2}
\begin{array}{ll}
\ds \q \mE \int_{0}^T \big\langle v_1(s),
\widehat Q^{(0)}(\xi_2,u_2,v_2)(s)
\big\rangle_{H}ds + \mE \int_{0}^T \big\langle
 Q^{(0)}(\xi_1,u_1,v_1) (s),
v_2(s) \big\rangle_{H}ds \\
\ns\ds =\mE \int_{0}^T
 \[\big\langle
\big(Q v_1\big)(s), x_2 (s)
\big\rangle_{H}+\big\langle x_1 (s),
\big(\widehat Q v_2\big)(s)
\big\rangle_{H}\]ds,
\end{array}
\end{equation}
holds for any $\xi_1,\xi_2\in
L^4_{\cF_0}(\O;H)$, $u_1(\cd), u_2(\cd)\in
L^2_{\dbF}(0,T;L^4(\O;H))$ and
$v_1(\cd),v_2(\cd)\in \cH$. Here, $x_1(\cd)$
and $x_2(\cd)$ solve accordingly
\eqref{op-fsystem2} and \eqref{op-fsystem3}
with $t=0$.
\end{theorem}

{\it Proof}\,: As in the proof of Theorem
\ref{op well th} (but with the set
$\{r_j\}_{j=1}^\infty$ (introduced at the very
beginning of Step 4) being replaced by
$\{r_j\}_{j=1}^\infty\cup
\{t_1,t_2,\cdots,t_{n_0}\}$), we introduce the
equation (\ref{op-bsystem4}) (approximating to
the equation \eqref{op-bsystem3}), and the
equations (\ref{op-fsystem5}) and
(\ref{op-fsystem6}) (which are accordingly
finite approximations of the equations
\eqref{op-fsystem2} and \eqref{op-fsystem3}),
and obtain the approximate variational equality
\eqref{eq2} for $P^{n,\l}(\cd)$ and
$Q^{n,\l}(\cd)$. Also, we fix a sequence
$\{\l_m\}_{m=1}^\infty\subset \rho(A)$ such
that $\l_m\to\infty$ as $m\to\infty$. We divide
the rest of proof into two steps.

\ms

{\bf Step 1}. Choose $\xi_2^n=0$ and $u_2^n=0$
in the equation \eqref{op-fsystem6}. Then,
there is a constant $C_1(\l_m)>0$, independent
of $t$ and $n$, such that
$|x_2^{n,\l}|_{L^\infty(t,T;L^4(\O;H))}\leq
C_1(\l_m)|v_2^n|_{L^2(t,T;L^4(\O;H))}$. Without
loss of generality, we may assume that
\begin{equation}\label{10.1l}
\frac{1}{\ds\max_{1\leq i\leq
n_0}(t_{i+1}-t_i)}>2^{12}|C_1(\l_m)|^4|K|^4_{L^4_\dbF(0,T;L^\infty(\O;\cL(D(A))))}.
\end{equation}
(Otherwise, we may choose a refined partition
of $[0,T]$ so that \eqref{10.1l} holds).
Letting $\xi_1\in
L^4_{\cF_{t_{n_0}}}(\O;D(A))$, $u_1^n =
-(A^{n,\l_m}+J_n)\xi_1$ and $v_1^n=-K_n\xi_1$
in \eqref{op-fsystem5}, and letting $\xi_2^n=0$
and $u_2^n=0$ in \eqref{op-fsystem6}, by
(\ref{eq2}) with $t=t_{n_0}$, we find that
\begin{equation}\label{6.16eq1}
\begin{array}{ll}
\ds \q\mE\big\langle P^{n}_T \xi_1,
x_2^{n,\l_m}(T) \big\rangle_{\dbR^n} - \mE
\int_{t_{n_0}}^T \big\langle
F_n(s) \xi_1, x_2^{n,\l_m}(s) \big\rangle_{\dbR^n}ds\\
\ns\ds = \mE \int_{t_{n_0}}^T \big\langle
P^{n,\l_m}(s)u_1^n(s),
x_2^{n,\l_m}(s)\big\rangle_{\dbR^n}ds  + \mE
\int_{t_{n_0}}^T \big\langle
P^{n,\l_m}(s)K_n(s)\xi_1,
v_2^n(s)\big\rangle_{\dbR^n}ds \\
\ns\ds \q   + \mE \int_{t_{n_0}}^T \big\langle
P^{n,\l_m}(s)v_1^{n}(s),
K_n(s)x_2^{n,\l_m} + v_2^n(s) \big\rangle_{\dbR^n}ds  \\
\ns\ds \q  - \mE \int_{t_{n_0}}^T \big\langle
Q^{n,\l_m}(s)K_n(s) \xi_1, x_2^{n,\l_m}(s)
\big\rangle_{\dbR^n}ds + \mE \int_{t_{n_0}}^T
\big\langle Q^{n,\l_m}(s) \xi_1^n, v_2^n(s)
\big\rangle_{\dbR^n}ds.
\end{array}
\end{equation}
First, we find a $\xi_1\in
L^4_{\cF_{t_{n_0}}}(\O;D(A))$ with
$|\xi_1|_{L^4_{\cF_{t_{n_0}}}(\O;D(A))}=1$ such
that
$$
\big|Q^{n,\l_m}(\cd)
\xi_1^n\big|_{L^{2}_\dbF(t_{n_0},T;L^{\frac{4}{3}}(\O;H))}
\geq
\frac{1}{2}\big|\!\big|Q^{n,\l_m}(\cd)\big|\!\big|_{\cL(L^4_{\cF_{t_{n_0}}}(\O;D(A)),\;
L^{2}_\dbF({t_{n_0}},T;L^{\frac{4}{3}}(\O;H)))}.
$$
Next, we find a $v_2\in
L^2_\dbF({t_{n_0}},T;L^4(\O;H))$ with
$|v_2|_{L^2_\dbF({t_{n_0}},T;L^4(\O;H))}=1$ so
that
$$
\mE \int_{t_{n_0}}^T \big\langle Q^{n,\l_m}(s)
\xi_1^n, v_2^n(s) \big\rangle_{\dbR^n}ds \geq
\frac{1}{2}\big|Q^{n,\l_m}(\cd)
\xi_1^n\big|_{L^{2}_\dbF({t_{n_0}},T;L^{\frac{4}{3}}(\O;H))}.
$$
Hence,
\begin{equation}\label{6.16eq2}
\mE \int_{t_{n_0}}^T \big\langle Q^{n,\l_m}(s)
\xi_1^n, v_2^n(s) \big\rangle_{\dbR^n}ds \geq
\frac{1}{4}\big|\!\big|Q^{n,\l_m}(\cd)\big|\!\big|_{\cL(L^4_{\cF_{t_{n_0}}}(\O;D(A)),\;
L^{2}_\dbF({t_{n_0}},T;L^{\frac{4}{3}}(\O;H)))}.
\end{equation}
On the other hand, by \eqref{10.1l}, it follows
that
$$
\begin{array}{ll}\ds
\q\Big|\mE \int_{t_{n_0}}^T \big\langle
Q^{n,\l_m}(s)K_n
\xi_1, x_2^{n,\l_m}(s) \big\rangle_{\dbR^n}ds\Big|\\
\ns\ds \leq
|x_2^{n,\l_m}|_{L^\infty_\dbF(t_{n_0},T;L^4(\O;H))}
\int_{t_{n_0}}^T
\big|Q^{n,\l_m}(s)K_n \xi_1 \big|_{L^{\frac{4}{3}}_{\cF_s}(\O;H)} ds\\
\ns\ds\leq
\sqrt{T-t_{n_0}}|x_2^{n,\l_m}|_{L^\infty_\dbF(t_{n_0},T;L^4(\O;H))}|K_{n_0}|_{
L^\infty_{\cF_{t_{n_0}}}\!\!(\O;\cL(D(A)))}\!\big|\!\big|Q^{n,\l_m}(\cd)\big|\!\big|_{\cL(L^4_{\cF_{t_{n_0}}}(\O;D(A)),\,
L^{2}_\dbF(t_{n_0},T;L^{\frac{4}{3}}(\O;H)))}\\
\ns\ds\leq
(t_{n_0+1}-t_{n_0})^{\frac14}C_1(\l_m)|K|_{L^4_\dbF(0,T;L^\infty(\O;\cL(D(A))))}\big|\!\big|Q^{n,\l_m}(\cd)\big|\!\big|_{\cL(L^4_{\cF_{t_{n_0}}}(\O;D(A)),\;
L^{2}_\dbF(t_{n_0},T;L^{\frac{4}{3}}(\O;H)))}\\
\ns\ds\leq
\frac{1}{8}\big|\!\big|Q^{n,\l_m}(\cd)\big|\!\big|_{\cL(L^4_{\cF_{t_{n_0}}}(\O;D(A)),\;
L^{2}_\dbF({t_{n_0}},T;L^{\frac{4}{3}}(\O;H)))}.
\end{array}
$$
This, together with \eqref{6.16eq2}, implies
that
 \bel{pp1}
\begin{array}{ll}\ds
\q\mE \int_{t_{n_0}}^T \big\langle
Q^{n,\l_m}(s)K_n(s) \xi_1, x_2^{n,\l_m}(s)
\big\rangle_{\dbR^n}ds + \mE \int_{t_{n_0}}^T
\big\langle Q^{n,\l_m}(s)
\xi_1^n, v_2^n(s) \big\rangle_{\dbR^n}ds \\
\ns\ds \geq
\frac{1}{8}\big|\!\big|Q^{n,\l_m}(\cd)\big|\!\big|_{\cL(L^4_{\cF_{t_{n_0}}}(\O;D(A)),\;
L^{2}_\dbF({t_{n_0}},T;L^{\frac{4}{3}}(\O;H)))}.
\end{array}
 \ee
On the other hand, from the choice of $\xi_1$,
$u_1^n$ and $v_1^n$, we find that
$$
|u_1^n|_{L^2_\dbF(0,T;L^4(\O;H))}\leq
C(\l_m)\big(|A|_{\cL(D(A),\;H)} +
|J|_{L^4_\dbF(0,T;L^\infty(\O;\cL(H)))}\big).
$$
Hence, by the estimate (\ref{eq205}), it
follows that
\begin{equation}\label{6.16eq4}
\begin{array}{ll}\ds
\q\Big|\mE\big\langle P^{n}_T \xi_1,
x_2^{n,\l_m}(T) \big\rangle_{\dbR^n} - \mE
\int_{t_{n_0}}^T \big\langle F_n(s) \xi_1,
x_2^{n,\l_m}(s) \big\rangle_{\dbR^n}ds - \mE
\int_{t_{n_0}}^T \big\langle
P^{n,\l_m}(s)u_1^n(s),
x_2^{n,\l_m}(s)\big\rangle_{\dbR^n}ds \\
\ns\ds \q - \mE \int_{t_{n_0}}^T \!\big\langle
K_n(s)\xi_1,
P^{n,\l_m}(s)v_2^n(s)\big\rangle_{\dbR^n}ds
-\!\mE \int_{t_{n_0}}^T \!\big\langle
P^{n,\l_m}(s)v_1^{n}(s), K_n(s)x_2^{n,\l_m} +
v_2^n(s) \big\rangle_{\dbR^n}ds  \Big|\\
\ns\ds \leq\!
C(\l_m)\big(|P_T|_{L^4_{\cF_T}(\O;\cL(H))} +
|F|_{L^1_\dbF(0,T;L^2(\O;\cL(H)))}\big)\big(1 +
|A|_{\cL(D(A),\;H)} +
|(J,K)|_{(L^4_\dbF(0,T;L^\infty(\O;\cL(H))))^2}\big).
\end{array}
\end{equation}
Combining \eqref{6.16eq1} and
\eqref{pp1}--\eqref{6.16eq4}, we find that
 \bel{pp2}
\begin{array}{ll}\ds
\q\big|\!\big|Q^{n,\l_m}(\cd)\big|\!\big|_{\cL(L^4_{\cF_{t_{n_0}}}(\O;D(A)),\;
L^{2}_\dbF({t_{n_0}},T;L^{\frac{4}{3}}(\O;H)))}\\
\ns\ds \leq\!
C(\l_m)\big(|P_T|_{L^4_{\cF_T}(\O;\cL(H))} +
|F|_{L^1_\dbF(0,T;L^2(\O;\cL(H)))}\big)\big(1 +
|A|_{\cL(D(A),\;H)} +
|(J,K)|_{(L^4_\dbF(0,T;L^\infty(\O;\cL(H))))^2}\big).
\end{array}
 \ee
By \eqref{pp2} and Corollary \ref{cor1}, there
exist a bounded, pointwisely defined linear
operator $ Q^{\l_m}_{t_{n_0}}$ from
$L^4_{\cF_{t_{n_0}}}(\O;D(A))$ to
$L^{2}_\dbF(t_{n_0},T;L^{\frac{4}{3}}(\O;H))$,
and a subsequence $\{n_k^{(5)}\}_{k=1}^\infty$
of $\{n_k^{(4)}\}_{k=1}^\infty$ (defined
between (\ref{eq405x}) and (\ref{eq406})) such
that
\begin{equation}\label{6.17eq1}
\mbox{(w)}\mbox{-}\lim_{k\to\infty}Q^{n_k^{(5)},\l_{m}}\xi
= Q^{\l_m}_{t_{n_0}}\xi \q\mbox{ in }
L^{2}_\dbF(t_{n_0},T;L^{\frac{4}{3}}(\O;H)),\qq
\forall\; \xi\in L^4_{\cF_{t_{n_0}}}(\O;D(A)).
\end{equation}

Next, letting $\xi_1\in
L^4_{\cF_{t_{n_0-1}}}(\O;D(A))$, $u_1^n =
-(A^{n,\l_m}+J_n)\xi_1$ and $v_1^n=-K_n\xi_1$
in \eqref{op-fsystem5}, and letting $\xi_2=0$
and $u_2^n=0$ in \eqref{op-fsystem6}, by
\eqref{eq2} with $t=t_{n_0-1}$, we find that
\begin{equation}\label{6.16eq1x}
\begin{array}{ll}
\ds \q\mE\big\langle P^{n}_T \xi_1,
x_2^{n,\l}(T) \big\rangle_{\dbR^n} - \mE
\int_{t_{n_0-1}}^T \big\langle
F_n(s) \xi_1, x_2^{n,\l_m}(s) \big\rangle_{\dbR^n}ds\\
\ns\ds = \mE \int_{t_{n_0-1}}^T\big\langle
P^{n,\l}(s)u_1^n(s),
x_2^{n,\l_m}(s)\big\rangle_{\dbR^n}ds  + \mE
\int_{t_{n_0-1}}^T\big\langle
P^{n,\l_m}(s)K_n(s)\xi_1,
v_2^n(s)\big\rangle_{\dbR^n}ds \\
\ns\ds \q   + \mE \int_{t_{n_0-1}}^T
\big\langle P^{n,\l_m}(s)v_1^{n}(s),
K_n(s)x_2^{n,\l_m} +
v_2^n(s) \big\rangle_{\dbR^n}ds  \\
\ns\ds \q  + \mE \int_{t_{n_0-1}}^T \big\langle
 Q^{n,\l_m}(s) K_n(s) \xi_1, x_2^{n,\l_m}(s)
\big\rangle_{\dbR^n}ds + \mE \int_{t_{n_0-1}}^T
\big\langle Q^{n,\l_m}(s) \xi_1^n, v_2^n(s)
\big\rangle_{\dbR^n}ds.
\end{array}
\end{equation}
On the other hand, for these data $\xi_1$,
$u_1^n$, $v_1^n$, $\xi_2$, $u_2^n$ and $v_2^n$,
from the variational equality \eqref{eq2} with
$t=t_{n_0}$, we obtain that
\begin{equation}\label{6.16eq1xx}
\begin{array}{ll}
\ds \q\mE\big\langle P^{n}_T \xi_1,
x_2^{n,\l_m}(T) \big\rangle_{\dbR^n} - \mE
\int_{t_{n_0}}^T \big\langle
F_n(s) \xi_1, x_2^{n,\l_m}(s) \big\rangle_{\dbR^n}ds\\
\ns\ds = \mE \big\langle
P^{n,\l_m}(t_{n_0})\xi_1,
x_2^{n,\l_m}(t_{n_0})\big\rangle_{\dbR^n}+\mE
\int_{t_{n_0}}^T \big\langle
P^{n,\l_m}(s)u_1^n(s),
x_2^{n,\l_m}(s)\big\rangle_{\dbR^n}ds  \\
\ns\ds \q + \mE \!\int_{t_{n_0}}^T \!
\!\big\langle P^{n,\l_m}(s)K_n(s)\xi_1,
v_2^n(s)\big\rangle_{\dbR^n}ds + \mE
\!\int_{t_{n_0}}^T \! \!\big\langle
P^{n,\l_m}(s)v_1^{n}(s), K_n(s)x_2^{n,\l_m} +
v_2^n(s) \big\rangle_{\dbR^n}ds   \\
\ns\ds \q  + \mE \int_{t_{n_0}}^T \big\langle
 Q^{n,\l_m}(s) K_n(s) \xi_1, x_2^{n,\l_m}(s)
\big\rangle_{\dbR^n}ds + \mE \int_{t_{n_0}}^T
\big\langle Q^{n,\l_m}(s) \xi_1^n, v_2^n(s)
\big\rangle_{\dbR^n}ds.
\end{array}
\end{equation}
From \eqref{6.16eq1x} and \eqref{6.16eq1xx}, it
follows that
\begin{equation}\label{6.16eq1xxx}
\3n\begin{array}{ll} \ds \q \mE \big\langle
P^{n,\l_m}(t_{n_0})\xi_1,
x_2^{n,\l_m}(t_{n_0})\big\rangle_{\dbR^n} -\mE
\int_{t_{n_0-1}}^{t_{n_0}} \big\langle
F_n(s) \xi_1, x_2^{n,\l_m}(s) \big\rangle_{\dbR^n}ds\\
\ns\ds = \mE \int_{t_{n_0-1}}^{t_{n_0}}
\big\langle P^{n,\l_m}(s)u_1^n(s),
x_2^{n,\l_m}(s)\big\rangle_{\dbR^n}ds + \mE
\int_{t_{n_0-1}}^{t_{n_0}} \big\langle
P^{n,\l_m}(s)K_n(s)\xi_1,
v_2^n(s)\big\rangle_{\dbR^n}ds\\
\ns\ds \q   + \mE
\int_{t_{n_0-1}}^{t_{n_0}}\big\langle
P^{n,\l_m}(s)v_1^{n}(s), K_n(s)x_2^{n,\l_m} +
v_2^n(s) \big\rangle_{\dbR^n}ds   \\
\ns\ds \q  -\mE
\int_{t_{n_0-1}}^{t_{n_0}}\big\langle
Q^{n,\l_m}(s)K_n(s) \xi_1, x_2^{n,\l_m}(s)
\big\rangle_{\dbR^n}ds + \mE
\int_{t_{n_0-1}}^{t_{n_0}} \big\langle
Q^{n,\l_m}(s) \xi_1^n, v_2^n(s)
\big\rangle_{\dbR^n}ds.
\end{array}
\end{equation}
Similar to \eqref{pp2}, from
\eqref{6.16eq1xxx}, one obtains that
 \bel{pp3}
\begin{array}{ll}\ds
\q\big|\!\big|Q^{n,\l_m}(\cd)\big|\!\big|_{\cL(L^4_{\cF_{t_{{n_0}-1}}}(\O;D(A)),\;
L^{2}_\dbF(t_{{n_0}-1},t_{n_0};L^{\frac{4}{3}}(\O;H)))}\\
\ns\ds \leq\!
C(\l_m)\big(|P_T|_{L^4_{\cF_T}(\O;\cL(H))} +
|F|_{L^1_\dbF(0,T;L^2(\O;\cL(H)))}\big)\big(1 +
|A|_{\cL(D(A),\;H)} +
|(J,K)|_{(L^4_\dbF(0,T;L^\infty(\O;\cL(H))))^2}\big).
\end{array}
 \ee
By \eqref{pp3} and utilizing  Corollary
\ref{cor1}, we conclude that there exist a
bounded, pointwisely defined linear operator $
Q^{\l_m}_{t_{n_0-1}}$ from
$L^4_{\cF_{t_{n_0-1}}}(\O;D(A))$ to
$L^{2}_\dbF(t_{n_0-1},t_{n_0};L^{\frac{4}{3}}(\O;H))$,
and a subsequence $\{n_k^{(6)}\}_{k=1}^\infty$
of $\{n_k^{(5)}\}_{k=1}^\infty$  such that
\begin{equation}\label{6.17eq1xx}
\mbox{(w)}\mbox{-}\lim_{k\to\infty}Q^{n_k^{(6)},\l_{m}}\xi
= Q^{\l_m}_{t_{{n_0}-1}} \xi \q\mbox{ in }
L^{2}_\dbF(t_{{n_0}-1},t_{n_0};L^{\frac{4}{3}}(\O;H)),\qq
\forall\; \xi\in L^4_{\cF_{t_{i}}}(\O;D(A)).
\end{equation}

Generally, for any $i=1,2,\cds,n_0$, there
exist a bounded, pointwisely defined linear
operator $Q^{\l_m}_{t_i}$ from
$L^4_{\cF_{t_{i}}}(\O;D(A))$ to
$L^{2}_\dbF(t_{i},t_{i+1};L^{\frac{4}{3}}(\O;H))$,
and a subsequence
$\{n_k^{(5+n_0-i)}\}_{k=1}^\infty$ of
$\{n_k^{(4+n_0-i)}\}_{k=1}^\infty$ such that
\begin{equation}\label{6.17eq1xxx}
\mbox{(w)}\mbox{-}\lim_{k\to\infty}Q^{n_k^{(5+n_0-i)},\l_{m}}\xi
= Q^{\l_m}_{t_i} \xi \q\mbox{ in }
L^{2}_\dbF(t_{i},t_{i+1};L^{\frac{4}{3}}(\O;H)),\qq
\forall\; \xi\in L^4_{\cF_{t_i}}(\O;D(A)).
\end{equation}
Since $ Q^{\l_m}_{t_i}$ is pointwisely defined,
for $\ae (t,\o)\in (t_{i},t_{i+1})\times\O$,
there is a $ q^{\l_m}_{t_i}(t,\o)\in
\cL(D(A),\;H)$ such that
 $$
 \big(Q^{\l_m}_{t_i} \xi\big)(t,\o)= q^{\l_m}_{t_i}(t,\o)\xi(\o),\qq \forall\; \xi\in
L^4_{\cF_{t_i}}(\O;D(A)).
 $$

For each $i=1,2,\cds,n_0$, write
 $$
\cH_i=\Big\{\sum_{j=1}^\ell
\chi_{O_j\cap([t_i,T]\times\O)}(\cd)h_j\;\Big|\;
\ell\in\dbN,\, O_j\in\cM, \, h_j\in
L^4_{\cF_{t_i}}(\O;D(A))\Big\}.
 $$
It is clear that $\cH_i$ is dense in
$L^{2}_\dbF(t_i,T;L^{4}(\O;H))$ and
$\cH\subset\cH_1$. Define an operator
$Q^{i,\l_m}$ from $\cH_i$ to
$L^2_\dbF(t_i,T;L^{\frac{4}{3}}(\O;H))$ as
follows: For any $\ds v=\sum_{j=1}^\ell
\chi_{O_j\cap([t_i,T]\times\O)}(\cd)h_j\in\cH_i$
with $\ell\in\dbN$, $O_j\in\cM$ and $h_j\in
L^4_{\cF_{t_i}}(\O;D(A))$,
$$
\big(Q^{i,\l_m}v\big)(t,\o)=\sum_{\g=i}^{n_0}\sum_{j=1}^\ell
\chi_{[t_\g,t_{\g+1})}(t)\chi_{O_j\cap([t_i,T]\times\O)}(t,\o)q^{\l_m}_{t_\g}
(t,\o)h_j, \ \ae (t,\o)\in (0,T)\times \O.
$$
It is easy to check that $ Q^{i,\l_m} v\in
L^2_\dbF(t_i,T;L^{\frac{4}{3}}(\O;H))$, and
$Q^{i,\l_m}$ is a pointwisely defined linear
operator from $\cH_i$ to
$L^2_\dbF(t_i,T;L^{\frac{4}{3}}(\O;H))$. Also,
for the above $ v$, we have
$$
Q^{{n_{k}^{(n_0+4)}},\l_m}(s)v^{{n_{k}^{(n_0+4)}}}(s)
=\sum_{j=1}^\ell
\chi_{O_j\cap([t_i,T]\times\O)}Q^{{n_{k}^{(n_0+4)}},\l_m}(s)\G_{n_k^{(n_0+4)}}h_j.
$$
Hence,
$$
\begin{array}{ll}\ds
\q
Q^{{n_{k}^{(n_0+4)}},\l_m}(\cd)v^{{n_{k}^{(n_0+4)}}}(\cd)-
\big(Q^{i,\l_m} v\big)(\cd)\\
\ns\ds = \sum_{j=1}^\ell
\chi_{O_j\cap([t_i,T]\times\O)}(\cd)\[
Q^{{n_{k}^{(n_0+4)}},\l_m}(\cd)\G_{n_k^{(n_0+4)}}h_j-\big(
Q^{i,\l_m} h_j\big)(\cd)\].
\end{array}
$$
This gives that
\begin{equation}\label{s6eq1}
\mbox{(w)-}\lim_{k\to\infty}Q^{{n_{k}^{(n_0+4)}},\l_m}(\cd)v^{{n_{k}^{(n_0+4)}}}(\cd)=
Q^{i,\l_m} v \q\mbox{ in }
L^{2}_\dbF(t_i,T;L^{\frac{4}{3}}(\O;H)),\q
\forall\; v\in \cH_i.
\end{equation}

Similarly, one can find a subsequence
$\{n_{k}^{(n_0+5)}\}_{k=1}^\infty\subset\{n_{k}^{(n_0+4)}\}_{k=1}^\infty$
and a pointwisely defined linear operator
$\widehat Q^{i,\l_m}$ from $\cH_i$ to
$L^2_\dbF(t_i,T;L^{\frac{4}{3}}(\O;H))$ such
that
\begin{equation}\label{s6eq1x}
\mbox{(w)-}\lim_{k\to\infty}Q^{{n_{k}^{(n_0+5)}},\l_{m}}(\cd)^*v^{{n_{k}^{(n_0+5)}}}(\cd)=
\widehat Q^{i,\l_m} v \q\mbox{ in }
L^{2}_\dbF(t_i,T;L^{\frac{4}{3}}(\O;H)),\q
\forall\; v\in \cH_i.
\end{equation}
For any $v_1,v_2\in\cH_i$, we have
\begin{equation}\label{10.8eq1}
\begin{array}{ll}\ds
\mE\int_{t_i}^T\big\langle Q^{i,\l_{m}}v_1,v_2
\big\rangle_H
dt\3n&\ds=\lim_{k\to\infty}\mE\int_{t_i}^T\big\langle
Q^{{n_{k}^{(n_0+5)}},\l_{m}}(t)v_1^{{n_{k}^{(n_0+5)}}}(t),v_2^{{n_{k}^{(n_0+5)}}}(t)
\big\rangle_H dt \\
\ns&\ds =
\lim_{k\to\infty}\mE\int_{t_i}^T\big\langle
v_1^{{n_{k}^{(n_0+5)}}}(t),Q^{{n_{k}^{(n_0+5)}},\l_{m}}(t)^*v_2^{{n_{k}^{(n_0+5)}}}(t)
\big\rangle_H dt\\
\ns&\ds = \mE\int_{t_i}^T\big\langle
v_1,\widehat Q^{i,\l_{m}}v_2 \big\rangle_H dt.
\end{array}
\end{equation}

For any $\xi_1,\xi_2\in L^4_{\cF_{t_i}}(\O;H)$,
$u_1(\cd), u_2(\cd)\in
L^2_{\dbF}(t_i,T;L^4(\O;H))$ and
$v_1(\cd),v_2(\cd)\in \cH_i$, by
\eqref{s6eq1}--\eqref{s6eq1x}, it is easy to
see that
\begin{equation}\label{6.17eq2}
\begin{array}{ll}\ds
\lim_{k\to\infty}  \mE \int_{t_i}^T
\[\big\langle Q^{{n_{k}^{(n_0+5)}},\l_{m}}(s)v_1^{{n_{k}^{(n_0+5)}}}(s),
x_2^{{n_{k}^{(n_0+5)}},\l_{m}}(s)
\big\rangle_{\dbR^{n_{k}^{(n_0+5)}}}\\\ns\ds\qq\qq\qq+\big\langle
Q^{{n_{k}^{(n_0+5)}},\l_{m}}(s)x_1^{{n_{k}^{(n_0+5)}},\l_{m}}(s),
v_2^{{n_{k}^{(n_0+5)}}}(s)
\big\rangle_{\dbR^{n_{k}^{(n_0+5)}}}\]ds \\
\ns\ds  =  \mE \int_{t_i}^T
 \[\big\langle
\big(Q^{i,\l_m} v_1\big)(s), x_2 (s)
\big\rangle_{H}+\big\langle x_1 (s),
\big(\widehat Q^{i,\l_m} v_2\big)(s)
\big\rangle_{H}\]ds.
\end{array}
\end{equation}
Therefore,
\begin{equation}\label{6.17eq4}
\begin{array}{ll}
\ds \q \mE \int_{t_i}^T \big\langle v_1(s),
\widehat Q^{(t_i,\l_m)}(\xi_2,u_2,v_2)(s)
\big\rangle_{H}ds + \mE \int_{t_i}^T
\big\langle
 Q^{(t_i,\l_m)}(\xi_1,u_1,v_1) (s),
v_2(s) \big\rangle_{H}ds \\
\ns\ds =\mE \int_{t_i}^T
 \[\big\langle
\big(Q^{i,\l_m} v_1\big)(s), x_2 (s)
\big\rangle_{H}+\big\langle x_1 (s),
\big(\widehat Q^{i,\l_m} v_2\big)(s)
\big\rangle_{H}\]ds.
\end{array}
\end{equation}
By \eqref{eq5} and \eqref{6.17eq4}, we find
that
\begin{equation}\label{6.17eq3}
\begin{array}{ll}
\ds \q\mE\big\langle P_T x_1^{\l_m}(T),
x_2^{\l_m}(T) \big\rangle_{H} - \mE
\int_{t_i}^T \big\langle
F(s) x_1^{\l_m}(s), x_2^{\l_m}(s) \big\rangle_{H}ds\\
\ns\ds =\mE\big\langle P^{\l_m}(0)\xi_1,\xi_2
\big\rangle_{H} + \mE \int_{t_i}^T \big\langle
P^{\l_m}(s)u_1(s),
x_2^{\l_m}(s)\big\rangle_{H}ds  + \mE
\int_{t_i}^T \big\langle
P^{\l_m}(s)x_1^{\l_m}(s),
u_2(s)\big\rangle_{H}ds\\
\ns\ds \q   + \mE \int_{t_i}^T \big\langle
P^{\l_m}(s)K(s)x_1^{\l_m}(s),
v_2(s)\big\rangle_{H}ds + \mE \int_{t_i}^T
\big\langle P^{\l_m}(s)v_1(s),
K(s)x_2^{\l_m}(s)+v_2(s)
\big\rangle_{H}ds  \\
\ns\ds \q   + \mE \int_{t_i}^T
 \[\big\langle
\big(Q^{i,\l_m} v_1\big)(s), x_2^{\l_m}(s)
\big\rangle_{H}+\big\langle x_1^{\l_m} (s),
\big(\widehat Q^{i,\l_m} v_2\big)(s)
\big\rangle_{H}\]ds,
\end{array}
\end{equation}
holds for any $\xi_1,\xi_2\in
L^4_{\cF_{t_i}}(\O;H)$, $u_1(\cd), u_2(\cd)\in
L^2_{\dbF}(t_i,T;L^4(\O;H))$ and
$v_1(\cd),v_2(\cd)\in \cH_i$
($i=1,2,\cds,n_0$).

\ms

{\bf Step 2}. We now take $m\to\infty$ in
\eqref{6.17eq3}. The argument below is very
similar to Step 1. Choose $\xi_2 =0$ and $u_2
=0$ in the equation \eqref{op-fsystem3.1}. By
\eqref{zx-s1} and similar to \eqref{op-fsystem5
enest x}, there is a constant $C_2
>0$, independent of $t$ and $m$, such that
$|x_2^{\l_m}|_{L^\infty(t,T;L^4(\O;H))}\leq
C_2|v_2|_{L^2(t,T;L^4(\O;H))}$. Without loss of
generality, we may assume that
\begin{equation}\label{10.8eq2}
\frac{1}{\ds\max_{1\leq i\leq
n_0}(t_{i+1}-t_i)}>2^{12}|C_2
|^4|K|^4_{L^4_\dbF(0,T;L^\infty(\O;\cL(D(A))))}.
\end{equation}
Letting $\xi_1\in
L^4_{\cF_{t_{n_0}}}(\O;D(A))$, $u_1 = -(A^{
\l_m}+J )\xi_1$ and $v_1 =-K \xi_1$ in
\eqref{op-fsystem2.1}, and letting $\xi_2 =0$
and $u_2 =0$ in \eqref{op-fsystem3.1}, by
\eqref{6.17eq3} and \eqref{10.8eq1}, we find
that
\begin{equation}\label{10.8eq2s}
\begin{array}{ll}
\ds \q\mE\big\langle P_T \xi_1, x_2^{\l_m}(T)
\big\rangle_{H} - \mE \int_{t_{n_0}}^T
\big\langle
F(s) \xi_1, x_2^{\l_m}(s) \big\rangle_{H}ds\\
\ns\ds = \mE \int_{t_{n_0}}^T \big\langle
P^{\l_m}(s)u_1(s),
x_2^{\l_m}(s)\big\rangle_{H}ds  + \mE
\int_{t_{n_0}}^T \big\langle
P^{\l_m}(s)K(s)\xi_1,
v_2(s)\big\rangle_{H}ds \\
\ns\ds \q   + \mE \int_{t_{n_0}}^T \big\langle
P^{\l_m}(s)v_1(s),
K(s)x_2^{\l_m} + v_2(s) \big\rangle_{H}ds  \\
\ns\ds \q  - \mE \int_{t_{n_0}}^T \big\langle
Q^{n_0,\l_m} K(s) \xi_1, x_2^{\l_m}(s)
\big\rangle_{H}ds + \mE \int_{t_{n_0}}^T
\big\langle Q^{n_0,\l_m} \xi_1, v_2(s)
\big\rangle_{H}ds.
\end{array}
\end{equation}
Similar to  \eqref{pp2}, by (\ref{10.8eq2s})
and noting \eqref{10.8eq2}, we have the
following estimate:
 \begin{equation}\label{10.8eq5}
\begin{array}{ll}\ds
\q\big|\!\big|Q^{n_0,\l_m}(\cd)\big|\!\big|_{\cL(L^4_{\cF_{t_{n_0}}}(\O;D(A)),\;
L^{2}_\dbF({t_{n_0}},T;L^{\frac{4}{3}}(\O;H)))}\\
\ns\ds \leq\!
C\big(|P_T|_{L^4_{\cF_T}(\O;\cL(H))} +
|F|_{L^1_\dbF(0,T;L^2(\O;\cL(H)))}\big)\big(1 +
|A|_{\cL(D(A),\;H)} +
|(J,K)|_{(L^4_\dbF(0,T;L^\infty(\O;\cL(H))))^2}\big).
\end{array}
\end{equation}
By \eqref{10.8eq5} and Corollary \ref{cor1},
there exist a bounded, pointwisely defined
linear operator $Q_{t_{n_0}}$ from
$L^4_{\cF_{t_{n_0}}}(\O;D(A))$ to
$L^{2}_\dbF(t_{n_0},T;L^{\frac{4}{3}}(\O;H))$,
and a subsequence
$\{\l_{m_j}^{(2)}\}_{j=1}^\infty$ of
$\{\l_{m_j}^{(1)}\}_{j=1}^\infty$ (defined at
the beginning of the Step 6 in the proof of
Theorem 6.1) such that
 $$
\mbox{(w)}\mbox{-}\lim_{j\to\infty}Q^{n_0,\l_{m_j}^{(2)}}\xi
= Q_{t_{n_0}}\xi \q\mbox{ in }
L^{2}_\dbF(t_{n_0},T;L^{\frac{4}{3}}(\O;H)),\qq
\forall\; \xi\in L^4_{\cF_{t_{n_0}}}(\O;D(A)).
 $$

Generally, for any $i=1,2,\cds,n_0$, there
exist a bounded, pointwisely defined linear
operator $Q_{t_i}$ from
$L^4_{\cF_{t_{i}}}(\O;D(A))$ to
$L^{2}_\dbF(t_{i},t_{i+1};L^{\frac{4}{3}}(\O;H))$,
and a subsequence
$\{\l_m^{(n_0-i+2)}\}_{m=1}^\infty$ of
$\{\l_m^{(n_0-i+1)}\}_{m=1}^\infty$ such that
 $$
\mbox{(w)}\mbox{-}\lim_{j\to\infty}Q^{i,\l_{m_j}^{(n_0-i+2)}}\xi
= Q_{t_i} \xi \q\mbox{ in }
L^{2}_\dbF(t_{i},t_{i+1};L^{\frac{4}{3}}(\O;H)),\qq
\forall\; \xi\in L^4_{\cF_{t_i}}(\O;D(A)).
 $$
Since $ Q_{t_i}$ is pointwisely defined, for
$\ae (t,\o)\in (t_{i},t_{i+1})\times\O$, there
is a $ q_{t_i}(t,\o)\in \cL(D(A),\;H)$ such
that
 $$
 \big(Q_{t_i} \xi\big)(t,\o)= q_{t_i}(t,\o)\xi(\o),\qq \forall\; \xi\in
L^4_{\cF_{t_i}}(\O;D(A)).
 $$
Define a linear operator $Q$ from $\cH$ to
$L^2_\dbF(0,T;L^{\frac{4}{3}}(\O;H))$ as
follows: For any $\ds v=\sum_{j=1}^\ell
\chi_{O_i}h_i\in\cH$ with $\ell\in\dbN$,
$O_i\in\cM$ and $h_i\in D(A)$,
$$
\big(Qv\big)(t,\o)=\sum_{\g=1}^{n_0}\sum_{j=1}^\ell
\chi_{[t_\g,t_{\g+1})}(t)\chi_{O_j}(t,\o)q_{t_\g}
(t,\o)h_j, \ \ \ae (t,\o)\in (0,T)\times \O.
$$
Then,
\begin{equation}\label{1-eq14}
\mbox{(w)-}\lim_{j\to\infty}Q^{1,\l^{(n_0+1)}_{m_j}}
v = Q v \q\mbox{ in }
L^{2}_\dbF(0,T;L^{\frac{4}{3}}(\O;H)),\q
\forall\; v\in \cH.
\end{equation}
By a similar argument, we see that there exist
a linear operator $\widehat Q$ from $\cH$ to
$L^{2}_\dbF(0,T;L^{\frac{4}{3}}(\O;H))$, and a
subsequence
$\{\l^{(n_0+2)}_{m_j}\}_{j=1}^\infty$ of
$\{\l^{(n_0+1)}_{m_j}\}_{j=1}^\infty$  such
that
\begin{equation}\label{10.8eq14}
\mbox{(w)-}\lim_{j\to\infty} \widehat
Q^{1,\l^{(n_0+2)}_{m_j}} v = \widehat Q v
\q\mbox{ in }
L^{2}_\dbF(0,T;L^{\frac{4}{3}}(\O;H)),\q
\forall\; v\in \cH.
\end{equation}

Now,  choosing arbitrarily $\xi_1,\xi_2\in H$,
$u_1(\cd), u_2(\cd)\in
L^2_{\dbF}(0,T;L^4(\O;H))$ and
$v_1(\cd),v_2(\cd)\in \cH$, by \eqref{6.17eq2}
and \eqref{1-eq14}--\eqref{10.8eq14}, we find
that
\begin{equation}\label{10.9eq1}
\begin{array}{ll}\ds
\lim_{j\to\infty}\lim_{k\to\infty}  \mE
\int_{0}^T
\[\big\langle Q^{{n_{k}^{(2n_0+4)}},\l^{(n_0+2)}_{m_j}}(s)v_1^{{n_{k}^{(2n_0+4)}}}(s),
x_2^{{n_{k}^{(2n_0+4)}},\l^{(n_0+2)}_{m_j}}(s)
\big\rangle_{\dbR^{n_{k}^{(2n_0+4)}}}\\\ns\ds\qq\qq\qq+\big\langle
Q^{{n_{k}^{(2n_0+4)}},\l^{(n_0+2)}_{m_j}}(s)x_1^{{n_{k}^{(2n_0+4)}},\l^{(n_0+2)}_{m_j}}(s),
v_2^{{n_{k}^{(2n_0+4)}}}(s)
\big\rangle_{\dbR^{n_{k}^{(2n_0+4)}}}\]ds \\
\ns\ds  =  \mE \int_{0}^T
 \[\big\langle
\big(Q v_1\big)(s), x_2 (s)
\big\rangle_{H}+\big\langle x_1 (s),
\big(\widehat Q v_2\big)(s)
\big\rangle_{H}\]ds.
\end{array}
\end{equation}
Combining  \eqref{eqbilinear1}, \eqref{6.6eq6},
\eqref{6.6eq11}, \eqref{3Q}, \eqref{6.18eq3}
and \eqref{10.9eq1}, we conclude that the
desired equality \eqref{10.9eq2} holds for any
$\xi_1,\xi_2\in L^4_{\cF_0}(\O;H)$, $u_1(\cd),
u_2(\cd)\in L^2_{\dbF}(0,T;L^4(\O;H))$ and
$v_1(\cd),v_2(\cd)\in \cH$. This completes the
proof of Theorem \ref{10.1th}.
\endpf

\br
1) We conjecture that the same conclusion in
Theorem \ref{10.1th} still holds for any $K\in
L^4_\dbF(0,T;L^\infty(\O;\cL(H)))$, or at least
for any $\ds K=\sum_{i=1}^{n_0}
\chi_{[t_i,t_{i+1})}(t)K_i$ with $n_0\in\dbN$,
$0=t_1<t_2<\cdots<t_{n_0}<t_{n_0+1}=T$, and
$K_i\in L^\infty_{\cF_{t_i}}(\O;\cL(H))$ (If
the later is true, then we may drop the
assumption  $b_x(\cd,\bar x(\cd),\bar
u(\cd))\in
L^4_\dbF(0,T;L^\infty(\O;\cL(D(A))))$ in
Theorem \ref{s7th max}). However, we cannot
prove it at this moment.

\ms

2) In some sense, the operators $Q$ and
$\widehat Q$ given in Theorem \ref{10.1th} play
similar roles as the operators $Q$ and $Q^*$,
where the later operator $Q$ is the second
component of the transposition solution
$(P(\cd),Q(\cd))$ to \eqref{op-bsystem3}.
\er

\section{Necessary condition for optimal controls, the case of convex
control domains} \label{s6}

For the sake of completeness, in this section,
we shall give a necessary condition for optimal
controls of the system \eqref{fsystem1} for the
case of special control domain $U$, i.e., $U$
is a convex subset of another separable Hilbert
space $H_1$, and the metric of $U$ is
introduced by the norm of $H_1$ (i.e.,
$d(u_1,u_2)=|u_1-u_2|_{H_1}$).

To begin with, we introduce the following
further assumptions for $a(\cd,\cd,\cd)$,
$b(\cd,\cd,\cd)$, $g(\cd,\cd,\cd)$ and
$h(\cd)$.

\ms

\no{\bf (A3)} {\it The maps $a(t,x,u)$ and
$b(t,x,u)$, and the functional $g(t,x,u)$ and
$h(x)$ are $C^1$ with respect to $x$ and $u$.
Moreover, there exists a constant $C_L>0$ such
that, for any $(t,x,u)\in [0,T]\times H\times
U$,
\begin{equation}\label{ab}
\left\{
\begin{array}{ll}\ds
 |\!|a_x(t,x,u)|\!|_{\cL(H)}+|\!|b_x(t,x,u)|\!|_{\cL(H)} + |g_x(t,x,u)|_H+|h_x(x)|_{H}\leq C_L,\\
\ns\ds |\!|a_u(t,x,u)|\!|_{\cL(H_1,H)}+
|\!|b_u(t,x,u)|\!|_{\cL(H_1,H)} +
|g_u(t,x,u)|_{H_1} \leq C_L.
\end{array}
\right.
\end{equation}}

\medskip

Our result in this section is as follows.

\begin{theorem}\label{th max}
Assume that $x_0\in L^2_{\cF_0}(\O;H)$. Let the
assumptions (A1), (A2) and (A3) hold, and let
$(\bar x(\cd),\bar u(\cd))$ be an optimal pair
of Problem (P). Let $(y(\cdot),Y(\cdot))$ be
the transposition solution of the equation
\eqref{bsystem1} with $p=2$, and $y_T$ and
$f(\cd,\cd,\cd)$ given by
 \bel{zv1}
 \left\{
 \ba{ll}
\ds y_T = -h_x\big(\bar x(T)\big),\\\ns
 \ds f(t,y_1,y_2)=-a_x(t,\bar x(t),\bar
u(t))^*y_1 - b_x\big(t,\bar x(t),\bar
u(t)\big)^*y_2 + g_x\big(t,\bar x(t),\bar
u(t)\big).
 \ea\right.
 \ee
Then,
\begin{equation}\label{maxth ine1}
\begin{array}{ll}\ds
\Re\big\langle a_u(t,\bar x(t),\bar u(t))^*
y(t) + b_u(t,\bar x(t),\bar u(t))^*Y(t) -
g_u(t,\bar u(t),\bar x(t)), u -\bar u(t)
\big\rangle_{H_1} \leq 0,\\
\ns\ds \hspace{9.8cm}\;\,\ae [0,T]\times \O,
\forall\; u \in U.
\end{array}
\end{equation}
\end{theorem}

\medskip

{\it Proof}\,: We use the convex perturbation
technique and divide the proof into several
steps.

\ms

{\bf Step 1}. For the optimal pair $(\bar
x(\cdot),\bar u(\cdot))$, we fix arbitrarily a
control $u(\cdot)\in \cU[0,T]$ satisfying $
u(\cd)-\bar u(\cd)\in
L^2_\dbF(0,T;L^2(\O;H_1))$. Since $U$ is
convex, we see that
$$
u^\e(\cdot) = \bar u(\cdot) + \e [u(\cdot) -
\bar u(\cdot)] = (1-\e)\bar u(\cdot) + \e
u(\cdot) \in \cU[0,T], \q\forall\;\e \in [0,1].
$$
Denote by $x^\e(\cdot)$ the state process of
\eqref{fsystem1} corresponding to the control
$u^\e(\cdot)$. By Lemma \ref{well lemma s1}, it
follows that
\begin{equation}\label{th max eq0.0xz}
|x^\e|_{C_\dbF(0,T;L^2(\O;H))}\leq C
\big(1+|x_0|_{L^2_{\cF_0}(\O;H)}\big),\q
\forall\; \e \in [0,1].
\end{equation}
Write $\ds x_1^\e(\cd) =
\frac{1}{\e}\big[x^\e(\cd)-\bar x(\cd)\big]$
and $\d u(\cd) = u(\cd) - \bar u(\cd)$. Since
$(\bar x(\cd),\bar u(\cd))$ satisfies
\eqref{fsystem1}, it is easy to see that
$x_1^\e(\cdot)$ satisfies the following
stochastic differential equation:
\begin{equation}\label{fsystem3x}
\left\{
\begin{array}{lll}\ds
dx_1^\e = \big(Ax_1^\e +  a_1^\e x^\e_1 +
a_2^\e\d u  \big)dt + \big( b_1^\e x^\e_1 +
b_2^\e\d u \big)dw(t) &\mbox{ in
}(0,T],\\
\ns\ds x_1^\e(0)=0,
\end{array}
\right.
\end{equation}
where
\begin{equation}\label{tatb}
\left\{
\begin{array}{ll}
\ds   a_1^\e (t)  = \int_0^1 a_x(t,\bar x(t) +
\si\e x_1^\e(t), u^\e(t))d\si,\qq  a_2^\e (t)
= \int_0^1 a_u(t,\bar x(t), \bar u(t)+\si\e\d
u(t))d\si,\\
\ns\ds  b_1^\e(t)  = \int_0^1 b_x(t,\bar x(t) +
\si\e x_1^\e(t), u^\e(t))d\si,\qq  b_2^\e(t)  =
\int_0^1 b_u(t,\bar x(t), \bar u(t)+\si\e\d
u(t))d\si.
\end{array}
\right.
\end{equation}
Consider the following stochastic differential
equation:
\begin{equation}\label{fsystem3.1}
\left\{
\begin{array}{lll}\ds
dx_2 = \big[Ax_2 + a_1(t)x_2 +  a_2(t)\d u
\big]dt + \big[ b_1(t) x_2 + b_2(t)\d u \big]
dw(t) &\mbox{ in
}(0,T],\\
\ns\ds x_2(0)=0,
\end{array}
\right.
\end{equation}
where
\begin{equation}\label{barab}
\left\{
\begin{array}{ll}\ds
a_1(t) = a_x(t,\bar x(t),\bar u(t)),\q a_2(t) =
a_u(t,\bar
x(t),\bar u(t)),\\
\ns\ds  b_1(t) = b_x(t,\bar x(t),\bar u(t)),\q
b_2(t) = b_u(t,\bar x(t),\bar u(t)).
\end{array}
\right.
\end{equation}

\ms

{\bf Step 2}. In this step, we shall show that
\begin{equation}\label{th max eq0.2}
\lim_{\e\to 0+} \Big|x_1^\e -
x_2\Big|_{L_\dbF^\infty(0,T;L^2(\O;H))}=0.
\end{equation}

First, using Lemma \ref{BDG} and by the
assumption (A1), we find that
\begin{equation}\label{th max eq0}
\begin{array}{ll}\ds
\mE|x_1^\e(t)|^2_H  = \mE\Big| \int_0^t S(t-s)
a_1^\e(s) x_1^\e(s) ds + \int_0^t S(t-s)
a_2^\e(s)\d u(s) ds  \\
\ns\ds \hspace{2.4cm} +
\int_0^t S(t-s) b_1^\e(s) x_1^\e(s) dw(s) +  \int_0^t S(t-s)b_2^\e(s)\d u(s) dw(s)\Big|_H^2\\
\ns\ds \hspace{1.7cm}\leq C \mE\[\Big| \int_0^t
S(t-s) a_1^\e(s) x_1^\e(s) ds \Big|_H^2
+  \Big|\int_0^t S(t-s) b_1^\e(s) x_1^\e(s) dw(s) \Big|_H^2 \\
\ns\ds \hspace{2.4cm}  + \Big| \int_0^t S(t-s)a_2^\e(s)\d u(s)ds\Big|_H^2 + \Big| \int_0^t S(t-s)b_2^\e(s)\d u(s)dw(s)\Big|_H^2\] \\
\ns\ds \hspace{1.7cm} \leq C\[ \int_0^t
\mE|x_1^\e(s)|_H^2ds + \int_0^T\mE |\d
u(s)|_{H_1}^2 dt\].
\end{array}
\end{equation}
It follows from \eqref{th max eq0} and
Gronwall's inequality that
\begin{equation}\label{th max eq0.1}
\begin{array}{ll}\ds
\mE|x_1^\e(t)|^2_H \leq  C|\bar u -
u|^2_{L^2_\dbF(0,T;L^2(\O;H_1))},\q\forall\;t\in
[0,T].
\end{array}
\end{equation}
By a similar computation, we see that
\begin{equation}\label{th max eq2}
\begin{array}{ll}\ds
\mE|x_2(t)|^2_H \leq  C |\bar u -
u|^2_{L^2_\dbF(0,T;L^2(\O;H_1))},\q\forall\;t\in
[0,T].
\end{array}
\end{equation}

On the other hand, put $x_3^\e = x_1^\e - x_2$.
Then, $x_3^\e$ solves the following equation:
\begin{equation}\label{fsystem3x2}
\left\{
\begin{array}{lll}\ds
dx_3^\e = \big[Ax_3^\e +  a_1^\e(t) x^\e_3 +
\big(a_1^\e(t)- a_1(t)\big)x_2 +
\big( a_2^\e(t)-a_2(t)\big)\d u \big]dt \\
\ns\ds \hspace{1.3cm} + \big[  b_1^\e(t) x^\e_3
+ \big(b^\e_1(t) - b_1(t)\big) x_2 + \big(
b_2^\e(t)-b_2(t)\big)\d u \big] dw(t) &\mbox{
in
}(0,T],\\
\ns\ds x_3^\e(0)=0.
\end{array}
\right.
\end{equation}
It follows from \eqref{th max
eq2}--\eqref{fsystem3x2} that
\begin{equation}
\begin{array}{ll}\ds
\q\mE|x_3^\e(t)|^2_{H} \\
\ns\ds = \mE\Big| \int_0^t S(t-s)
a_1^\e(s)x_3^\e(s)ds
+ \int_0^t S(t-s) b_1^\e(s)x_3^\e(s) dw(s) \\
\ns\ds\q +\int_0^t S(t-s)\big[
a^\e_1(s)-a_1(s)\big] x_2(s) ds
+\int_0^t S(t-s)\big[ b_1^\e(s) - b_1(s)\big] x_2(s) dw(s) \nonumber\\
\ns\ds\q+\int_0^t S(t-s)\big[
a^\e_2(s)-a_2(s)\big]\d u(s) ds
+\int_0^t S(t-s)\big[ b_2^\e(s) - b_2(s)\big]\d u(s)dw(s)\Big|_H^2\\
\ns\ds  \leq C\[\mE\int_0^t |x_3^\e(s)|_H^2 ds \\
\ns\ds \q+
|x_2(\cd)|^2_{L^\infty_\dbF(0,T;L^2(\O;H))}
\int_0^T \mE\big( |\!|a^\e_1(s)-a_1(s)|\!|_{\cL(H)}^2  + |\!|b^\e_1(s)-b_1(s)|\!|_{\cL(H)}^2 \big) dt \\
\ns\ds \q+ |u-\bar
u|^2_{L^2_\dbF(0,T;L^2(\O;H_1))}\int_0^T
\mE\big(|\!|a_2^\e(s) -
a_2(s)|\!|_{\cL(H_1,H)}^2 + |\!|b_2^\e(s) -
b_2(s)|\!|_{\cL(H_1,H)}^2 \big)dt
\]\\
\ns\ds \leq C(1+|u-\bar
u|^2_{L^2_\dbF(0,T;L^2(\O;H_1))})\Big\{\mE\int_0^t
|x_3^\e(s)|_H^2
ds + \int_0^T \mE\[ |\!|a^\e_1(s)-a_1(s)|\!|_{\cL(H)}^2  \\
\ns\ds \q +
|\!|b^\e_1(s)-b_1(s)|\!|_{\cL(H)}^2+
|\!|a_2^\e(s) - a_2(s)|\!|_{\cL(H_1,H)}^2 +
|\!|b_2^\e(s) - b_2(s)|\!|_{\cL(H_1,H)}^2
\]dt\Big\}.
\end{array}
\end{equation}
This, together with Gronwall's inequality,
implies that
\begin{equation}\label{th max eq3}
\begin{array}{ll}\ds
\mE|x_3^\e(t)|^2_{H} \leq Ce^{C|u-\bar
u|_{L^2_\dbF(0,T;L^2(\O;H_1))}}\int_0^T \mE\[ |\!|a^\e_1(s)-a_1(s)|\!|_{\cL(H)}^2  + |\!|b^\e_1(s)-b_1(s)|\!|_{\cL(H)}^2\\
\ns\ds \qq\qq\q\;\, + |\!|a_2^\e(s) -
a_2(s)|\!|_{\cL(H_1,H)}^2 + |\!|b_2^\e(s) -
b_2(s)|\!|_{\cL(H_1,H)}^2 \]ds,\q \forall\,
t\in [0,T].
\end{array}
\end{equation}

Note that \eqref{th max eq0.1} implies
$x^\e(\cd)\to \bar x(\cd)$ (in $H$) in
probability, as $\e\to0$. Hence, by
\eqref{tatb}, \eqref{barab} and the continuity
of $a_x(t,\cd,\cd)$, $b_x(t,\cd,\cd)$,
$a_u(t,\cd,\cd)$ and $b_u(t,\cd,\cd)$, we
deduce that
$$
\ba{ll}\ds \lim_{\e\to 0}\int_0^T\mE\[
|\!|a^\e_1(s)-a_1(s)|\!|_{\cL(H)}^2 +
|\!|b^\e_1(s)-b_1(s)|\!|_{\cL(H)}^2
\\\ns\ds\qq\qq+ |\!|a_2^\e(s) - a_2(s)|\!|_{\cL(H_1,H)}^2 +|\!|b_2^\e(s) -
b_2(s)|\!|_{\cL(H_1,H)}^2 \]ds=0. \ea
$$
This, combined with (\ref{th max eq3}), gives
\eqref{th max eq0.2}.

\ms

{\bf Step 3}. Since $(\bar x(\cdot),\bar
u(\cdot))$ is an optimal pair of Problem (P),
from \eqref{th max eq0.2}, we find that
\begin{equation}\label{var 1}
\begin{array}{ll}\ds
 0\leq \lim_{\e\to 0}\frac{\cJ(u^\e(\cdot)) - \cJ(\bar u(\cdot))}{\e} \\
\ns\ds\;\;\;= \Re\left\{\dbE\int_0^T
\Big[\big\langle g_1(t),x_2(t)\big\rangle_H +
\big\langle g_2(t),\d u(t)
\big\rangle_{H_1}\Big] dt + \dbE\big\langle
h_x(\bar x(T)),x_2(T)\big\rangle_H\right\},
\end{array}
\end{equation}
where
$$
g_1(t) = g_x(t,\bar x(t),\bar u(t)),\q g_2(t) =
g_u(t,\bar x(t),\bar u(t)).
$$

Now, by the definition of the transposition
solution to \eqref{bsystem1} (with $y_T$ and
$f(\cd,\cd,\cd)$ given by (\ref{zv1})), it
follows that
\begin{equation}\label{max eq1}
\begin{array}{ll}\ds
\q-\dbE \big\langle h_x(\bar
x(T)),x_2(T)\big\rangle_H -
\dbE\int_0^T \big\langle g_1(t),x_2(t)\big\rangle_H dt \\
\ns\ds = \dbE \int_0^T \Big[\big\langle
a_2(t)\d u(t), y(t)\big\rangle_H + \big\langle
b_2(t)\d u(t), Y(t)\big\rangle_H\Big]dt.
\end{array}
\end{equation}
Combining \eqref{var 1} and \eqref{max eq1}, we
find
\begin{eqnarray}\label{max ine2}
\Re\dbE\int_0^T\big\langle a_2(t)^* y(t) +
b_2(t)^*Y(t) - g_2(t), u(t)-\bar u(t)
\big\rangle_{H_1}dt\leq 0
\end{eqnarray}
holds for any $u(\cdot)\in \cU[0,T]$ satisfying
$u(\cd)-\bar u(\cd)\in
L^2_\dbF(0,T;L^2(\O;H_1))$. Hence, by means of
Lemma \ref{lemma4}, we conclude that
\begin{eqnarray}\label{max ine3}
\Re\big\langle a_2(t)^* y(t) + b_2(t)^*Y(t) -
g_2(t), u - \bar u(t) \big\rangle_{H_1} \leq 0,
\,\ae [0,T]\times \O,\;\forall \,u \in U.
\end{eqnarray}
This completes the proof of Theorem \ref{th
max}.\endpf

\section{Necessary condition for optimal controls, the case of nonconvex
control domains} \label{s7}

In this section, we shall derive a necessary
condition for optimal controls of the system
\eqref{fsystem1} with a general nonconvex
control domain. For such a case, the convex
perturbation technique does not work any more.
We need to adopt the spike variation technique
to derive the desired necessary condition.

We need the further following conditions on
$a(\cd,\cd,\cd)$, $b(\cd,\cd,\cd)$,
$g(\cd,\cd,\cd)$ and $h(\cd)$:

\ms

\no{\bf (A4)} {\it The maps $a(t,x,u)$ and
$b(t,x,u)$, and the functional $g(t,x,u)$ and
$h(x)$ are $C^2$ with respect to $x$, and
$a_x(t,x,u)$, $b_x(t,x,u)$, $g_x(t,x,u)$,
$a_{xx}(t,x,u)$, $b_{xx}(t,x,u)$ and
$g_{xx}(t,x,u)$ are continuous with respect to
$u$. Moreover, there exists a constant $C_L>0$
such that
\begin{equation}\label{ab1}
\left\{
\begin{array}{ll}\ds
|\!|a_x(t,x,u)|\!|_{\cL(H)}+|\!|b_x(t,x,u)|\!|_{\cL(H)}+|g_x(t,x,u)
|_H+|h_x(x) |_H
\leq C_L,\\
\ns\ds |\!|a_{xx}(t,x,u)|\!|_{\cL(H\times
H,\;H)}+|\!|b_{xx}(t,x,u)|\!|_{\cL(H\times
H,\;H)}+|\!|g_{xx}(t,x,u)|\!|_{\cL(H)}+|\!|h_{xx}(x)|\!|_{\cL(H)}
 \leq C_L,\\
\ns\ds \hspace{10.8cm} \forall\; (t,x,u)\in
[0,T]\times H\times U.
\end{array}
\right.
\end{equation}}

\ms

Let
\begin{equation}\label{H}
\begin{array}{ll}\ds
\dbH(t,x,u,k_1,k_2) \= \big\langle k_1,a(t,x,u)  \big\rangle_H + \big\langle k_2, b(t,x,u)  \big\rangle_H - g(t,x,u),\\
\ns\ds \hspace{4cm} (t,x,u,k_1,k_2)\in
[0,T]\times H \times U\times H\times H.
\end{array}
\end{equation}

We have the following result.

\begin{theorem}\label{s7th max}
Suppose that $H$ is a separable Hilbert space,
$L^p_{\cF_T}(\O;\dbC)$ ($1\leq p < \infty$) is
a separable Banach space, $U$ is a separable
metric space, and $x_0\in L^8_{\cF_0}(\O;H)$.
Let the assumptions (A1), (A2) and (A4) hold,
and let $(\bar x(\cd),\bar u(\cd))$ be an
optimal pair of Problem (P). Let
$(y(\cdot),Y(\cdot))$ be the transposition
solution to \eqref{bsystem1} with $p=2$, and
$y_T$ and $f(\cd,\cd,\cd)$ given by
(\ref{zv1}). Assume that $b_x(\cd,\bar
x(\cd),\bar u(\cd))\in
L^4_\dbF(0,T;L^\infty(\O;\cL(D(A))))$, and
$(P(\cd),Q^{(\cd)},\widehat Q^{(\cd)})$ is the
relaxed transposition solution to the equation
\eqref{op-bsystem3} in which $P_T$, $J(\cd)$,
$K(\cd)$ and $F(\cd)$ are given by
 \bel{zv2}
 \left\{
 \ba{ll}
\ds P_T = -h_{xx}\big(\bar x(T)\big),\\\ns
 \ds J(t) = a_x(t,\bar x(t),\bar
u(t)), \\
\ns \ds K(t) =b_x(t,\bar x(t),\bar u(t)), \\
\ns \ds F(t)= -\dbH_{xx}\big(t,\bar x(t),\bar
u(t),y(t),Y(t)\big).
 \ea\right.
 \ee
Then,
\begin{equation}\label{s7maxth ine1}
\begin{array}{ll}\ds
\Re\dbH\big(t,\bar x(t),\bar u(t),y(t),Y(t)\big) - \Re\dbH\big(t,\bar x(t),u,y(t),Y(t)\big) \\
\ns\ds  - \frac{1}{2}\Big\langle  P(t)\Big[
b\big(t,\bar x(t),\bar u(t)\big)-b\big(t,\bar
x(t),u\big)
\Big], b\big(t,\bar x(t),\bar u(t)\big)-b\big(t,\bar x(t),u\big) \Big\rangle_H \geq 0,\\
\ns\ds \hspace{9cm} \,\ae [0,T]\times
\O,\forall \,u \in U.
\end{array}
\end{equation}
\end{theorem}

{\it Proof:} We divide the proof into several
steps.

\ms

{\bf Step 1}. For each $\e>0$, let $E_\e\subset
[0,T]$ be a measurable set with measure $\e$.
Put
\begin{equation}\label{s7ue}
u^\e(\cdot) = \left\{
\begin{array}{ll}
\ds \bar u(t), & t\in [0,T]\setminus E_\e,\\
\ns\ds u(t), & t\in E_\e.
\end{array}
\right.
\end{equation}
where  $u(\cdot)$ is an arbitrary given element
in $\cU[0,T]$.

We introduce some notations which will be used
in the sequel.
\begin{equation}\label{s7tatb1}
\left\{
\begin{array}{ll}
\ds a_1(t) = a_x(t,\bar x(t),\bar u(t)), \q
b_1(t) = b_x(t,\bar
x(t),\bar u(t)),\q g_1(t) = g_x(t,\bar x(t),\bar u(t)),\\
\ns\ds a_{11}(t) = a_{xx}(t,\bar x(t),\bar
u(t)), \q b_{11}(t) = b_{xx}(t,\bar x(t),\bar
u(t)), \q g_{11}(t)
= g_{xx}(t,\bar x(t),\bar u(t)),\\
\ns\ds \tilde a_1^\e (t)  = \int_0^1 a_x\big(t,\bar x(t) + \si (x^\e (t) - \bar x(t)), u^\e (t)\big)d\si,\\
\ns\ds  \tilde b_1^\e(t)  = \int_0^1 b_x\big(t,\bar x(t) + \si (x^\e(t) - \bar x(t)), u^\e(t)\big)d\si,\\
\ns\ds \tilde a_{11}^\e (t)  = 2\int_0^1 (1-\si) a_{xx}\big(t,\bar x(t) + \si (x^\e(t) - \bar x(t)), u^\e (t)\big)d\si,\\
\ns\ds  \tilde b_{11}^\e (t)  = 2\int_0^1
(1-\si) b_{xx}\big(t,\bar x(t) + \si (x^\e(t) -
\bar x(t)), u^\e (t)\big)d\si,
\end{array}
\right.
\end{equation}
and
\begin{equation}\label{s7tatb2}
\left\{
\begin{array}{ll}
\ds \d a(t)  = a(t,\bar x(t), u(t)) - a(t,\bar
x(t),\bar u(t)),
\\
\ns\ds\d b(t)  = b(t,\bar x(t), u(t)) -
b(t,\bar x(t),\bar u(t)),\\
\ns\ds\d g(t)  = g(t,\bar x(t), u(t)) -
g(t,\bar x(t),\bar u(t)),\\
\ns\ds \d a_1(t) = a_x(t,\bar x(t), u(t)) -
a_x(t,\bar x(t),\bar u(t)), \\
\ns\ds \d b_1(t) = b_x(t,\bar x(t), u(t)) -
b_x(t,\bar x(t),\bar
u(t)),\\
\ns\ds \d a_{11}(t) = a_{xx}(t,\bar x(t), u(t))
-
a_{xx}(t,\bar x(t),\bar u(t)), \\
\ns\ds \d b_{11}(t) = b_{xx}(t,\bar x(t), u(t))
- b_{xx}(t,\bar
x(t),\bar u(t)),\\
\ns\ds \d g_1(t) = g_x(t,\bar x(t), u(t)) -
g_x(t,\bar x(t),\bar u(t)), \\
\ns\ds  \d g_{11}(t) = g_{xx}(t,\bar x(t),
u(t)) - g_{xx}(t,\bar x(t),\bar u(t)).
\end{array}
\right.
\end{equation}

Let $x^\e(\cdot)$ be the state process of the
system \eqref{fsystem1} corresponding to the
control $u^\e(\cdot)$. Then, $x^\e(\cdot)$
solves
\begin{eqnarray}\label{s7fsystem2}
\left\{
\begin{array}{lll}\ds
d x^\e = \big[Ax^\e +a(t,x^\e,u^\e)\big]dt + b(t,x^\e,u^\e)dw(t) &\mbox{ in }(0,T],\\
\ns\ds x^\e(0)=x_0.
\end{array}
\right.
\end{eqnarray}
 By Lemma \ref{well lemma s1}, we know that
\begin{equation}\label{th max eq0.0zx}
|x^\e|_{C_\dbF([0,T];L^8(\O;H))}\leq C
\big(1+|x_0|_{L^8_{\cF_0}(\O;H)}\big),\q
\forall\; \e
>0.
\end{equation}
 Let
$x_1^\e(\cd) = x^\e(\cd)-\bar x(\cd)$. Then, by
(\ref{th max eq0.0zx}) and noting that the
optimal pair $(\bar x(\cd),\bar u(\cd))$ solves
the equation \eqref{fsystem1}, we see that
$x_1^\e(\cdot)$ satisfies the following
stochastic differential equation:
\begin{equation}\label{fsystem3zx}
\left\{
\begin{array}{lll}\ds
dx_1^\e = \big[Ax_1^\e + \tilde a_1^\e(t)
x^\e_1 + \chi_{E_\e} (t)\d a(t) \big]dt + \big[
\tilde b_1^\e(t) x^\e_1 + \chi_{E_\e} (t)\d
b(t) \big]dw(t) &\mbox{ in
}(0,T],\\
\ns\ds x_1^\e(0)=0.
\end{array}
\right.
\end{equation}

Consider the following two stochastic
differential equations:
\begin{equation}\label{s7fsystem3.1}
\left\{
\begin{array}{lll}\ds
dx_2^\e = \big[Ax_2^\e + a_1(t)x_2^\e \big]dt +
\big[ b_1(t) x_2^\e + \chi_{E_\e} (t)\d b(t)
\big] dw(t) &\mbox{ in
}(0,T],\\
\ns\ds x_2^\e(0)=0
\end{array}
\right.
\end{equation}
and\footnote{Recall that, for any
$C^2$-function $f(\cd)$ defined on a Banach
space $X$ and $x_0\in X$, $f_{xx}(x_0)\in \cL
(X\times X,X)$. This means that, for any
$x_1,x_2\in X$, $f_{xx}(x_0)(x_1,x_2)\in X$.
Hence, by (\ref{s7tatb1}),
$a_{11}(t)\big(x_2^\e,x_2^\e\big)$ (in
(\ref{s7fsystem3.2})) stands for $a_{xx}(t,\bar
x(t),\bar u(t))\big(x_2^\e(t),x_2^\e(t)\big)$.
One has a similar meaning for
$b_{11}(t)\big(x_2^\e,x_2^\e\big)$ and so on.}
\begin{equation}\label{s7fsystem3.2}
\left\{
\begin{array}{lll}\ds
dx_3^\e = \[Ax_3^\e + a_1(t)x_3^\e +
\chi_{E_\e}(t)\d a(t) +
\frac{1}{2}a_{11}(t)\big(x_2^\e,x_2^\e\big) \]dt \\
\ns\ds \hspace{1cm} + \[ b_1(t) x_3^\e +
\chi_{E_\e} (t)\d b_1(t)x_2^\e +
\frac{1}{2}b_{11}(t)\big(x_2^\e,x_2^\e\big)\]
dw(t) &\mbox{ in
}(0,T],\\
\ns\ds x_3^\e(0)=0.
\end{array}
\right.
\end{equation}
In the following Steps 2--4, we shall show that
\begin{equation}\label{s7th max eq0.2}
|x_1^\e - x_2^\e -
x_3^\e|_{L^\infty_\dbF(0,T;L^2(\O;H))}=o(\e),\qq\hb{as
} \e\to0.
\end{equation}

\ms

{\bf Step 2}. In this step, we provide some
estimates on $x_i^\e$ ($i=1,2,3$).

First of all, by a direct computation, we find
\begin{equation}\label{s7th max eq0}
\begin{array}{ll}\ds
 \mE|x_1^\e(t)|^8_H  = \mE\Big| \int_0^t S(t-s)
 \tilde a_1^\e(s) x_1^\e(s) ds
 +  \int_0^t S(t-s)\chi_{E_\e}(s)\d a(s) ds \\
\ns\ds \hspace{2.4cm}  +
  \int_0^t S(t-s)\tilde b_1^\e(s) x_1^\e(s) dw(s) +  \int_0^t S(t-s)\chi_{E_\e}(s)\d b(s) dw(s)\Big|_H^8\\
\ns\ds \hspace{1.7cm}\leq C \bigg\{\mE\Big|
\int_0^t S(t-s)\tilde a_1^\e(s) x_1^\e(s) ds
\Big|_H^8
 +  \mE\Big| \int_0^t S(t-s)\chi_{E_\e}(s)\d a(s)
ds\Big|_H^8\\
\ns\ds \hspace{2.4cm}  + \mE\Big|\int_0^t
S(t-s)\tilde b_1^\e(s) x_1^\e(s) dw(s)
\Big|_H^8 + \mE\Big| \int_0^t
S(t-s)\chi_{E_\e}(s)\d b(s)
dw(s)\Big|_H^8\bigg\}.
\end{array}
\end{equation}
Now, we estimate the terms in the right hand
side of the inequality \eqref{s7th max eq0} one
by one. For the first term, we have
 \bel{1z1}
\begin{array}{ll}\ds
\mE\Big| \int_0^t S(t-s)\tilde a_1^\e(s)
x_1^\e(s) ds \Big|_H^8 \leq C \int_0^t
\mE\big|\tilde a_1^\e(s) x_1^\e(s)\big|_H^8 ds
\leq C \mE \int_{0}^t \big| x_1^\e(s) \big|_H^8
ds.
\end{array}
\ee
By the last condition in (\ref{ab0}), it
follows that
 \bel{1z1-1}
\begin{array}{ll}\ds
\q|\d a(s)|_H=\big| a\big(s,\bar x(s),u(s)
\big)-a\big(s,\bar x(s),\bar u(s)\big) \big|_H \\
\ns\ds \le \Big|a\big(s,\bar x(s),u(s)
\big)-a\big(s,0,u(s)
\big)\Big|_H+\Big|a\big(s,0,\bar
u(s)\big)-a\big(s,\bar x(s),\bar
u(s)\big)\Big|_H\\
\ns\ds\q+\Big|a\big(s,0,u(s)
\big)-a\big(s,0,\bar u(s)\big)\Big|_H \\
\ns\ds \le \Big|\int_0^1a_x\big(s,\si\bar
x(s),u(s) \big)\bar x(s)d\si
\Big|_H+\Big|\int_0^1a_x\big(s,\si\bar
x(s),\bar u(s) \big)\bar
x(s)d\si \Big|_H+C_L\\
\ns\ds \le C\big[|\bar x(s) |_H+1\big],\qq \ae
s\in [0,T].
\end{array}
 \ee
Hence, using Lemma \ref{well lemma s1}, we have
the following estimate:
 \bel{1z1-2}
\begin{array}{ll}\ds
\q\mE\Big| \int_0^t S(t-s)\chi_{E_\e}(s)\d a(s)
ds\Big|_H^8\leq C\mE\Big\{\int_0^t
\chi_{E_\e}(s)\big|\d
a(s)\big|_H ds\Big\}^8\\
\ns\ds \le C\mE\Big\{\int_0^T
\chi_{E_\e}(s)\big[|\bar x(s)
|_H+1\big]  ds\Big\}^8\\
\ns\ds\le C\mE\Big\{\Big[\int_0^T
\chi_{E_\e}(s)ds\Big]^{7/8}\Big[\int_0^T
\chi_{E_\e}(s)\big(|\bar x(s)
|_H^8+1\big)ds \Big]^{1/8} \Big\}^8 \\
\ns\ds\le C\e^7\int_0^T
\chi_{E_\e}(s)\big(\mE|\bar x(s) |_H^8+1\big)ds
\le
C\big(|x(\cd)|_{C_\dbF([0,T];L^8(\O;H))}^8+1\big)\e^7\int_0^T
\chi_{E_\e}(s)ds \\
\ns\ds \leq C(x_0) \e^8.
\end{array}
 \ee
Here and henceforth, $ C(x_0)$ is a generic
constant (depending on $x_0$, $T$, $A$ and
$C_L$), which may be different from line to
line. By Lemma \ref{BDG}, similar to
(\ref{1z1}), it follows that
 \bel{1z3}
\begin{array}{ll}\ds
\mE\Big|\int_0^t S(t-s)\tilde b_1^\e(s)
x_1^\e(s) dw(s) \Big|_H^8 \leq C \mE \[\int_0^t
\Big| S(t-s)\tilde b_1^\e(s)x_1^\e(s) \Big|^2_H
ds\]^4 \leq C\mE \int_{0}^t \big| x_1^\e(s)
\big|^8ds.
\end{array}
 \ee
Similar to (\ref{1z1-1}), we have
 \bel{1z3-1}
|\d b(s)|_H\le C\big[|\bar x(s) |_H+1\big],\qq
\ae s\in [0,T].
 \ee
Hence, by Lemma \ref{BDG} again, similar to
(\ref{1z1-2}), one has
 \bel{1z4}
 \ba{ll}\ds
\q\mE\Big| \int_0^t S(t-s)\chi_{E_\e}(s)\d b(s)
dw(s)\Big|_H^8 \leq C
\[\mE\int_0^t  \chi_{E_\e}(s)|\d b(s)|_H^2 ds\]^4 \\
\ns\ds \le C\mE\Big\{\int_0^T
\chi_{E_\e}(s)\big[|\bar x(s)
|_H^2+1\big]  ds\Big\}^4\\
\ns\ds\le C\mE\Big\{\Big[\int_0^T
\chi_{E_\e}(s)ds\Big]^{3/4}\Big[\int_0^T
\chi_{E_\e}(s)\big(|\bar x(s)
|_H^8+1\big)ds \Big]^{1/4} \Big\}^4 \\
\ns\ds\le C\e^3\int_0^T
\chi_{E_\e}(s)\big(\mE|\bar x(s) |_H^8+1\big)ds
\le
C\big(|x(\cd)|_{C_\dbF([0,T];L^8(\O;H))}^8+1\big)\e^3\int_0^T
\chi_{E_\e}(s)ds \\
\ns\ds \leq C(x_0) \e^4.
\end{array}
 \ee
Therefore, combining (\ref{s7th max eq0}),
(\ref{1z1}), (\ref{1z1-2}), (\ref{1z3}) and
(\ref{1z4}), we end up with
  $$
 \mE|x_1^\e(t)|^8_H   \leq C(x_0)\[ \int_0^t \mE|x_1^\e(s)|_H^8ds +
 \e^8 + \e^4  \],\qq \ae t\in [0,T].
 $$
This, together with Gronwall's inequality,
implies that
\begin{equation}\label{s7th max eq0.1}
 |x_1^\e(\cd)|^8_{L^\infty_\dbF(0,T;L^8(\O;H))} \leq
C(x_0) \e^4.
\end{equation}
From the inequality \eqref{s7th max eq0.1} and
H\"{o}lder's inequality, we find that
\begin{equation}\label{s7th max eq0.1b}
\ba{ll}\ds
|x_1^\e(\cd)|^4_{L^\infty_\dbF(0,T;L^4(\O;H))}
\leq C(x_0) \e^2,\\\ns\ds
|x_1^\e(\cd)|^2_{L^\infty_\dbF(0,T;L^2(\O;H))}
\leq C(x_0) \e.
 \ea
\end{equation}

By a similar computation, we have
\begin{equation}\label{s7th max eq1}
\begin{array}{ll}\ds
\q\mE|x_2^\e(t)|^8_H \\
\ns\ds = \mE\Big| \int_0^t
S(t-s)a_1(s)x_2^\e(s) ds + \int_0^t
S(t-s)b_1(s)x_2^\e(s) dw(s)  +
\int_0^t S(t-s)\chi_{E_\e}(s)\d b(s) dw(s)\Big|_H^8\\
\ns\ds  \leq C \bigg[\mE\Big| \int_0^t
S(t-s)a_1(s)x_2^\e(s) ds \Big|_H^8
+ \mE\Big|\int_0^t S(t-s)b_1(s)x_2^\e(s) dw(s) \Big|_H^8 \\
\ns\ds \q  + \mE\Big| \int_0^t S(t-s)\chi_{E_\e}(s)\d b(s) dw(s)\Big|_H^8\bigg] \\
\ns\ds   \leq C(x_0)\(\int_0^t
\mE|x_2^\e(s)|_H^8ds  + \e^4 \).
\end{array}
\end{equation}
By means of Gronwall's inequality once more,
(\ref{s7th max eq1}) leads to
\begin{equation}\label{s7th max eq2}
|x_2^\e(\cd)|^8_{L^\infty_\dbF(0,T;L^8(\O;H))}
\leq C(x_0)\e^4.
\end{equation}
From inequality \eqref{s7th max eq2} and
utilizing H\"{o}lder's inequality again, we get
\begin{equation}\label{s7th max eq2b}
 \ba{ll}\ds |x_2^\e(\cd)|^4_{L^\infty_\dbF(0,T;L^4(\O;H))} \leq
C(x_0)\e^2,\\\ns\ds
|x_2^\e(\cd)|^2_{L^\infty_\dbF(0,T;L^2(\O;H))}
\leq C(x_0)\e. \ea
\end{equation}

Similar to (\ref{1z1-2}), we have
$$\mE\Big| \int_0^t
S(t-s)\chi_{E_\e}(s)\d a(s) ds\Big|_H^4\leq
C(x_0) \e^4.
 $$
Hence, it follows from Lemma \ref{BDG} and
\eqref{s7th max eq2}--\eqref{s7th max eq2b}
that
\begin{equation}
\begin{array}{ll}\ds
\q
|x_3^\e(t)|^4_{ L^4_{\cF_t}(\O;H)}\nonumber\\
\ns\ds = \mE\Big| \int_0^t S(t-s) a_1(s)x_3^\e
(s)ds + \int_0^t S(t-s)\chi_{E_\e}(s)\d a(s) ds
+ \frac{1}{2}\int_0^t S(t-s)
a_{11}(s)\big(x_2^\e(s),x_2^\e(s)\big) ds \\
\ns\ds \q + \int_0^t S(t-s) b_1(s) x_3^\e(s)
dw(s) +  \int_0^t
S(t-s)\chi_{E_\e}(s)\d b_1(s)x_2^\e(s) dw(s) \\
\ns\ds \q  + \frac{1}{2}\int_0^t S(t-s)
b_{11}(s)\big(x_2^\e(s),x_2^\e
(s)\big)dw(s)\Big|_H^4
\\
\ns\ds \leq C(x_0)\mE\[\!
\int_0^t\!|x_3^\e(s)|_H^4ds \!+\! \Big|
\!\int_0^t S(t-s)\chi_{E_\e}(s)\d a(s)
ds\Big|_H^4 \!+\! \int_0^T\!|x_2^\e(t)|_H^8 ds
\!+ \! \Big|\!\int_0^T\!|\chi_{E_\e}
x_2^\e(t)|_H^2dt\Big|^2 \]\\
\ns\ds \leq C(x_0)\[\mE
\int_0^t\!|x_3^\e(s)|_H^4ds \!+\e^4+
\mE\Big(\,\Big|\int_0^T
\chi_{E_\e}(s)ds\Big|^{1/2}\Big|\int_0^T
\chi_{E_\e}(s)|x_2^\e(s)
|_H^4ds \Big|^{1/2} \,\Big)^2 \]\\
\ns\ds \leq C(x_0)\[\mE \int_0^t
|x_3^\e(s)|_H^4ds  +\e^4+\e \int_0^T
\chi_{E_\e}(s)\mE|x_2^\e(s)
|_H^4ds  \]\\
\ns\ds \leq C(x_0)\[ \mE\int_0^t
|x_3^\e(s)|_H^4ds + \e^4 \],
\end{array}
\end{equation}
this, together with Gronwall's inequality,
implies that
\begin{equation}\label{s7th max eq2.1}
|x_3^\e(\cd)|^4_{L^\infty_\dbF(0,T;L^4(\O;H))}\leq
C(x_0)\e^4.
\end{equation}
Then, by H\"{o}lder's inequality, we conclude
that
\begin{equation}\label{s7th max eq2.1b}
|x_3^\e(\cd)|^2_{L^\infty_\dbF(0,T;L^2(\O;H))}
\leq C(x_0) \e^2.
\end{equation}

\ms

{\bf Step 3}. We now estimate $x_4^\e \= x_1^\e
- x_2^\e$. Clearly, $x_4^\e$ solves the
following equation:
\begin{equation}\label{fsystem3z2}
\left\{
\begin{array}{lll}\ds
dx_4^\e = \big\{Ax_4^\e + \tilde a_1^\e(t)
x^\e_4 + \big[\tilde a_1^\e(t) -
a_1(t)\big]x_2^\e + \chi_{E_\e}(t)\d
a(t) \big\}dt \\
\ns\ds\hspace{1.2cm} + \big\{\tilde b_1^\e(t)
x^\e_4 + [\tilde b_1^\e(t) - b_1(t)] x_2^\e
\big\} dw(t) &\mbox{ in
}(0,T],\\
\ns\ds x_4^\e(0)=0.
\end{array}
\right.
\end{equation}
Hence,
\begin{equation}\label{s7th max eq2.2}
\begin{array}{ll}\ds
\q\mE|x_4^\e(t)|^2_{H} \\
\ns\ds  = \mE\Big| \int_0^t S(t-s)\tilde
a_1^\e(s) x_4^\e(s) ds +\int_0^t
S(t-s)\big[\tilde a_1^\e(s) - a_1(s)\big]
x_2^\e(s) ds+ \int_0^t S(t-s)\chi_{E_\e}(s)\d
a(s)ds
  \\
\ns\ds \q+ \int_0^t S(t-s)\tilde b_1^\e(s)
x_4^\e(s) dw(s)   + \int_0^t S(t-s)\big[\tilde
b_1^\e(s) - b_1(s)\big]
 x_2^\e(s) dw(s)\Big|_H^2\\
\ns\ds \leq C\bigg[\mE\Big| \int_0^t
S(t-s)\tilde a_1^\e(s) x_4^\e(s) ds \Big|_H^2
+\mE\Big|\int_0^t S(t-s)\big[\tilde a_1^\e(s) -
a_1(s)\big] x_2^\e(s) ds \Big|_H^2
  \\
\ns\ds \q +  \mE\Big|\int_0^t
S(t-s)\chi_{E_\e}(s)\d a(s)ds \Big|_H^2
+\mE\Big| \int_0^t S(t-s)\tilde b_1^\e(s)
x_4^\e(s) dw(s) \Big|_H^2
  \\
\ns\ds \q +  \mE\Big|\int_0^t S(t-s)\big[\tilde
b_1^\e(s) - b_1(s)\big]
 x_2^\e(s) dw(s)\Big|_H^2\bigg].
\end{array}
\end{equation}

We now estimate the terms in the right hand
side of the inequality \eqref{s7th max eq2.2}
one by one. It is easy to see that
 \bel{2z1}
\mE\Big| \int_0^t S(t-s)\tilde a_1^\e(s)
x_4^\e(s) ds  \Big|_H^2 \leq C\mE \int_0^t
\big| \tilde a_1^\e(s)x_4^\e(s) \big|^2_H ds
\leq C(x_0) \mE \int_{0}^t \big| x_4^\e(s)
\big|^2ds.
 \ee
By (\ref{s7th max eq2b}), we have
 \bel{2z6}
\begin{array}{ll}\ds
\q \big|\!\big|\tilde a_1^\e(s) -
a_1(s)\big|\!\big|_{\cL(H)}
  \\
\ns\ds =\Big|\!\Big| \int_0^1
\big[a_x\big(s,\bar x(s) + \si x_1^\e (w) ,
u^\e (s)\big)- a_x(s,\bar x(s),\bar u(s))\big]d\si \Big|\!\Big|_{\cL(H)}\\
\ns\ds = \Big|\!\Big|\int_0^1
\big[a_x\big(s,\bar x(s) + \si x_1^\e (s) ,
u^\e (s)\big)- a_x(s,\bar x(s),u^\e (s))\\
\ns\ds\qq\qq+a_x(s,\bar x(s),u^\e (s))- a_x(s,\bar x(s),\bar u(s))\big]d\si \Big|\!\Big|_{\cL(H)} \\
\ns\ds = \Big|\!\Big|\int_0^1 \Big[\si
\int_0^1a_{xx}\big(s,\bar x(s) + \eta\si x_1^\e
(s) ,
u^\e (s)\big)x_1^\e (s)d\eta+\chi_{E_e}(s)\d a_1(s)\Big]d\si \Big|\!\Big|_{\cL(H)}  \\
\ns\ds \leq C \left[\big| x_1^\e(s) \big|_H+
\chi_{E_\e}(s)\right],\qq \ae s\in [0,T].
\end{array}
 \ee
Hence,
  \bel{2z2}
\begin{array}{ll}\ds
\q \mE\Big|\int_0^t S(t-s)\big[\tilde a_1^\e(s)
- a_1(s)\big] x_2^\e(s) ds \Big|_H^2
  \\
\ns\ds \le C\mE\int_0^t \big|\!\big|\tilde
a_1^\e(s) - a_1(s)\big|\!\big|_{\cL(H)}^2
\big|x_2^\e(s)\big|_H^2 ds
  \\
\ns\ds \leq C
|x_2^\e(\cd)|^2_{L^\infty_\dbF(0,T;L^4(\O;H))}
\int_0^T \[\mE\big|\!\big|\tilde a_1^\e(s) -
a_1(s)\big|\!\big|_{\cL(H)}^4\]^{1/2} ds\\
\ns\ds \leq C(x_0)\e
\int_{0}^T\left[\chi_{E_\e}(t)+ \mE\big|
x_1^\e(t) \big|^4_H\right]^{1/2}dt\\
\ns\ds\leq C(x_0)\e \int_{0}^T\[\chi_{E_\e}(t)+
\(\mE\big| x_1^\e(t) \big|^4_H\)^{1/2}\]dt \leq
C(x_0)\e^2.
\end{array}
 \ee
Similar to (\ref{1z1-2}), we have
 \bel{2z3} \mE\Big| \int_0^t S(t-s)\chi_{E_\e}(s)\d a(s)
ds\Big|_H^2 \leq C(x_0) \e^2.
 \ee
By Lemma \ref{BDG} and similar to (\ref{2z1}),
it follows that
 \bel{2z4}
\begin{array}{ll}\ds
 \mE\Big| \int_0^t S(t-s)\tilde b_1^\e(s)
x_4^\e(s) dw(s) \Big|_H^2  \leq C \mE \int_0^t
\Big|
S(t-s)\tilde b_1^\e(s)x_4^\e(s) \Big|^2_H ds  \\
\ns\ds \leq C \mE \int_{0}^t \big| \tilde
b_1^\e(s)x_4^\e(s) \big|^2ds \leq C \mE
\int_{0}^t \big| x_4^\e(s) \big|^2ds.
\end{array}
 \ee
Similar to (\ref{2z6}), we have
 $$
\big|\!\big|\tilde b_1^\e(s) -
b_1(s)\big|\!\big|_{\cL(H)} \leq C \left[\big|
x_1^\e(s) \big|_H+ \chi_{E_\e}(s)\right],\qq
\ae s\in [0,T].
 $$
Hence, similar to (\ref{2z2}), one obtains that
 \bel{2z5}
\mE\Big|\int_0^t S(t-s)\big[\tilde b_1^\e(s) -
b_1(s)\big] x_2^\e(s) dw(s) \Big|_H^2  \leq
C\mE\int_0^t \big|\!\big|\tilde b_1^\e(s) -
b_1(s)\big|\!\big|_{\cL(H)}^2\big|
x_2^\e(s)\big|_H^2 ds \leq C(x_0)\e^2.
 \ee

Combining (\ref{s7th max eq2.2}), (\ref{2z1}),
(\ref{2z2}), (\ref{2z3}), (\ref{2z4}) and
(\ref{2z5}), we obtain that
$$
\mE|x_4^\e(t)|^2_{H} \leq  C(x_0) \(\int_0^t
\mE|x_4^\e(s)|_H^2 ds + \e^2\).
$$
Utilizing Gronwall's inequality again, we find
that
\begin{equation}\label{s7th max eq3}
|x_4^\e(\cd)|^2_{L^\infty_\dbF(0,T;L^2(\O;H))}
\leq C(x_0) \e^2,\q \forall t\in [0,T].
\end{equation}

\ms

{\bf Step 4}.  We are now in a position to
estimate $\mE|x_1^\e(t) - x_2^\e(t) -
x_3^\e(t)|^2_{H}=\mE|x_4^\e(t) - x_3^\e(t)
|^2_H$.

Let $x_5^\e(\cd) = x_4^\e(\cd) - x_3^\e(\cd)$.
It is clear that $x_5^\e(\cd)=x_1^\e(t) -
x_2^\e(t) - x_3^\e(t)=x^\e(\cd)-\bar x(\cd) -
x_2^\e(t) - x_3^\e(t)$. We claim that
$x_5^\e(\cd)$ solves the following equation
(Recall (\ref{s7tatb1})--(\ref{s7tatb2}) for
the notations):
\begin{equation}\label{s7fsystem4}
\left\{
\begin{array}{ll}
\ds dx_5^\e = \Big[Ax_5^\e + a_1(t) x_5^\e +
\chi_{E_\e}(t)\d a_1(t)x_1^\e +
\frac{1}{2}\big(\tilde a_{11}^\e(t) -
a_{xx}(t,\bar
x, u^\e) \big)\big(x_1^\e,x_1^\e\big)  \\
\ns\ds \hspace{1.3cm} +
\frac{1}{2}\chi_{E_\e}(t) \d
a_{11}(t)\big(x_1^\e,x_1^\e\big) +
\frac{1}{2}a_{11}(t)\big(x_1^\e,x_1^\e\big) -
\frac{1}{2}a_{11}(t)\big(x_2^\e,x_2^\e\big)\Big]dt\\
\ns\ds \hspace{1.3cm} + \Big[b_1(t) x_5^\e +
\chi_{E_\e}(t)\d b_1(t) x_4^\e +
\frac{1}{2}\big(\,\tilde b_{11}^\e(t) -
b_{xx}(t,\bar x,u^\e)
\big)\big(x_1^\e,x_1^\e\big)  \\
\ns\ds\hspace{1.3cm} +
\frac{1}{2}\chi_{E_\e}(t)\d b_{11}(t)\big(
x_1^\e,x_1^\e\big) +
\frac{1}{2}b_{11}(t)\big(x_1^\e,x_1^\e\big) -
\frac{1}{2}b_{11}(t)\big(x_2^\e,x_2^\e\big)
\Big]dw(t), &\mbox{ in }(0,T],\\
\ns\ds x_5^\e(0) = 0.
\end{array}
\right.
\end{equation}
Indeed, by (\ref{s7fsystem2}),
\eqref{fsystem1}, (\ref{s7fsystem3.1}) and
(\ref{s7fsystem3.2}), it is easy to see that
the drift term for the equation solved by
$x_5^\e(\cd)$ is as follows:
\begin{equation}\label{zeq1}
\begin{array}{ll}\ds
Ax^\e +a(t,x^\e,u^\e)-A\bar x -a(t,\bar x,\bar
u)-Ax_2^\e - a_1(t)x_2^\e- Ax_3^\e -
a_1(t)x_3^\e  \\
\ns\ds \q- \chi_{E_\e}(t)\d a(t) -
\frac{1}{2}a_{11}(t)\big(x_2^\e,x_2^\e\big)\\
\ns\ds = Ax_5^\e+a(t,x^\e,u^\e)-a(t,\bar
x,u^\e) - a_1(t)(x_2^\e + x_3^\e)-
\frac{1}{2}a_{11}(t)\big(x_2^\e,x_2^\e\big).
\end{array}
\end{equation}
For $\si\in [0,1]$, write $f(\si )=a(t,\bar
x+\si  x_1^\e,u^\e)$. Then, by Taylor's formula
with the integral type remainder, we see that
$$
f(1)-f(0)=f'(0) + \int_0^1(1-\si ) f''(\si
)d\si .
$$
Since $f'(\si )=a_x(t,\bar x+\si
x_1^\e,u^\e)x_1^\e$ and $f''(\si
)=a_{xx}(t,\bar x+\si
x_1^\e,u^\e)(x_1^\e,x_1^\e)$,
 we obtain that
\begin{equation}\label{zeq2}
\begin{array}{ll}\ds
a(t,x^\e,u^\e)-a(t,\bar x,u^\e) = a_x(t,\bar
x,u^\e)x_1^\e +
\int_0^1 (1-\si )a_{xx}(t,\bar x+\si  x_1^\e,u^\e)(x_1^\e,x_1^\e)d\si \\
\ns\ds \hspace{3.89cm} = a_x(t,\bar x,u^\e)
x_1^\e + \frac{1}{2}\tilde
a_{11}^\e(t)(x_1^\e,x_1^\e).
\end{array}
\end{equation}
Next,
\begin{equation}\label{zeq3}
\begin{array}{ll}\ds
\q a_x(t,\bar x,u^\e) x_1^\e -  a_1(t)(x_2^\e + x_3^\e)\\
\ns\ds = a_x(t,\bar x,u^\e) x_1^\e -a_x(t,\bar
x,\bar u)x_1^\e +
a_1(t)(x_1^\e - x_2^\e - x_3^\e)\\
\ns\ds = \chi_{E_\e}\big[a_x(t,\bar
x,u)-a_x(t,\bar x,\bar
u)\big]x_1^\e  + a_1(t) x_5^\e \\
\ns\ds = \chi_{E_\e}\d a_1(t)x_1^\e + a_1(t)
x_5^\e.
\end{array}
\end{equation}
Further,
\begin{equation}\label{zeq4}
\begin{array}{ll}\ds
\q\frac{1}{2}\tilde a_{11}^\e(t)(x_1^\e,x_1^\e)
-\frac{1}{2} a_{11}(t)(x_2^\e,x_2^\e)\\
\ns\ds = \frac{1}{2} \tilde
a_{11}^\e(x_1^\e,x_1^\e)-\frac{1}{2}
a_{xx}(t,\bar x,u^\e)
(x_1^\e,x_1^\e)+\frac{1}{2}a_{xx}(t,\bar
x,u^\e) (x_1^\e,x_1^\e)-\frac{1}{2}
a_{11}(t)(x_1^\e,x_1^\e)\\
\ns\ds \q +\frac{1}{2}
a_{11}(t)(x_1^\e,x_1^\e)-\frac{1}{2}
a_{11}(t)(x_2^\e,x_2^\e)\\
\ns\ds = \frac{1}{2} \(\tilde
a_{11}^\e(t)-a_{xx}(t,\bar
x,u^\e)\)(x_1^\e,x_1^\e) +
\frac{1}{2}\chi_{E_\e}\d
a_{11}(t)(x_1^\e,x_1^\e) +\frac{1}{2}
a_{11}(t)(x_1^\e,x_1^\e)-\frac{1}{2}
a_{11}(t)(x_2^\e,x_2^\e).
\end{array}
\end{equation}
By \eqref{zeq1}--\eqref{zeq4}, we conclude that
$$
\begin{array}{ll}\ds
\q Ax^\e +a(t,x^\e,u^\e)-A\bar x -a(t,\bar
x,\bar u)-Ax_2^\e - a_1(t)x_2^\e- Ax_3^\e -
a_1(t)x_3^\e  \\
\ns\ds \q- \chi_{E_\e}(t)\d a(t) -
\frac{1}{2}a_{11}(t)\big(x_2^\e,x_2^\e\big)-\chi_{E_\e}(t)\d
b_1(t)
x_2^\e\\
\ns\ds = Ax_5^\e +a_1(t) x_5^\e +
\chi_{E_\e}(t)\d a_1(t)x_1^\e +
\frac{1}{2}\big(\tilde a_{11}^\e(t) -
a_{xx}(t,\bar
x, u^\e) \big)\big(x_1^\e,x_1^\e\big)  \\
\ns\ds \q + \frac{1}{2}\chi_{E_\e}(t) \d
a_{11}(t)\big(x_1^\e,x_1^\e\big) +
\frac{1}{2}a_{11}(t)\big(x_1^\e,x_1^\e\big) -
\frac{1}{2}a_{11}(t)\big(x_2^\e,x_2^\e\big).\end{array}
$$
Similarly, the diffusion term (for the equation
solved by $x_5^\e(\cd)$) is as follows:
 $$
\begin{array}{ll}\ds
\q b(t,x^\e,u^\e) -b(t,\bar x,\bar u) -
b_1(t)x_2^\e- b_1(t)x_3^\e - \chi_{E_\e}(t)\d
b(t) x_2^\e-
\frac{1}{2}b_{11}(t)\big(x_2^\e,x_2^\e\big)\\
\ns\ds = b_1(t) x_5^\e + \chi_{E_\e}(t)\d
b_1(t) x_4^\e + \frac{1}{2}\big(\,\tilde
b_{11}^\e(t) - b_{xx}(t,\bar x,u^\e)
\big)\big(x_1^\e,x_1^\e\big)  \\
\ns\ds\q+ \frac{1}{2}\chi_{E_\e}(t)\d
b_{11}(t)\big( x_1^\e,x_1^\e\big) +
\frac{1}{2}b_{11}(t)\big(x_1^\e,x_1^\e\big) -
\frac{1}{2}b_{11}(t)\big(x_2^\e,x_2^\e\big).\end{array}
$$
This verifies that $x_5^\e(\cd)$ satisfies the
equation (\ref{s7fsystem4}).

From (\ref{s7fsystem4}), we see that, for any
$t\in [0,T]$,
\begin{equation}\label{s7fem4}
\begin{array}{ll}
\ds \mE\big| x_5^\e(t)\big|_H^2  \le C
\Big[\mE\Big| \int_0^t S(t-s)
\big[a_1(s)x_5^\e(s) +\chi_{E_\e}(s)\d a_1(s)
x_1^\e(s)\big]ds \Big|_H^2\\\ns\ds
 \qq\qq\q\q\;\;+\mE\Big| \int_0^t S(t-s) \big(\tilde a_{11}^\e(s)
- a_{xx}(s,\bar x(s), u^\e(s))
\big)\big(x_1^\e(s),x_1^\e(s)\big)ds
\Big|_H^2\\\ns\ds
 \qq\qq\q\q\;\;+\mE\Big| \int_0^t S(t-s) \chi_{E_\e}(s) \d
a_{11}(s)\big(x_1^\e(s),x_1^\e(s)\big)ds
\Big|_H^2\\\ns\ds
 \qq\qq\q\q\;\;+\mE\Big| \int_0^t S(t-s) \big[a_{11}(s)\big(x_1^\e(s),x_1^\e(s)\big) -
a_{11}(s)\big(x_2^\e(s),x_2^\e(s)\big)\big]ds
\Big|_H^2
\\\ns\ds
 \qq\qq\q\q\;\;+\mE\Big| \int_0^t S(t-s) \big[b_1(s)x_5^\e(s) +\chi_{E_\e}(s)\d
b_1(s) x_4^\e(s)\big]dw(s) \Big|_H^2\\\ns\ds
 \qq\qq\q\q\;\;+\mE\Big| \int_0^t S(t-s) \big(\tilde b_{11}^\e(s)
- b_{xx}(s,\bar x(s), u^\e(s))
\big)\big(x_1^\e(s),x_1^\e(s)\big)dw(s)
\Big|_H^2\\\ns\ds
 \qq\qq\q\q\;\;+\mE\Big| \int_0^t S(t-s) \chi_{E_\e}(s) \d
b_{11}(s)\big(x_1^\e(s),x_1^\e(s)\big)dw(s)
\Big|_H^2\\\ns\ds
 \qq\qq\q\q\;\;+\mE\Big| \int_0^t S(t-s) \big[b_{11}(s)\big(x_1^\e(s),x_1^\e(s)\big) -
b_{11}(s)\big(x_2^\e(s),x_2^\e(s)\big)\big]dw(s)
\Big|_H^2\Big].
 \ea
\end{equation}

We now estimate the ``drift" terms in the right
hand side of \eqref{s7fem4}. By \eqref{s7th max
eq0.1b}, we have the following estimate:
\begin{equation}\label{s7eq0.1}
\begin{array}{ll}\ds
\q\mE\Big| \int_0^t S(t-s) \big[a_1(s)x_5^\e(s)
+\chi_{E_\e}(s)\d a_1(s)
x_1^\e(s)\big]ds \Big|_H^2\\
\ns\ds \leq C\Big[\mE\int_0^t \big| x_5^\e(s)
\big|_H^2 ds +\mE\Big|\int_0^T
 \chi_{E_\e}(s)|x_1^\e(\cd)|_Hds \Big|^2\Big]\\
\ns\ds \leq C\Big[\mE\int_0^t \big| x_5^\e(s)
\big|_H^2 ds +\mE\Big|\int_0^T
 \chi_{E_\e}(s)ds \int_0^T
 \chi_{E_\e}(s)|x_1^\e(s)|_H^2ds \Big|\Big] \\
\ns\ds \leq C \Big[\mE\int_0^t \big| x_5^\e(s)
\big|_H^2 ds +\e\int_0^T
 \chi_{E_\e}(s)\mE|x_1^\e(s)|_H^2ds \Big]\\
\ns\ds \leq C\Big[\mE\int_0^t \big| x_5^\e(s)
\big|_H^2 ds +\e |x_1^\e(\cd)
|_{L^\infty_\dbF(0,T;L^2(\O;H))}^2\int_0^t
\chi_{E_\e}(s)ds \Big]\\\ns\ds\leq
C(x_0)\Big[\mE\int_0^t \big| x_5^\e(s)
\big|_H^2 ds + \e^3\Big].
\end{array}
\end{equation}
By (\ref{s7tatb1}) and recalling that
$x_1^\e(\cd) = x^\e(\cd)-\bar x(\cd)$, we see
that, for $\ae s\in[0,T]$,
 \bel{2z9}
 \ba{ll}
 \ds\q\big|\!\big| \,\tilde a_{11}^\e(s) - a_{xx}(s,\bar x(s),
u^\e(s)) \big|\!\big|_{\cL(H\times H,
\,H)}\\\ns\ds =\Big|\!\Big| 2\int_0^1 (1-\si)
a_{xx}\big(s,\bar x(s) + \si x_1^\e(s), u^\e
(s)\big)d\si- a_{xx}(s,\bar x(s), u^\e(s))
\Big|\!\Big|_{\cL(H\times H, \,H)}\\\ns\ds
=\Big|\!\Big| 2\int_0^1 (1-\si)
\[a_{xx}\big(s,\bar x(s) + \si x_1^\e(s), \bar
u (s)\big)- a_{xx}(s,\bar x(s), \bar
u(s))\]d\si \\\ns\ds
 \q+2\int_0^1 (1-\si) \chi_{E_\e}(s)a_{xx}\big(s,\bar
x(s) + \si x_1^\e(s),  u
(s)\big)d\si+\chi_{E_\e}(s)a_{xx}(s,\bar x(s),
u(s)) \Big|\!\Big|_{\cL(H\times H,
\,H)}\\\ns\ds \le C\[\int_0^1
\big|\!\big|a_{xx}\big(s,\bar x(s) + \si
x_1^\e(s), \bar u (s)\big)- a_{xx}(s,\bar x(s),
\bar u(s))\big|\!\big|_{\cL(H\times H,
\,H)}d\si +\chi_{E_\e}(s)\].
 \ea
 \ee
Hence, by \eqref{s7th max eq0.1} and noting the
continuity of $a_{xx}(t,x,u)$ with respect to
$x$, we have
\begin{equation}\label{s7eq0.2}
\begin{array}{ll}\ds
\q \mE\Big| \int_0^t S(t-s) \big(\,\tilde
a_{11}^\e(s) - a_{xx}(s,\bar x(s),
u^\e(s)) \big)\big( x_1^\e(s), x_1^\e(s)\big) ds \Big|_H^2 \\
\ns\ds  \leq C\mE\int_0^T \big|\!\big| \,\tilde
a_{11}^\e(s) - a_{xx}(s,\bar x(s),
u^\e(s)) \big|\!\big|_{\cL(H\times H, \,H)}^2|x_1^\e(s)|^4_Hdt\\
\ns\ds  \leq
C|x_1^\e(\cd)|^4_{L^\infty_\dbF(0,T;L^8(\O;H))}
\int_0^T \[\mE\big|\!\big| \,\tilde
a_{11}^\e(s) - a_{xx}(s,\bar x(s),
u^\e(s)) \big|\!\big|_{\cL(H\times H, \,H)}^4\]^{1/2}ds\\
\ns\ds  \leq C(x_0)\e^2
\int_0^T\!\!\[\mE\int_0^1
\big|\!\big|a_{xx}\big(s,\bar x(s) \!+\! \si
x_1^\e(s), \bar u (s)\big)\!-\! a_{xx}(s,\bar
x(s), \bar u(s))\big|\!\big|_{\cL(H\times H,
\,H)}^4d\si \!+\!\chi_{E_\e}(s)\]^{1/2}\!ds
 \\
\ns\ds =o(\e^2),\qq\hb{as }\e\to0.
\end{array}
\end{equation}
Also, it holds that
\begin{equation}\label{s7eq0.3}
\begin{array}{ll}
\ds \q\mE\Big|\int_0^t S(t-s)\chi_{E_\e}(s) \d
a_{11}(s)\big(x_1^\e(s),
 x_1^\e(s) \big) ds \Big|_H^2 \\
\ns\ds \leq C
|x_1^\e(\cd)|^4_{L^\infty_\dbF(0,T;L^8(\O;H))}
\int_0^T \chi_{E_\e}(s)
\[\mE\big|\!\big| \d
a_{11}(s)\big|\!\big|_{\cL(H\times H,\;H)}^4\]^{1/2} ds \\
\ns\ds  \leq C(x_0)\e^2 \int_0^T
\chi_{E_\e}(t)dt = C(x_0)\e^3.
\end{array}
\end{equation}
By means of \eqref{s7th max eq0.1b},
\eqref{s7th max eq2b} and \eqref{s7th max eq3},
and noting that $x_4^\e = x_1^\e - x_2^\e$, we
obtain that
\begin{equation}\label{s7th max eq4}
\begin{array}{ll}
\ds \q\frac{1}{2}\mE\Big| \int_0^t
S(t-s)\[a_{11}(s)\big(x_1^\e(s),x_1^\e(s)\big) - a_{11}(s)\big(x_2^\e(s),x_2^\e(s)\big) \] ds \Big|^2_H \\
\ns\ds = \frac{1}{2}\mE\Big| \int_0^t
S(t-s)\[a_{11}(s)\big(x_4^\e(s),x_1^\e(s)\big) + a_{11}(s)\big(x_2^\e(s),x_4^\e(s)\big) \] ds \Big|^2_H \\
\ns\ds  \leq
C\big(|x_1^\e(\cd)|^2_{L^\infty_\dbF(0,T;L^2(\O;H))}
+|x_2^\e(\cd)|^2_{L^\infty_\dbF(0,T;L^2(\O;H))}\big)|x_4^\e(\cd)|^2_{L^\infty_\dbF(0,T;L^2(\O;H))}\\
\ns\ds \leq C(x_0)\e^3.
\end{array}
\end{equation}

Next, we estimate the ``diffusion" terms in the
right hand side of \eqref{s7fem4}. Similar to
(\ref{s7eq0.1}) and noting \eqref{s7th max
eq3}, we obtain that
\begin{equation}\label{s7eq0.4}
\begin{array}{ll}\ds
\q\mE\Big| \int_0^t S(t-s) \big[b_1(s)x_5^\e(s)
+\chi_{E_\e}(s)\d
b_1(s) x_4^\e(s)\big]dw(s) \Big|^2_H \\
\ns\ds\leq \mE\int_0^t \big|S(t-s)
\big[b_1(s)x_5^\e(s) +\chi_{E_\e}(s)\d b_1(s)
x_4^\e(s)\big]\big|_H^2 ds  \leq
C(x_0)\Big[\mE\int_0^t \big| x_5^\e(s)
\big|_H^2 ds + \e^3\Big].
\end{array}
\end{equation}
By virtue of \eqref{s7th max eq0.1} again,
similar to (\ref{s7eq0.2}), we find that
\begin{equation}\label{s7eq0.5}
\begin{array}{ll}\ds
\q\mE\Big|\int_0^t S(t-s) \big( \tilde
b_{11}^\e(s) - b_{xx}(s,\bar
x(s),u^\e(s)) \big)\big(x_1^\e(s),x_1^\e(s)\big) dw(s)\Big|_H^2 \\
\ns\ds =  \mE\int_0^t \Big|S(t-s) \big( \tilde
b_{11}^\e(s) - b_{xx}(s,\bar x(s),u^\e(s))
\big)\big(x_1^\e(s),x_1^\e(s)\big)\Big|^2_H ds\\
\ns\ds  \leq
C|x_1^\e(\cd)|^4_{L^\infty_\dbF(0,T;L^8(\O;H))}
\int_0^T \[\mE\big|\!\big| \,\tilde
b_{11}^\e(s) - b_{xx}(s,\bar x(s),
u^\e(s)) \big|\!\big|_{\cL(H\times H, \,H)}^4\]^{1/2}ds\\
\ns\ds  \leq C(x_0)\e^2
\int_0^T\!\!\[\mE\!\int_0^1
\big|\!\big|b_{xx}\big(s,\bar x(s) \!+\! \si
x_1^\e(s), \bar u (s)\big)\!-\! b_{xx}(s,\bar
x(s), \bar u(s))\big|\!\big|_{\cL(H\times H,
\,H)}^4d\si +\chi_{E_\e}(s)\]^{1/2}\!ds
 \\
\ns\ds =o(\e^2),\qq\hb{as }\e\to0.
 \ea
 \ee
Similar to (\ref{s7eq0.3}), we have
\begin{equation}\label{s7eq0.6}
\begin{array}{ll}\ds
\q\mE\Big|\int_0^t S(t-s)\chi_{E_\e}(s) \d
b_{11}(s)\big(x_1^\e(s),
 x_1^\e(s) \big) dw(s) \Big|_H^2 \\
\ns\ds = \mE\int_0^t \big|S(t-s)\chi_{E_\e}(s)
\d b_{11}(s)\big(x_1^\e(s),
 x_1^\e(s) \big) \big|_H^2 ds  \leq C(x_0)\e^3.
\end{array}
\end{equation}
Similar to (\ref{s7th max eq4}), it holds that
\begin{equation}\label{s7eq0.7}
\begin{array}{ll}\ds
\q\mE \Big| \int_0^t S(t-s)
\[b_{11}(s)\big(x_1^\e(s),x_1^\e(s)\big) - b_{11}(s)\big(x_2^\e(s),x_2^\e(s)\big) \] dw(s) \Big|_H^2
\\ \ns\ds = \mE  \int_0^t \Big|S(t-s)
\[b_{11}(s)\big(x_1^\e(s),x_1^\e(s)\big) - b_{11}(s)\big(x_2^\e(s),x_2^\e(s)\big) \]\Big|_H^2 ds \leq C(x_0) \e^3.
\end{array}
\end{equation}

From \eqref{s7fem4}--(\ref{s7eq0.1}) and
(\ref{s7eq0.2})--\eqref{s7eq0.7}, we conclude
that
\begin{equation}\label{s7eq1}
\mE|x_5^\e(t)|_H^2 \leq  C(x_0) \mE\int_0^t
\big| x_5^\e(s) \big|_H^2 ds +o(\e^2),\qq
\hb{as } t\to0.
\end{equation}
By means of Gronwall's inequality again, we get
 \bel{sqqeq0.2}
|x_5^\e(\cd)|_{L^\infty_\dbF(0,T;L^2(\O;H))}^2
=o(\e^2),\qq\hb{as } t\to0.
 \ee
This gives \eqref{s7th max eq0.2}.

\ms

{\bf Step 5}. We are now in a position to
complete the proof.

We need to compute the value of $\cJ(u^\e(\cd))
- \cJ(\bar u(\cd))$.
\begin{equation}\label{s7eq3}
\begin{array}{ll}\ds
\q  \cJ(u^\e(\cd)) -  \cJ(\bar u(\cd))\\
\ns\ds = \mE\int_0^T \big[g(t,x^\e(t),u^\e(t))
- g(t,\bar x(t),\bar u(t))\big]
dt + \mE h\big(x^\e(T)\big) - \mE h\big(\bar x(T)\big)\\
\ns\ds  = \Re\mE \int_0^T \Big\{
\chi_{E_\e}(t)\d g(t) + \big\langle g_x(t,\bar
x(t),
u^\e(t)),x_1^\e(t) \big\rangle_H  \\
\ns\ds \q + \int_0^1 \big\langle(1-\si)
g_{xx}\big(t, \bar x(t) + \si
x_1^\e(t),u^\e(t)\big)x_1^\e(t),x_1^\e(t)\big\rangle_H
d\si \Big\}dt \\
\ns\ds  \q +\Re\mE\big\langle h_x(\bar x(T)),
x_1^\e(T)\big\rangle_H + \Re\mE\int_0^1
\big\langle (1-\si) h_{xx}\big( \bar x(T)+ \si
x_1^\e(T)\big)x_1^\e(T), x_1^\e(T)\big\rangle_H
d\si.
\end{array}
\end{equation}
This, together with the definition of
$x_1^\e(\cd)$, $x_2^\e(\cd)$, $x_3 ^\e(\cd)$,
$x_4 ^\e(\cd)$  and $x_5^\e(\cd)$, yields that
\begin{equation}\label{s7eq4}
\begin{array}{ll}\ds
\q  \cJ(u^\e(\cd)) -  \cJ(\bar u(\cd))\\
\ns\ds = \Re\mE\int_0^T \Big\{ \chi_{E_\e}(t)\d
g(t) + \big\langle \d g_1(t),x_1^\e(t)
\big\rangle_H \chi_{E_\e}(t) + \big\langle
g_1(t), x_2^\e(t) + x_3^\e(t) \big\rangle_H +
\big\langle g_1(t),
x_5^\e(t) \big\rangle_H \\
\ns\ds\q  + \int_0^1\big\langle (1-\si) \big[
g_{xx}\big(t, \bar x(t) + \si x_1^\e(t),
u^\e(t)\big) - g_{xx}\big(t,\bar
x(t),u^\e(t)\big)\big]x_1^\e(t), x_1^\e(t)
\big\rangle_H d \si \\
\ns\ds \q + \frac{1}{2}\big\langle \d
g_{11}(t)x_1^\e(t), x_1^\e(t) \big\rangle_H
\chi_{E_\e}(t) + \frac{1}{2}\big\langle
g_{11}(t)x_2^\e(t),x_2^\e(t) \big\rangle_H +
\frac{1}{2} \big\langle g_{11}(t)x_4^\e(t),
x_1^\e(t) + x_2^\e(t) \big\rangle_H \Big\}dt \\
\ns\ds \q + \Re\mE\big\langle h_x\big(\bar
x(T)\big), x_2^\e(t) \!+\! x_3^\e(t)
\big\rangle_H \!+\! \Re\mE \big\langle
h_x\big(\bar x(T)\big), x_5^\e(t) \big\rangle_H
\!+\! \frac{1}{2}\Re\mE\big\langle
h_{xx}\big(\bar x(T)\big) x_2^\e(t),x_2^\e(t)
\big\rangle_H \\
\ns\ds \q + \frac{1}{2}\Re\mE \big\langle
h_{xx}\big( \bar x(T) \big)x_4^\e(T), x_1^\e(T)
+
x_2^\e(T)  \big\rangle_H \\
\ns\ds \q + \Re\mE\int_0^1 \big\langle
 (1-\si)\big[ h_{xx}\big( \bar x(T) + \si x_1^\e(T) \big) - h_{xx}\big(\bar x(T)\big)
\big]x_1^\e(T),x_1^\e(T) \big\rangle_H d \si.
\end{array}
\end{equation}
Similar to (\ref{2z9}), for $\ae t\in[0,T]$, we
find that
 \bel{2z10}
 \ba{ll}
 \ds\q\big|\!\big|\int_0^1 (1-\si) \big[ g_{xx}\big(t, \bar
x(t) + \si x_1^\e(t), u^\e(t)\big) -
g_{xx}\big(t,\bar
x(t),u^\e(t)\big)\big]d\si\big|\!\big|_{\cL(H\times
H,\;H)}\\\ns\ds =\big|\!\big| \int_0^1 (1-\si)
\[g_{xx}\big(t,\bar x(t) + \si x_1^\e(t), \bar u (t)\big)-
g_{xx}\big(t,\bar x(t), \bar u(t)\big)\]d\si
\\\ns\ds
 \q+\int_0^1 (1-\si) \chi_{E_\e}(t)g_{xx}\big(t,\bar
x(t) + \si x_1^\e(t),  u
(t)\big)d\si+\chi_{E_\e}(t)g_{xx}\big(t,\bar
x(t), u(t)\big) \big|\!\big|_{\cL(H\times
H,\;H)} d\si
\\\ns\ds \le C\[\int_0^1
\big|\!\big|g_{xx}\big(t,\bar x(t) + \si
x_1^\e(t), \bar u (t)\big)- g_{xx}\big(t,\bar
x(t), \bar u(t)\big)\big|\!\big|_{\cL(H\times
H,\;H)}d\si +\chi_{E_\e}(t)\].
 \ea
 \ee
By (\ref{s7eq4}), noting \eqref{s7th max
eq0.1}, \eqref{s7th max eq2}, \eqref{s7th max
eq2.1}, \eqref{s7th max eq3}, \eqref{sqqeq0.2}
and (\ref{2z10}), and using the continuity of
both $h_{xx}(x)$ and $g_{xx}(x)$ with respect
to $x$, we end up with
\begin{equation}\label{s7eq5}
\begin{array}{ll}
\ds
\q  \cJ(u^\e(\cd)) -  \cJ(\bar u(\cd))\\
\ns\ds = \Re\mE \int_0^T\[ \big\langle
g_1(t),x_2^\e(t) + x_3^\e(t) \big\rangle_H +
\frac{1}{2}\big\langle
g_{11}(t)x_2^\e(t),x_2^\e(t) \big\rangle_H +
\chi_{E_\e}(t)\d g(t)
\] dt \\
\ns\ds\q + \Re\mE \big\langle h_x\big(\bar
x(T)\big), x_2^\e(T)+x_3^\e(T) \big\rangle_H +
\frac{1}{2}\Re\mE\big\langle h_{xx}\big(\bar
x(T)\big)x_2^\e(t),x_2^\e(t) \big\rangle_H +
o(\e).
\end{array}
\end{equation}

In the sequel, we shall get rid of
$x_2^\e(\cd)$ and $x_3^\e(\cd)$ in
\eqref{s7eq5} by solutions to the equations
\eqref{bsystem1} and \eqref{op-bsystem3}.

First, by the definition of the transposition
solution to \eqref{bsystem1} (with $y_T$ and
$f(\cd,\cd,\cd)$ given by (\ref{zv1})), we
obtain that
\begin{equation}\label{s7eq6}
-\mE\big\langle h_x(\bar
x_T)),x_2^\e(T)\big\rangle_H - \mE \int_0^T
\big\langle g_1(t),x_2^\e(t)\big\rangle_H dt =
\mE \int_0^T\big\langle Y(t), \d
b(t)\big\rangle_H\chi_{E_\e}(t) dt
\end{equation}
and
\begin{equation}\label{s7eq7}
\begin{array}{ll}
\ds \q-\mE\big\langle h_x(\bar
x_T)),x_3^\e(T)\big\rangle_H - \mE \int_0^T
\big\langle g_1(t),x_3^\e(t)\big\rangle_H dt
\\
\ns\ds = \mE \int_0^T \Big\{  \frac{1}{2}\[
\big\langle y(t),a_{11}(t)\big(x_2^\e(t),
x_2^\e(t)\big) \big\rangle_H + \big\langle
Y(t), b_{11}(t)\big(
x_2^\e(t), x_2^\e(t)\big) \big\rangle_H \] \\
\ns\ds \hspace{1.8cm} + \chi_{E_\e}(t)
\[ \big\langle y(t),\d a(t) \big\rangle_H +
\big\langle Y,\d b_1(t)x_2^\e(t) \big\rangle_H
 \]\Big\}dt.
\end{array}
\end{equation}
According to \eqref{s7eq5}--\eqref{s7eq7},  we
conclude that
\begin{equation}\label{s7eq8}
\begin{array}{ll} \ds
\q\cJ(u^\e(\cd)) -  \cJ(\bar u(\cd))\\
\ns\ds = \frac{1}{2}\Re\mE\int_0^T\[
\big\langle g_{11}(t)x_2^\e(t),
x_2^\e(t)\big\rangle_H - \big\langle
y(t),a_{11}(t)\big(x_2^\e(t),
x_2^\e(t)\big) \big\rangle_H \\
\ns\ds \q - \big\langle Y,
b_{11}(t)\big(x_2^\e(t),
x_2^\e(t)\big)\big\rangle_H
\]dt + \Re\mE\int_0^T \chi_{E_\e}(t)\[ \d g(t) - \big\langle
y(t),\d a(t)\big\rangle_H \\
\ns\ds \q-\big\langle Y(t),\d b(t)
\big\rangle_H
\]dt + \frac{1}{2}\Re\mE \big\langle h_{xx}\big(\bar
x(T)\big)x_2^\e(T), x_2^\e(T) \big\rangle_H +
o(\e).
\end{array}
\end{equation}

Next, by the definition of the relaxed
transposition solution to \eqref{op-bsystem3}
(with $P_T$, $J(\cd)$, $K(\cd)$ and $F(\cd)$
given by (\ref{zv2})), and noting (\ref{s7th
max eq2b}), we obtain that
\begin{equation}\label{s7eq9}
 \ba{ll}\ds
-\mE\big\langle h_{xx}\big(\bar x(T)\big)
x_2^\e(T), x_2^\e(T) \big\rangle_H +
\mE\int_0^T \big\langle \dbH_{xx}\big(t,\bar
x(t),\bar u(t),y(t),Y(t)\big) x_2^\e(t), x_2^\e(t) \big\rangle_H dt\\
\ns\ds = \dbE\int_0^T \chi_{E_\e}(t)\big\langle
b_1(t)x_2^\e(t), P(t)^*\d b(t)\big\rangle_{H}
dt + \dbE\int_0^T \chi_{E_\e}(t)\big\langle
P(t)\d b(t),
b_1(t)x_2^\e(t)\big\rangle_{H} dt\\
\ns\ds \q  + \dbE\int_0^T
\chi_{E_\e}(t)\big\langle P(t)\d b(t), \d
b(t)\big\rangle_{H} dt + \dbE\int_0^T
\chi_{E_\e}(t)\big\langle
\d b(t),\widehat Q^{(0)}(0,0,\chi_{E_\e}\d b)(t)\big\rangle_{H} dt \\
\ns\ds \q + \dbE\int_0^T
\chi_{E_\e}(t)\big\langle Q^{(0)}(0,0,\d
b)(t),\d b(t)\big\rangle_{H} dt.
 \ea
\end{equation}

Now, we estimate the terms in the right hand
side of \eqref{s7eq9}. It is clear that
$P(t)^*=P(t)$ for $t\in (0,T)$, and hence
\begin{equation}\label{s7eq9.1}
\begin{array}{ll}\ds
\q\Big|\dbE\int_0^T \chi_{E_\e}(t)\big\langle
b_1(t)x_2^\e(t), P(t)^*\d b(t)\big\rangle_{H}
dt\Big|\\
\ns\ds \leq
|x_2^\e(\cd)|_{L^\infty_\dbF(0,T;L^4(\O;H))}|b_1|_{L^\infty_\dbF(0,T;
\cL(H))}\int_{E_\e}|P(t)^*\d
b(t)|_{L_{\cF_t}^{\frac{4}{3}}(\O;H)}dt
\\
\ns\ds \leq C(x_0)\sqrt{\e}\int_{E_\e} |P(t)\d
b(t)|_{L_{\cF_t}^{\frac{4}{3}}(\O;H)}dt= o(\e).
\end{array}
\end{equation}
Similarly,
\begin{equation}\label{s7eq9.2}
\Big|\dbE\int_0^T \chi_{E_\e}(t)\big\langle
P(t)\d b(t), b_1(t)x_2^\e(t)\big\rangle_{H}
dt\Big| =o(\e).
\end{equation}

In what follows, for any $\tau\in [0,T)$, we
choose $E_{\e}=[\tau,\tau+\e]\subset [0,T]$. We
find a sequence
$\{\beta_n\}_{n=1}^\infty\subset \cH$ (recall
\eqref{cH} for the definition of $\cH$) such
that
$$
\lim_{n\to\infty}\beta_n = \d b \q\mbox{ in }
L^4_\dbF(0,T;H).
$$
Hence,
 \begin{equation}\label{10.9qqq3}
|\beta_n|_{L^4_\dbF(0,T;H)}\le
C(x_0)<\infty,\qq\forall\;n\in\dbN,
\end{equation}
and there is a subsequence
$\{n_k\}_{k=1}^\infty\subset \{
n\}_{n=1}^\infty$ such that
\begin{equation}\label{s7eq9.2-11}
\lim_{k\to\infty} |\b_ {n_k}(t)-\d
b(t)|_{L^4_{\cF_t}(\O;H)} = 0\q\mbox{ for }\ae
t\in [0,T].
\end{equation}
For any $\ell\in\dbN$, let
$t_j=\frac{j-1}{\ell}T$ for $j=1,\cds,\ell+1$.
Since the set of simple processes is dense in
$L^4_\dbF(0,T;L^\infty(\O;\cL(D(A))))$, we can
find a $\ds b_1^\ell\equiv b_1^\ell(t,\o) =
\sum_{j=1}^\ell\chi_{[t_j,t_{j+1})}(t)
f_j(\o)$, where $ f_j\in
L^\infty_{\cF_{t_j}}(\O;\cL(D(A)))$, such that
 \bel{jkl}
\lim_{\ell\to\infty}|b_1^{\ell}-b_1|_{L^4_\dbF(0,T;L^\infty(\O;\cL(D(A))))}=0.
 \ee
It follows that
 \begin{equation}\label{10.9q3}
|b_1^{\ell}|_{
L^4_\dbF(0,T;L^\infty(\O;\cL(H)))}\le
C(x_0)<\infty,\qq\forall\;\ell\in\dbN.
\end{equation}

Denote by
$(P^{\ell}(\cd),Q^{(\cd,\ell)},\widehat
Q^{(\cd,\ell)})$ the relaxed transposition
solution to the equation \eqref{op-bsystem3}
with $K$ replaced by $b_1^{\ell}$, and $P_T$,
$J$ and $F$ given as in \eqref{zv2}. Also,
denote by $Q^{\ell}$ and $\widehat Q^{\ell}$
the corresponding pointwisely defined linear
operators from $\cH$ to
$L^2_\dbF(0,T;L^{\frac{4}{3}}(\O;H))$, given in
Theorem \ref{10.1th}. By Theorem \ref{10.2th1}
and noting \eqref{jkl}--\eqref{10.9q3}, we see
that
\begin{equation}\label{10.9eq6}
\left\{
\begin{array}{ll}\ds
\lim_{\ell\to\infty}\big|\!\big|Q^{(0,\ell)}(0,0,\cd)-Q^{(0)}(0,0,\cd)\big|\!\big|_{\cL(L^2_\dbF(0,T;L^4(\O;H)),\;L^2_\dbF(0,T;L^{\frac{4}{3}}(\O;H)))}=0,\\
\ns\ds \lim_{\ell\to\infty}\big|\!\big|\widehat
Q^{(0,\ell)}(0,0,\cd)-\widehat
Q^{(0)}(0,0,\cd)\big|\!\big|_{\cL(L^2_\dbF(0,T;L^4(\O;H)),\;L^2_\dbF(0,T;L^{\frac{4}{3}}(\O;H)))}=0.
\end{array}
\right.
\end{equation}

Consider the following equation:
\begin{equation}\label{s7fsystem3.1x}
\left\{
\begin{array}{lll}\ds
dx_{2,n_k}^{\e,\ell} =
\big[Ax_{2,n_k}^{\e,\ell} +
a_1(t)x_{2,n_k}^{\e,\ell} \big]dt + \big[
b_1^{\ell}(t) x_{2,n_k}^{\e,\ell} +
\chi_{E_{\e}} (t)\b_ {n_k}(t) \big] dw(t)
&\mbox{ in
}(0,T],\\
\ns\ds x_{2,n_k}^{\e,\ell}(0)=0.
\end{array}
\right.
\end{equation}
We have
\begin{equation}\label{s7th max eq1x}
\begin{array}{ll}\ds
\q\mE|x_{2,n_k}^{\e,\ell}(t)|^4_H \\
\ns\ds = \mE\Big| \int_0^t
S(t-s)a_1(s)x_{2,n_k}^{\e,\ell}(s) ds +
\int_0^t
S(t-s)b_1^{\ell}(s)x_{2,n_k}^{\e,\ell}(s)
dw(s) \\
\ns\ds \qq +
\int_0^t S(t-s)\chi_{E_{\e}}(s)\b_ {n_k}(s) dw(s)\Big|_H^4\\
\ns\ds  \leq C \bigg[\mE\Big| \int_0^t
S(t-s)a_1(s)x_{2,n_k}^{\e,\ell}(s) ds \Big|_H^4
+ \mE\Big|\int_0^t S(t-s)b_1^{\ell}(s)x_{2,n_k}^{\e,\ell}(s) dw(s) \Big|_H^4 \\
\ns\ds \q  + \mE\Big| \int_0^t S(t-s)\chi_{E_{\e}}(s)\b_ {n_k}(s) dw(s)\Big|_H^4\bigg] \\
\ns\ds   \leq
C\[\int_0^t|b_1^{\ell}(s)|_{L^\infty(\O;\cL(H))}^4
\mE|x_{2,n_k}^{\e,\ell}(s)|_H^4 ds  + \e
\int_{E_{\e}}\mE|\beta_{n_k}(s)|_{H}^4ds \].
\end{array}
\end{equation}
By \eqref{10.9qqq3} and (\ref{10.9q3}), thanks
to Gronwall's inequality, \eqref{s7th max eq1x}
leads to
\begin{equation}\label{s7th max eq2x}
|x_{2,n_k}^{\e,\ell}(\cd)|^4_{L^\infty_\dbF(0,T;L^4(\O;H))}
\leq C(x_0,\ell,k)\e^2.
\end{equation}
Here and henceforth, $C(x_0,\ell,k)$ is a
generic constant (depending on $x_0$, $\ell$,
$k$, $T$, $A$ and $C_L$), which may be
different from line to line.
 For any fixed
$i,k\in\dbN$, since $Q^{\ell}\b_ {n_k}\in
L^2_{\dbF}(0,T;L^{\frac{4}{3}}(\O;H))$, by
\eqref{s7th max eq2x}, we find that
\begin{equation}\label{s7eq9.3}
\begin{array}{ll}\ds
\q\Big|\dbE\int_0^T \chi_{E_{\e}}(t)\big\langle
\big(Q^{\ell} \b_
{n_k}\big)(t),x_{2,n_k}^{\e,\ell}(t)\big\rangle_{H}
dt \Big|\leq
|x_{2,n_k}^{\e,\ell}(\cd)|_{L^\infty_\dbF(0,T;L^4(\O;H))}
\int_{E_{\e}}\big|\big(Q^{\ell} \b_ {n_k}\big)(t)\big|_{L^{\frac{4}{3}}_{\cF_t}(\O;H)}dt\\
\ns\ds \leq C(x_0,\ell,k)\sqrt{{\e}}
\int_{E_{\e}}\big|\big(Q^{\ell} \b_
{n_k}\big)(t)\big|_{L^{\frac{4}{3}}_{\cF_t}(\O;H)}dt
= o({\e}), \qq\hbox{as }\e\to0.
\end{array}
\end{equation}
Similarly,
\begin{equation}\label{s7eq9.3x}
\Big|\dbE\int_0^T \chi_{E_{\e}}(t)\big\langle
x_{2,n_k}^{\e,\ell}(t),\big(\widehat
Q^{\ell}\b_ {n_k}\big)(t)\big\rangle_{H} dt
\Big| = o({\e}), \qq\hbox{as }\e\to0.
\end{equation}

From \eqref{10.9eq2} in Theorem \ref{10.1th},
and noting that both $Q^{\ell}$ and $\widehat
Q^{\ell}$ are pointwisely defined, we arrive at
the following equality:
 \bel{wws1}
\begin{array}{ll}
\ds \q  \mE \int_{0}^T \big\langle
\chi_{E_{\e}}(t)\b_ {n_k}(t), \widehat
Q^{(0,\ell )}(0,0,\chi_{E_{\e}} \b_ {n_k})(t)
\big\rangle_{H}dt +  \mE
 \int_{0}^T \big\langle
 Q^{(0,\ell )}(0,0,\chi_{E_{\e}} \b_ {n_k}) (t), \chi_{E_{\e}}\b_ {n_k}(t) \big\rangle_{H}dt \\
\ns\ds =\mE \int_{0}^T
\chi_{E_{\e}}\[\big\langle \big(Q^{\ell} \b_
{n_k}\big)(t), x_{2,n_k}^{\e,\ell} (t)
\big\rangle_{H}+\big\langle x_{2,n_k}^{\e,\ell}
(t), \big(\widehat Q^{\ell} \b_ {n_k}\big)(t)
\big\rangle_{H}\]dt.
\end{array}
 \ee
Hence,
 \bel{wws2}
\begin{array}{ll}\ds
\q\mE \int_{0}^T \big\langle \chi_{E_{\e}}(t)\d
b(t), \widehat Q^{(0)}(0,0,\chi_{E_{\e}}\d
b)(t) \big\rangle_{H}dt + \mE \int_{0}^T
\big\langle
 Q^{(0)}(0,0,\chi_{E_{\e}}\d
b) (t), \chi_{E_{\e}}(t)\d
b(t) \big\rangle_{H}dt \\
\ns\ds \q-\mE \int_{0}^T \chi_{E_{\e}}(t)
\[\big\langle \big(Q^{\ell} \b_ {n_k}\big)(t),
x_{2,n_k}^{\e,\ell} (t)
\big\rangle_{H}+\big\langle x_{2,n_k}^{\e,\ell}
(t), \big(\widehat Q^{\ell} \b_ {n_k}\big)(t)
\big\rangle_{H}\]dt\\
\ns\ds = \mE \int_{0}^T \big\langle
\chi_{E_{\e}}(t)\d b(t), \widehat
Q^{(0)}(0,0,\chi_{E_{\e}} \d b)(t)
\big\rangle_{H}dt + \mE \int_{0}^T \big\langle
 Q^{(0)}(0,0,\chi_{E_{\e}} \d
b) (t), \chi_{E_{\e}}(t)\d
b(t) \big\rangle_{H}dt \\
\ns\ds \q -\mE \int_{0}^T  \big\langle
\chi_{E_{\e}}(t)\b_
{n_k}(t), \widehat Q^{(0,\ell )}(0,0,\chi_{E_{\e}} \b_ {n_k} )(t) \big\rangle_{H}dt \\
\ns\ds \q - \mE \int_{0}^T \big\langle
 Q^{(0,\ell )}(0,0,\chi_{E_{\e}} \b_ {n_k} ) (t), \chi_{E_{\e}}(t)\b_ {n_k}(t) \big\rangle_{H}dt.
\end{array}
 \ee
It is easy to see that
\begin{equation}\label{s7eq9.2xx}
\begin{array}{ll}\ds \q\Big|\mE \int_{0}^T
\big\langle \chi_{E_{\e}}(t)\d b(t), \widehat
Q^{(0)}(0,0,\chi_{E_{\e}} \d b)(t)
\big\rangle_{H}dt  - \mE \int_{0}^T \big\langle
\chi_{E_{\e}}(t)\b_ {n_k}(t), \widehat
Q^{(0,\ell )}(0,0,\chi_{E_{\e}} \b_ {n_k})(t)
\big\rangle_{H}dt \Big|\\
\ns\ds \leq \Big|\mE \int_{0}^T \big\langle
\chi_{E_{\e}}(t)\d b(t), \widehat
Q^{(0)}(0,0,\chi_{E_{\e}} \d b)(t)
\big\rangle_{H}dt - \mE \int_{0}^T \big\langle
\chi_{E_{\e}}(t)\d b(t), \widehat
Q^{(0)}(0,0,\chi_{E_{\e}} \b_ {n_k})(t)
\big\rangle_{H}dt \Big|\\
\ns\ds \q + \Big|\mE \int_{0}^T\!\! \big\langle
\chi_{E_{\e}}(t)\d b(t), \widehat
Q^{(0)}(0,0,\chi_{E_{\e}} \b_ {n_k})(t)
\big\rangle_{H}dt - \mE  \int_{0}^T\!\!
\big\langle \chi_{E_{\e}}(t)\b_ {n_k}(t),
\widehat Q^{(0)}(0,0,\chi_{E_{\e}} \b_
{n_k})(t)
\big\rangle_{H}dt \Big|\\
\ns\ds \q +\!\Big|\mE\! \int_{0}^T\!\!\!
\big\langle \chi_{E_{\e}}(t)\b_ {n_k}(t),
\widehat Q^{(0)}(0,0,\chi_{E_{\e}} \b_
{n_k})(t) \big\rangle_{H}dt - \mE \!\int_{0}^T
\!\!\!\big\langle \chi_{E_{\e}}(t)\b_ {n_k}(t),
\widehat Q^{(0,\ell )}(0,0,\chi_{E_{\e}} \b_
{n_k})(t) \big\rangle_{H}dt \Big|.
\end{array}
\end{equation}
From \eqref{s7eq9.2-11} and the density of the
Lebesgue point, we find that for $\ae\tau\in
[0,T)$, it holds that
\begin{equation}\label{s7eq9.3-1}
\begin{array}{ll}\ds
\q\lim_{k\to\infty}\lim_{\e\to
0}\frac{1}{{\e}}\Big|\mE \int_{0}^T \big\langle
Q^{(0)}(0,0,\chi_{E_{\e}} \d
b)(t),\chi_{E_{\e}}(t)\d b(t) \big\rangle_{H}dt
\\
\ns\ds \qq\qq\q - \mE \int_{0}^T \big\langle
Q^{(0)}(0,0,\chi_{E_{\e}} \d b)(t),\chi_{E_{\e}}(t)\b_ {n_k}(t) \big\rangle_{H}dt \Big|\\
\ns\ds \leq \lim_{k\to\infty}\lim_{\e\to
0}\frac{1}{{\e}}|Q^{(0)}(0,0,\chi_{E_{\e}} \d
b) |_{L^2_\dbF(0,T;L^{\frac{4}{3}}(\O;H))}
\[\int_0^T \chi_{E_{\e}}(t) \(\mE|\d b(t) - \b_ {n_k}(t)|^4_{H}\)^{\frac{1}{2}} dt
\Big]^{\frac{1}{2}}\\
\ns\ds \leq C\lim_{k\to\infty}\lim_{\e\to
0}\frac{1}{{\e}}\big|\chi_{E_{\e}} \d
b\big|_{L^2_\dbF(0,T;L^4(\O;H))}
\[\int_0^T \chi_{E_{\e}}(t) \(\mE|\d b(t) - \b_ {n_k}(t)|^4_{H}\)^{\frac{1}{2}} dt
\Big]^{\frac{1}{2}}\\
\ns\ds \leq C\lim_{k\to\infty}\lim_{\e\to
0}\frac{|\d
b(\tau)|_{L^4_{\cF_\tau}(\O;H)}}{\sqrt{{\e}}}\[\int_0^T
\chi_{E_{\e}}(t) \(\mE|\d b(t) - \b_
{n_k}(t)|^4_{H}\)^{\frac{1}{2}} dt
\Big]^{\frac{1}{2}}\\
\ns\ds = C\lim_{k\to\infty}\lim_{\e\to 0}|\d
b(\tau)|_{L^4_{\cF_\tau}(\O;H)}\[\frac{1}{{\e}}\int_\tau^{\tau+{\e}}
|\d b(t) - \b_ {n_k}(t)|_{L^4_{\cF_t}(\O;H)}^2
dt
\Big]^{\frac{1}{2}}\\
\ns\ds = C\lim_{k\to\infty}|\d
b(\tau)|_{L^4_{\cF_\tau}(\O;H)}|\d b(\tau) -
\b_ {n_k}(\tau)|_{L^4_{\cF_\tau}(\O;H)}
\\\ns\ds= 0.
\end{array}
\end{equation}
Similarly,
\begin{equation}\label{s7eq9.3-1x}
\begin{array}{ll}\ds
\q\lim_{k\to\infty}\lim_{\e\to
0}\frac{1}{{\e}}\Big|\mE \int_{0}^T \big\langle
\chi_{E_{\e}}(t)\d b(t), \widehat
Q^{(0)}(0,0,\chi_{E_{\e}} \b_ {n_k})(t)
\big\rangle_{H}dt \\
\ns\ds \qq\qq\q - \mE  \int_{0}^T \big\langle
\chi_{E_{\e}}(t)\b_ {n_k}(t), \widehat
Q^{(0)}(0,0,\chi_{E_{\e}} \b_ {n_k})(t)
\big\rangle_{H}dt \Big|\\
\ns\ds \leq \lim_{k\to\infty}\lim_{\e\to
0}\frac{1}{{\e}}\big|\widehat
Q^{(0)}(0,0,\chi_{E_{\e}} \b_
{n_k})\big|_{L^2_\dbF(0,T;L^{\frac{4}{3}}(\O;H))}
\[\int_0^T \chi_{E_{\e}}(t) \(\mE|\d b(t) - \b_ {n_k}(t)|^4_{H}\)^{\frac{1}{2}} dt
\Big]^{\frac{1}{2}}\\
\ns\ds \leq C\lim_{k\to\infty}\lim_{\e\to
0}\frac{1}{{\e}}\big|\chi_{E_{\e}} \b_
{n_k}\big|_{L^2_\dbF(0,T;L^4(\O;H))}
\[\int_0^T \chi_{E_{\e}}(t) \(\mE|\d b(t) - \b_ {n_k}(t)|^4_{H}\)^{\frac{1}{2}} dt
\Big]^{\frac{1}{2}}\\
\ns\ds \leq C\lim_{k\to\infty}\lim_{\e\to
0}\frac{1}{{\e}}\Big\{\big|\chi_{E_{\e}} \d
b\big|_{L^2_\dbF(0,T;L^4(\O;H))}
\[\int_0^T \chi_{E_{\e}}(t) \(\mE|\d b(t) - \b_ {n_k}(t)|^4_{H}\)^{\frac{1}{2}} dt
\Big]^{\frac{1}{2}}\\\ns\ds\qq\qq\qq\qq\q+\int_0^T \chi_{E_{\e}}(t) \(\mE|\d b(t) - \b_ {n_k}(t)|^4_{H}\)^{\frac{1}{2}} dt\Big\}\\
\ns\ds \leq C\lim_{k\to\infty}\lim_{\e\to
0}\Big\{\frac{|\d
b(\tau)|_{L^4_{\cF_\tau}(\O;H)}}{\sqrt{{\e}}}\[\int_0^T
\chi_{E_{\e}}(t) \(\mE|\d b(t) - \b_
{n_k}(t)|^4_{H}\)^{\frac{1}{2}} dt
\Big]^{\frac{1}{2}}\\\ns\ds\qq\qq\qq\q+\frac{1}{{\e}}\int_0^T
\chi_{E_{\e}}(t) \(\mE|\d b(t) - \b_
{n_k}(t)|^4_{H}\)^{\frac{1}{2}} dt
\Big\}\\
\ns\ds = C\lim_{k\to\infty}\lim_{\e\to
0}\Big\{|\d
b(\tau)|_{L^4_{\cF_\tau}(\O;H)}\[\frac{1}{{\e}}\int_\tau^{\tau+{\e}}
|\d b(t) - \b_ {n_k}(t)|_{L^4_{\cF_t}(\O;H)}^2
dt
\Big]^{\frac{1}{2}}\\\ns\ds\qq\qq\qq\q+\frac{1}{{\e}}\int_\tau^{\tau+{\e}}
|\d b(t) - \b_ {n_k}(t)|_{L^4_{\cF_t}(\O;H)}^2
dt
\Big\}\\
\ns\ds = C\lim_{k\to\infty}[|\d
b(\tau)|_{L^4_{\cF_\tau}(\O;H)}|\d b(\tau) -
\b_ {n_k}(\tau)|_{L^4_{\cF_\tau}(\O;H)}+|\d
b(\tau) - \b_
{n_k}(\tau)|_{L^4_{\cF_\tau}(\O;H)}^2]
\\\ns\ds= 0.
\end{array}
\end{equation}

From \eqref{10.9eq6} and the density of the
Lebesgue point, we find that for $\ae\tau\in
[0,T)$, it holds that
\begin{equation}\label{s7eq9.2xxx}
\begin{array}{ll}\ds
\q\lim_{k\to\infty}\lim_{\ell\to\infty}\lim_{\e\to
0}\frac{1}{{\e}}\Big|\mE \int_{0}^T \big\langle
\chi_{E_{\e}}(t)\b_ {n_k}(t), \widehat
Q^{(0)}(0,0,\chi_{E_{\e}} \b_ {n_k})(t)
\big\rangle_{H}dt \\
\ns\ds \qq\qq\qq\qq\q - \mE \int_{0}^T
\big\langle \chi_{E_{\e}}(t)\b_
{n_k}(t),\widehat Q^{(0,\ell
)}(0,0,\chi_{E_{\e}}
\b_ {n_k})(t) \big\rangle_{H}dt \Big|\\
\ns\ds \leq
\lim_{k\to\infty}\lim_{\ell\to\infty}\lim_{\e\to
0}\frac{1}{{\e}}\big|\chi_{E_{\e}} \b_
{n_k}\big|_{L^2_\dbF(0,T;L^{4}(\O;H))}
\(\int_0^T\big|\widehat Q^{(0)}(0,0,\chi_{E_{\e}}\b_ {n_k})\\
\ns\ds \hspace{3.5cm}-\widehat Q^{(0,\ell
)}(0,0,\chi_{E_{\e}}\b_
{n_k})\big|^2_{L_{\cF_s}^{\frac{4}{3}}(\O;H))}ds\)^{\frac{1}{2}}
\\
\ns\ds \leq
\!\lim_{k\to\infty}\lim_{\ell\to\infty}\lim_{\e\to
0}\frac{1}{{\e}}\big|\chi_{E_{\e}} \b_
{n_k}\big|^2_{L^2_\dbF(0,T;L^{4}(\O;H))}\!
\big|\!\big|\widehat
Q^{(0,\ell)}(0,0,\cd)\!-\!\widehat
Q^{(0)}(0,0,\cd)\!\big|\!\big|_{\cL(L^2_\dbF(0,T;L^4(\O;H)),\,L^2_\dbF(0,T;L^{\frac{4}{3}}(\O;H)))}
\\
\ns\ds = \lim_{k\to\infty}\lim_{\ell\to\infty}
\big|\b_ {n_k}(\tau,\cd)\big|^2_{L^{4}(\O;H)}
\big|\!\big|\widehat Q^{(0,\ell)}(0,0,\cd)-\widehat Q^{(0)}(0,0,\cd)\big|\!\big|_{\cL(L^2_\dbF(0,T;L^4(\O;H)),\;L^2_\dbF(0,T;L^{\frac{4}{3}}(\O;H)))}\\
\ns\ds = 0.
\end{array}
\end{equation}
From \eqref{s7eq9.2xx}--\eqref{s7eq9.2xxx}, we
find that
\begin{equation}\label{s7eq9.3-1xx}
\begin{array}{ll}\ds
\q\lim_{k\to\infty}\lim_{\ell\to\infty}\lim_{\e\to
0}\frac{1}{{\e}}\Big|\mE \int_{0}^T \big\langle
\chi_{E_{\e}}(t)\d b(t), \widehat
Q^{(0)}(0,0,\chi_{E_{\e}} \d b)(t)
\big\rangle_{H}dt \\
\ns\ds \qq - \mE  \int_{0}^T \big\langle
\chi_{E_{\e}}(t)\b_ {n_k}(t), \widehat
Q^{(0,\ell )}(0,0,\chi_{E_{\e}} \b_ {n_k})(t)
\big\rangle_{H}dt \Big| = 0.
\end{array}
\end{equation}
By a similar argument, we obtain that
\begin{equation}\label{s7eq9.3-1xxx}
\begin{array}{ll}\ds
\q\lim_{k\to\infty}\lim_{\ell\to\infty}\lim_{\e\to
0}\frac{1}{{\e}}\Big|\mE \int_{0}^T \big\langle
Q^{(0)}(0,0,\chi_{E_{\e}} \d
b)(t),\chi_{E_{\e}}(t)\d b(t)
\big\rangle_{H}dt \\
\ns\ds \qq - \mE  \int_{0}^T \big\langle
Q^{(0,\ell )}(0,0,\chi_{E_{\e}} \b_
{n_k})(t),\chi_{E_{\e}}(t)\b_ {n_k}(t)
\big\rangle_{H}dt \Big| = 0.
\end{array}
\end{equation}

From \eqref{s7eq9.3}--\eqref{s7eq9.3x},
\eqref{wws1}--\eqref{wws2} and
\eqref{s7eq9.3-1xx}--\eqref{s7eq9.3-1xxx}, we
obtain that
\begin{equation}\label{s7eq9.2-111}
 \ba{ll}\ds
\Big|\dbE\int_0^T \chi_{E_{\e}}(t)\big\langle
\d b(t),\widehat Q^{(0)}(0,0,\chi_{E_{\e}}\d
b)(t)\big\rangle_{H} dt + \dbE\int_0^T
\chi_{E_{\e}}(t)\big\langle Q^{(0)}(0,0,\d
b)(t),\d b(t)\big\rangle_{H} dt\Big|
\\\ns\ds=o({\e}), \qq\hbox{as }\e\to0.
 \ea
\end{equation}

Therefore, we have
\begin{equation}\label{s7eq10}
 \ba{ll}\ds
 \q\cJ(u^{\e}(\cd)) -  \cJ(\bar u(\cd))\\\ns\ds = \Re\mE\int_0^T \[
 \d g(t) - \big\langle
y(t),\d a(t)\big\rangle_H -\big\langle Y(t),\d
b(t) \big\rangle_H - \frac{1}{2}\big\langle
P(t)\d b(t), \d b(t) \big\rangle_H
 \]\chi_{E_{\e}}(t)dt + o({\e}).
 \ea
\end{equation}
Since $\bar u(\cd)$ is the optimal control,
$\cJ(u^{\e}(\cd)) - \cJ(\bar u(\cd))\geq 0$.
Thus,
\begin{equation}\label{s7eq11}
\Re\mE\int_0^T \chi_{E_{\e}}(t)\[
  \big\langle
y(t),\d a(t)\big\rangle_H +\big\langle Y(t),\d
b(t) \big\rangle_H -\d g(t)+
\frac{1}{2}\big\langle P(t)\d b(t), \d b(t)
\big\rangle_H
 \]dt \leq o({\e}),
\end{equation}
as $\e\to0$.

Finally,  similar to \cite{Kushner, Zhou}, from
\eqref{s7eq11}, we obtain \eqref{s7maxth ine1}.
This completes the proof of Theorem \ref{s7th
max}.
\endpf

\br
1) We believe that  $b_x(\cd,\bar x(\cd),\bar
u(\cd))\in
L^4_\dbF(0,T;L^\infty(\O;\cL(D(A))))$ is a
technical condition in Theorem \ref{s7th max}
but we cannot drop it at this moment (because
we need to use Theorem \ref{10.1th}). It is
easy to see that this condition is satisfied
for one of the following cases:

i) The operator $A$ is a bounded linear
operator on $H$;

ii) The diffusion term $b(t,x,u)$ is
independent of the state variable $x$; or

iii) Some further regularities for $x_0$,
$a(\cd, \cd,\cd)$ and $b(\cd, \cd,\cd)$ are
imposed, say $x_0\in L^8_{\cF_0}(\O;D(A))$, and
the Assumption (A1) holds also when the space
$H$ is replaced by $D(A)$.

\ms

2) If the equation \eqref{op-bsystem3}, with
$P_T$, $J(\cd)$, $K(\cd)$ and $F(\cd)$ given by
\eqref{zv2}, admits a transposition solution
$(P(\cd),Q(\cd))$, then the assumption
$b_x(\cd,\bar x(\cd),\bar u(\cd))\in
L^4_\dbF(0,T;L^\infty(\O;\cL(D(A))))$ is not
needed (for the same conclusion in Theorem
\ref{s7th max}). Indeed, in this case, by
Definition \ref{op-definition2}, instead of
\eqref{s7eq9}, we have
\begin{equation}\label{s7eq9zx}
 \ba{ll}\ds
-\mE\big\langle h_{xx}\big(\bar x(T)\big)
x_2^\e(T), x_2^\e(T) \big\rangle_H +
\mE\int_0^T \big\langle \dbH_{xx}\big(t,\bar
x(t),\bar u(t),y(t),Y(t)\big) x_2^\e(t), x_2^\e(t) \big\rangle_H dt\\
\ns\ds = \dbE\int_0^T \chi_{E_\e}(s)\big\langle
b_1(s)x_2^\e(s), P(s)^*\d b(s)\big\rangle_{H}
ds + \dbE\int_0^T \chi_{E_\e}(s)\big\langle
P(s)\d b(s),
b_1(s)x_2^\e(s)\big\rangle_{H} ds\\
\ns\ds \q  + \dbE\int_0^T
\chi_{E_\e}(s)\big\langle P(s)\d b(s), \d
b(s)\big\rangle_{H} ds + \dbE\int_0^T
\chi_{E_\e}(s)\big\langle
Q(s)\d b(s),x_2^\e(s)\big\rangle_{H} ds \\
\ns\ds \q + \dbE\int_0^T
\chi_{E_\e}(s)\big\langle Q(s)x_2^\e(s),\d
b(s)\big\rangle_{H} ds.
 \ea
\end{equation}
The estimates \eqref{s7eq9.1}-\eqref{s7eq9.2}
are still valid. On the other hand, by
$Q(\cd)\d b(\cd)\in
L^1_{\dbF}(0,T;L^{\frac{4}{3}}(\O;H))$, it
holds that
\begin{equation}\label{s7eqqq9.3}
\begin{array}{ll}\ds
\q\dbE\int_0^T \chi_{E_\e}(s)\big\langle
Q(s)\d b(s),x_2^\e(s)\big\rangle_{H} ds \\
\ns\ds \leq
|x_2^\e(s)|_{L^\infty_\dbF(0,T;L^4(\O;H))}
\int_{E_\e}|Q(s)\d b(s)|_{L^{\frac{4}{3}}(\O;H)}ds\\
\ns\ds \leq C\sqrt{\e} \int_{E_\e}|Q(s)\d
b(s)|_{L^{\frac{4}{3}}(\O;H)}ds= o(\e).
\end{array}
\end{equation}
Similarly, noting that $Q(t)^*=Q(t)$ for $\ae
t\in (0,T)$, we obtain that
\begin{equation}\label{s7eqqq9.2}
\Big|\dbE\int_0^T \chi_{E_\e}(s)\big\langle
Q(s)x_2^\e(s),\d b(s)\big\rangle_{H} ds\Big|
=\Big|\dbE\int_0^T \chi_{E_\e}(s)\big\langle
x_2^\e(s),Q(s)^*\d b(s)\big\rangle_{H} ds\Big|
=o(\e).
\end{equation}
Combining \eqref{s7eq9zx},
\eqref{s7eq9.1}-\eqref{s7eq9.2} and
\eqref{s7eqqq9.3}--\eqref{s7eqqq9.2}, we still
have \eqref{s7eq10}, which leads to the desired
result.

\ms

3) For concrete equations, say for the
controlled stochastic heat equations, one may
obtain better results than that of Theorem
\ref{s7th max}. Related work will be presented
elsewhere.
\er

\section*{Acknowledgement}

This work is supported by the NSF of China
under grants 10831007 and 11101070. Part of
this work was finished when the first ({\it
resp.},  second) author worked at the ``Basque
Center for Applied Mathematics, Basque Country,
Spain" ({\it resp.}, ``Key Laboratory of
Systems and Control, Academy of Mathematics and
System Science, Chinese Academy of Sciences,
China"). The second author would like to
express his appreciation to Professor Gengsheng
Wang for his helpful comments on the first two
versions of this paper.


\begin{thebibliography}{99}


\bibitem{Al-H1}A. Al-Hussein. \it Backward stochastic partial differential equations
driven by infinite dimensional martingales and
applications. \sl Stochastics. \rm {\bf 81}
(2009), 601--626.

\bibitem{Al-H2}A. Al-Hussein. \it Necessary conditions for optimal control
of stochastic evolution equations in Hilbert
spaces. \sl Appl. Math. Optim. \rm {\bf 63}
(2011),  385--400.

\bibitem{AGY}V.V. Anh, W. Grecksch and J. Yong. \it Regularity of
backward stochastic Volterra integral equations
in Hilbert spaces. \sl Stoch. Anal. Appl. \rm
{\bf 29} (2010), 146--168.

\bibitem{Bensoussan1} A.~Bensoussan. \it Lecture on Stochastic
Control. \rm In: \sl Nonlinear Filtering and
Stochastic Control. \rm Edited by S. K. Mitter
and A. Moro. Lecture Notes in Mathematics, vol.
{\bf 972}. Springer-Verlag, Berlin, 1982,
1--62.

\bibitem{Bensoussan2}A. Bensoussan. \it Stochastic maximum principle for distributed
parameter systems. \sl J. Franklin Inst. \rm
{\bf 315} (1983), 387--406.

\bibitem{Bismut} J.-M. Bismut. \sl Analyse Convexe et Probabiliti\'es. \rm Ph D Thesis, Facult\'e des
Sciences de Paris, Paris, France, 1973.

\bibitem{B1}J.-M. Bismut. \it An introductory approach to duality in optimal stochastic
control. \sl SIAM Rev. \rm {\bf 20} (1978),
62--78.

\bibitem{C} J. B. Conway. \sl A Course in Functional Analysis (Second Edition). \rm Springer-Verlag, New York,
1994.

\bibitem{Prato} G.~Da Prato and J.~Zabczyk. \sl Stochastic Equations in Infinite Dimensions. \rm Cambridge University Press,
Cambridge, 1992.

\bibitem{Halmos} P. R. Halmos. \sl A Hilbert Space Problem Book (Second Edition, Revised and Enlarged). \rm Springer-Verlag, New York, 1982.

\bibitem{Haussmann} U.~G.~Haussmann. \it General necessary conditions for optimal control of stochastic system. \sl
Math. Prog. Study. \rm {\bf 6} (1976), 34--48.

\bibitem{HP1} Y. Hu and S. Peng. \it Maximum principle for semilinear stochastic evolution control
systems. \sl Stoch. \& Stoch. Rep. \rm {\bf 33}
(1990), 159--180.

\bibitem{HP2} Y. Hu and S. Peng. \it Adapted solution of backward semilinear stochastic evolution
equations. \sl Stoch. Anal. \& Appl. \rm {\bf
9} (1991), 445--459.

\bibitem{Kushner} H.~J.~Kushner. \it Necessary conditions for continuous parameter stochastic optimization
problems. \sl SIAM J. Control. \rm {\bf 10}
(1972), 550--565.

\bibitem{LY} X.~Li and J.~Yong. \sl Optimal Control Theory for Infinite-Dimensional
Systems. \rm Systems \& Control: Foundations \&
Applications. Birkh\"auser Boston, Inc.,
Boston, MA, \rm 1995.

\bibitem{LK} K.~Liu. \sl Stability of Infinite Dimensional Stochastic Differential Equations with Applications.
\rm Pitman Monographs and Surveys in Pure and
Applied Mathematics, vol. {\bf 135}. Chapman \&
Hall/CRC, 2006.

\bibitem{LYZ}Q. L\"{u},  J.~Yong and X.~Zhang. \it Representation of It\^o integrals by Lebesgue/Bochner
integrals. \sl J. Eur. Math. Soc. \rm {\bf 14}
(2012), 1795--1823.

\bibitem{LZ} Q. L\"{u} and X.~Zhang. \it Well-posedness of backward stochastic differential equations with
general filtration. \rm In submission. (see
http://arxiv.org/abs/1010.0026)

\bibitem{MY}J. Ma and J. Yong. \sl Forward-Backward Stochastic Differential
Equations and Their Applications. \rm Lecture
Notes in Math. vol. {\bf 1702}.
Springer-Verlag, New York, 1999.

\bibitem{MM}N. I. Mahmudova and M. A. McKibben. \it On backward stochastic evolution equations in Hilbert spaces and
optimal control. \sl Nonlinear Anal. \rm {\bf
67} (2007), 1260--1274.

\bibitem{NVW1}J. M. A. M. van Neerven, M. C. Veraar and L. W. Weis. \it Stochastic
integration in UMD Banach spaces. \sl Ann.
Probab. \rm {\bf 35} (2007), 1438--1478.

\bibitem{NVW2}J. M. A. M. van Neerven, M. C. Veraar and L. W. Weis. \it Stochastic
evolution equations in UMD Banach spaces. \sl
J. Funct. Anal. \rm {\bf 255} (2008), 940--993.

\bibitem{PP} E.~Pardoux and S.~Peng. \it Adapted solution of backward
stochastic equation. \sl Systems Control Lett.
\rm {\bf 14} (1990), 55--61.

\bibitem{Peng1} S.~Peng. \it A general stochastic maximum principle for optimal
control problems. \sl SIAM J. Control Optim.
\rm {\bf 28} (1990), 966--979.

\bibitem{PC} L. S. Pontryagin, V. G. Boltyanskii, R. V. Gamkrelidze and
E. F. Mischenko. \sl Mathematical Theory of
Optimal Processes. \rm Wiley, New York, 1962.

\bibitem{TL} S. Tang and X. Li. \it Maximum principle for optimal
control of distributed parameter stochastic
systems with random jumps. \rm In: \sl
Differential Equations, Dynamical Systems, and
Control Science. \rm Lecture Notes in Pure and
Appl. Math. vol. {\bf 152}. Dekker, New York,
1994, 867--890.

\bibitem{Tu} C. Tudor. \it Optimal control for semilinear stochastic evolution
equations. \sl Appl. Math. Optim. \rm {\bf 20}
(1989), 319--331.

\bibitem{WZ}  P.~Wang and X.~Zhang. \it Numerical solutions of
backward stochastic differential equations: a
finite transposition method. \sl C. R. Math.
Acad. Sci. Paris. \rm {\bf 349} (2011),
901--903.

\bibitem{YZ} J.~Yong and X. Y.~Zhou. \sl Stochastic Controls: Hamiltonian Systems and
HJB Equations. \rm Springer-Verlag, New York,
1999.

\bibitem{Zhou} X. Y.~Zhou. \it On the necessary conditions of optimal controls for
stochastic partial differential equations. \sl
SIAM J. Control Optim. \rm {\bf 31} (1993),
1462--1478.

\end{thebibliography}
\end{document}